\documentclass[12pt]{article}
\usepackage{mathrsfs}
\usepackage{fullpage}
\usepackage{url}
\usepackage{amsfonts}
\usepackage{graphicx}
\usepackage{mathrsfs}
\usepackage{amsmath}
\usepackage{amssymb}
\usepackage{float}
\usepackage{rotating,color}
\usepackage{xcolor}

\usepackage{natbib}
\setcitestyle{open={(},close={)}}
\usepackage[colorlinks=true,citecolor=blue]{hyperref}

\newcommand{\be}{\begin{equation}}
\newcommand{\ee}{\end{equation}}
\newcommand{\beaa}{\begin{eqnarray*}}
\newcommand{\eeaa}{\end{eqnarray*}}
\newcommand{\bea}{\begin{eqnarray}}
\newcommand{\eea}{\end{eqnarray}}
\newcommand{\lbl}{\label}

\newcommand{\bd}{\bold}

\newcommand{\red}[1]{\color{red}}

\newtheorem{theorem}{ \noindent T{\footnotesize HEOREM}}
\newtheorem{prop}{ \noindent P{\footnotesize ROPOSITION}}
\newtheorem{lemma}{ \noindent L{\footnotesize EMMA}}[section]

\newtheorem{example}{Example}[section]
\newtheorem{result}{Result}[section]
\newtheorem{remark}{Remark}[section]

\newcommand{\bm}{\boldsymbol}

\def\tr{\mathrm {tr}}

\def\B{{\bf B}}

\def\S{{\bf S}}
\def\R{{\bf R}}

\def\bmu{{\bm \mu}}

\def\X{{\bm X}}

\def\D{{\bf D}}

\def\tr{\mathrm {tr}}

\def\bms{{\bm\Sigma}}

\usepackage{mathrsfs}
\usepackage{fullpage}
\usepackage{amsfonts}
\usepackage{graphicx}
\usepackage{mathrsfs}
\usepackage{amsmath}
\usepackage{amssymb}
\usepackage{float}
\usepackage{rotating,color}
\usepackage{xcolor}
\begin{document}
 \title{Mean Test with Fewer Observation than Dimension and Ratio Unbiased Estimator for  Correlation Matrix}

\author{\textbf{Tiefeng Jiang} \\
\small{School of Statistics}\\
\small{University of Minnesota}\\
\small{224 Church St SE
Minneapolis, MN 55455}\\
  \texttt{\small jiang040@umn.edu}
\and
\textbf{Ping Li} \\
\small{Cognitive Computing Lab}\\
\small{Baidu Research}\\
\small{10900 NE 8th St. Bellevue, WA 98004}\\
  \texttt{\small liping11@baidu.com}
}

\date{}

\maketitle

\begin{abstract}
Hotelling's T-squared test is a classical tool to test if the normal mean of a multivariate normal distribution is a specified one or the means of two multivariate normal means are equal. When the population dimension is higher than the sample size, the test is no longer applicable. Under this situation, in this paper we revisit the tests proposed by \cite{srivastava2008test}, who  revise the Hotelling's statistics by replacing Wishart matrices with their diagonal matrices. They show the revised statistics are asymptotically normal. We use the random matrix theory to examine their statistics again and find that their discovery is just part of the big picture. In fact, we prove that their statistics, decided by the Euclidean norm of the population correlation matrix, can go to normal, mixing chi-squared distributions and a convolution of both. Examples are provided to show the phase transition phenomenon between the normal and mixing chi-squared distributions. The second contribution of ours is a rigorous derivation of an asymptotic ratio-unbiased-estimator of the squared Euclidean norm of the correlation matrix.

\end{abstract}

\section{Introduction}\label{sec1}

Among many statistical hypothesis testing problems,  the Hotelling's $T^2$ tests \citep{hotelling1931generalization} are classic ones to study if the mean of a multivariate normal distribution is equal to the given one, or if the means of two normal distributions are equal. Let us quickly review the two problems. First, assume $\{\X_1,\cdots,\X_n\}$ is a random sample from a  $p$-dimensional  normal
distribution $N_p(\bmu, \bms)$, where $\bmu$ is the mean vector and $\bms$ is the covariance matrix. Consider the test $H_0: \bmu=\bm 0 ~~\textrm{versus}~~H_1: \bmu\not=\bm 0.$ The Hotelling $T^2$ test statistic is defined by $T^2=(\bar{\X}-\mu)^T\hat{\S}^{-1}(\bar{\X}-\mu)$, where
$\bar{\X}$ and $\hat{\S}$ are the sample mean and the sample covariance matrix.  For the two population test, assume that $\{\X_{i1},\cdots,\X_{in_i}\}$ for $i=1,2$ are two independent
random samples  from  $N_p(\bmu_1,\bms)$ and $N_p(\bmu_2,\bms)$, respectively.  We aim to test $H_0:\bmu_1=\bmu_2\ \ \mbox{versus}\ \  H_1:\bmu_1\neq\bmu_2.$ The corresponding Hotelling's $T^2$ test statistic is given  by $\frac{n_1n_2}{n_1+n_2}(\bar{\X}_1-\bar{\X}_2)^{T} \hat{\S}^{-1}(\bar{\X}_1-\bar{\X}_2)$
 where $\bar{\X}_i$ is the sample mean vector of the $i$-th sample and $\hat{\S}$ is the pooled sample covariance matrix.

As a consequence of the likelihood ratio test, the above two Hotelling $T^2$ tests have very nice properties. For example, for the one-sample test, $T^2$ is invariant under linear transformations.  Also, $T^2$ is a uniformly most powerful test under the general linear group; see, for example, p. 211 from \cite{muirhead1982aspects} or p. 190 from \cite{anderson2003introduction}.

Despite its important role in classical statistics, Hotelling's $T^2$  have  some limitations. First, in order to guarantee that the sample covariance matrix is invertible, the sample size have to be larger than the population dimension. These are not true for some modern data in which the population dimension is larger or even much larger than the sample size. For example, for DNA microarray data, thousands of gene expression levels are often measured on a small number of subjects. From
\url{http://genomics-pubs.princeton.edu/oncology/}, which is a popular colon dataset,  one can see $p=2000$, $n_1=22$ and $n_2=40$.
Second, as \cite{bai1996effect} have showed,  the Hotelling two-sample test is  inconsistent as the population dimension and the sample size are comparable.

%\url{http://microarray.princeton.edu/oncology/affydata/index.html}

\newpage

To accommodate data with the feature of large $p$ small $n$ for the above one-sample and two-sample testing problems, \cite{dempster1958high,dempster1960significance} study the so-called ``non-exact" test. \cite{bai1996effect} construct a new test statistic by removing the inverse matrix from the definition of $T^2$ given earlier. Starting from this century, investigators begin to conceive new test statistics.
\cite{srivastava2008test}, \cite{srivastava2009test}, \cite{chen2010two},  and \cite{dong2016shrinkage} replace the inverse of the sample covariance matrix by mild quantities. \cite{lopes2011more} and \cite{srivastava2016raptt} use a sort of dimension reduction method to lower the population dimension and then use Hotelling's $T^2$ tests. For data with a certain of sparse nature, \cite{cai2014two} establish tests by measuring the maximum componentwise mean difference of appropriately transformed observations. Other research related  to sparsity  can be seen from \cite{zhong2013tests}, \cite{chen2014two}, \cite{wang2015high}, \cite{gregory2015two}  and \cite{guo2016tests}. The researchers \cite{biswas2014nonparametric}, \cite{chang2017simulation}, \cite{chakraborty2017tests}, and \cite{xue2020distribution}  consider the above test under non-normal assumptions.  Other contributions include \cite{park2013test} and  \cite{wu2006multivariate} for the consideration of scale-invariant tests and \cite{gretton2012optimal}  for  a kernel-based discrepancy measure. \cite{feng2016multivariate} work on a two-sample location problem via a multivariate-sign-based high-dimensional test. \cite{xu2016adaptive} and \cite{zhang2020simple_JASA} construct test statistics for the  two-sample test problem by $l_q$-norm. Recently, built upon \cite{chen2011regularized}, the authors  \citep{li2020adaptable} investigate the Hotelling's $T^2$ by a ridge-regularized method. \cite{zhang2020simple} extend the tests by \cite{wu2006multivariate} and  \cite{zhang2020simple_JASA} through a modification of the test  by \cite{srivastava2008test}.

In this paper, we will revisit the procedure by \cite{srivastava2008test} for the one-sample and two-sample problems of testing means.  We will briefly review them next and state our motivation.

\newpage

Let $\{\X_1,\cdots,\X_n\}$ be a random sample from a  $p$-dimensional  normal
distribution $N_p(\bmu, \bms)$ with correlation matrix $\bd{R}$. The sample mean and the sample covariance matrix defined by
\bea\lbl{wuoi0}
\bar{\X}=\frac{1}{n}\sum_{i=1}^n\X_i\ \ \mbox{and}\ \ \hat{\S}=\frac{1}{n}\sum_{i=1}^n (\X_i-\bar{\X})(\X_i-\bar{\X})^T.
\eea
Let $\hat{\D}$ be the diagonal matrix of $\hat{\S}$ and $\hat{\R}$ be the sample correlation matrix defined by $\hat{\R}=\hat{\D}^{-1/2}\hat{\S}\hat{\D}^{-1/2}.$
For the testing problem $H_0: \bmu=\bm 0$ vs $H_1: \bmu\not=\bm 0$, \cite{srivastava2008test} propose the following test statistic.
\begin{align}\lbl{jin_wuzu}
T_{SD}=\frac{n\bar{\X}^{T} \hat{\D}^{-1}\bar{\X}-p(n-1)(n-3)^{-1}}{\sqrt{2[\tr(\hat{\R}^2)-p^2(n-1)^{-1}]}}.
\end{align}
Under certain conditions, they show $T_{SD}$ converges to a normal distribution.

For two population case, the test is $H_0:\bmu_1=\bmu_2$ vs $H_1:\bmu_1\neq\bmu_2$. Assume that
$\{\X_{i1},\cdots,\X_{in_i}\}$ for $i=1,2$ are two independent
random samples from  $N_p(\bmu_1,\bms)$ and $N_p(\bmu_2,\bms)$, respectively. Let $\bar{\X}_i$ be the sample mean for the $i$-th sample and $\hat{\S}$ be the pooled sample covariance matrix defined by
\bea\lbl{duck_yazi}
\hat{\S}=\frac{1}{n_1+n_2}\Big[\sum_{j=1}^{n_1}
(\bd{X}_{1j}-\bar{\bd{X}}_1)(\bd{X}_{1j}-\bar{\bd{X}}_1)^T+\sum_{j=1}^{n_2}
(\bd{X}_{2j}-\bar{\bd{X}}_2)(\bd{X}_{2j}-\bar{\bd{X}}_2)^T\Big].
\eea
Assume $\hat{\D}$ is the diagonal matrix of $\hat{\S}$ and
$\hat{\R}=\hat{\D}^{-1/2}\hat{\S}\hat{\D}^{-1/2}$ is the pooled sample correlation matrix.  \cite{srivastava2008test} consider the following statistic defined by
\begin{align}\lbl{u82su}
T'_{SD}=\frac{\frac{n_1n_2}{n_1+n_2}(\bar{\X}_1-\bar{\X}_2)^{T} \hat{\D}^{-1}(\bar{\X}_1-\bar{\X}_2)-\frac{(n_1+n_2-2)p}{n_1+n_2-4}}
{\sqrt{2\big[\tr(\hat{\R}^2)-\frac{p^2}{n_1+n_2-2}\big]}}.
\end{align}
Under certain conditions, they show $T'_{SD}$ converges to a normal distribution. In fact, to improve the convergence speeds, they actually add a term $c_{p,n}=1+\frac{\tr(\hat{\R}^2)}{p^{3/2}}$ under the squared roots in the denominators of $T_{SD}$ and $T_{SD}'$, respectively. Based on their assumptions, $c_{p,n}$ goes to one. However, we will study a more general case in which $c_{p,n}$ may not go to one, and sometimes it even goes to infinity. This is the reason we dump the term $c_{p,n}$ from both $T_{SD}$ and $T_{SD}'$, respectively.

Evidently,  $T_{SD}$ and $T'_{SD}$ have a nice  property of scale-invariance, i.e., they are not changed if data are multiplied by a constant. Also, they can be directly computed. Our motivation in this paper to reexamine the tests by \cite{srivastava2008test} has three folds. (a) The conditions to guarantee the central limit theorems of $T_{SD}$ and $T'_{SD}$ are somehow stringent; see \eqref{SD08a}. Also the sample size and population dimension have to  satisfy that $p=o(n^2)$.  To make the method more applicable, we hope to relax the condition imposed on $\bms$ as well as that $p$ on $n$. (b) We would like to understand an interesting  observation by \cite{zhang2020simple} from their simulation: the distribution of $T'_{SD}$ sometimes looks like a normal curve, and other times it looks like a chi-square curve. (c) We plan to give a rigorous proof of the major ingredient of this theory, that is, an asymptotic ratio-unbiased-estimator of $\mbox{tr}(\bd{R}^2)$, where $\bd{R}$ is the population correlation matrix of $N_p(\bmu, \bms)$ aforementioned.

Now we state our findings. For (a), we have obtained the asymptotic distributions of $T_{SD}$ and $T'_{SD}$ for arbitrary $\bms$ in Theorems~\ref{Theorem1} and~\ref{Theorem2}, respectively. In particular, Theorem~\ref{Theorem1} holds for two extreme cases: the independent case with $\bd{R}=\bd{I}_p$ and the most dependent case, i.e.,  all entries of the population vector are identical. Our theory says the restriction on $p$ and $n$ will be $p=o(n^2)$ if the entries of the population vector are not far from independent. However, for dependent or very dependent case, our conclusion holds as long as $p$ is not more than a  polynomial order of $n$. For (b), we successfully understand the observation by \cite{zhang2020simple}. In fact, there are  indeed transition phenomena of the limiting distributions of $T_{SD}$ and $T'_{SD}$. They are sometimes normal, mixing chi-squared or the sum of two independent random variables, one has normal distribution and the other has a mixing chi-squared distribution. This can be quickly seen from \eqref{aiy327} and \eqref{o4398gh}. As for (c),  for arbitrary $\bms$ we have proved rigorously that an asymptotic ratio-unbiased-estimator of $\mbox{tr}(\bd{R}^2)$ is $\mbox{tr}(\hat{\bd{R}}^2)-\frac{p(p-1)}{n-1}$ in Theorem~\ref{Theorem3}. Especially the theorem is true for two extreme cases: $\bms$ is diagonal or proportional to a matrix whose entries are all equal to $1$. In addition, some of our partial calculations and heuristics indicate that the above phase transition from the Gaussian to mixing chi-squared distributions is possibly a universal phenomenon. We will present three testing procedures to justify this claim.

The above solution is conducted through the understanding of the sample correlation matrix $\hat{\R}$, a special random matrix, defined below \eqref{wuoi0}. Unlike the Gaussian orhtogonal/unitary/symplectic ensemble or Wishart matrices, $\hat{\R}$ lacks the orthogonal-invariant property. As a consequence, investigating the sample correlation matrix always cost extra energy than  working on  other popular matrices.   We employ the machinery for $\hat{\R}$ developed by  \cite{jiang2004limiting}, \cite{jiang2009variance}, \cite{cai2013distributions}, \cite{jiang2013central}, and \cite{fan2019largest}. In particular,  we extend the method conceived by \cite{jiang2004limiting} and  \cite{fan2019largest} to prove a weak law of large numbers for $\mbox{tr}(\hat{\bd{R}}^2)$ en route to obtain an asymptotic ratio-unbiased-estimator. More elaboration are provided in Section~\ref{Sec_wlln}.

% {\red IIn addition to simulation on powers of tests and limiting curves, we give applications to problems from three fields:  securities, neuroscience and genetics.}
The remaining of the paper is organized as follows. The one-sample mean and two-sample mean problems are studied in Sections~\ref{Sec_one_mean} and~\ref{Sec_two_mean}, respectively. In Section~\ref{Sec_wlln}, we give  an  asymptotic ratio-unbiased-estimator of $\mbox{tr}(\bd{R}^2)$, which is more applicable than the one by \cite{srivastava2008test} and \cite{srivastava2009test} (their result lacks a rigorous proof although it has been used in literature). The concluding remarks and discussions are presented in Section~\ref{Sec_concluding}, in which we particularly point out our findings on phase transitions between normal and mixing chi-squared distributions may also exist for some other testing procedures.   Finally, the proofs are given in Section~\ref{Sec_proofs}.

\section{One Sample Mean Test for Large $p$ and Small $n$}\lbl{Sec_one_mean}

Let $\{\X_1,\cdots,\X_n\}$ be a random sample from $N_p(\bmu, \bms)$. Consider the test that
\begin{align}\label{one}
H_0: \bmu=\bm 0 ~~\textrm{versus}~~H_1: \bmu\not=\bm 0.
\end{align}
Let $\bd{D}$ be the diagonal matrix of $\bd{\Sigma}.$ Then the $p\times p$ population correlation matrix is
\bea\lbl{38488}
\bd{R}=\bd{D}^{-1/2}\bd{\Sigma}\bd{D}^{-1/2}.
\eea
Similarly, let $\hat{\D}$ be the diagonal matrix of $\hat{\S}$ from~\eqref{wuoi0}. Then the sample correlation matrix $\hat{\R}$ is defined by
\bea\lbl{sample_corr_ma}
\hat{\R}=\hat{\D}^{-1/2}\hat{\S}\hat{\D}^{-1/2}.
\eea
\cite{srivastava2008test} and \cite{srivastava2009test} obtain a result on $T_{SD}$ from \eqref{jin_wuzu} as follows.

\begin{result}\lbl{AD2008}
Assume $n=O(p^{\zeta})$, $\frac{1}{2}<\zeta\leq 1$ and
\bea
&& 0<\lim_{p\to\infty}\frac{\mbox{tr}(\bd{R}^i)}{p}<\infty,~~~~i=1,2,3,4.\lbl{SD08a}
\eea
 If $\bmu=\bm 0$ then $T_{SD}\to N(0, 1)$ in distribution as $p\to\infty.$
\end{result}

A quick comment is that \eqref{SD08a} holds automatically for $i=1$ since all of the diagonal entries of $\bd{R}$ are equal to $1$, and hence $\mbox{tr}(\bd{R})=p$. In order to make approximation better, we now revise the statistic $T_{SD}$ slightly. Set
\bea\lbl{Statistics1}
T_{p,1}=\frac{n\bar{\X}^{T} \hat{\D}^{-1}\bar{\X}-pn(n-3)^{-1}}{\sqrt{ 2 \big|\tr(\hat{\R}^2)-p(p-1)(n-1)^{-1}\big|}}.
\eea
We are doing so because, by Lemma~\ref{sdu329}, the major contribution of the mean of $\mbox{tr}(\hat{\bd{R}}^2)$ is  $\mbox{tr}(\bd{R}^2)+p(p-1)(n-1)^{-1}$. Also, by Theorem~\ref{Theorem3}, $\tr(\hat{\R}^2)-p(p-1)(n-1)^{-1}>0$ as $n$ and $p$ are large enough. On the other hand, with probability one, $\tr(\hat{\R}^2)-p(p-1)(n-1)^{-1} \ne 0$ because $\tr(\hat{\R}^2)$ is a continuous function of Gaussian random variables. A discussion between $T_{p,1}$ and $T_{SD}$ will be elaborated shortly. Review the Frobenius norm, sometimes also called the Euclidean norm, $\|\bd{A}\|_F:=[\mbox{tr}(\bd{A}^T\bd{A})]^{1/2}$ for any matrix $\bd{A}$. For mathematical rigor, we assume the sample sizes $n$ depends on $p$. According to \cite{hu2016review}, the behavior of the statistic $T_{SD}$ is not known if $\bms$ has spikes. The following gives a complete characterization of the properties of $T_{SD}$ and $T_{p,1}$ in terms of the spikes of the correlation matrix $\bd{R}.$

\begin{theorem}\lbl{Theorem1} Suppose $\X_1,\cdots,\X_n$ is a random sample from $N_p(\bmu, \bms)$. Let $\bd{R}$ be the correlation matrix as in \eqref{38488} with eigenvalues  $\lambda_1\geq  \cdots \geq \lambda_p \geq 0$. Assume

(a)\ $\lim_{p\to\infty}\frac{\lambda_i}{\|\bd{R}\|_F}=\rho_i \geq 0$ for all $i\geq 1$;

(b)\ $\lim_{p\to\infty}\frac{p}{n\|\bd{R}\|_F}= 0$ and $\lim_{p\to\infty}\frac{p}{n^{a}}=0$ for some constant $a>0$.\\
If $\bmu=\bm 0$, then
$T_{SD}\to  b\xi_0+\frac{1}{\sqrt{2}}\sum_{i=1}^{\infty}\rho_i(\xi_i^2-1)$ in distribution, where $\xi_0, \xi_1, \xi_2, \cdots$ are i.i.d. $N(0, 1)$ and $b=(1-\sum_{i=1}^{\infty}\rho_i^2)^{1/2}$. The same conclusion also holds for  $T_{p,1}$.
\end{theorem}

Condition (a) considers the possibility that $\bd{R}$ may have spikes. If $\rho_1=0$, then  $\rho_i=0$ for every $i\geq 1$ due to the monotonicity of $\rho_i$, and we say there are no spikes in $\bd{R}$. In this situation we have $T_{p,1}\to N(0, 1)$ by Theorem~\ref{Theorem1}. Using  essentially the Fatou lemma, we have checked the given conditions actually imply $\sum_{i=1}^{\infty}\rho_i^2\leq 1$; see the proof of Lemma~\ref{vdu349j}. If $\sum_{i=1}^{\infty}\rho_i^2=1$, then $T_{p,1}$ converges to a mixing chi-square distribution. If  $\sum_{i=1}^{\infty}\rho_i^2\in (0,1)$, then the asymptotic distribution is a blend of a Gaussian distribution and a mixing chi-squared distribution.  Condition (b) characterizes the restriction between the sample size $n$ and population dimension $p$.

Observe that \eqref{SD08a} with  $i=3$ implies that $\lambda_1^3\leq \mbox{tr}(\bd{R}^3)=O(p)$. Thus, $\lambda_1=O(p^{1/3})$. Since all of the diagonal entries of $\bd{R}$ are identical to $1$ then $\mbox{tr}(\bd{R}^2)\geq p$. This shows that $\rho_1=0$ where $\rho_1$ is from Theorem~\ref{Theorem1}. Also, \eqref{SD08a} with  $i=2$ implies that $\|\bd{R}\|$ and $\sqrt{p}$ have the same order, thus $\frac{p}{n\|\bd{R}\|_F}\to 0$ as $p\to\infty$. This together with the condition ``$n=O(p^{\zeta})$, $\frac{1}{2}<\zeta\leq 1$" implies (b) from Theorem~\ref{Theorem1}.  So our theorem is more general than Result~\ref{AD2008} by \cite{srivastava2008test}.

Let $\X_1,\cdots,\X_n$ be a random sample from an $AR(1)$ model  with $\bd{R}=(\gamma^{|i-j|})$ and the absolute values of  $\gamma=\gamma_p$ staying away from $1$. By using the Gersgorin disc theorem [see, e.g., p. 344 from \cite{horn2012matrix}], the largest eigenvalue or $\bd{R}$ is of order $O(1)$. Hence  condition (a) of Theorem~\ref{Theorem1} holds with $\rho_1=0$. If condition (b) also holds, then both $T_{SD}$ and $T_{p,1}$ go to the standard normal distribution. The same conclusion is also valid for a banded correlation matrix $\bd{R}=(r_{ij})$ with $r_{ij}=0$ for $|j-i|\geq t$ where $t=t_p=o(\sqrt{p})$. In this case, the largest eigenvalue or $\bd{R}$ is of order $o(\sqrt{p})$. Similar results can be obtained for other patterned matrices including Toeplitz matrices, Hankel matrices and symmetric circulant matrices; see, e.g., \cite{brockwell2016introduction}.

Now let us look at some special features of Theorem~\ref{Theorem1}. First, we do not need the population matrix $\bms$ to be invertible as required in Hotelling's $T^2$ test \citep{hotelling1931generalization}. The largest discrepancy between $n$ and $p$ from \cite{srivastava2008test} is that  $p=o(n^2)$. Our range is that $p$ can be at any polynomial order of $n$ provided $p=o(n\|\bd{R}\|_F)$. If  the entries of the population vector are not very far from  independence, in the sense that  $\|\bd{R}\|_F$ is in the order of $\sqrt{p}$, then the restriction  $p=o(n\|\bd{R}\|_F)$ is reduced to $p=o(n^2)$. In the case that the entries of the population vector are very dependent such that $\|\bd{R}\|_F$ is in the order of $p$, then $p$ is allowed to take any polynomial order of $n$. To convince our readers for the dependent case and to make further discussions, we next study two extreme cases: independence and most dependence. The derivation of the results below does not use any techniques and steps from the proof of Theorem~\ref{Theorem1}.

\begin{prop}\lbl{Lemma_Remark_1}

Let $\X_1,\cdots,\X_n$ be a random sample from $N_p(\bmu, \bms)$ with $\bmu=\bd{0}$. Assume $n=n_p\to\infty$ as $p\to\infty$. The following hold.

(i)  Let all of the $p^2$ entries of $\bms$ be identical.  Then both  $T_{SD}$ and $T_{p,1}$  converge to $\frac{1}{\sqrt{2}}[\chi^2(1)-1]$ in distribution as $p\to\infty$ regardless of the relative speeds of $n$ and $p$.

(ii) Let  $\bms$ be a diagonal matrix whose diagonal entries are all positive, equivalently, $\bd{R}=\bd{I}_p$. Then
\bea\lbl{Theorem_extreme_cases}
T_{SD} \to
\begin{cases}
\eta, & \text{if $p/n^2\to 0$};\\
\eta+\sqrt{h/2}, & \text{if $p/n^2\to h$};\\
\infty, & \text{if $p/n^2\to \infty$}
\end{cases}
\eea
in distribution, where $\eta\sim N(0, 1)$. However,  $T_{p,1} \to N(0, 1)$
as $p\to\infty$ regardless of the speeds of $n$ and $p$ going to infinity.
\end{prop}

\begin{figure}[h!]
\mbox{
    \includegraphics[width=2.2in]{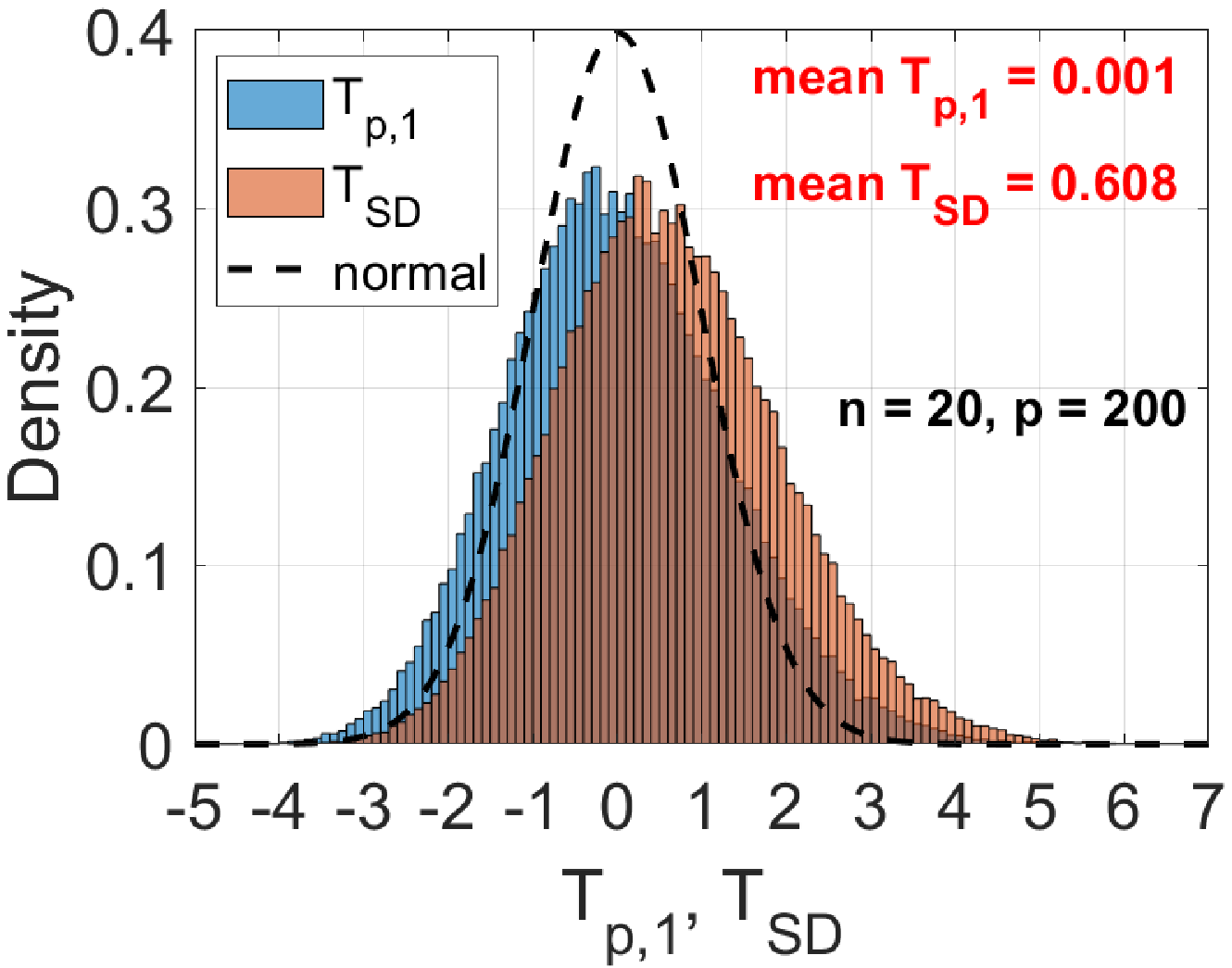}
    \includegraphics[width=2.2in]{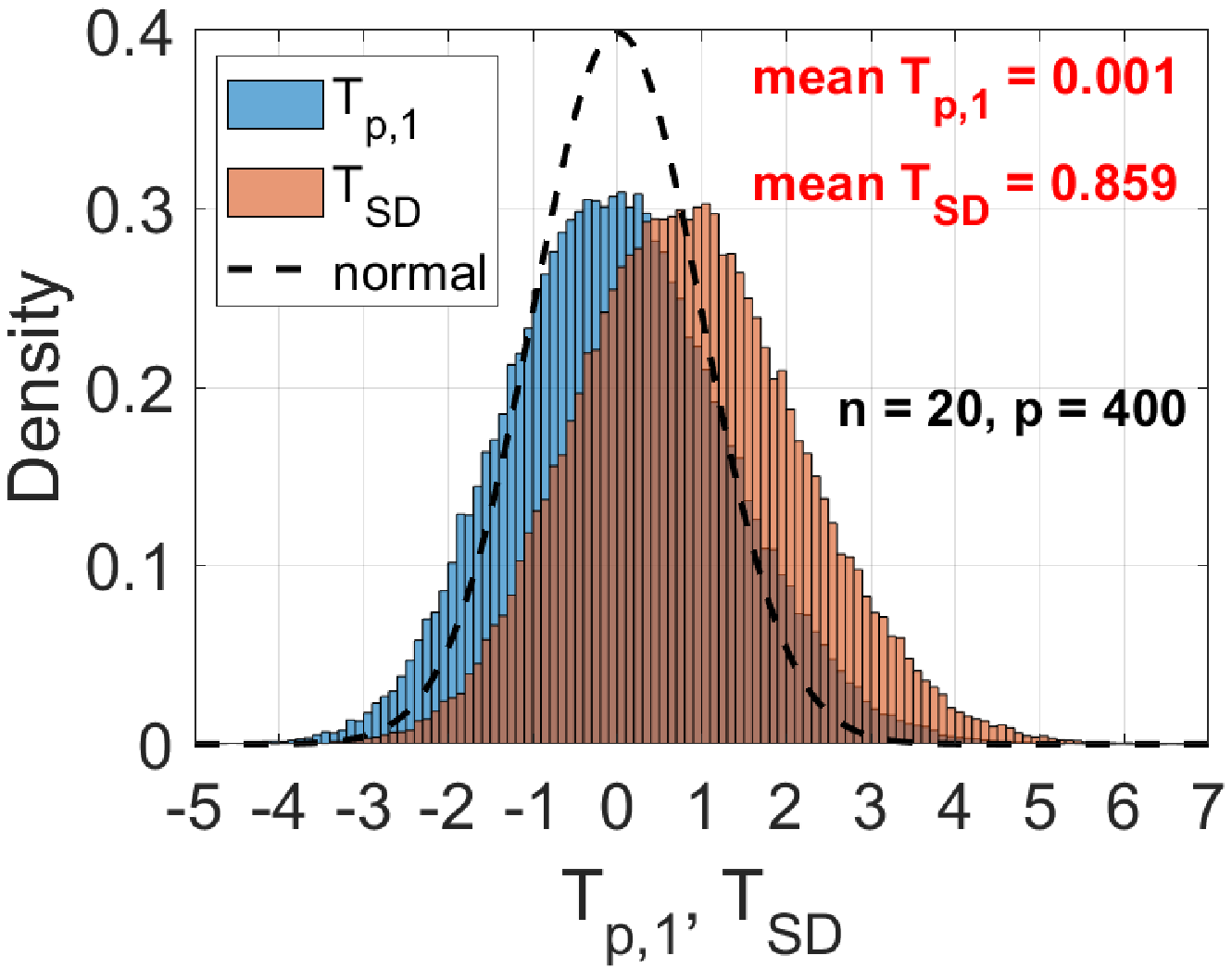}
    \includegraphics[width=2.2in]{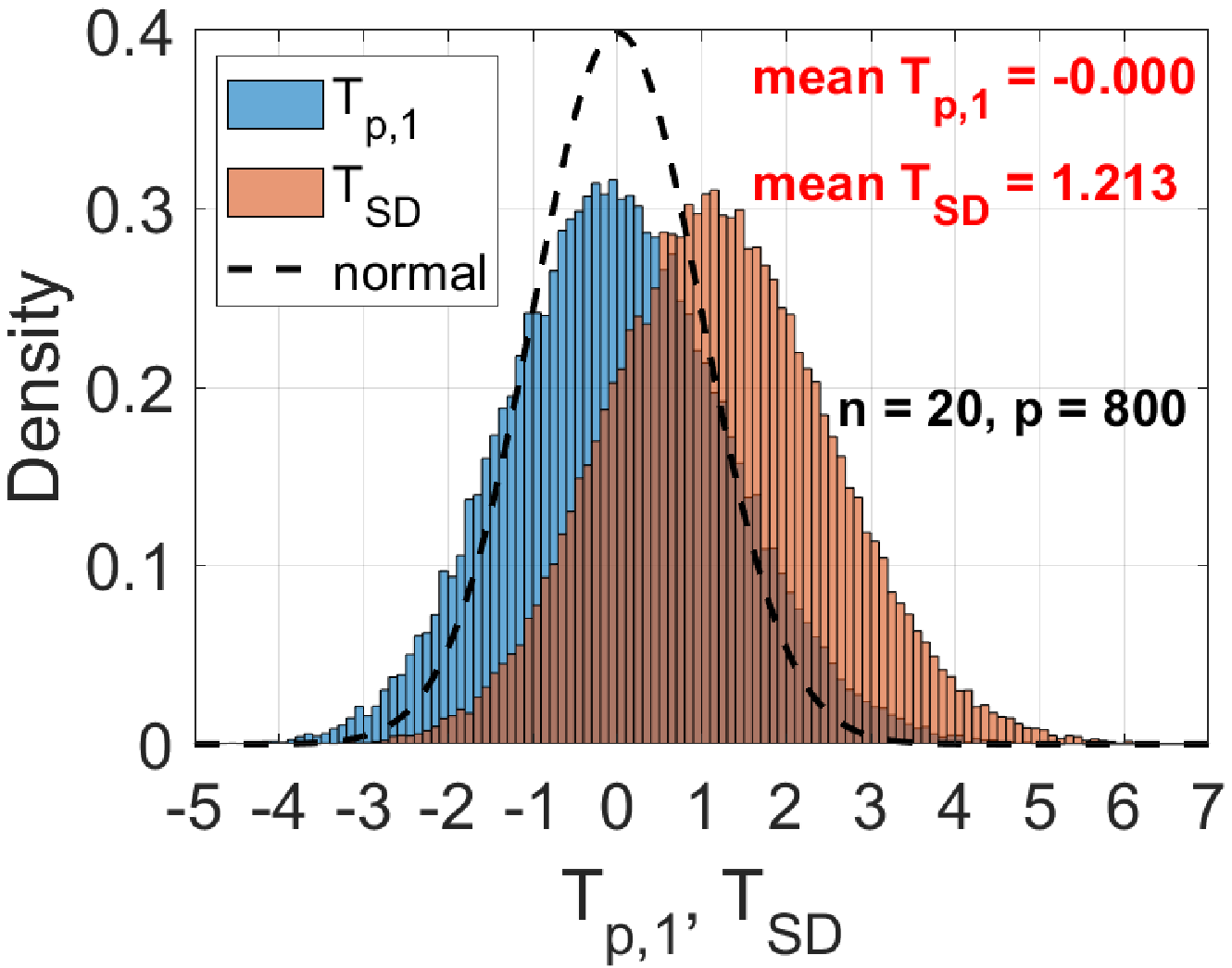}
}

\mbox{
    \includegraphics[width=2.2in]{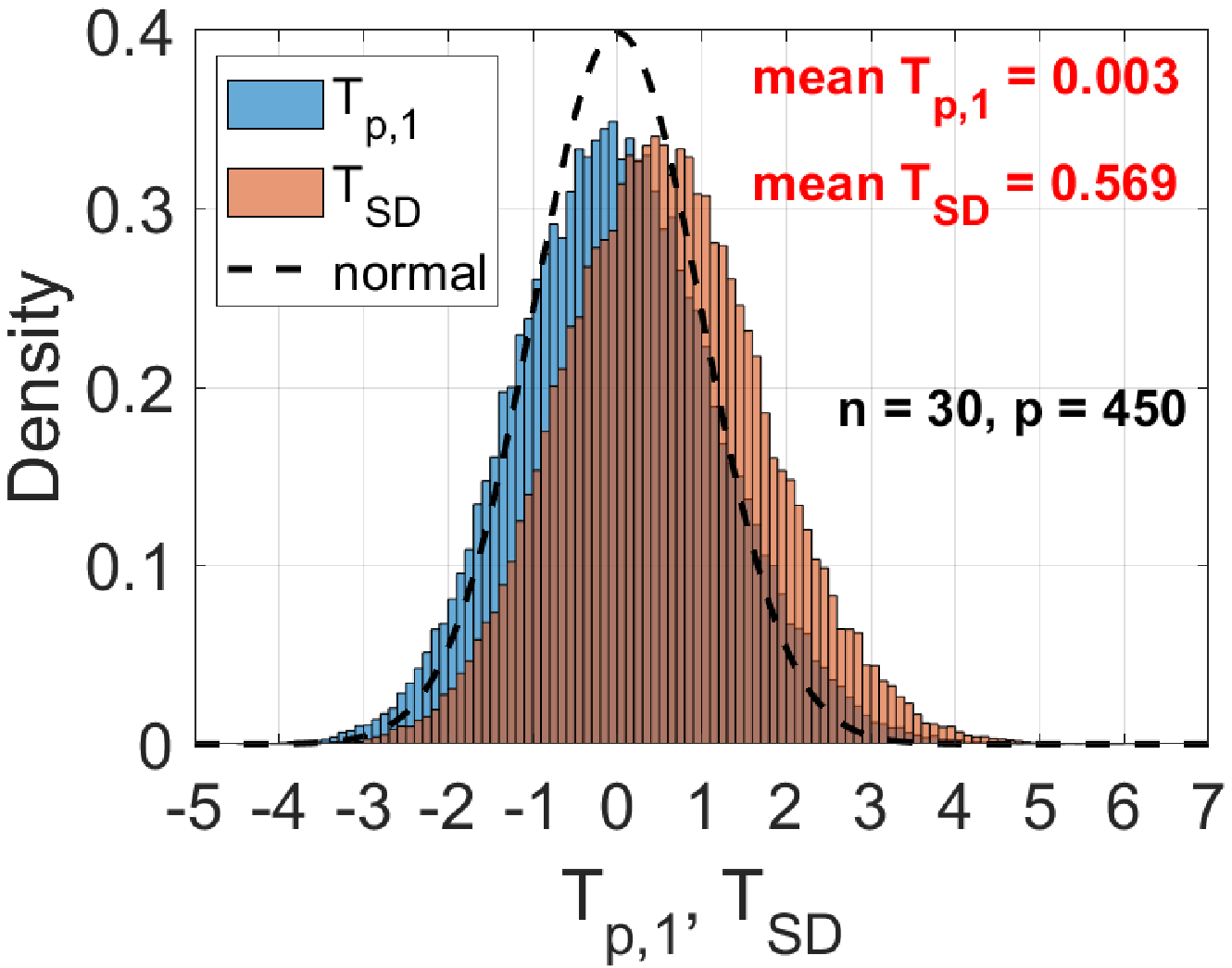}
    \includegraphics[width=2.2in]{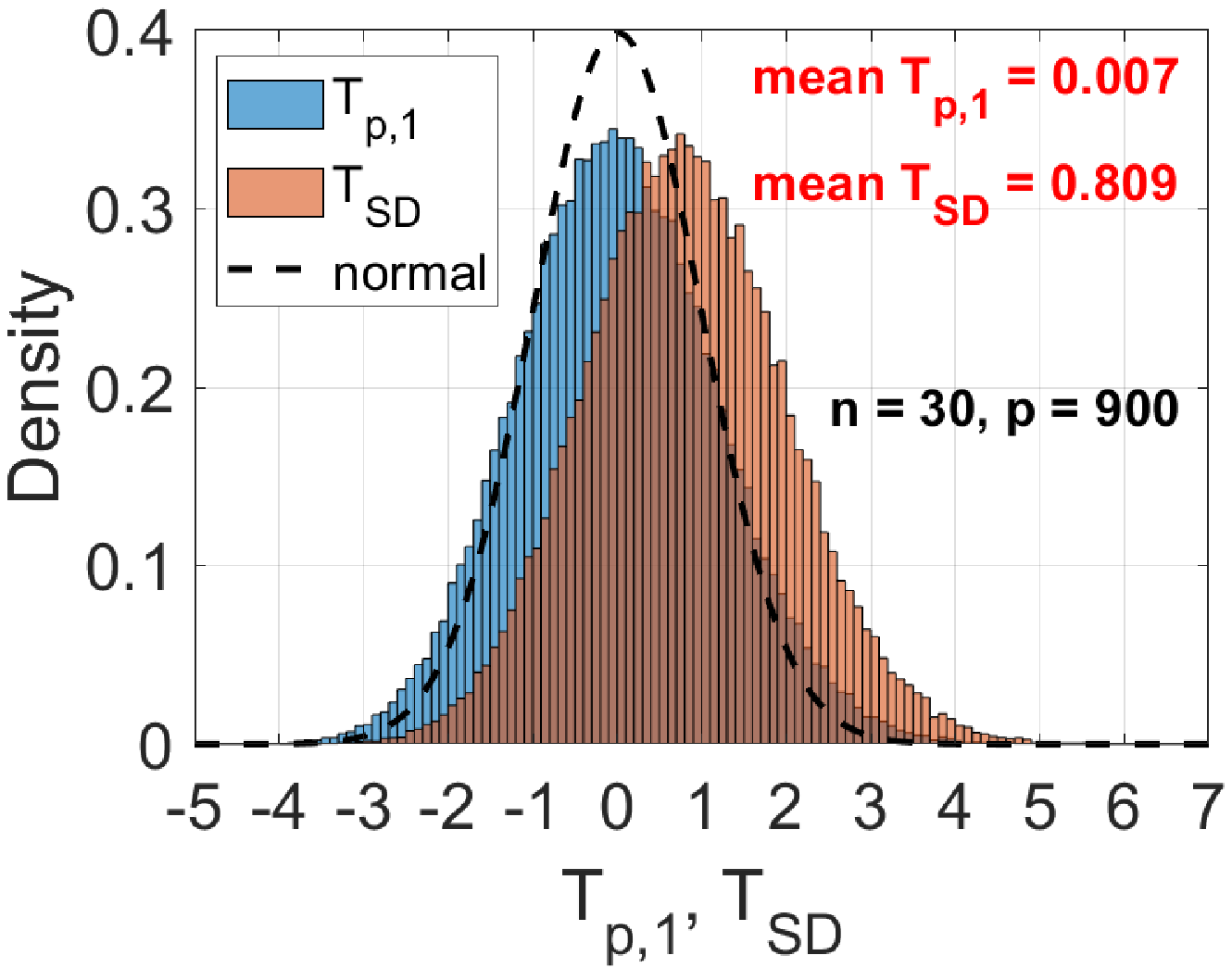}
    \includegraphics[width=2.2in]{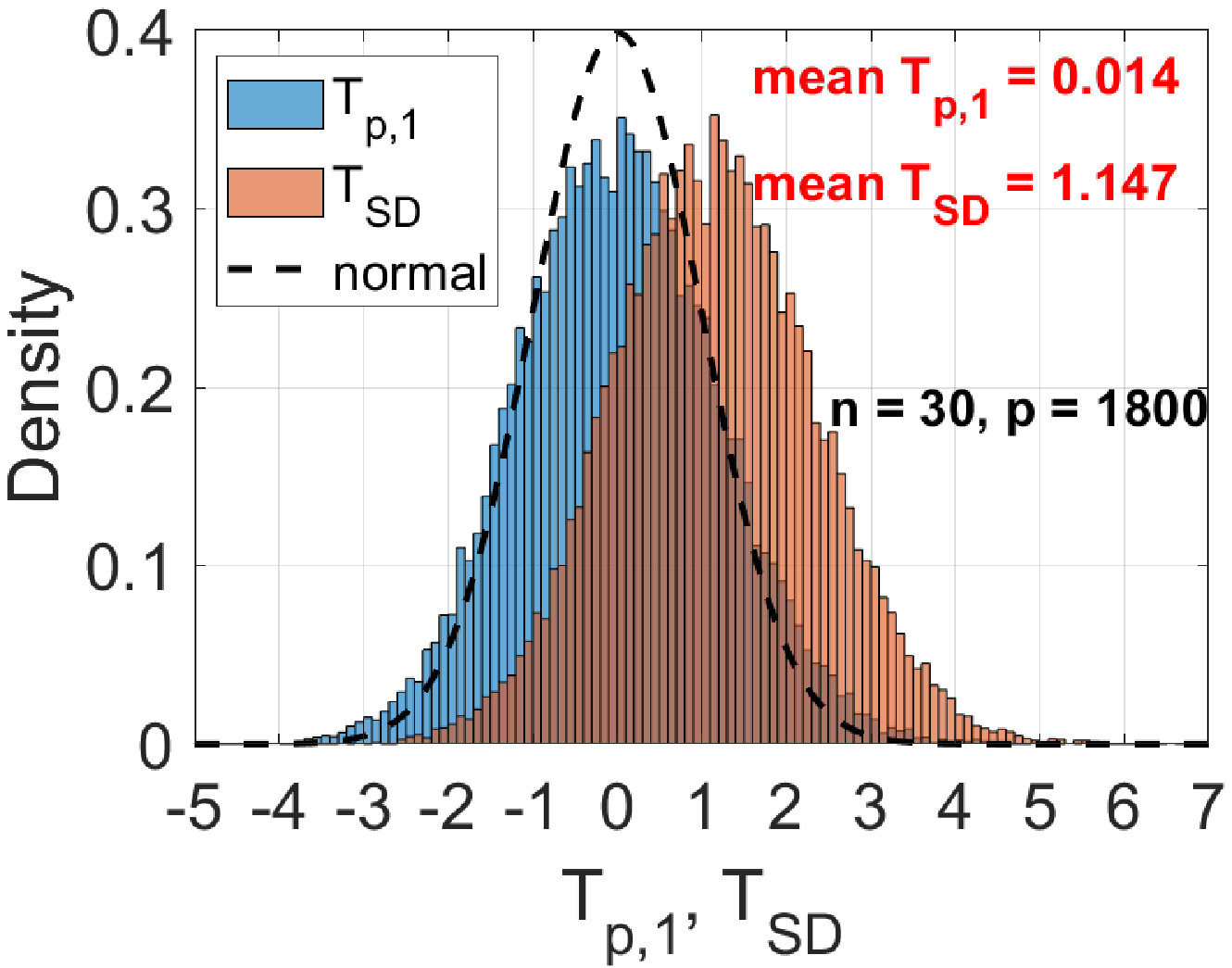}
}
\mbox{
    \includegraphics[width=2.2in]{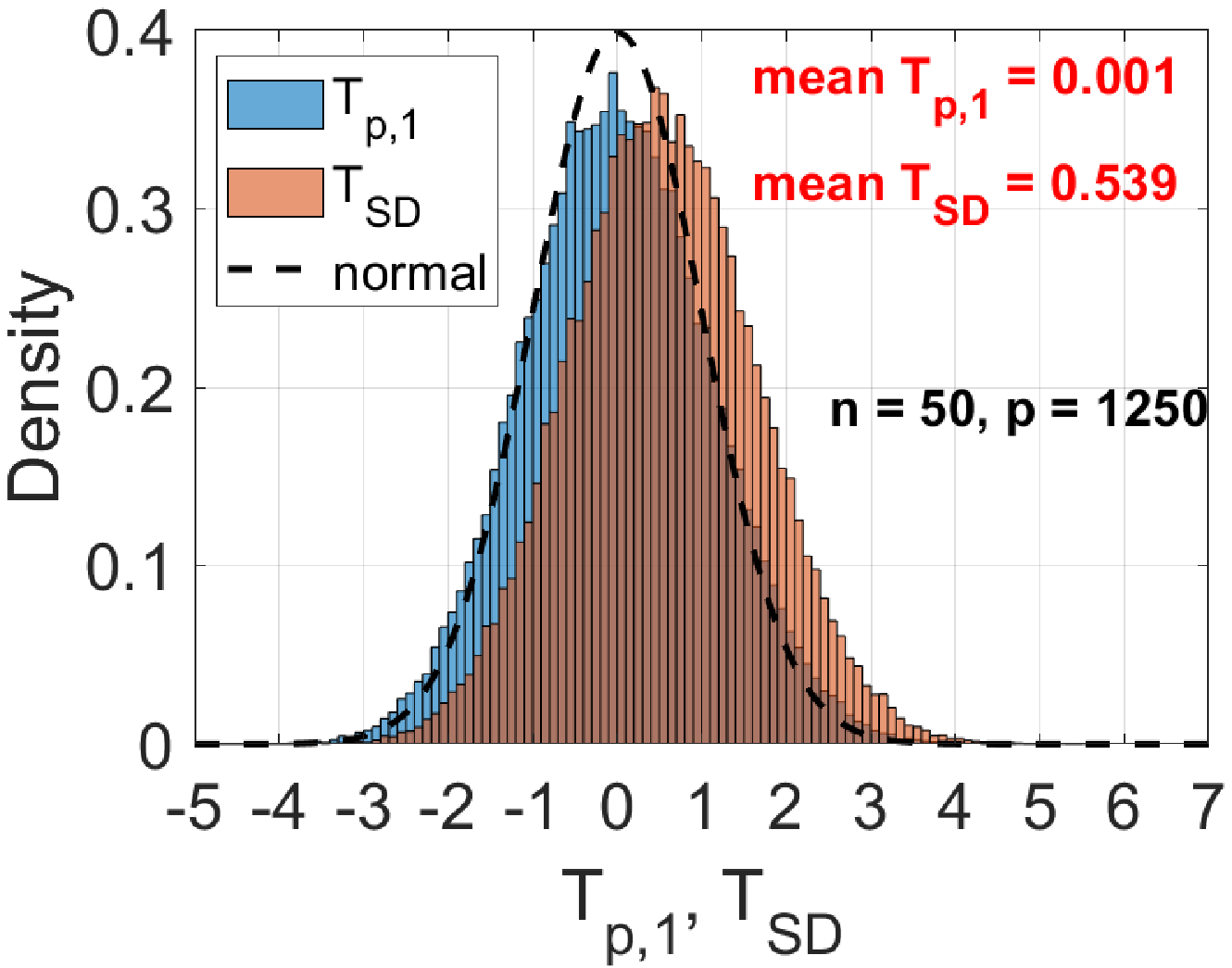}
    \includegraphics[width=2.2in]{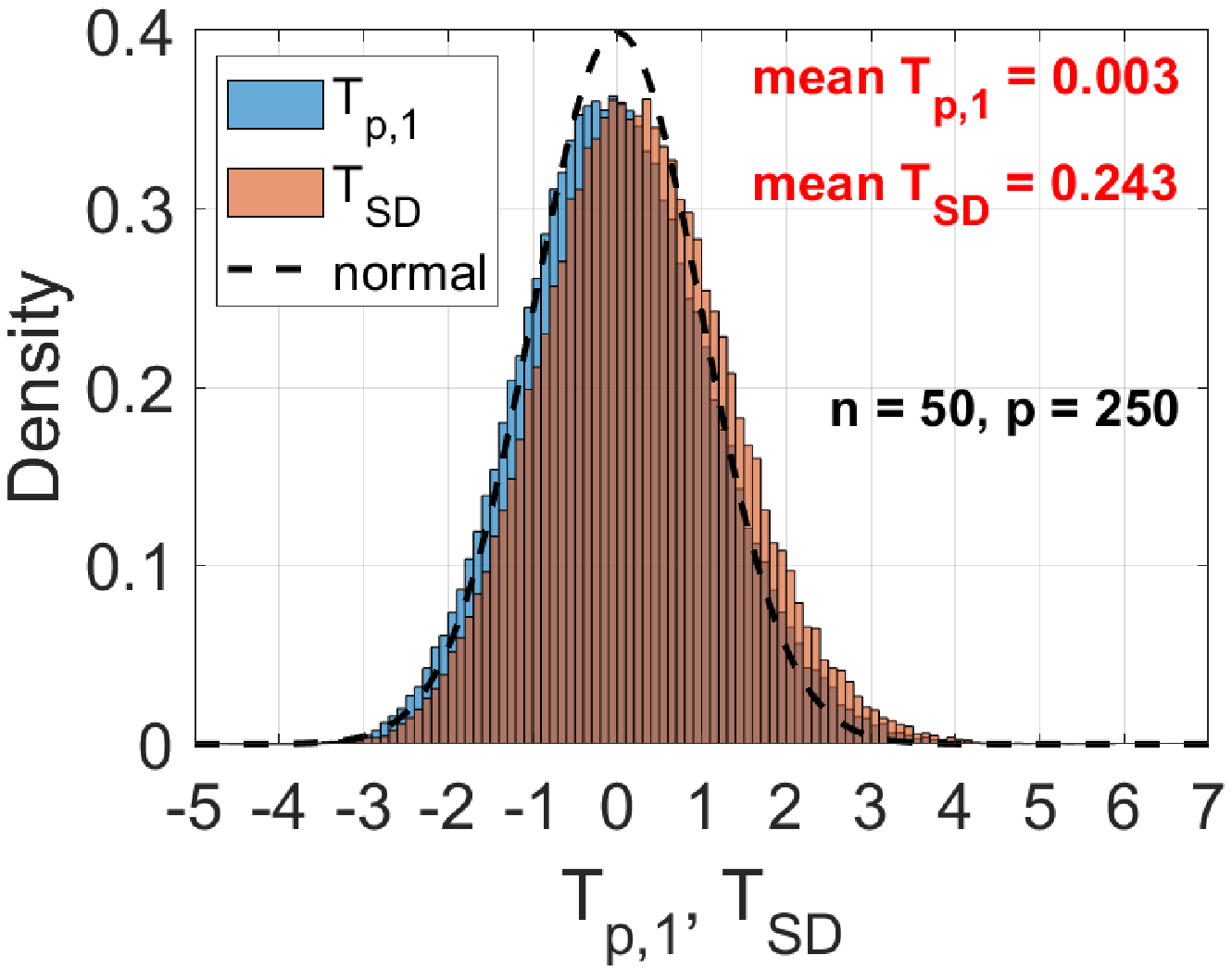}
    \includegraphics[width=2.2in]{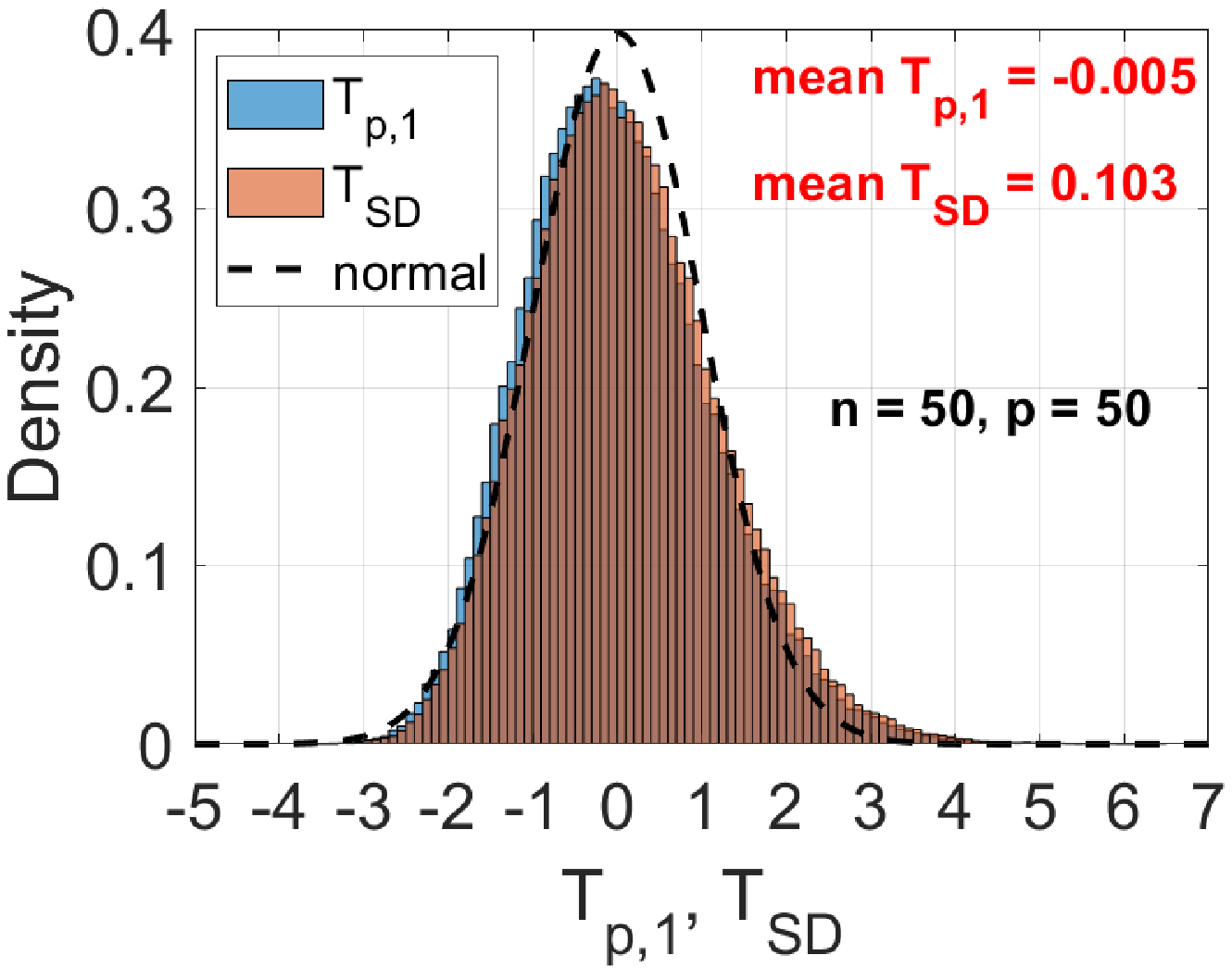}
}
    \caption{\small These pictures are based on Proposition~\ref{Lemma_Remark_1}, from $10^5$ simulations. It shows the centers of $T_{p,1}$ and $N(0,1)$ are always close, and the existence of shifts  between  density curves of  $T_{SD}$ and  $N(0, 1)$. The shifts are smaller as $p/n^2$ become smaller.}
    \label{fig:Prop1}
\end{figure}

In case (i) we see Theorem~\ref{Theorem1} holds without the assumption that $p=o(n^{a})$ for some constant $a>0$.
The conclusion for case (ii) says that the condition $\lim_{p\to\infty}\frac{p}{n\|\bd{R}\|_F}= 0$ in Theorem~\ref{Theorem1} is sharp, and $T_{p,1}$ is better than $T_{SD}$. In particular, we do not need the assumption $\lim_{p\to\infty}\frac{p}{n\|\bd{R}\|_F}= 0$ to assure $T_{p,1} \to N(0, 1)$. This is simply a coincidence because of the special structure of $\bms$ or $\bd{R}$. The condition  $\lim_{p\to\infty}\frac{p}{n\|\bd{R}\|_F}= 0$ is essentially required at handling the denominator of $T_{p,1}$ for arbitrary $\bd{R}$; see \eqref{Statistics1} and Theorem~\ref{Theorem3}. Figure~\ref{fig:Prop1} clearly shows the existence of shifts between the density curves of $T_{SD}$ and $N(0,1)$. However, with the ratio $p/n^2$ becoming smaller, the shift diminishes gradually.

\newpage

Practically, given the population correlation $\bd{R}$, we need to determine if it has spikes and how many spikes. Interested readers are referred to \cite{fan2020estimating} and \cite{morales2021asymptotics} for procedures to do so.  We next present some examples. Their verification are presented in Section~\ref{duhwoi9}.

\begin{example}\lbl{Ex1} Given $r\in [0,1]$, define $\bd{A}_m: = (1-r)\bd{I}_m + r \bd{J}_m$
for any $m\geq 1$, where  $\bd{I}_p$ and $\bd{J}_p$ denote the $p \times p$ identity matrix and the $p \times p$  matrix of ones.
Set $\bd{R}=\bd{A}_p$. Then $\lambda_1=1+(p-1)r$ and  $\lambda_2=\cdots = \lambda_p=1-r$. Assume $r=r_p$ and $\lim_{p\to\infty}\sqrt{p}\cdot r=c\geq 0$. Then
\bea\lbl{aiy327}
T_p\to
\begin{cases}
\xi_0, & \text{if $c=0$};\\
\frac{1}{\sqrt{c^2+1}}\xi_0+\frac{c}{\sqrt{2(c^2+1)}}(\xi_1^2-1),      & \text{if $c \in (0, \infty)$};\\
\frac{1}{\sqrt{2}}\cdot (\xi_1^2-1), & \text{if $c=\infty$}
\end{cases}
\eea
under condition $p=o(n^2)$ for the case $c\in [0,\infty)$ and under the conditions $nr\to\infty$ and  $p=o(n^a)$ for some constant $a>0$ for the case  $c=\infty$. Here $\xi_0$ and $\xi_1$ are i.i.d. $N(0, 1)$ and $T_p$ stands for $T_{SD}$ or $T_{p,1}$. The phase transition occurs between the normal and mixing chi-squared distributions as the entries of the correlation matrix change their values.
\end{example}

\begin{figure}
\mbox{    \includegraphics[width=2.2in]{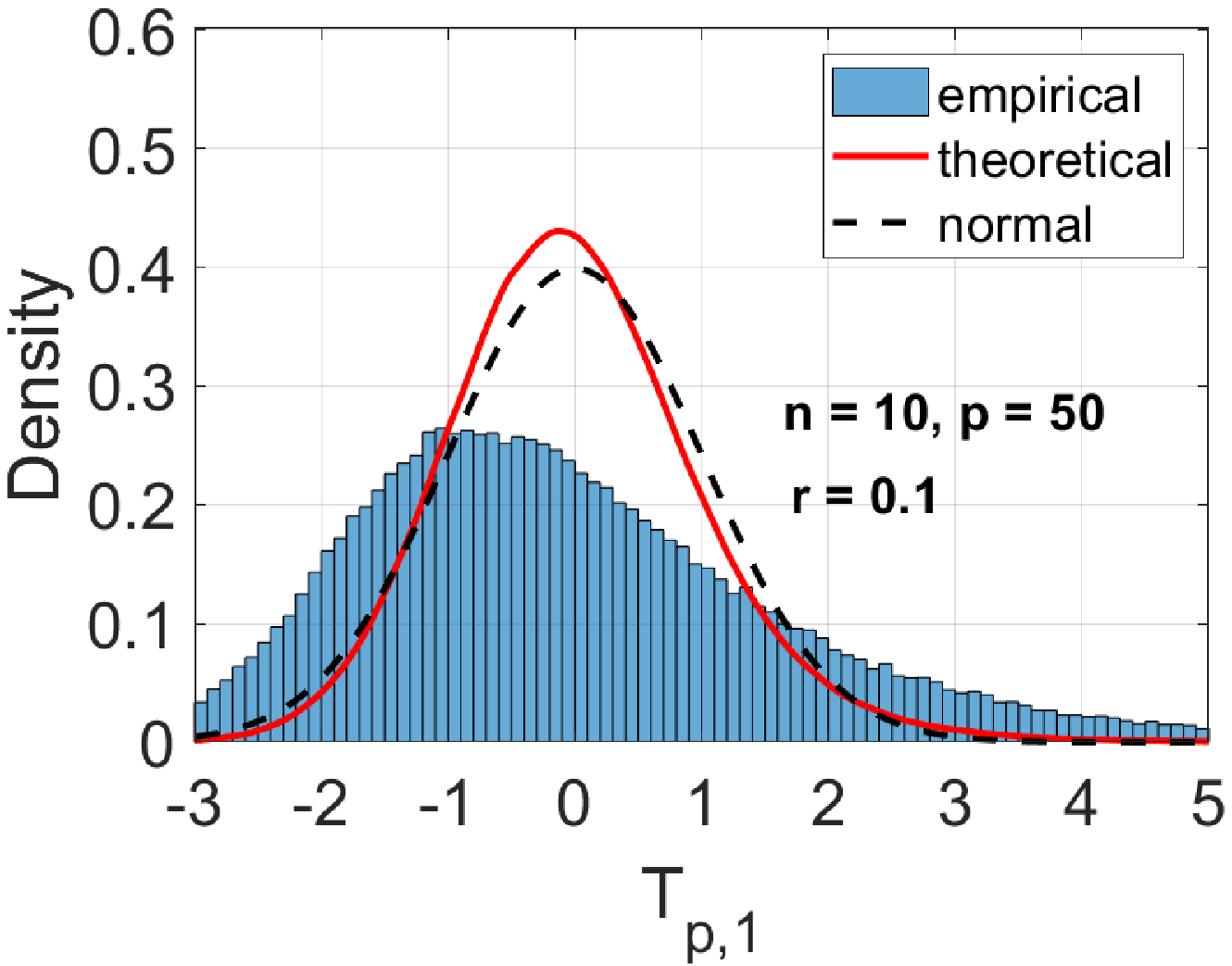}
    \includegraphics[width=2.2in]{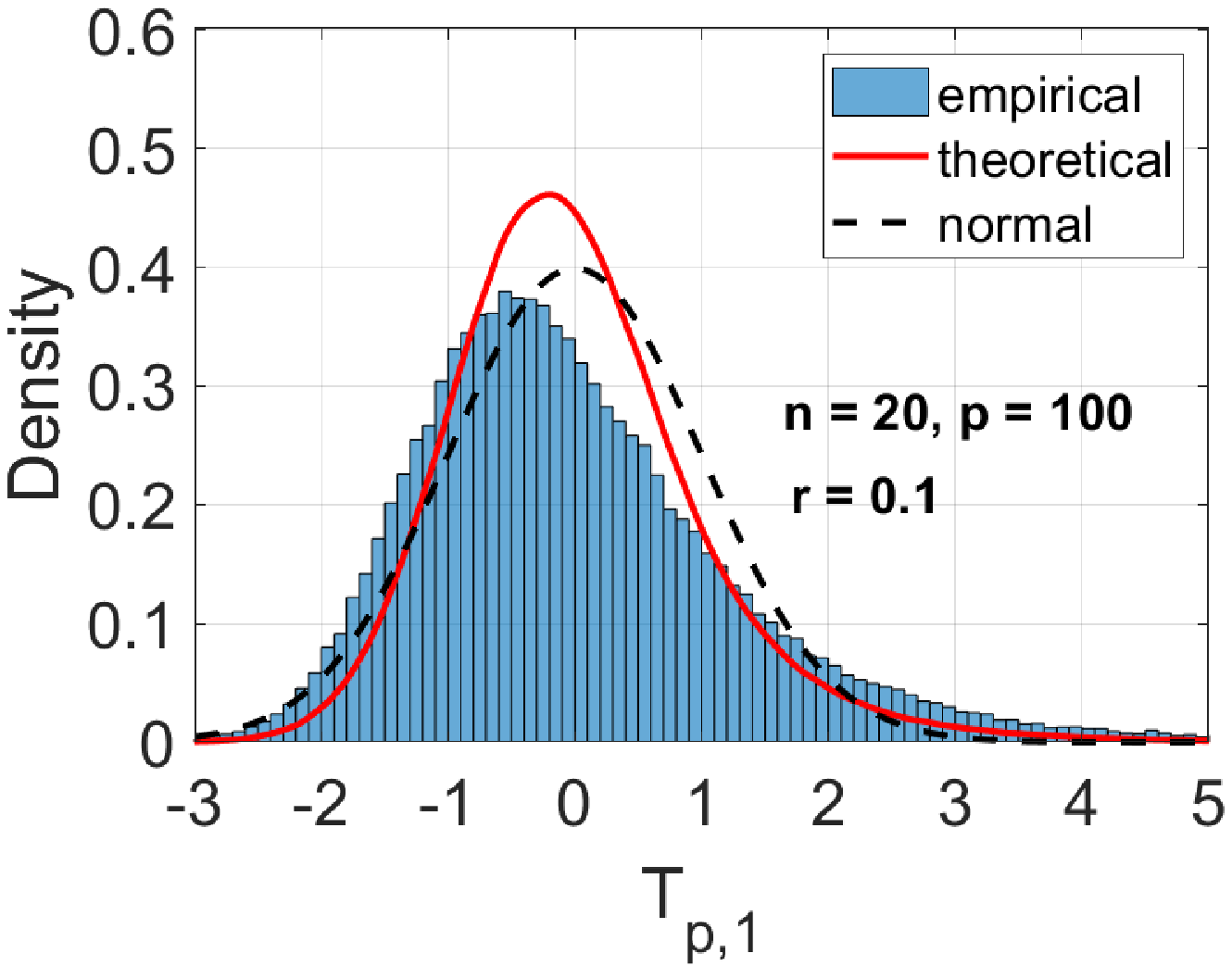}
    \includegraphics[width=2.2in]{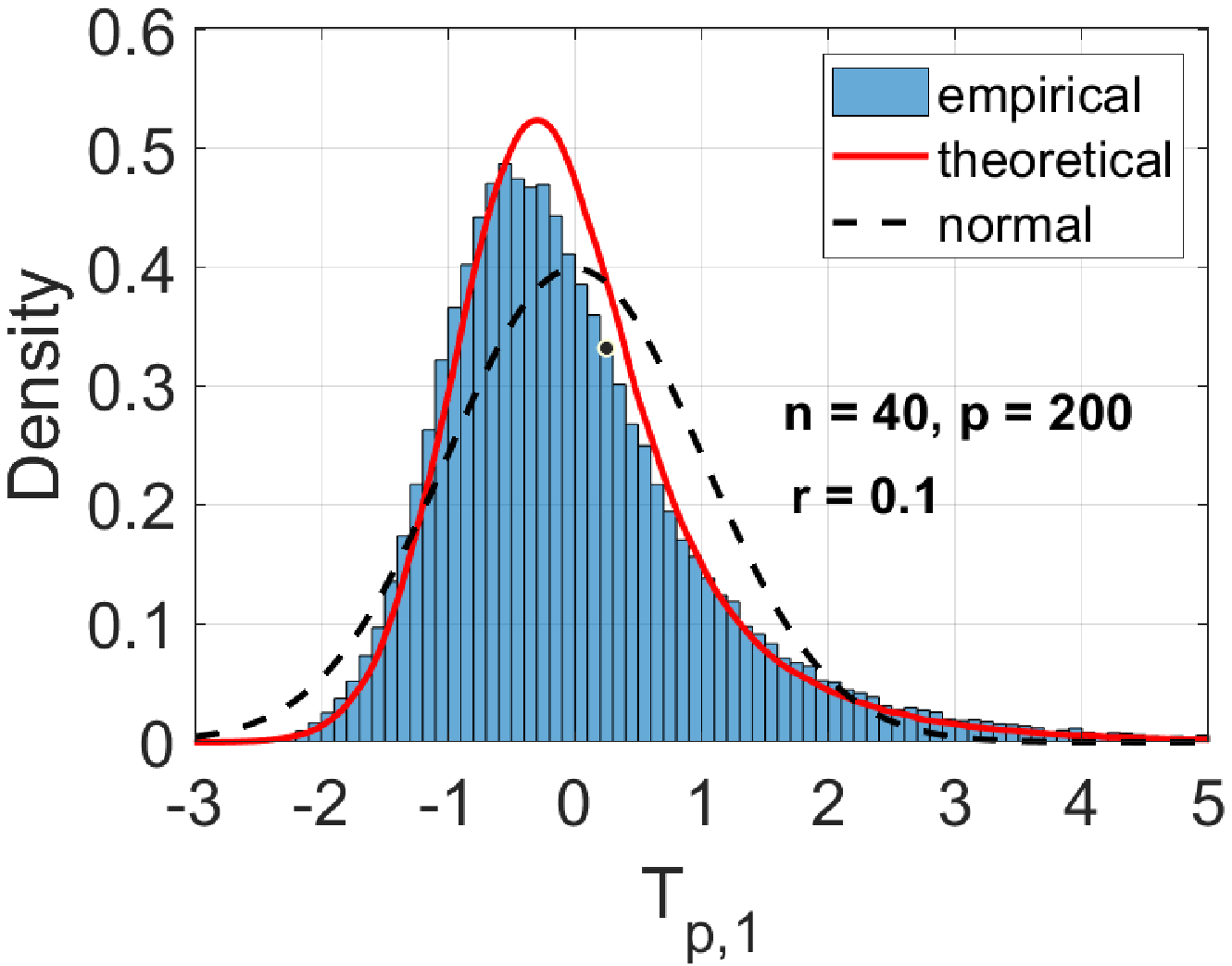}
}

\mbox{    \includegraphics[width=2.2in]{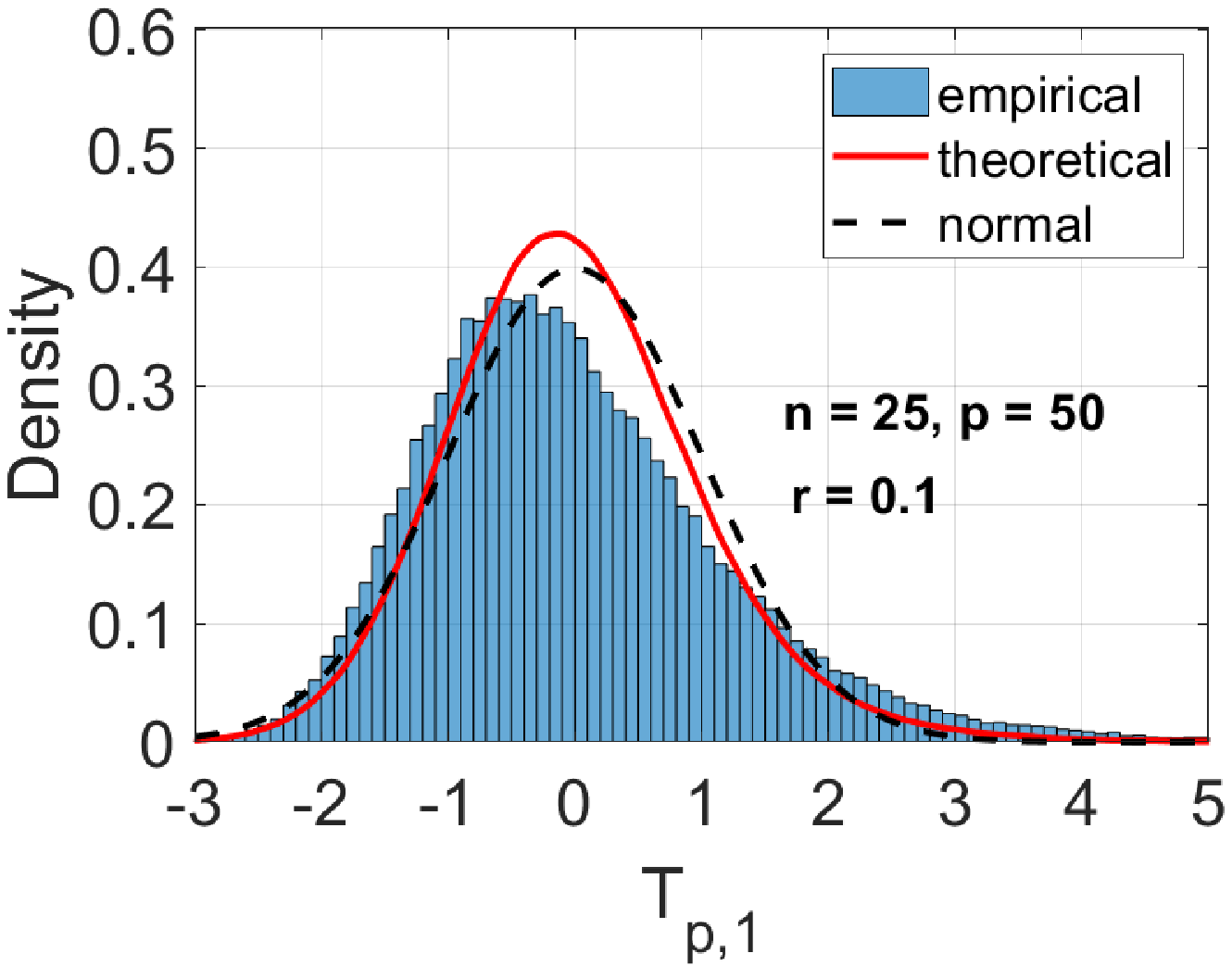}
    \includegraphics[width=2.2in]{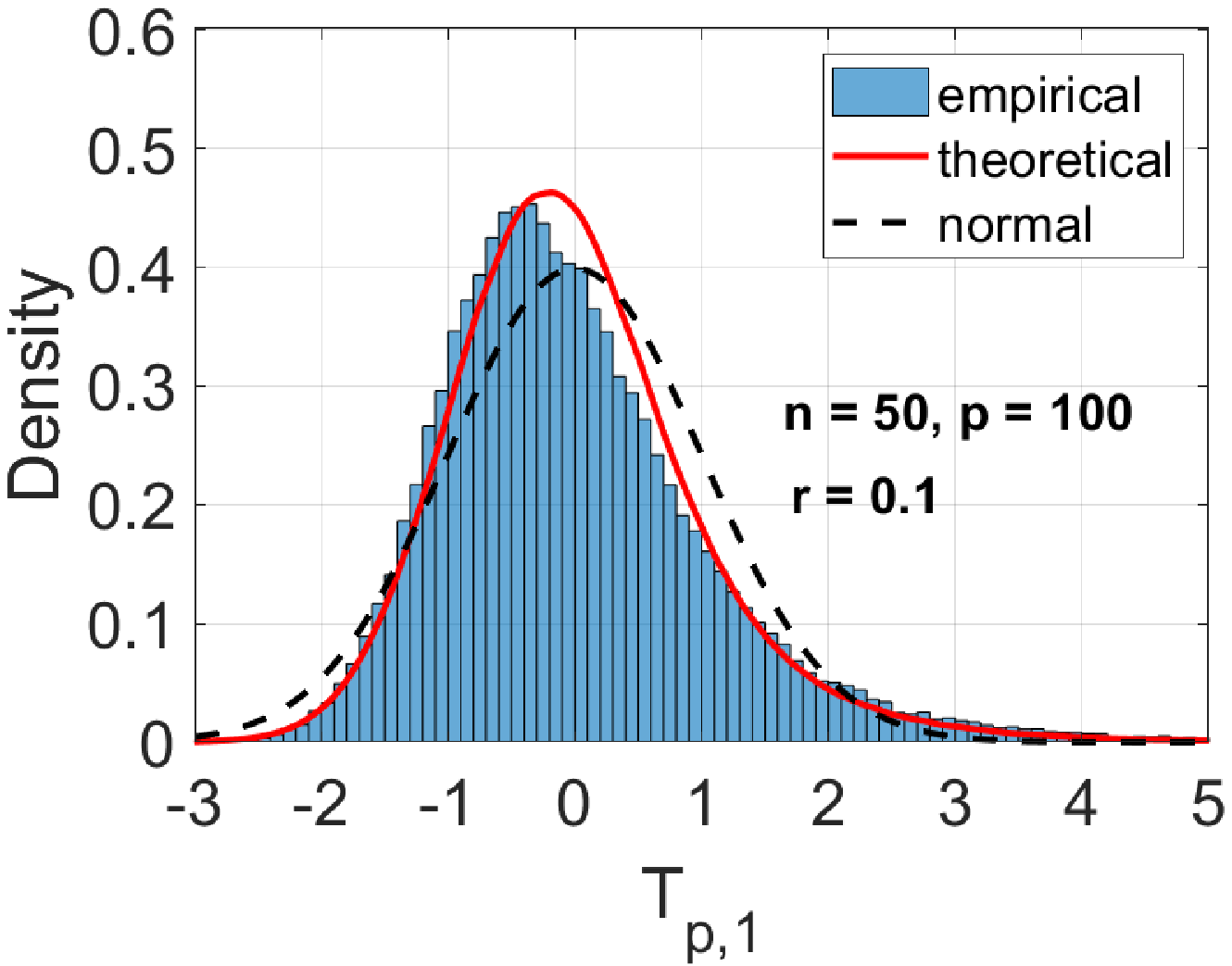}
    \includegraphics[width=2.2in]{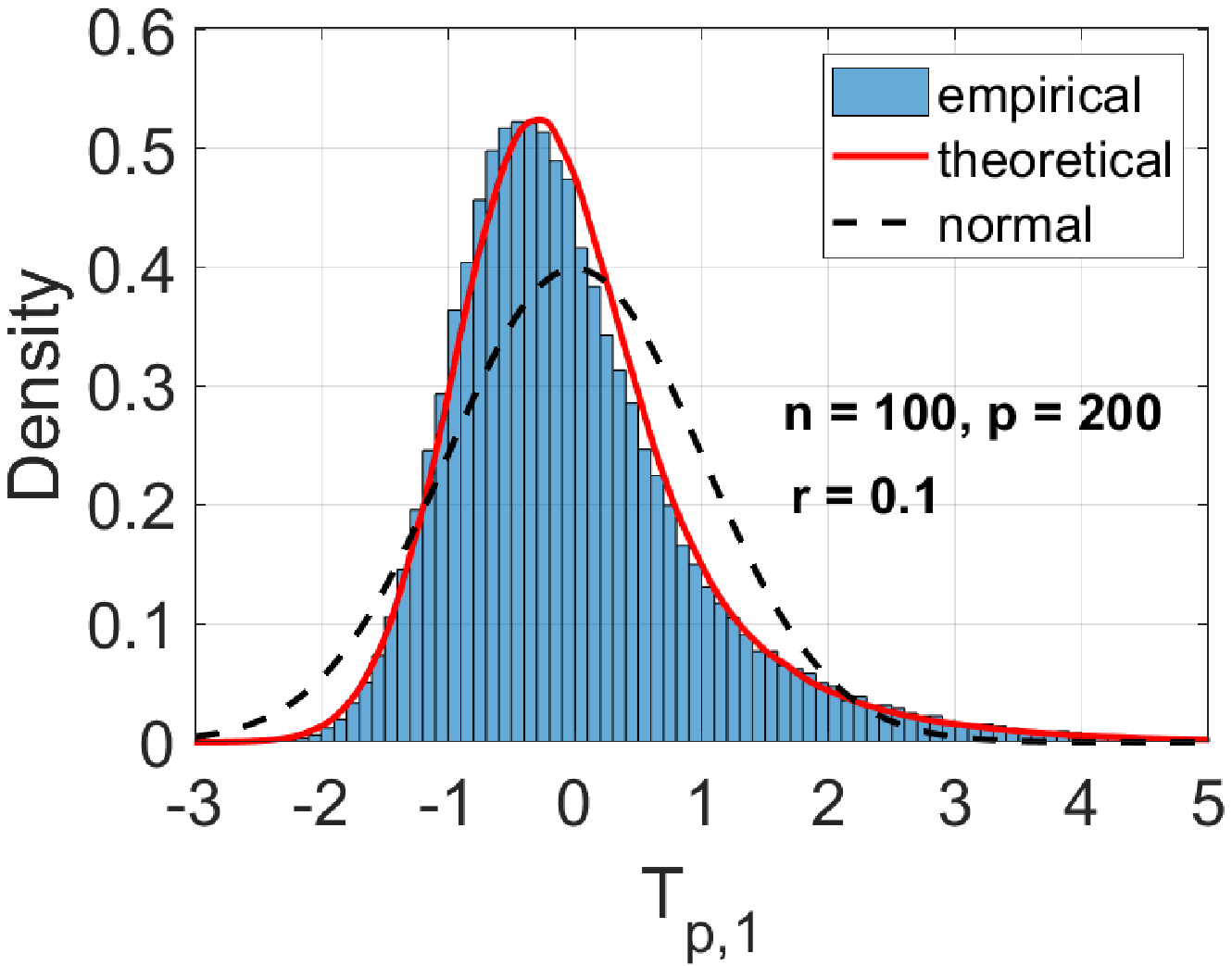}
}

\mbox{    \includegraphics[width=2.2in]{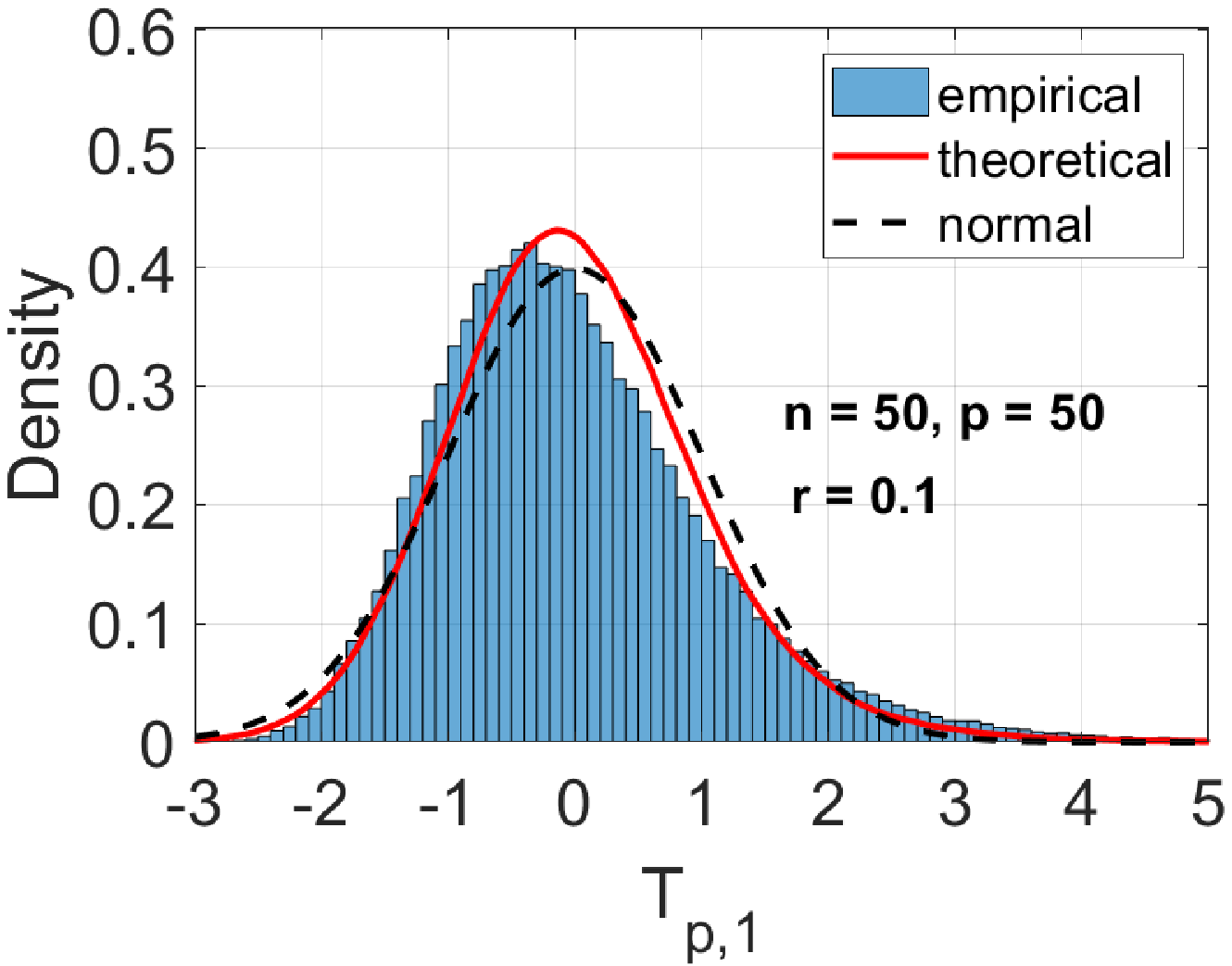}
    \includegraphics[width=2.2in]{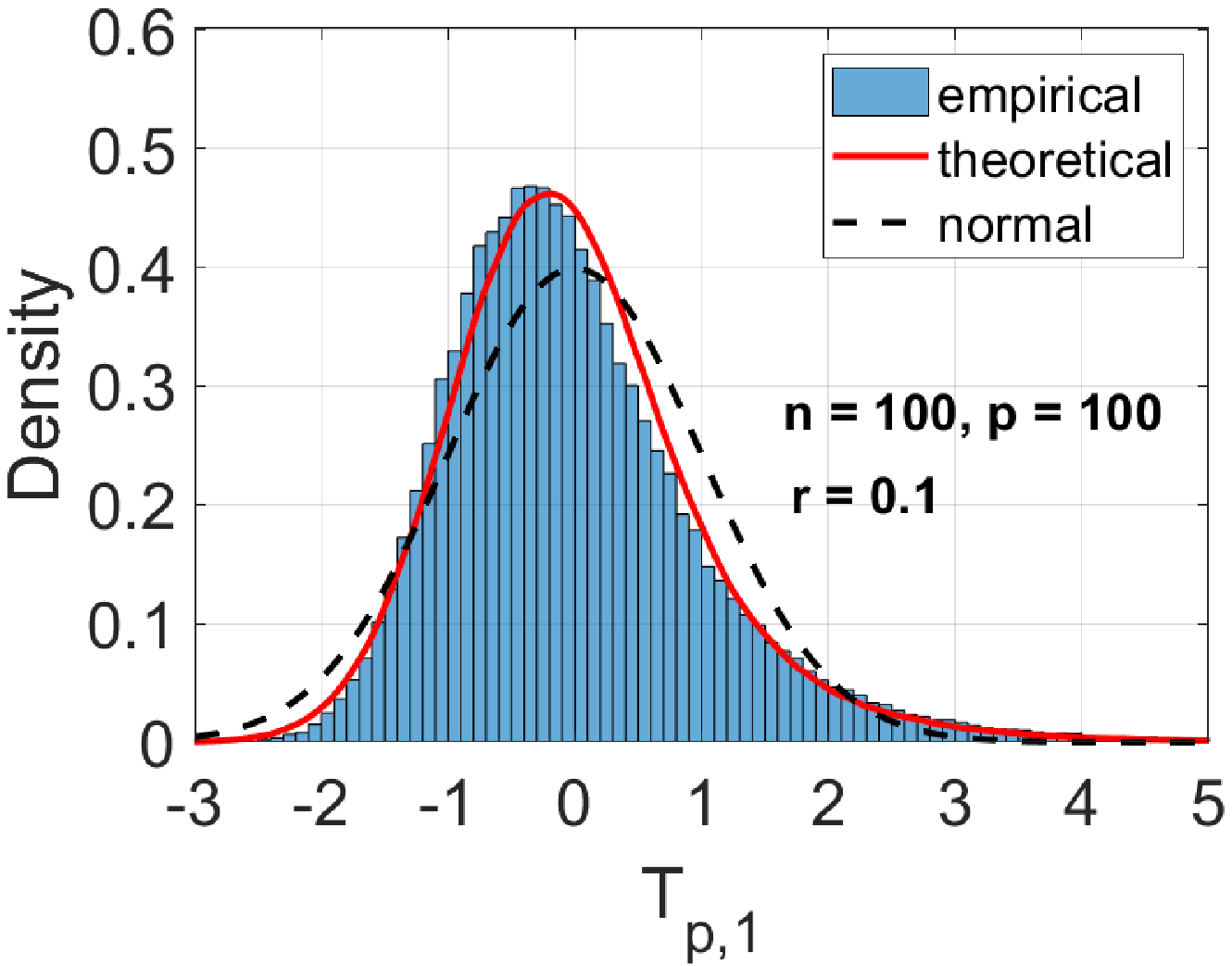}
    \includegraphics[width=2.2in]{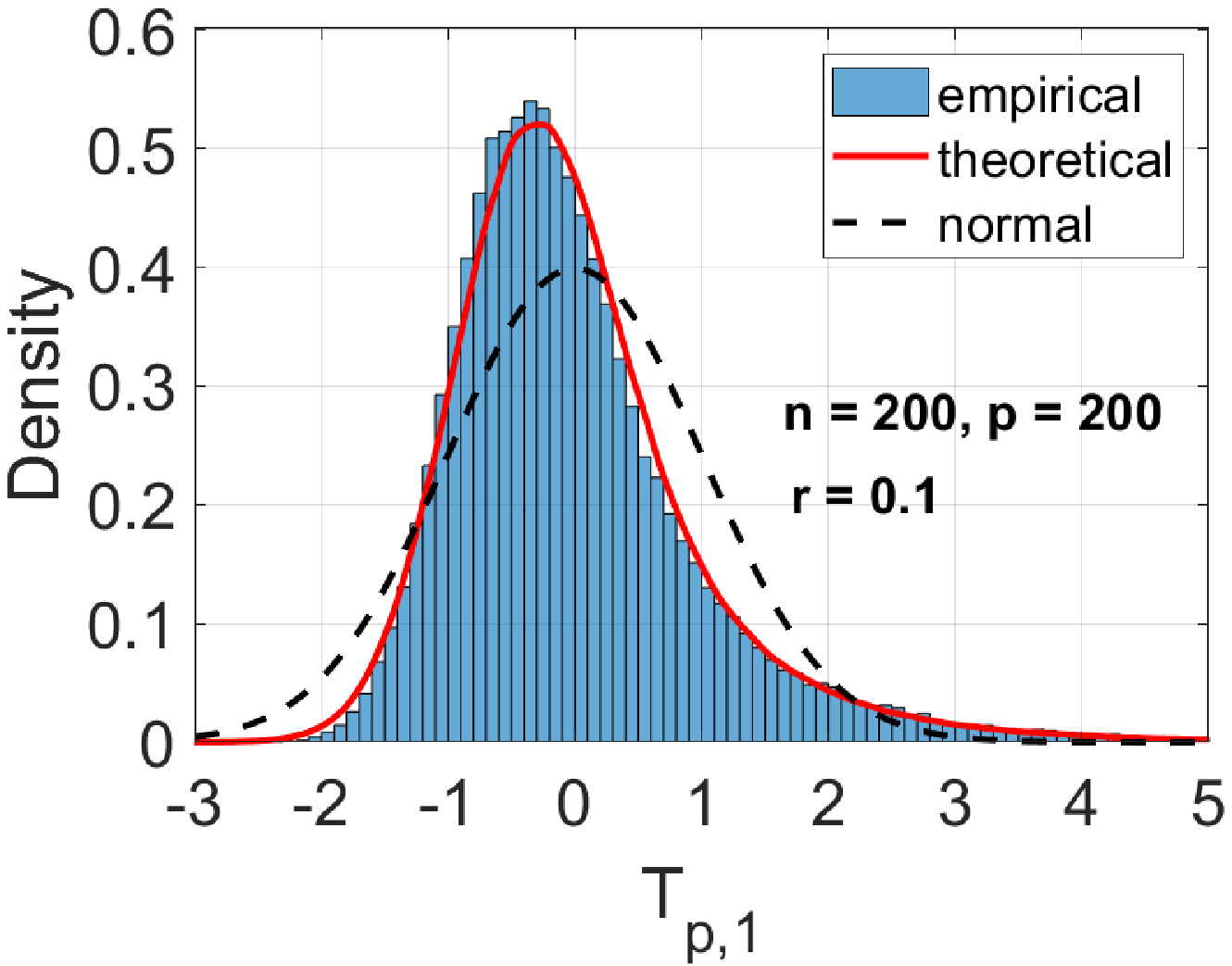}
}

\mbox{    \includegraphics[width=2.2in]{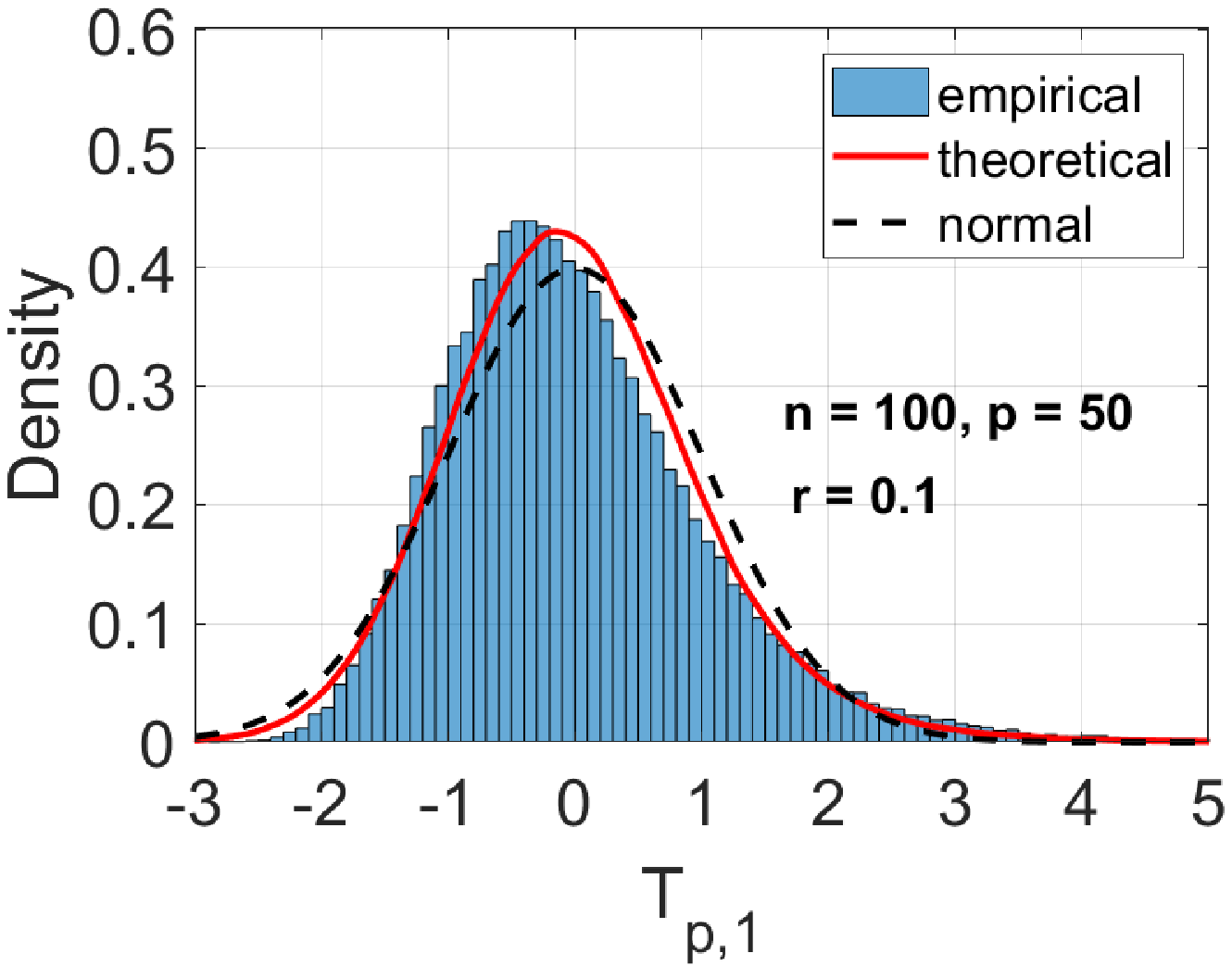}
    \includegraphics[width=2.2in]{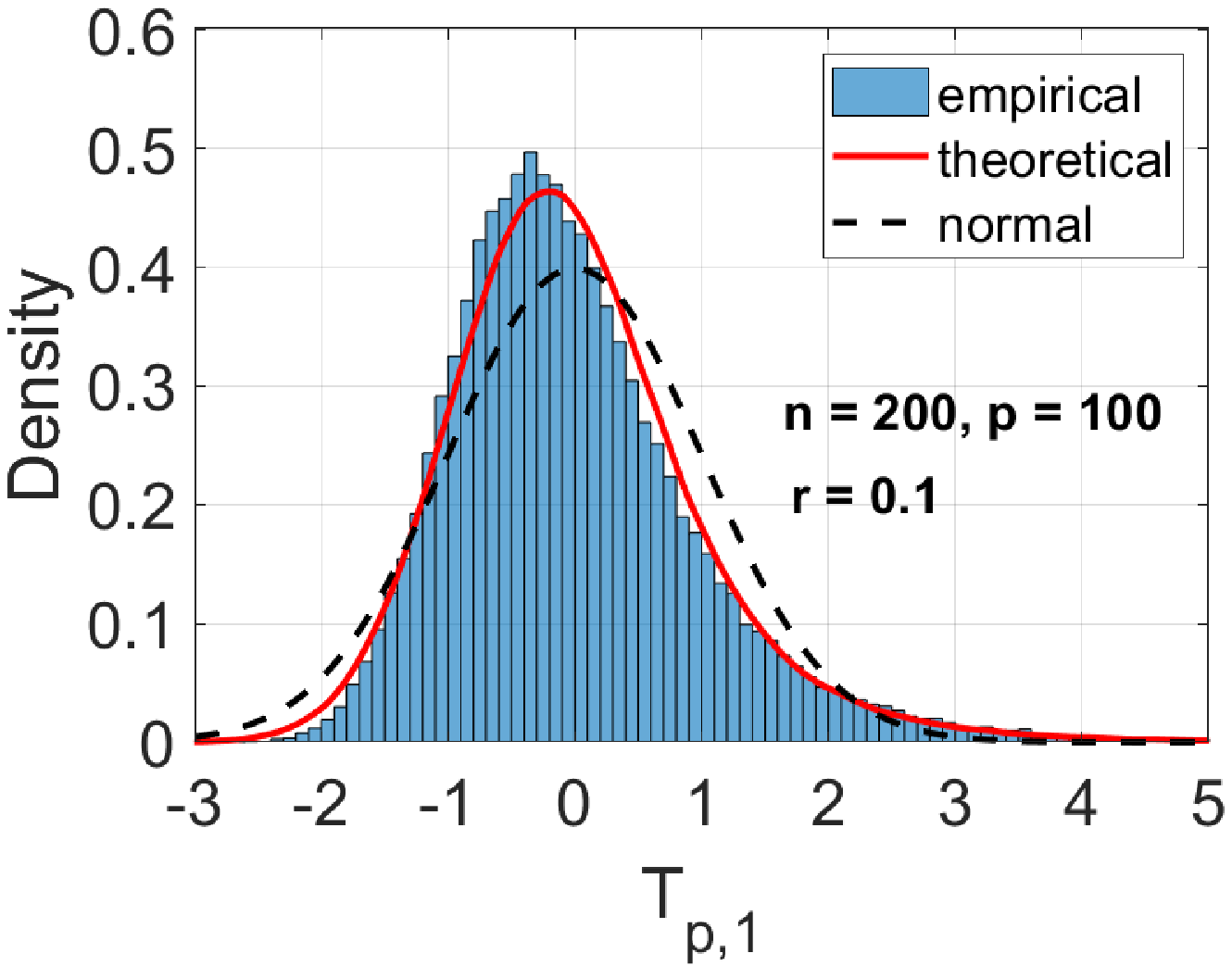}
    \includegraphics[width=2.2in]{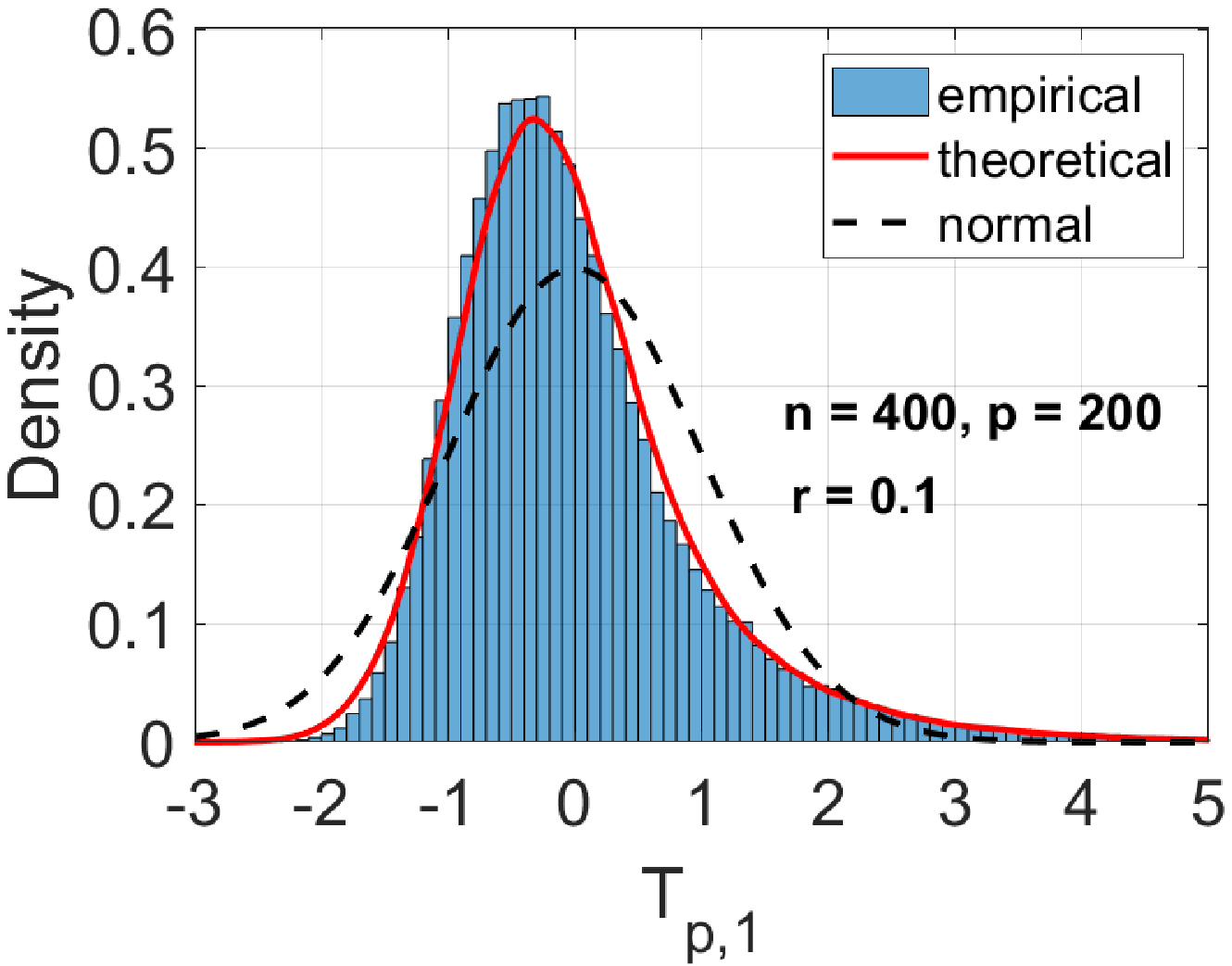}
}

    \caption{Example~\ref{Ex1}. For small $r=0.1$, the empirical curves in blue are close to our theoretical curves in red. As sample size $n$ becomes larger, our approximation becomes better; the normal approximation (black curve) is no longer valid.}
    \label{fig:Ex1_r01}
\end{figure}

\begin{figure}
\mbox{    \includegraphics[width=2.2in]{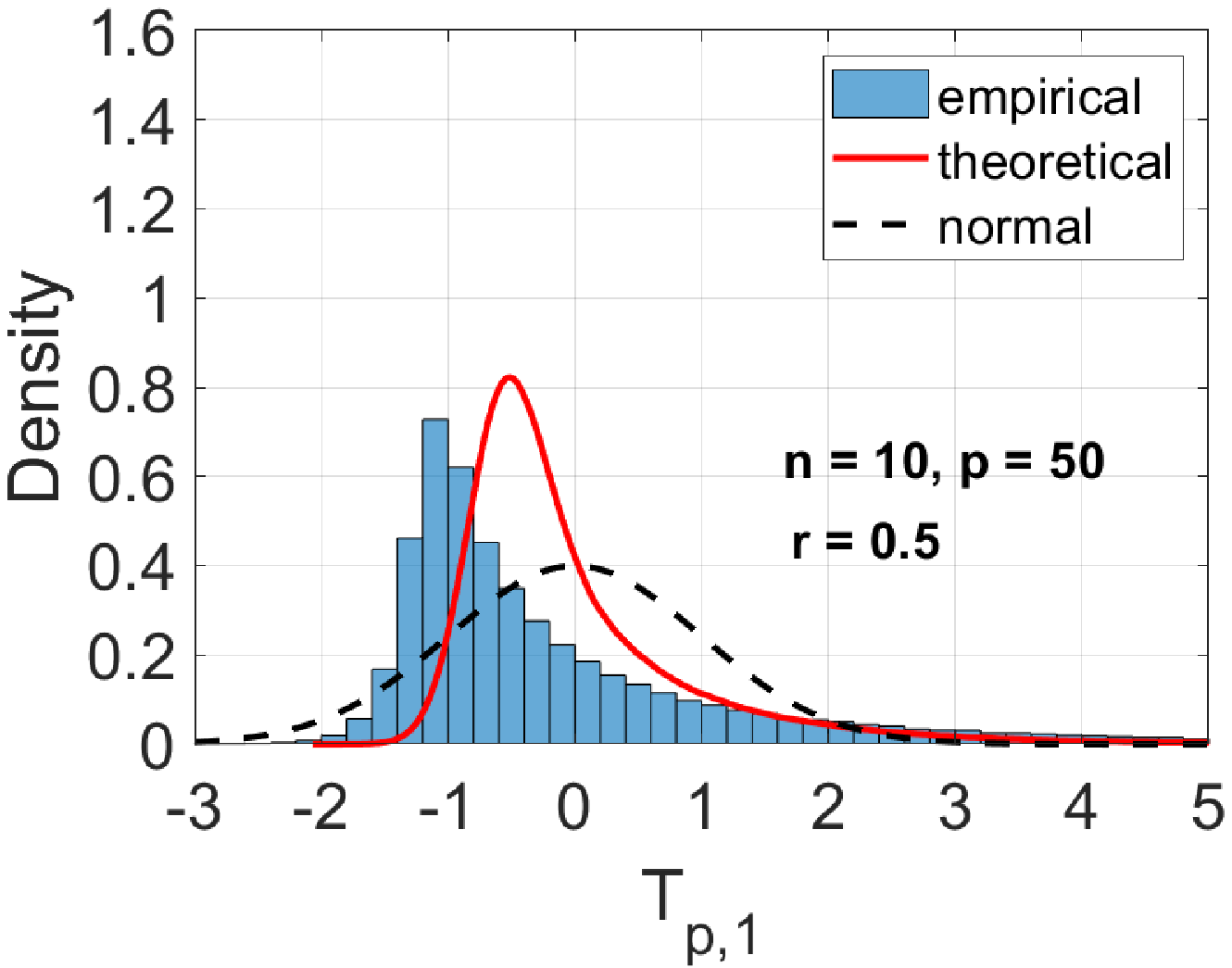}
    \includegraphics[width=2.2in]{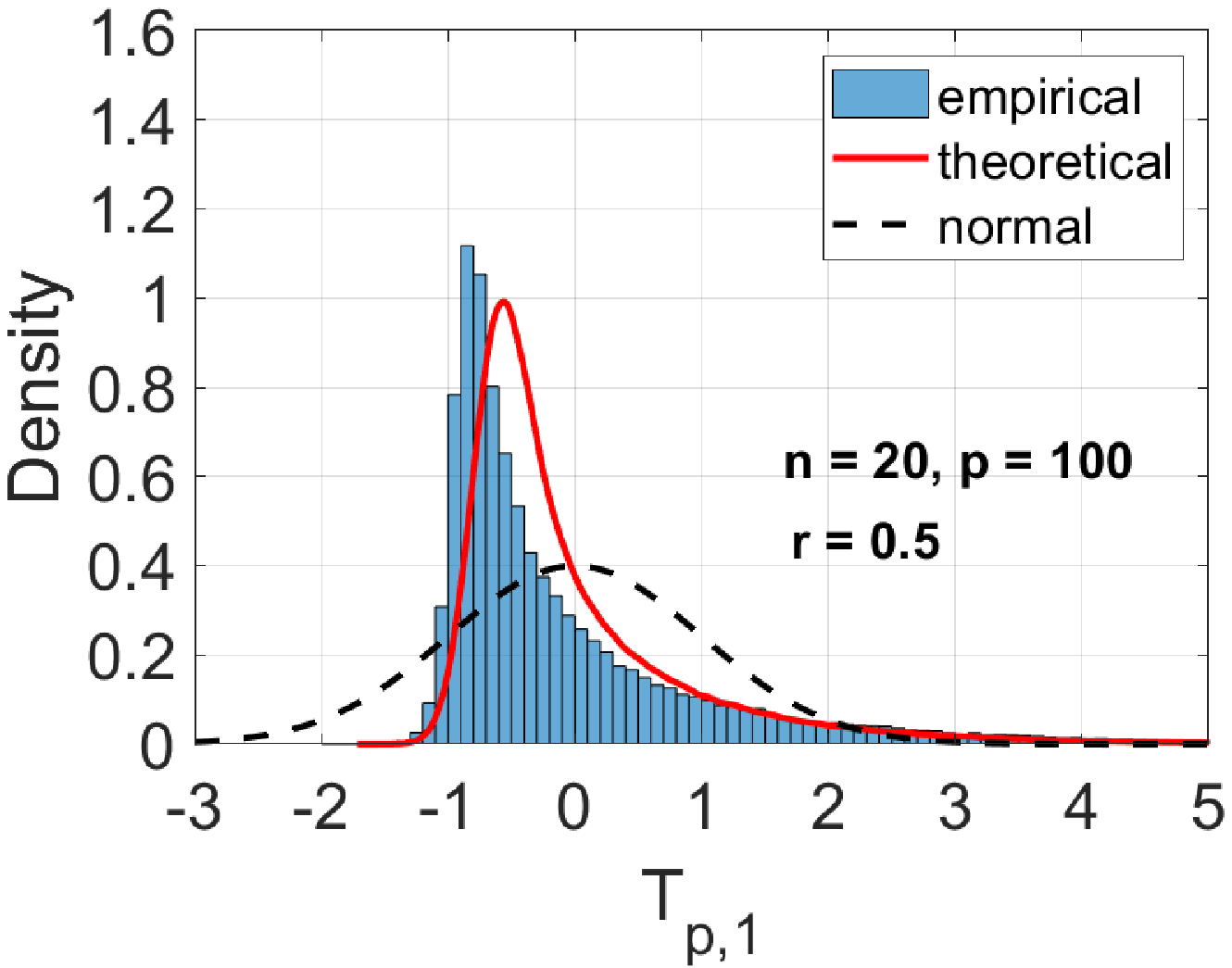}
    \includegraphics[width=2.2in]{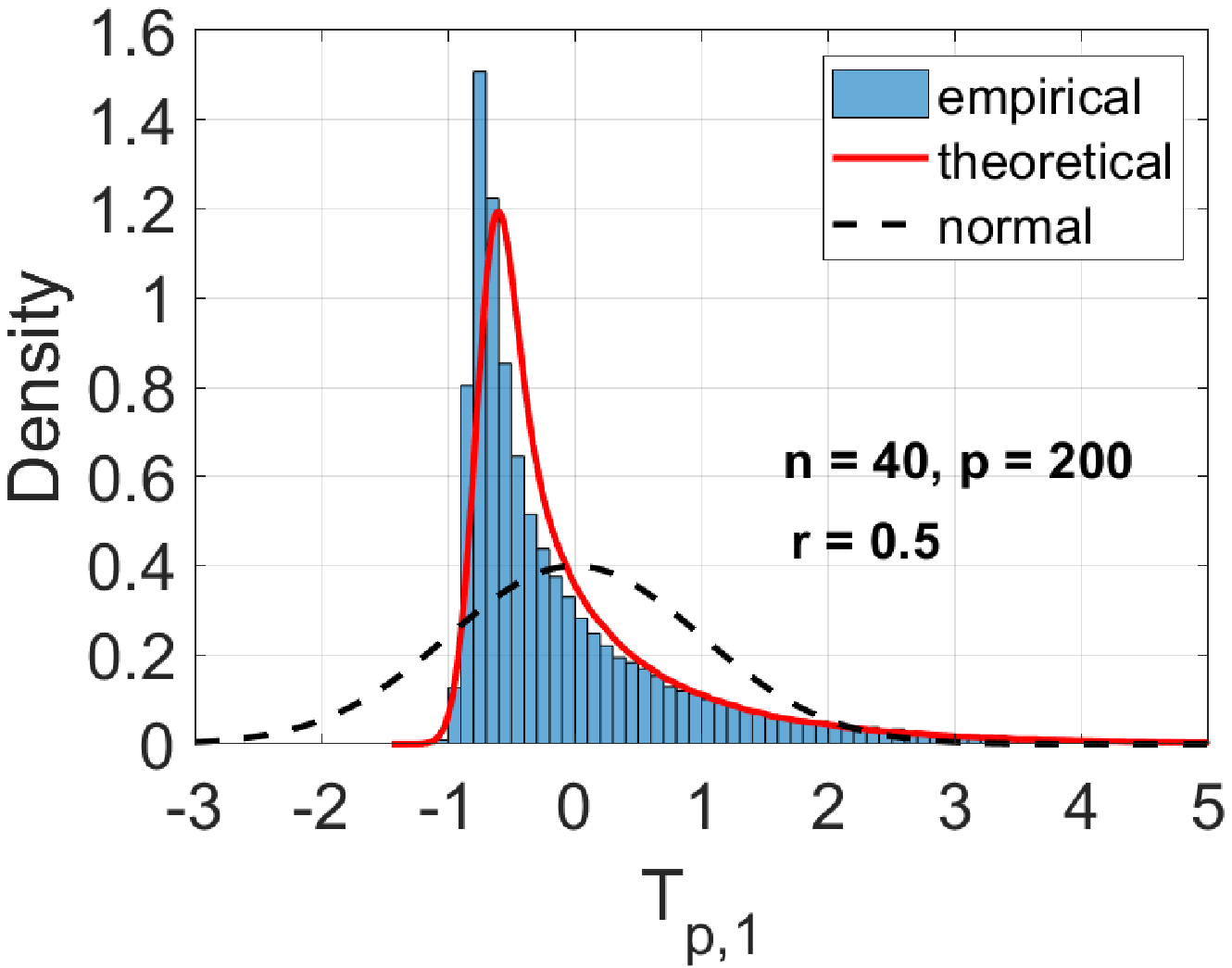}
}

\mbox{    \includegraphics[width=2.2in]{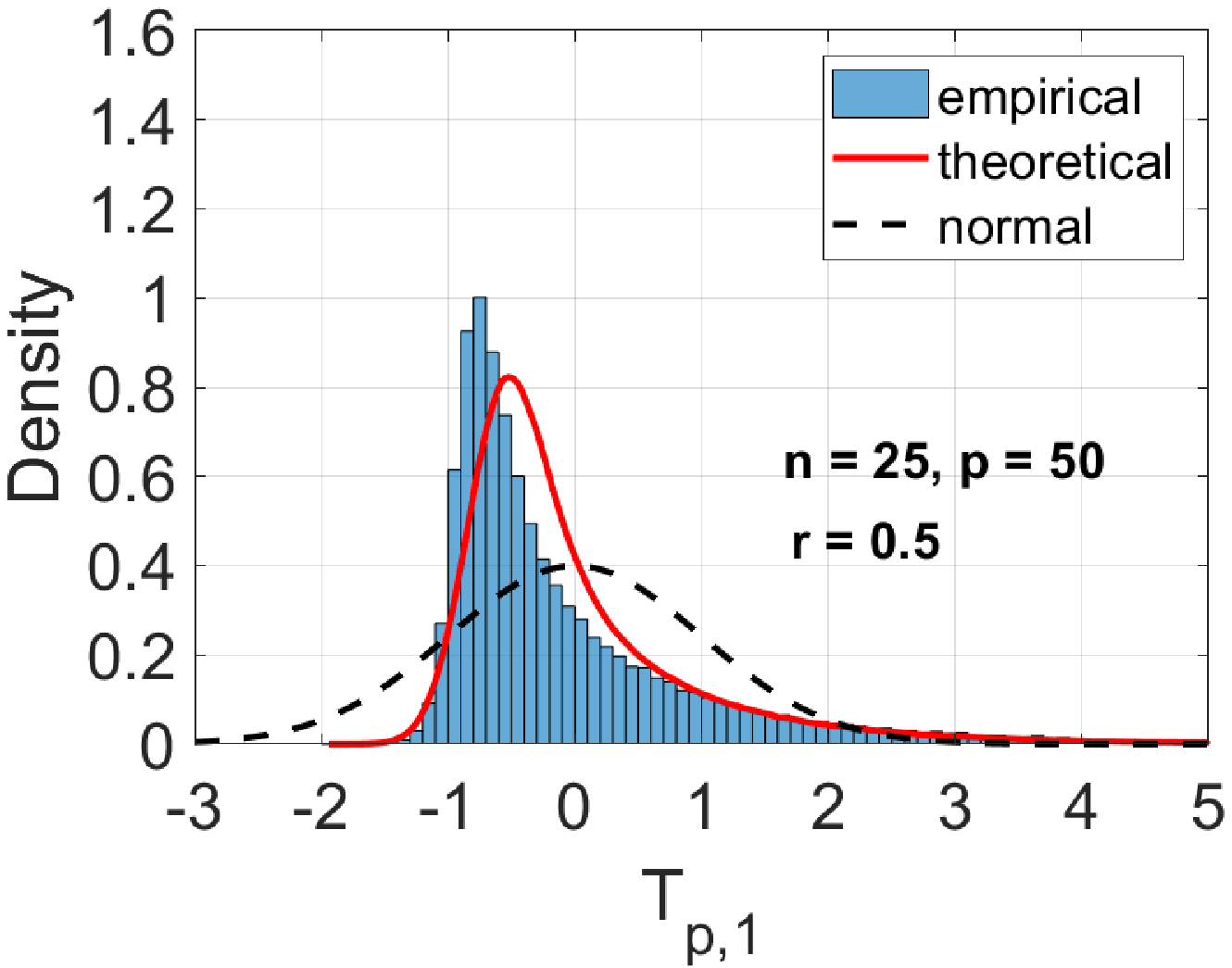}
    \includegraphics[width=2.2in]{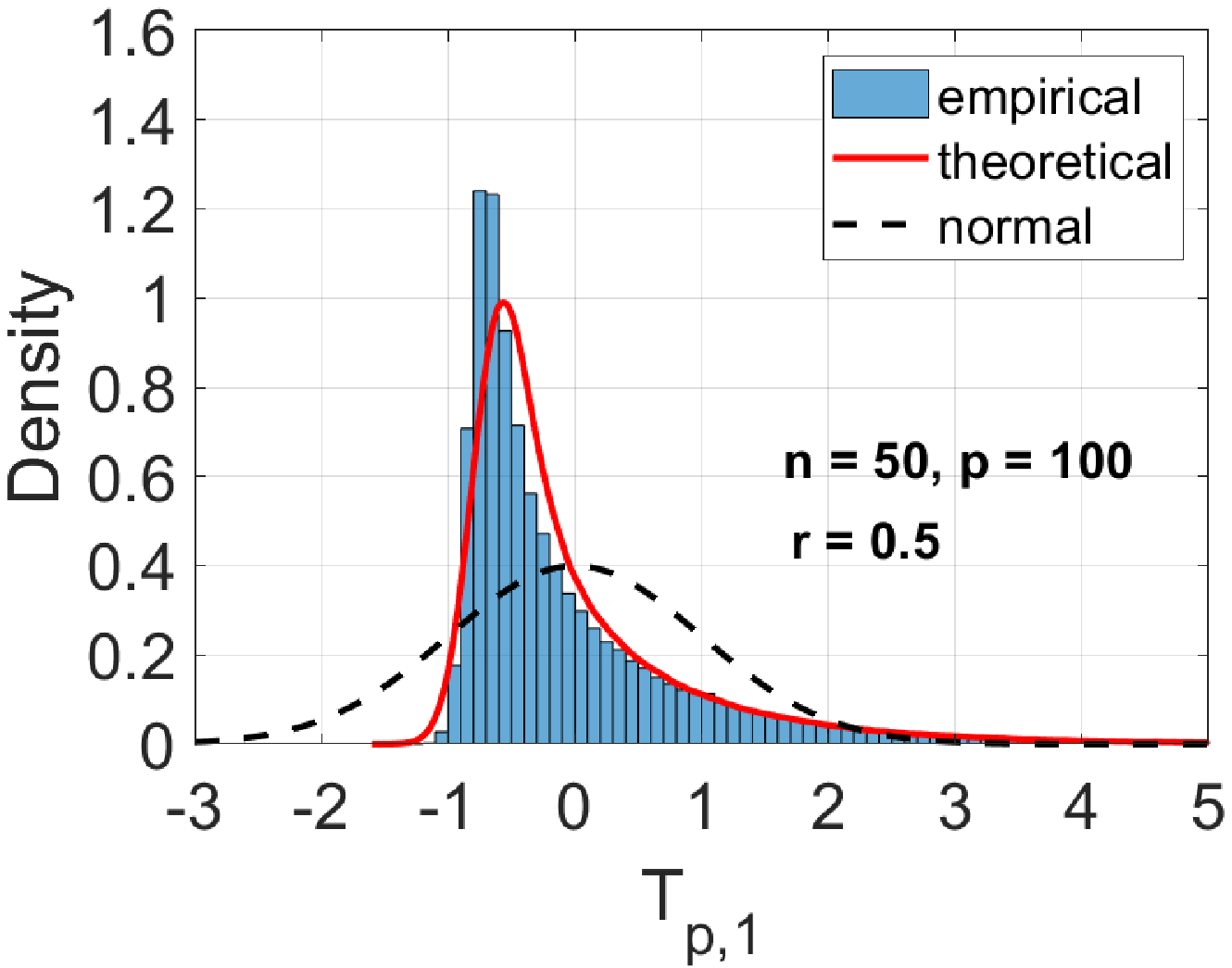}
    \includegraphics[width=2.2in]{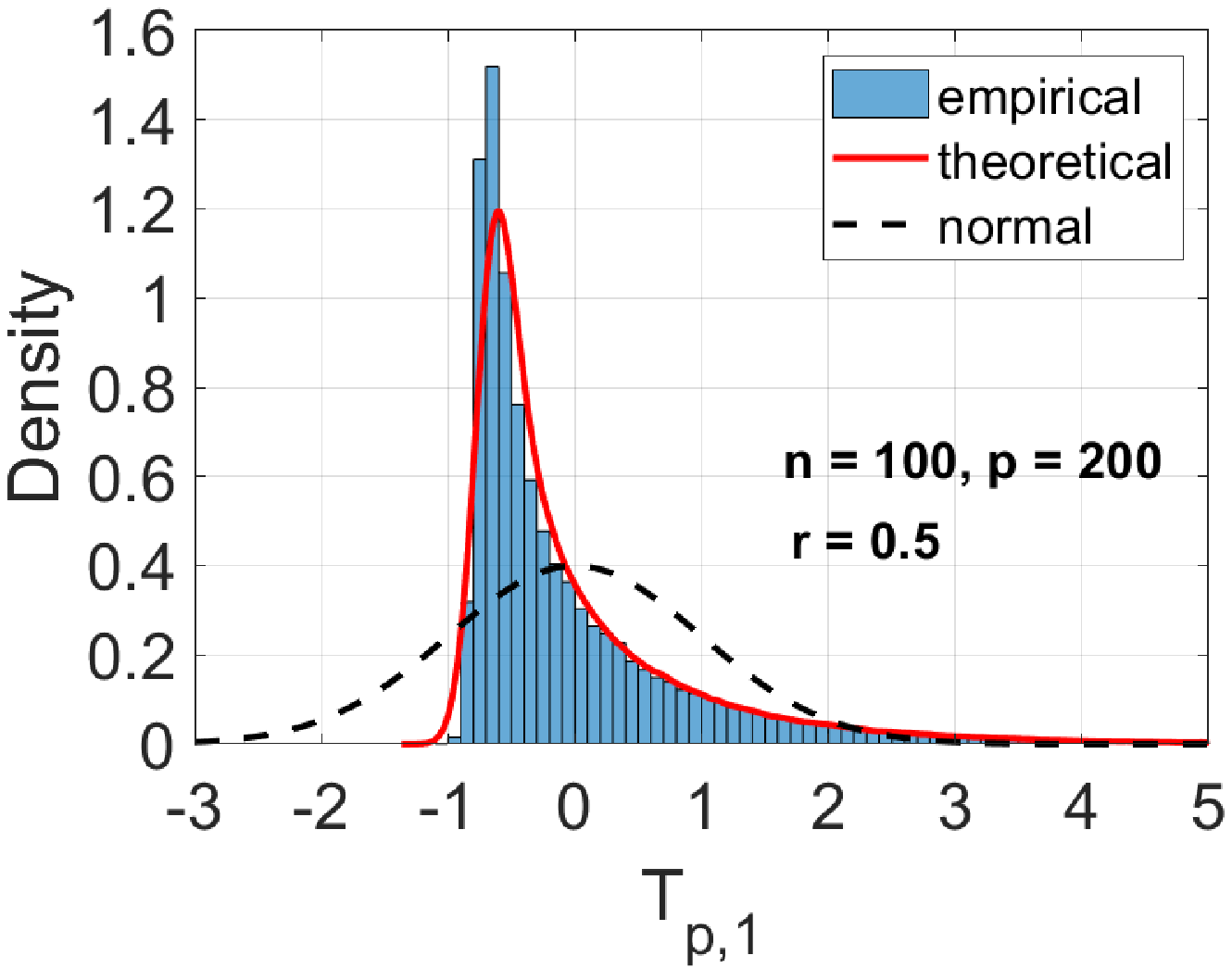}
}

\mbox{    \includegraphics[width=2.2in]{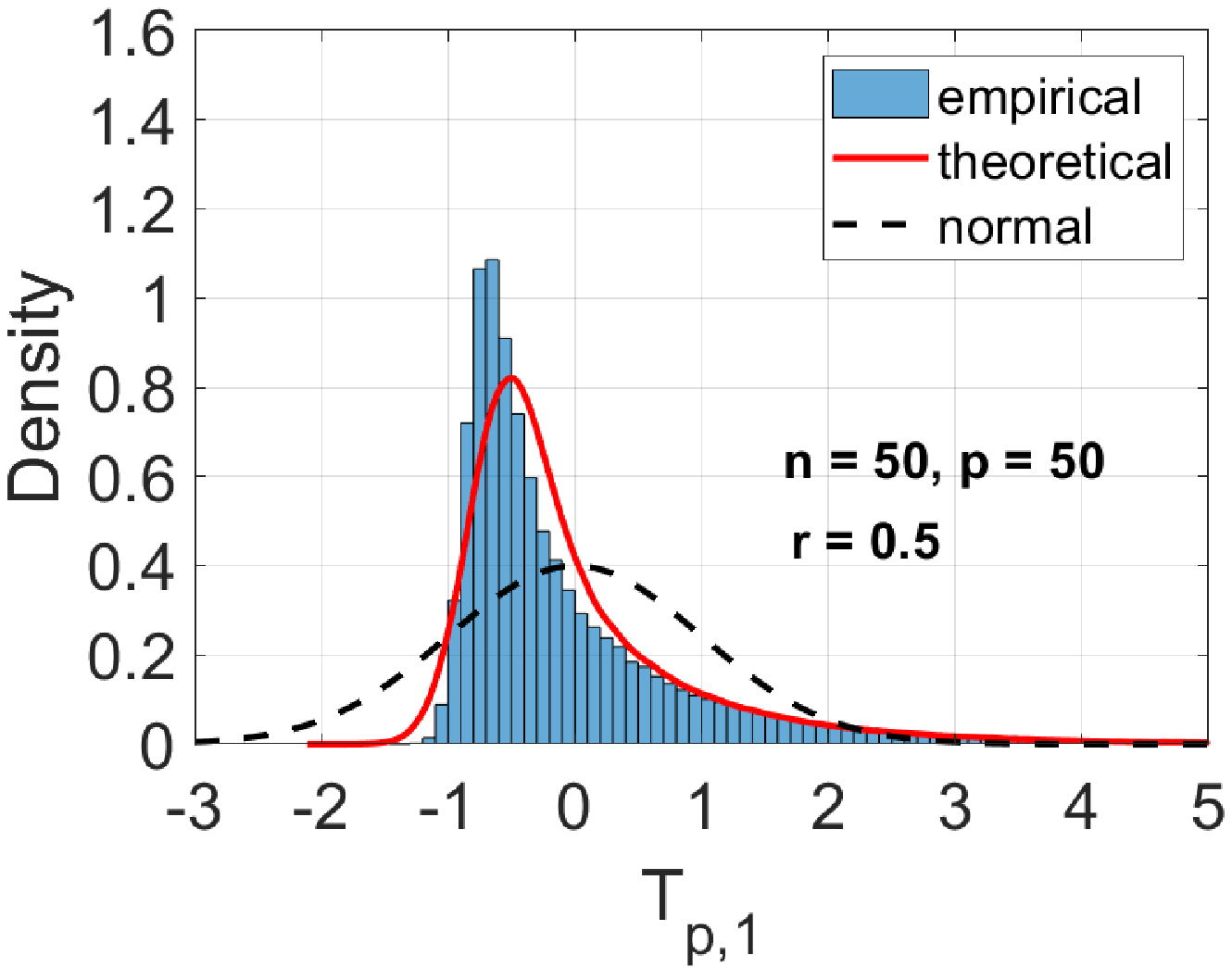}
    \includegraphics[width=2.2in]{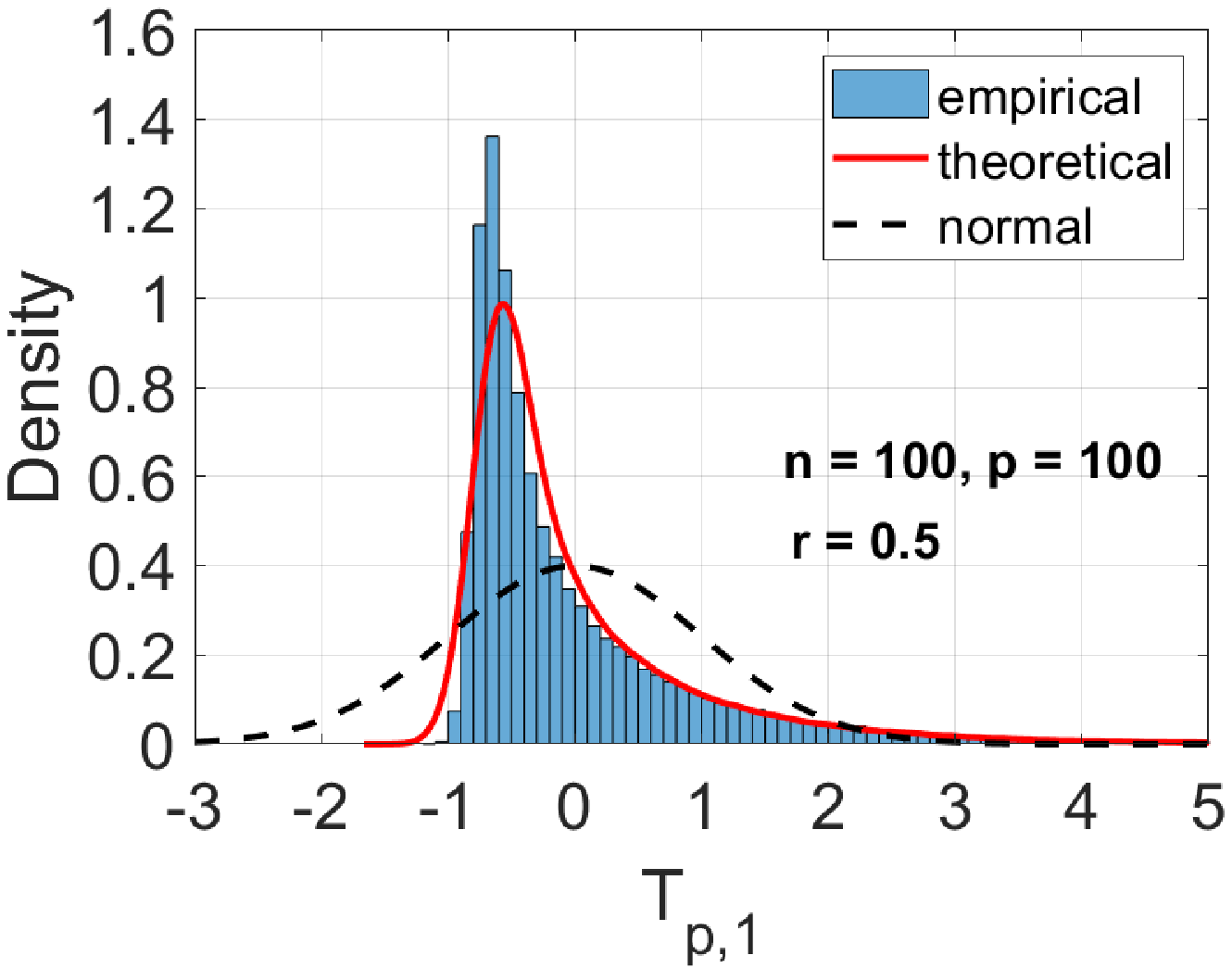}
    \includegraphics[width=2.2in]{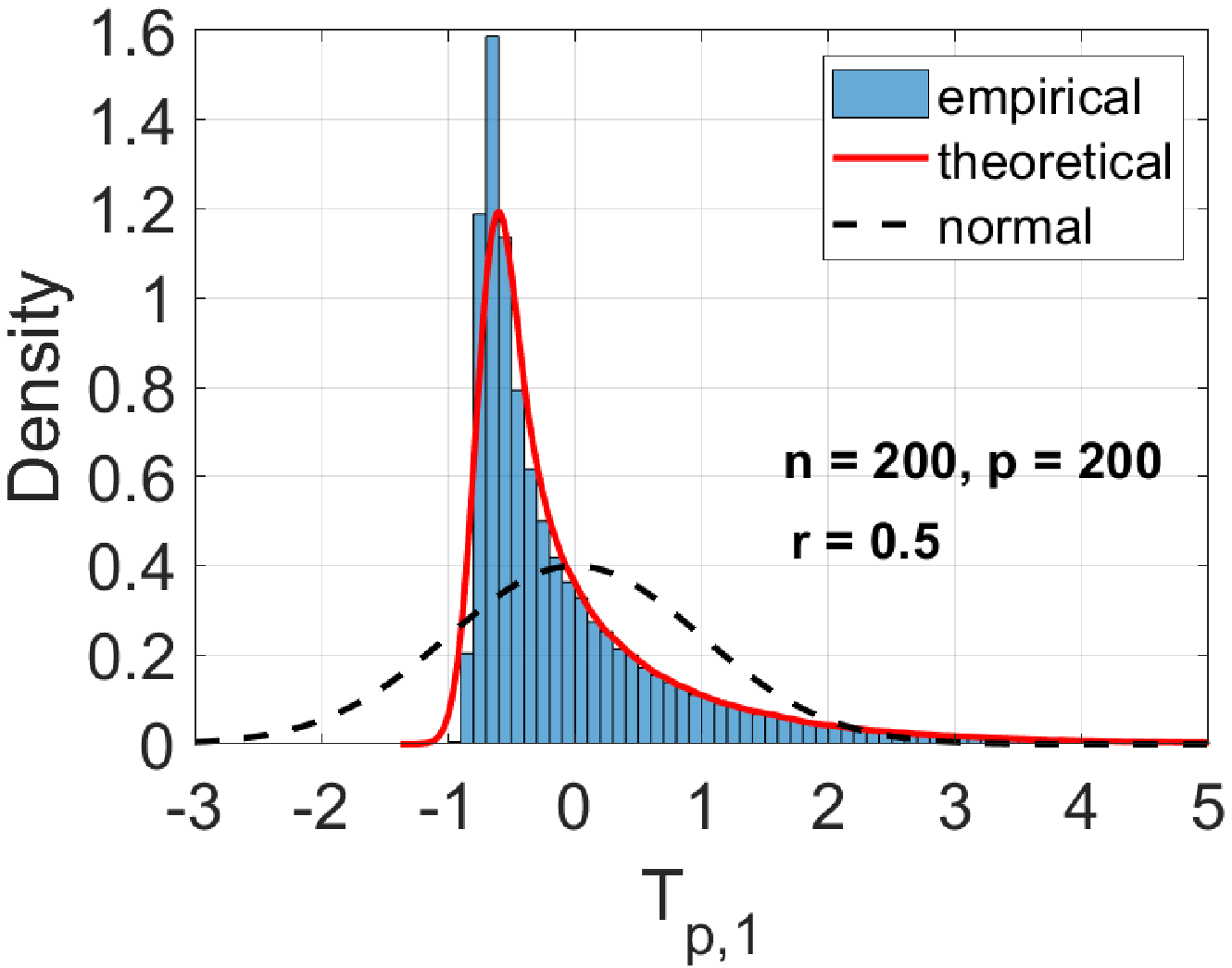}
}

\mbox{    \includegraphics[width=2.2in]{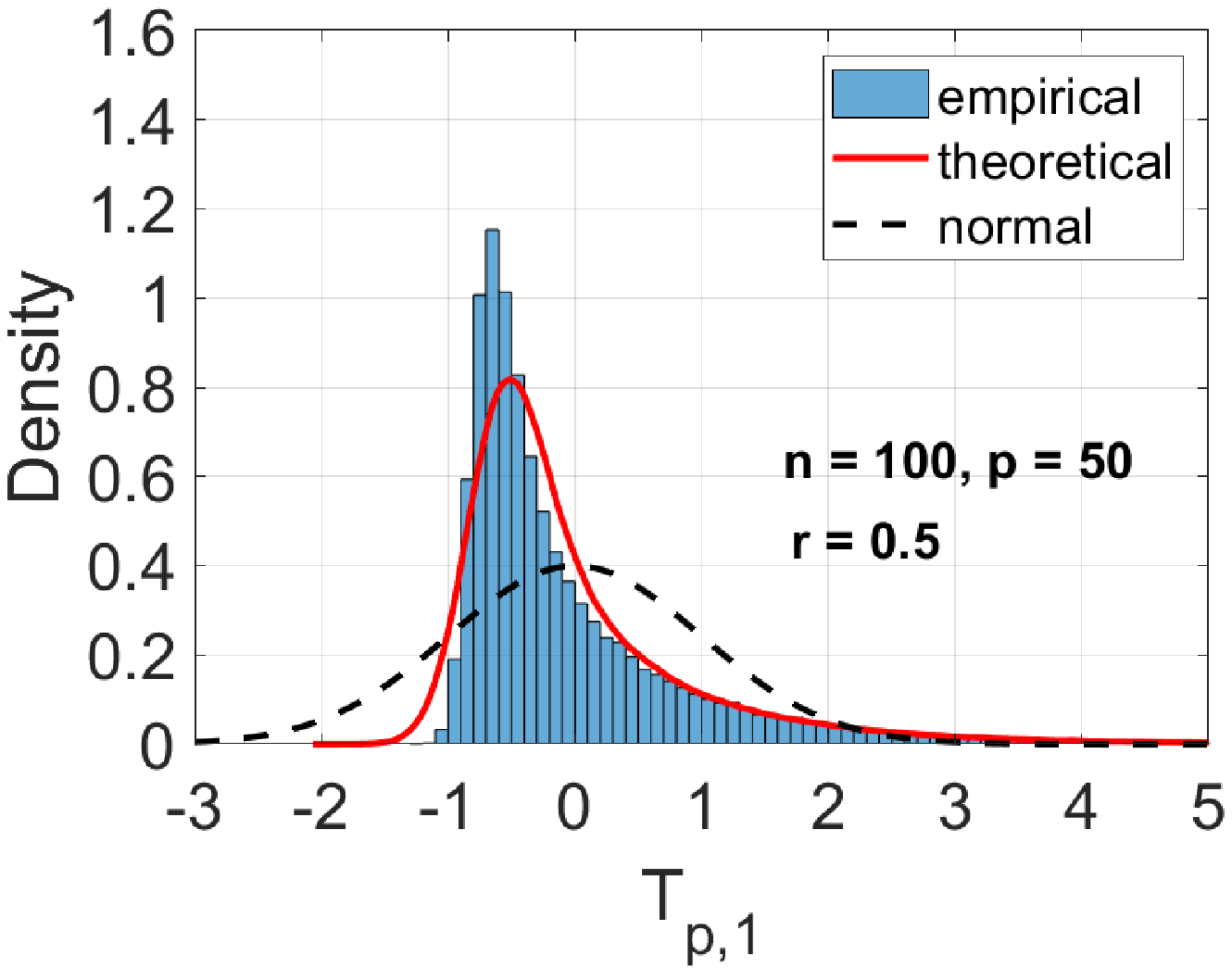}
    \includegraphics[width=2.2in]{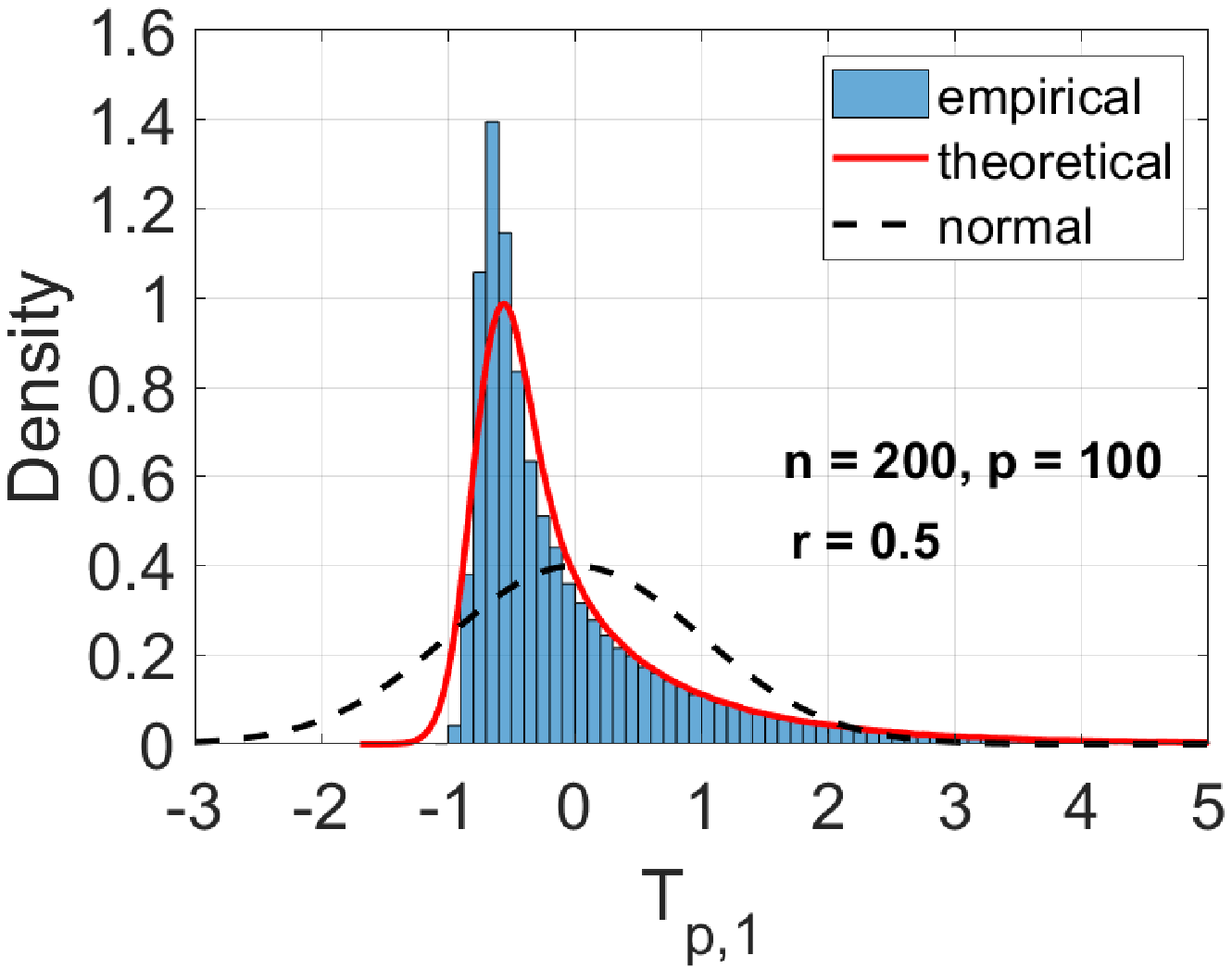}
    \includegraphics[width=2.2in]{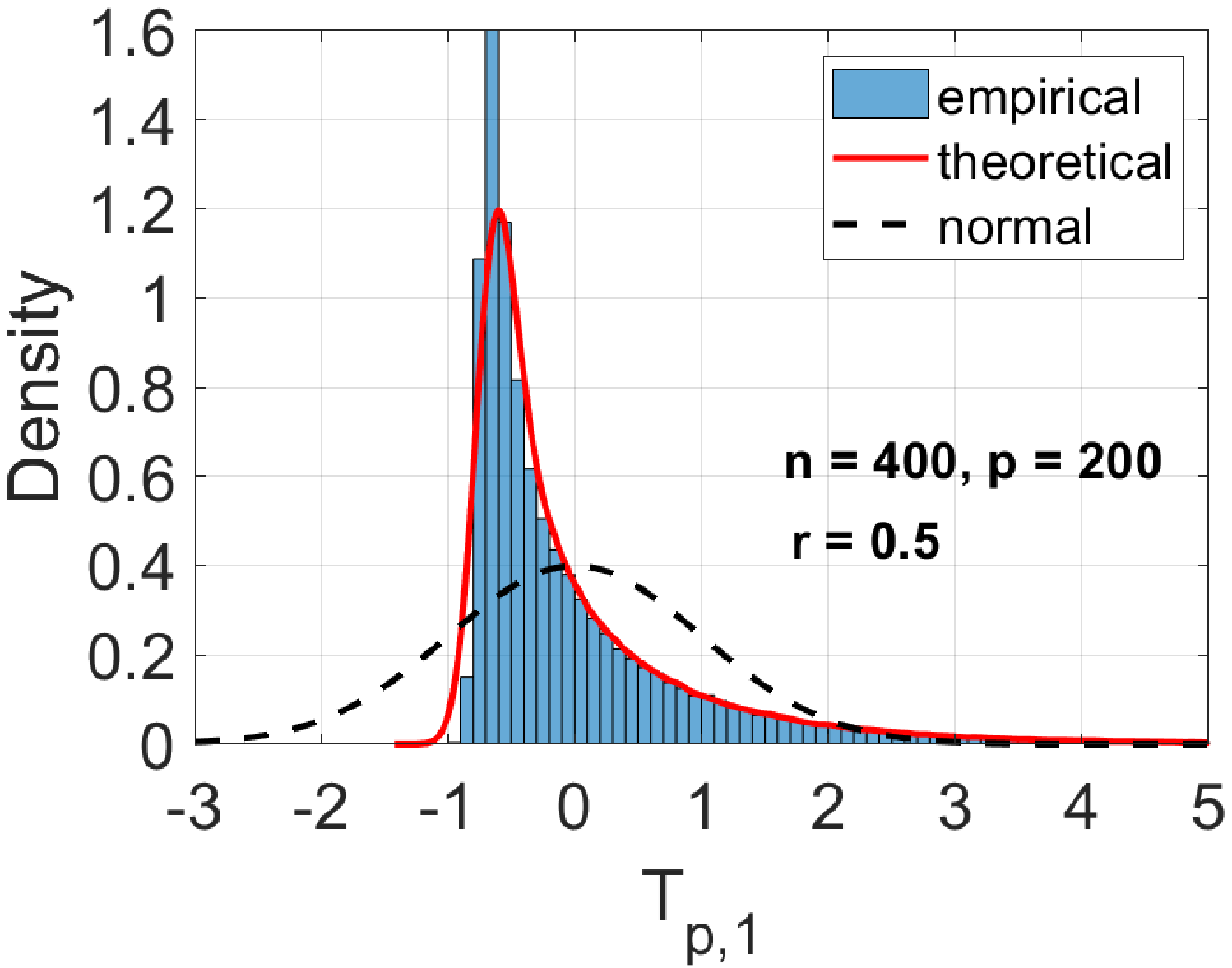}
}

\caption{Example~\ref{Ex1}. Compared to Figure~\ref{fig:Ex1_r01} with $r=0.1$, we take $r=0.5$ here. The empirical curves in blue are close to our theoretical curves in red. The normal approximation (black curve) from {\it Result~\ref{AD2008}} is no longer valid.}
\label{fig:Ex1_r05}
\end{figure}

We make Figures \ref{fig:Ex1_r01} and \ref{fig:Ex1_r05} based on Example~\ref{Ex1} to compare  the performances of our $T_{p,1}$ and $T_{SD}$ by  \cite{srivastava2008test}. For different values of $n,p,r$, the statistics $T_{p,1}$ and $T_{SD}$ are simulated for $10^5$ times. The elaboration is given below.

Figure \ref{fig:Ex1_r01} is designed for a small value of $r$ with $r=0.1$. Look at the first two pictures on the top row, that is, the ones with $(n,p)=(10, 50)$ and $(n,p)=(20, 100)$, respectively.  Neither $T_{p,1}$ nor $T_{SD}$ perform well. But this is expected  because the sample sizes are very small. However, by fixing the rate of $p/n$ and let $n$ and $p$ grow gradually, the approximations become better and better. In particular our theoretical curves are always close to the empirical ones, whereas the normal ones stay farther or much farther from the empirical ones.

Figure \ref{fig:Ex1_r05} is designed for a big value of $r$ with $r=0.5$. Except the first two pictures on the top row, which correspond to  $(n,p)=(10, 50)$ and $(n,p)=(20, 100)$, respectively, our theoretical curves match the empirical ones well. As the explanation aforementioned, that the two simulations do not behave well is understandable, simply because the sample sizes are very small. Our theoretical curves are constantly close to the empirical ones, the normal ones are nowhere close to the empirical ones.

In summary, the simulation indicates  our approximation stated in  Theorem~\ref{Theorem1} that  $(1-\sum_{i=1}^{\infty}\rho_i^2)^{1/2}\xi_0+\frac{1}{\sqrt{2}}\sum_{i=1}^{\infty}\rho_i(\xi_i^2-1)$  to $T_{SD}$ and $T_{p,1}$, respectively, outperforms the normal approximation from {\it Result~\ref{AD2008}}.

\newpage

\newpage\clearpage

\begin{example}\lbl{Ex2} For $r\in (0, 1)$, set $m=[p^{r}]$. Recall $\bd{A}_m$ from Example ~\ref{Ex1}. Define  $\bd{R}$ by
\beaa
\bd{R}=
\begin{pmatrix}
\bd{A}_m & \bd{0}\\
\bd{0} & \bd{I}_{p-m}
\end{pmatrix}
_{p\times p}
,
\eeaa
where the two ``$r$" from ``$m=[p^{r}]$" and ``$\bd{A}_m$" are the same one.
Assuming $\bmu=\bm 0$, we have
\bea\lbl{sfyi8}
T_{p}\to
\begin{cases}
\xi_0 & \text{if $0< r < \frac{1}{2}$};\\
\frac{2}{\sqrt{5}}\xi_0+\frac{1}{\sqrt{10}}\cdot (\xi_1^2-1),& \text{if $r = \frac{1}{2}$};\\
\frac{1}{\sqrt{2}}\cdot (\xi_1^2-1), & \text{if $\frac{1}{2}< r <1$}
\end{cases}
\eea
under condition $p=o(n^2)$ for $0< r \leq  \frac{1}{2}$ and $p=o(n^{1/(1-r)})$ for $\frac{1}{2}< r <1$, where $\xi_0$ and $\xi_1$ are i.i.d. $N(0, 1)$ as $r = \frac{1}{2}$. Here $T_p$ stands for $T_{SD}$ or $T_{p,1}$.  Obviously, there is a phase transition at $r=\frac{1}{2}$ as $r$ runs between $0$ and $1$.  The phase transition appears at both relative sizes of $n$ and $p$ together with the change of the values of matrix entries.
\end{example}

\begin{example}\lbl{Ex3} Set $m=[\log (p+2)]. $ Given $\tau\geq 0$, define integer $p'=p-[\tau\sqrt{p}]$, $\lambda_i=1+\tau 2^{-i}\big(1-2^{-m}\big)^{-1}\sqrt{p}$ for $1\leq i \leq m$ and $\lambda_i=1$ for $i=m+1, \cdots, p'-1$, $\lambda_{p'}=1+[\tau\sqrt{p}]-\tau\sqrt{p} \in [0, 1]$, and the rest of $\lambda_i$ are identical to $0$.
By Lemma~\ref{dsih328}, there exists a correlation matrix $\bd{R}$ such that $\bd{R}$ has eigenvalues $\lambda_i,\, 1\leq i \leq p.$ Assuming $\bmu=\bm 0$,  we have
\bea\lbl{sdispi0sa}
T_{p}\to
\begin{cases}
\xi_0, & \text{if $\tau=0$};\\
\sqrt{\frac{3}{\tau^2+3}}\xi_0+\sqrt{\frac{3\tau^2}{2(\tau^2+3)}}
\sum_{i=1}^{\infty}\frac{1}{2^i}(\xi_i^2-1),& \text{if $\tau\in (0, \infty)$};\\
\sqrt{\frac{3}{2}}\sum_{i=1}^{\infty}\frac{1}{2^i}(\xi_i^2-1), & \text{if $\tau\to \infty$}
\end{cases}
\eea
in distribution as $p\to\infty$ under condition $p=o(n^2)$.  Here $\xi_0, \xi_1, \cdots, \xi_d$ are i.i.d. $N(0, 1).$ Here $T_p$ stands for $T_{SD}$ or $T_{p,1}$. The statistic $T_{p}$ behaves like a rainbow, which has a Gaussian distribution at one end, and a mixing chi-squared distribution at the horizon, and a hybrid in between.
\end{example}

\newpage

An interesting remark is that, unlike the eigenvalues of a covariance matrix, the eigenvalues $\lambda_1\geq \lambda_2\geq \cdots\geq \lambda_p\geq 0$ of the $p\times p$  correlation matrix $\bd{R}$ cannot be arbitrary. By Lemma~\ref{dsih328}, they have to satisfy the so-called ``majorization" property: $\lambda_1+\cdots +\lambda_k \geq k$ for every $1\leq k \leq p$ and $\lambda_1+\cdots +\lambda_p=p$. Conversely, for any $\lambda_1\geq \lambda_2\geq \cdots\geq \lambda_p\geq 0$ with the ``majorization" property,  there is always a correlation matrix with  eigenvalues $\lambda_i,\, 1\leq i \leq p$.

\section{Two Sample Mean Test for Large $p$ and Small $n$}\lbl{Sec_two_mean}

We will study the two-sample mean testing problem
in the high-dimensional setting in this section. As before, assume that
$\{\X_{i1},\cdots,\X_{in_i}\}$ for $i=1,2$ are two independent
random samples with sizes $n_1$ and $n_2$, and from $p$-variate
normal distributions $N_p(\bmu_1,\bms)$ and $N_p(\bmu_2,\bms)$, respectively.  We
wish to test
\begin{align}\label{ht}
H_0:\bmu_1=\bmu_2\ \ \mbox{versus}\ \  H_1:\bmu_1\neq\bmu_2.
\end{align}
Let $\bar{\X}_i$ be  the sample mean  of the $i$-th sample  with $i=1, 2$ and $\hat{\S}$ be the pooled sample covariance matrix defined in \eqref{duck_yazi}. Set $n=n_1+n_2$. \cite{srivastava2008test} consider the  statistic $T'_{SD}$ from \eqref{u82su}.
Recall $\hat{\D}$ is the diagonal matrix of $\hat{\S}$ in \eqref{duck_yazi} and
\bea\lbl{efwei0}
\hat{\R}=\hat{\R}_p=\hat{\D}^{-1/2}\hat{\S}\hat{\D}^{-1/2}
\eea
is the pooled sample correlation matrix.
Similar to the discussion before \eqref{Statistics1} on the modification of $T_{SD}$, we make a little change for $T_{SD}'$ as follows. Define
\bea\lbl{Statistics2}
T_{p, 2}=\frac{\frac{n_1n_2}{n_1+n_2-1}(\bar{\X}_1-\bar{\X}_2)^{T} \hat{\D}^{-1}(\bar{\X}_1-\bar{\X}_2)-\frac{(n_1+n_2-1)p}{n_1+n_2-4}}
{\sqrt{2\big|\tr(\hat{\R}^2)-\frac{p(p-1)}{n_1+n_2-2}\big|}}.
\eea
\cite{srivastava2008test} derive the limiting distribution of $T_{SD}'$ under conditions similar to those from {\it Result~\ref{AD2008}} and an extra assumption ``$\frac{n_1}{n_1+n_2}\to c\in (0,1)$". On the other hand, \cite{hu2016review} ask for the properties of $T_{SD}'$ as the variances of the entries of the population vector are quite different. We will provide a  solution for $T'_{SD}$ and $T_{p, 2}$ next. For mathematical rigor, we assume that both sample sizes $n_1$ and $n_2$  depend on $p$.

\begin{theorem}\lbl{Theorem2} Let $\{\X_{i1},\cdots,\X_{in_i}\}$ for $i=1,2$ be two independent random samples from  $N_p(\bmu_1,\bms)$ and $N_p(\bmu_2,\bms)$, respectively. Let $\bd{R}$ be as in \eqref{38488} with eigenvalues $\lambda_1\geq  \cdots \geq \lambda_p \geq 0$. Assume

(a)~ $\lim_{p\to\infty}\frac{\lambda_i}{\|\bd{R}\|_F}=\rho_i \geq 0$ for all $i\geq 1$;

(b)~ $\lim_{p\to\infty}\frac{p}{(n_1+n_2)\|\bd{R}\|_F}= 0$ and $\lim_{p\to\infty}\frac{p}{(n_1+n_2)^{a}}=0$ for some constant $a>0$.\\
If $\bmu_1=\bmu_2$, then
$T_{p,2}\to  b\xi_0+\frac{1}{\sqrt{2}}\sum_{i=1}^{\infty}\rho_i(\xi_i^2-1)$
in distribution, where $\xi_0, \xi_1, \xi_2, \cdots$ are i.i.d. $N(0, 1)$ and $b=(1-\sum_{i=1}^{\infty}\rho_i^2)^{1/2}$. The same conclusion  also holds for ``\,$T_{p, 2}$".
\end{theorem}

The above result does not require condition ``$\min\{n_1, n_2\}\to \infty$" or  ``$\frac{n_1}{n_1+n_2}\to \gamma\in (0, 1)$", which are commonly assumed in literature for the two-sample testing problem.

As discussed below Theorem~\ref{Theorem1}, condition
 \eqref{SD08a} implies that $\rho_1=0$ and  $\frac{p}{n\|\bd{R}\|_F}\to 0$ as $p\to\infty$. Condition ``$n_1+n_2=O(p^{\zeta})$, $\frac{1}{2}<\zeta\leq 1$" from \eqref{SD08a} implies (b) from Theorem~\ref{Theorem2}.  So our theorem is more general than the conclusion for the two-sample test by \cite{srivastava2008test}.

If $\bd{R}$ is $AR(1)$ model, the banded model, the Toeplitz or Hankel matrices, please see  the discussion below Theorem~\ref{Theorem1}. They are also applied here.

Let $\X_1,\cdots,\X_n$ be a random sample from an $AR(1)$ model  with $\bd{R}=(\gamma^{|i-j|})$ and the absolute values of  $\gamma=\gamma_p$ staying away from $1$. By using the Gersgorin disc theorem [see, e.g., p. 344 from \cite{horn2012matrix}], the largest eigenvalue or $\bd{R}$ is of order $O(1)$. Hence  condition (a) of Theorem~\ref{Theorem1} holds with $\rho_1=0$. If condition (b) also holds, then both $T_{SD}$ and $T_{p,1}$ go to the standard normal distribution. The same conclusion is also valid for a banded correlation matrix $\bd{R}=(r_{ij})$ with $r_{ij}=0$ for $|j-i|\geq t$ where $t=t_p=o(\sqrt{p})$. In this case, the largest eigenvalue or $\bd{R}$ is of order $o(\sqrt{p})$. Similar results can be obtained for other patterned matrices including Toeplitz matrices, Hankel matrices and symmetric circulant matrices; see, e.g., \cite{brockwell2016introduction}.

Recall the definition of $T_{p, 2}$ from \eqref{Statistics2}. Let $\bd{D}$ be the diagonal matrix of $\bms$. Under some conditions including  $\|\bd{R}\|_F$ is of the order $p$, \cite{zhang2020simple} obtain the limiting distribution of a normalized $(\bar{\X}_1-\bar{\X}_2)^{T} \D^{-1}(\bar{\X}_1-\bar{\X}_2)$ as  $\frac{1}{\sqrt{2}}\sum_{i=1}^{\infty}\rho_i(\xi_i^2-1)$.  Although this quantity is not directly applicable because the unknown parameter matrix $\D$ is involved, it is indeed suggestive. Their result seems to be a special case of Theorem~\ref{Theorem2} with $\sum_{i=1}^{\infty}\rho_i^2=1$. In fact, under their assumptions, we have confirmed that   $\sum_{i=1}^{\infty}\rho_i^2=1$ in Lemma
\ref{Tuan_student} from Section~\ref{duhwoi9}.
Notice $\sqrt{p}\leq \|\bd{R}\|_F\leq p$ for any $\bd{R}$. \cite{srivastava2008test} consider the extreme case that $\|\bd{R}\|_F$ is of the scale of $\sqrt{p}$. This is essentially a weakly dependent situation with $\rho_i=0$ for each $i\geq 1$. According to Theorem~\ref{Theorem2}, the limiting distribution is a normal. On the other hand, roughly speaking, \cite{zhang2020simple} intend to study another extreme case in which $\|\bd{R}\|_F$ is of the scale $p$, and hence $\bd{R}$ is completely a spiked model. As explained above, the corresponding limiting distribution is a mixing chi-squared distribution. Our Theorem~\ref{Theorem2} handles the case for any $\bd{R}$ with $\|\bd{R}\|_F$ running everywhere between $\sqrt{p}$ and $p$, and the limiting distribution  turns out to be, interestingly enough, a convolution of both.

Now we demonstrate our results by giving some examples.
\begin{example}\lbl{Ex4} Recall Example~\ref{Ex1} and
$\bd{R}=\bd{A}_p$. For this example, \cite{zhang2020simple} observe that ``When $r = 0.01$, the histograms are quite symmetric and bell-shaped, indicating that a normal approximation as suggested by the theory of \cite{srivastava2008test} may be applied for approximating the null distribution of $T_{SD}$. However, when $r = 0.5$ and $r = 0.9$, the histograms are quite skewed, indicating that a normal approximation is no longer adequate." In fact, our Theorem~\ref{Theorem2} explains their inspection very accurately as follows. Assume $r=r_p$ and $\lim_{p\to\infty}\sqrt{p}\cdot r=c\geq 0$. By changing ``$n$" in Example~\ref{Ex1} to ``$n_1+n_2$", we have
\bea\lbl{o4398gh}
T_{p}\to
\begin{cases}
\xi_0, & \text{if $c=0$};\\
\frac{1}{\sqrt{c^2+1}}\xi_0+\frac{c}{\sqrt{2(c^2+1)}}(\xi_1^2-1),      & \text{if $c \in (0, \infty)$};\\
\frac{1}{\sqrt{2}}\cdot (\xi_1^2-1), & \text{if $c=\infty$}
\end{cases}
\eea
under condition $p=o((n_1+n_2)^2)$ for the case $c\in [0,\infty)$ and under the condition $p=o((n_1+n_2)^a)$ for some constant $a>0$ for the case  $c=\infty$, where $\xi_0$ and $\xi_1$ are i.i.d. $N(0, 1).$ Here $T_p$ stands for $T_{SD}'$ or $T_{p,2}$. Obviously, there is a phase transition between $c=0$ and $c=\infty$.
\end{example}
\newpage

\begin{figure}[h!]
\mbox{
\includegraphics[width=2.2in]{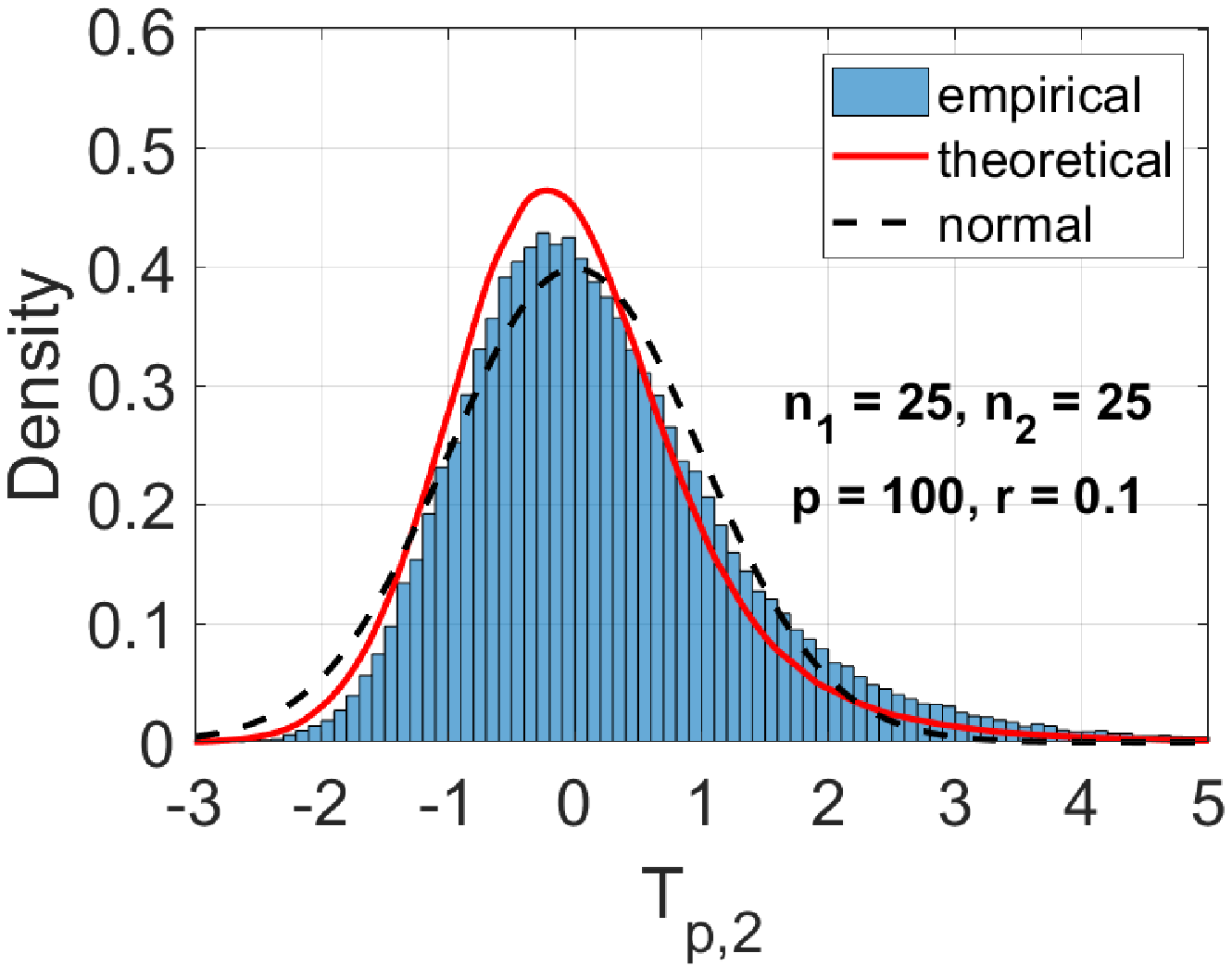}
\includegraphics[width=2.2in]{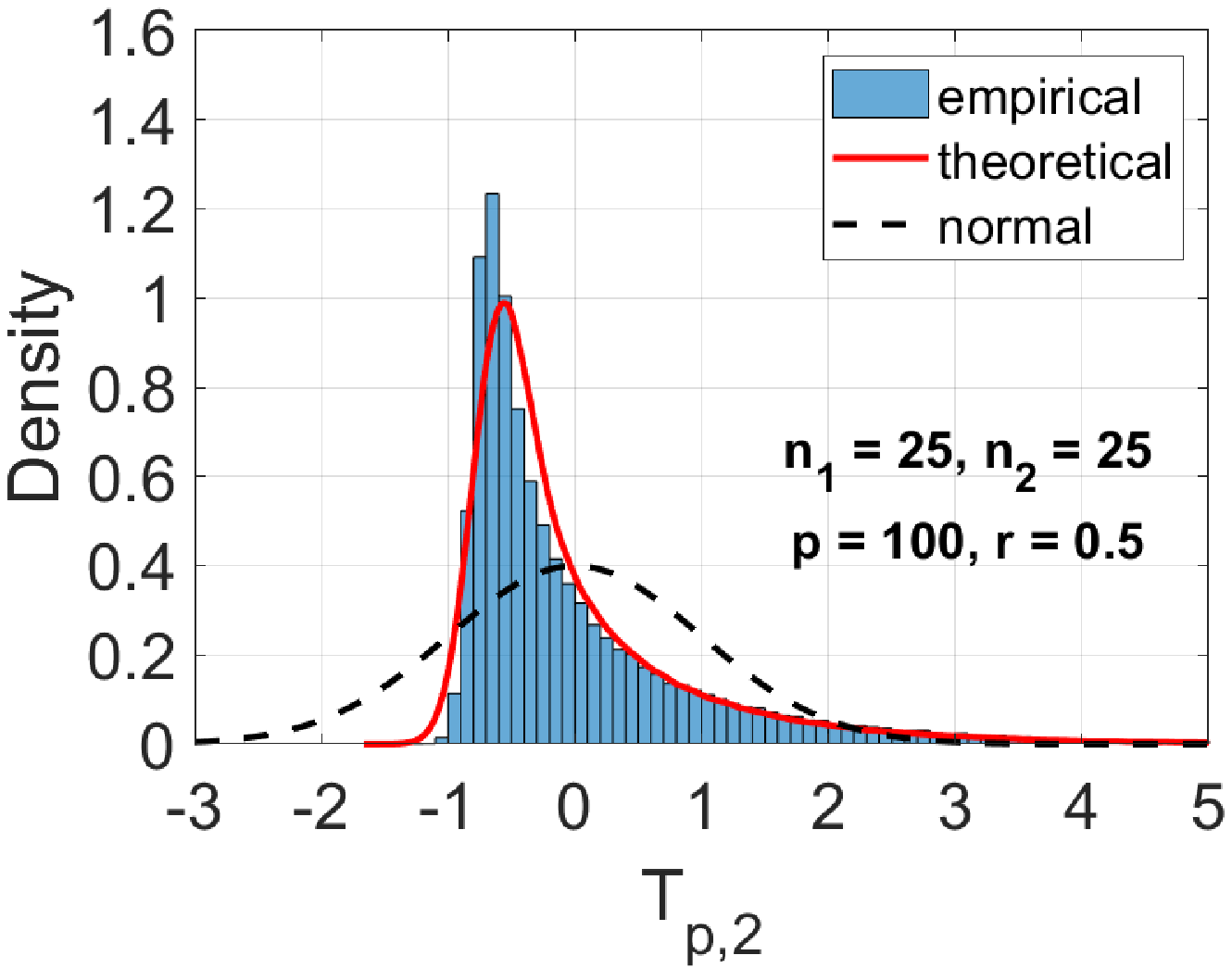}
\includegraphics[width=2.2in]{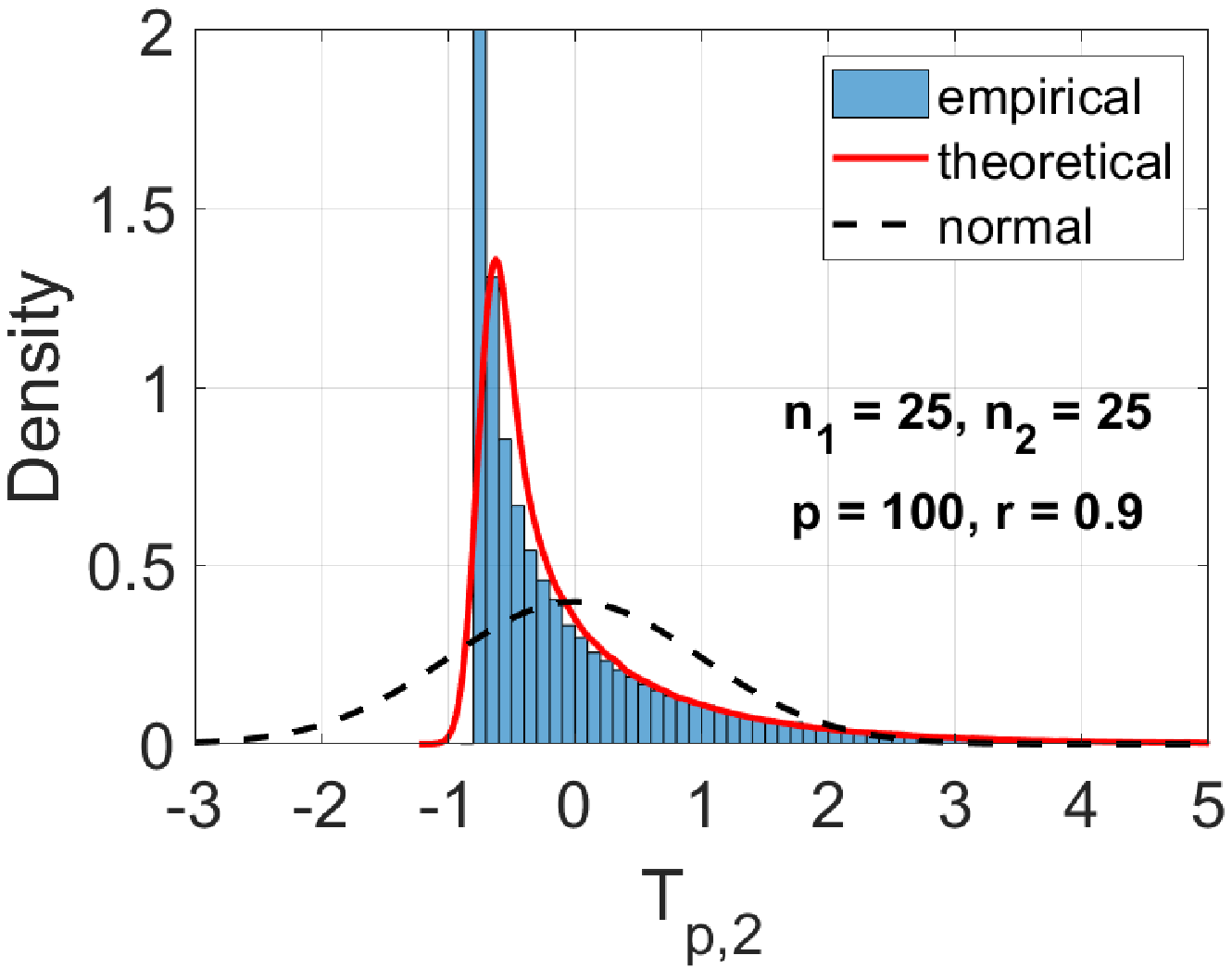}
}

\mbox{
\includegraphics[width=2.2in]{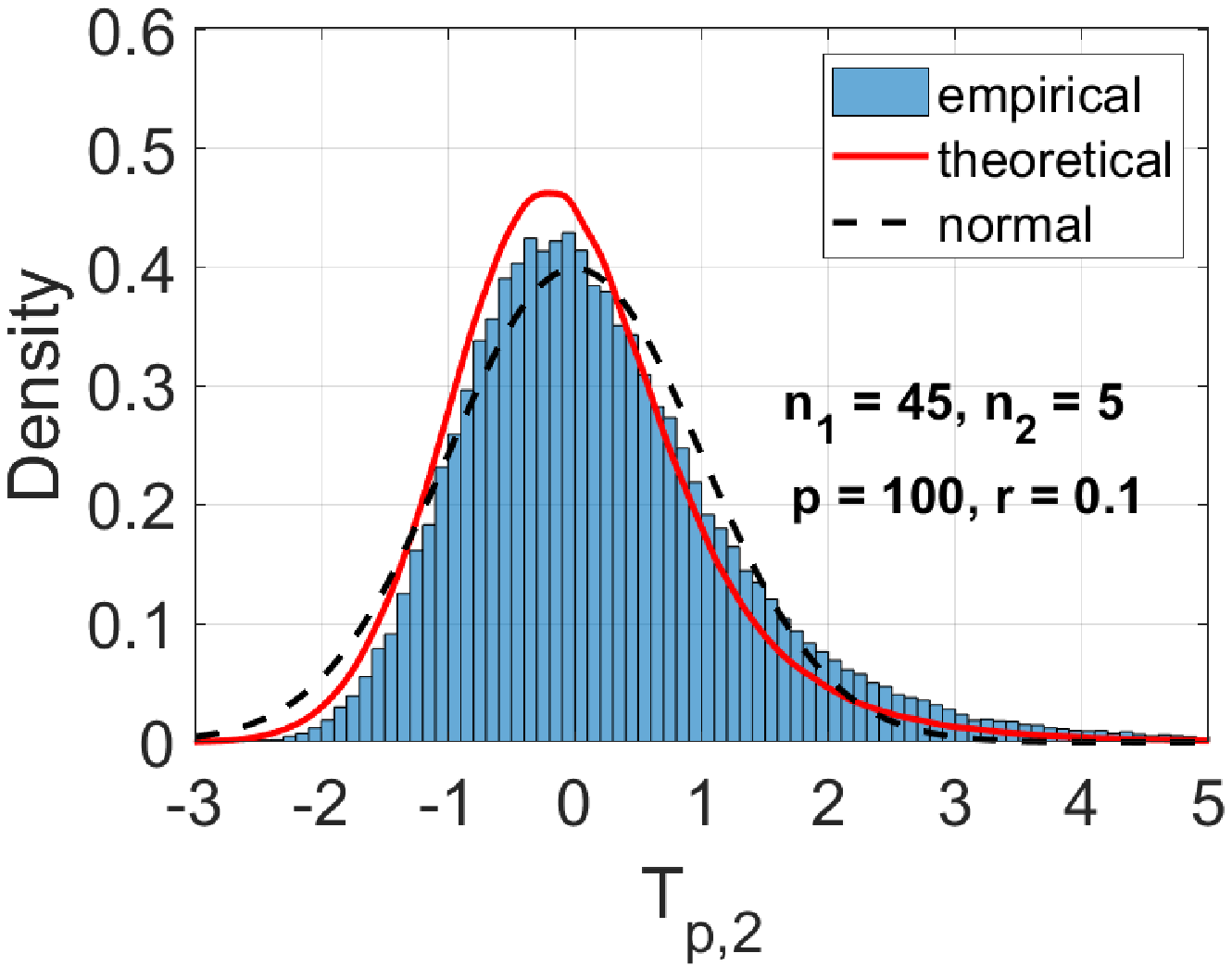}
\includegraphics[width=2.2in]{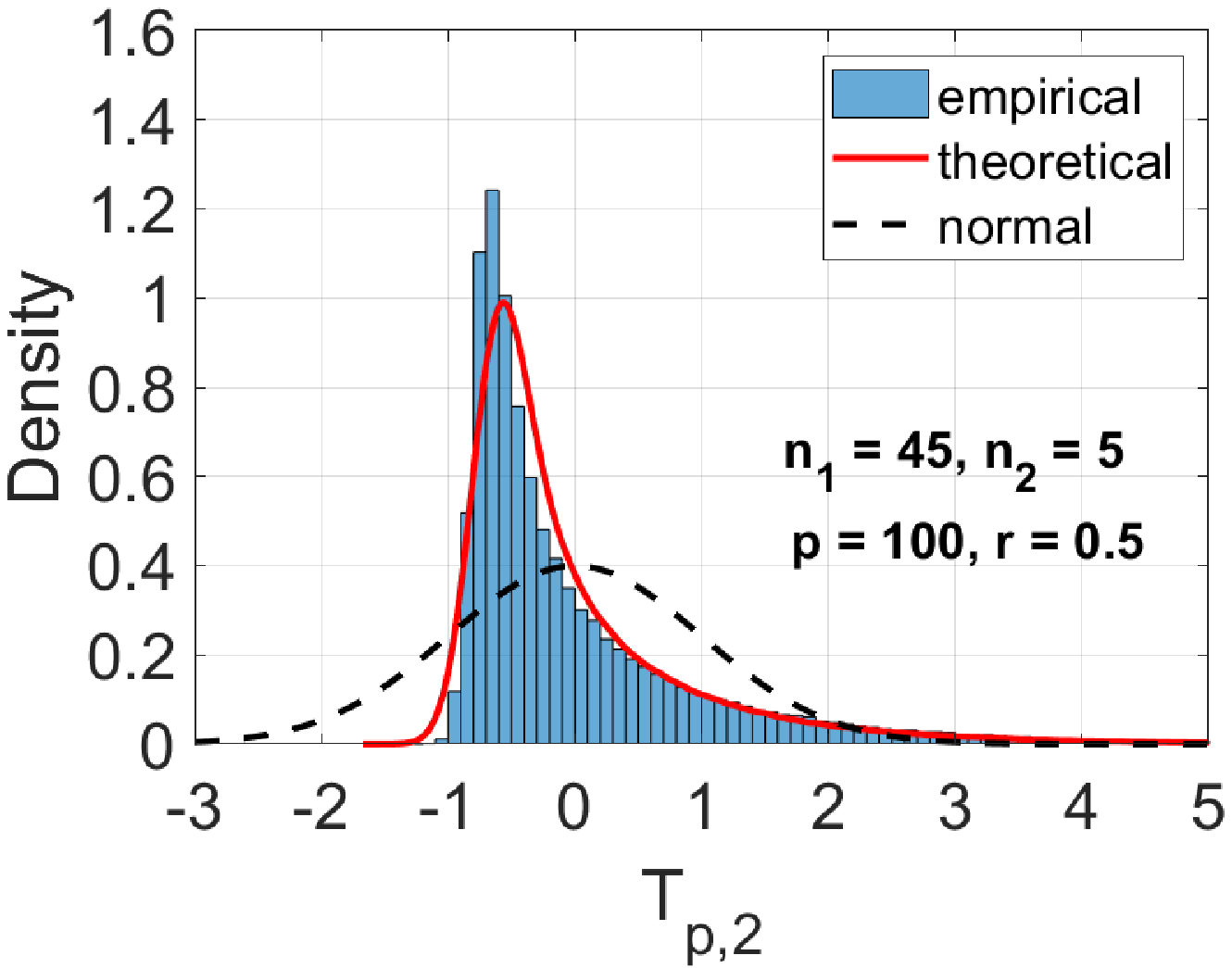}
\includegraphics[width=2.2in]{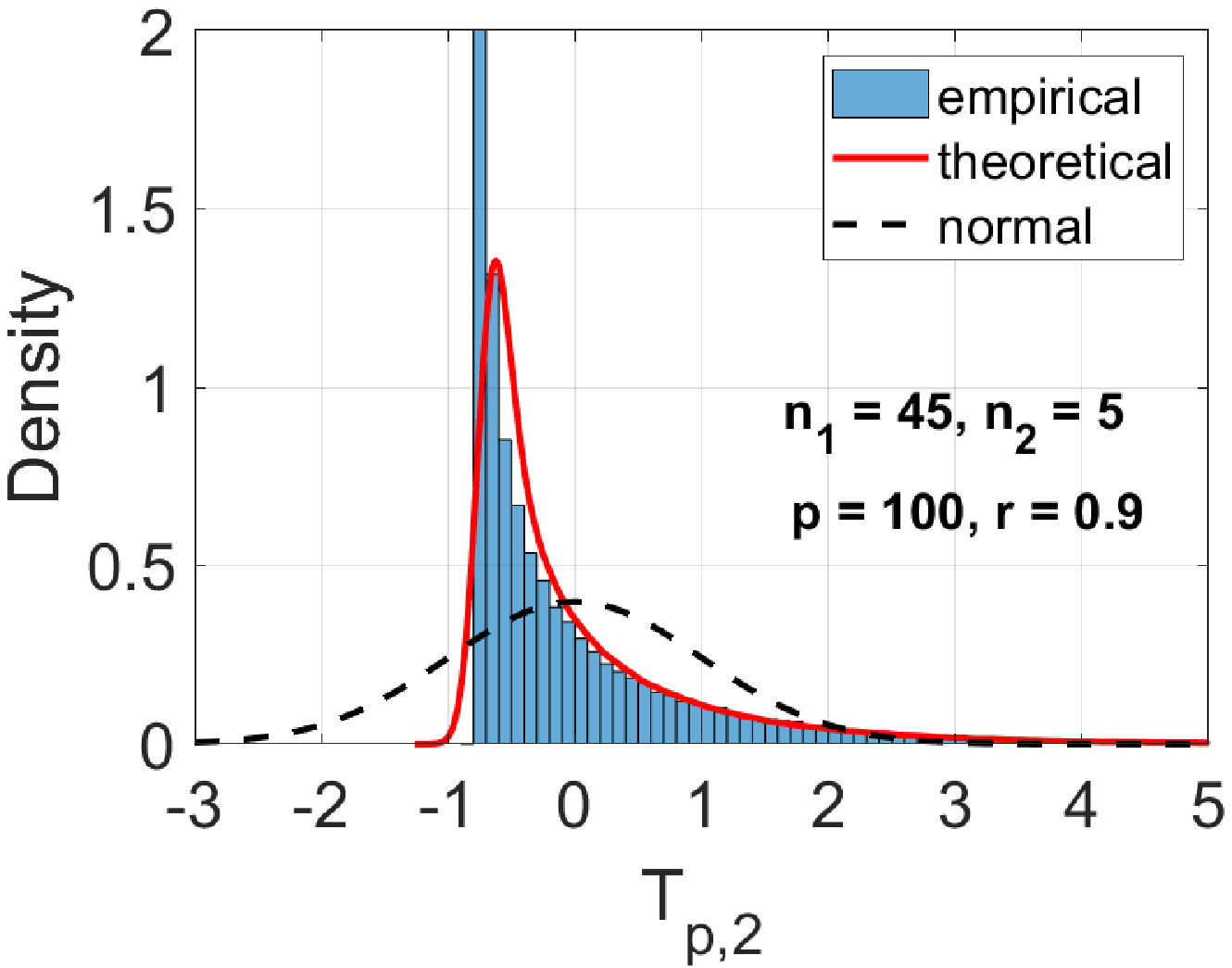}
}

\mbox{
\includegraphics[width=2.2in]{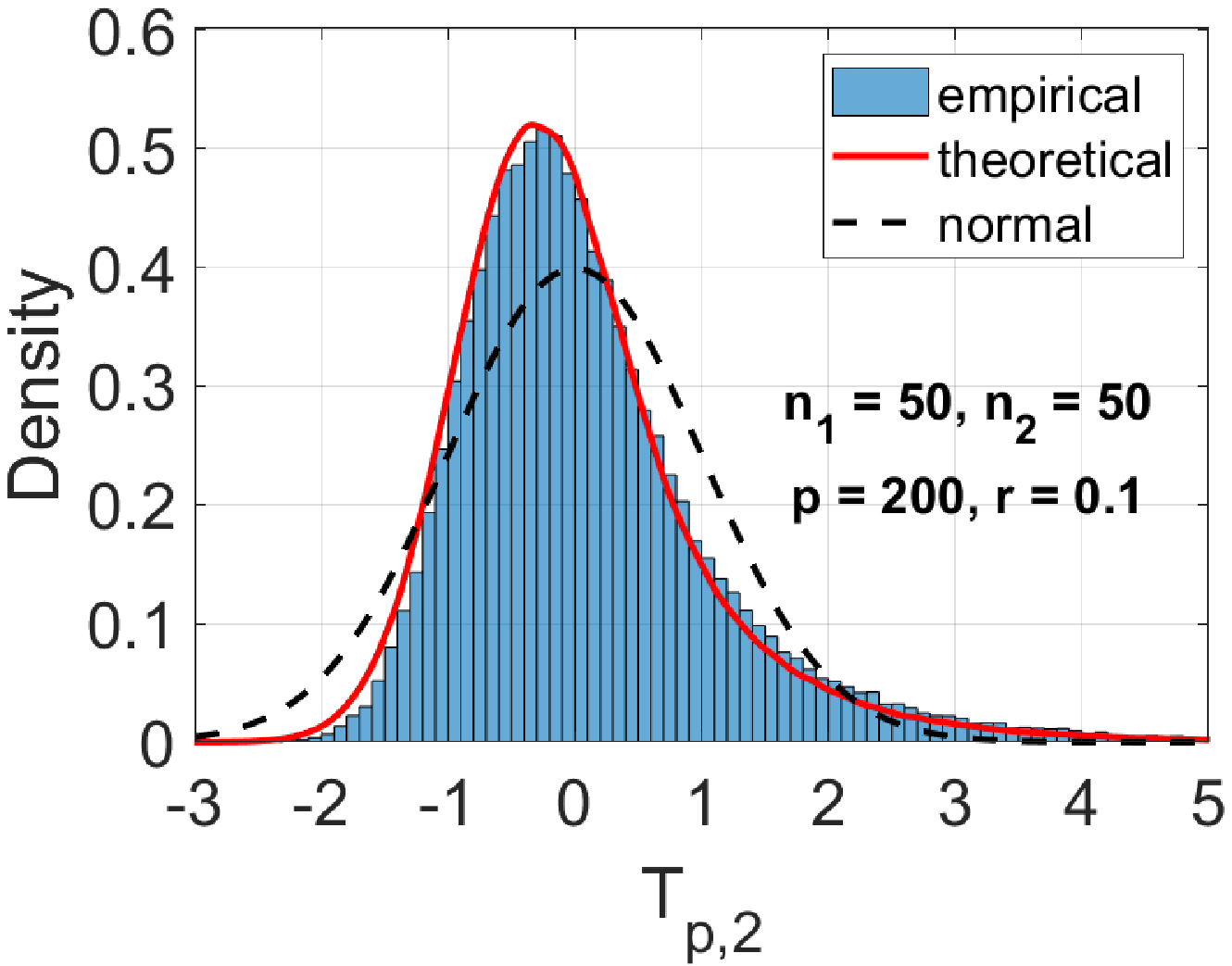}
\includegraphics[width=2.2in]{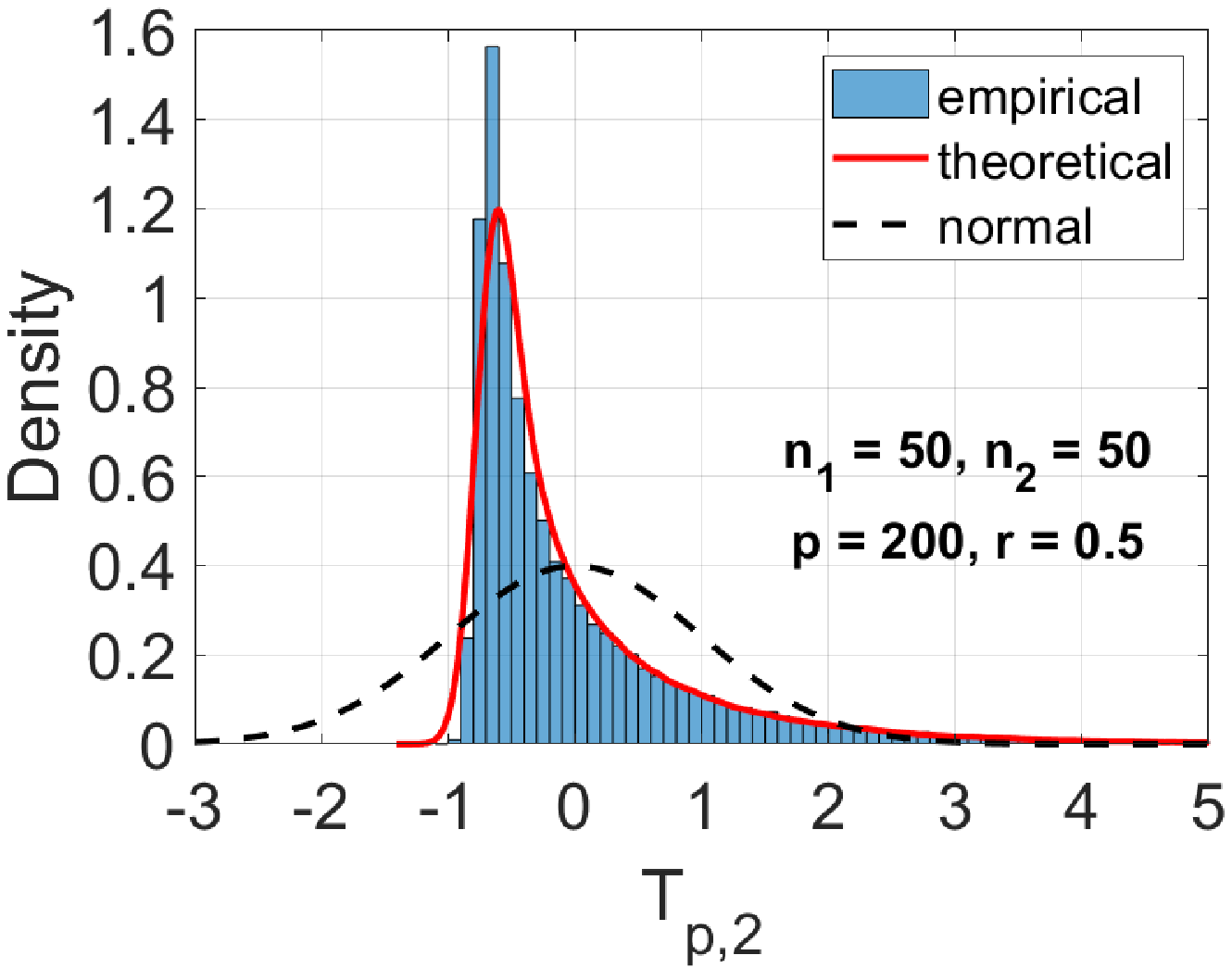}
\includegraphics[width=2.2in]{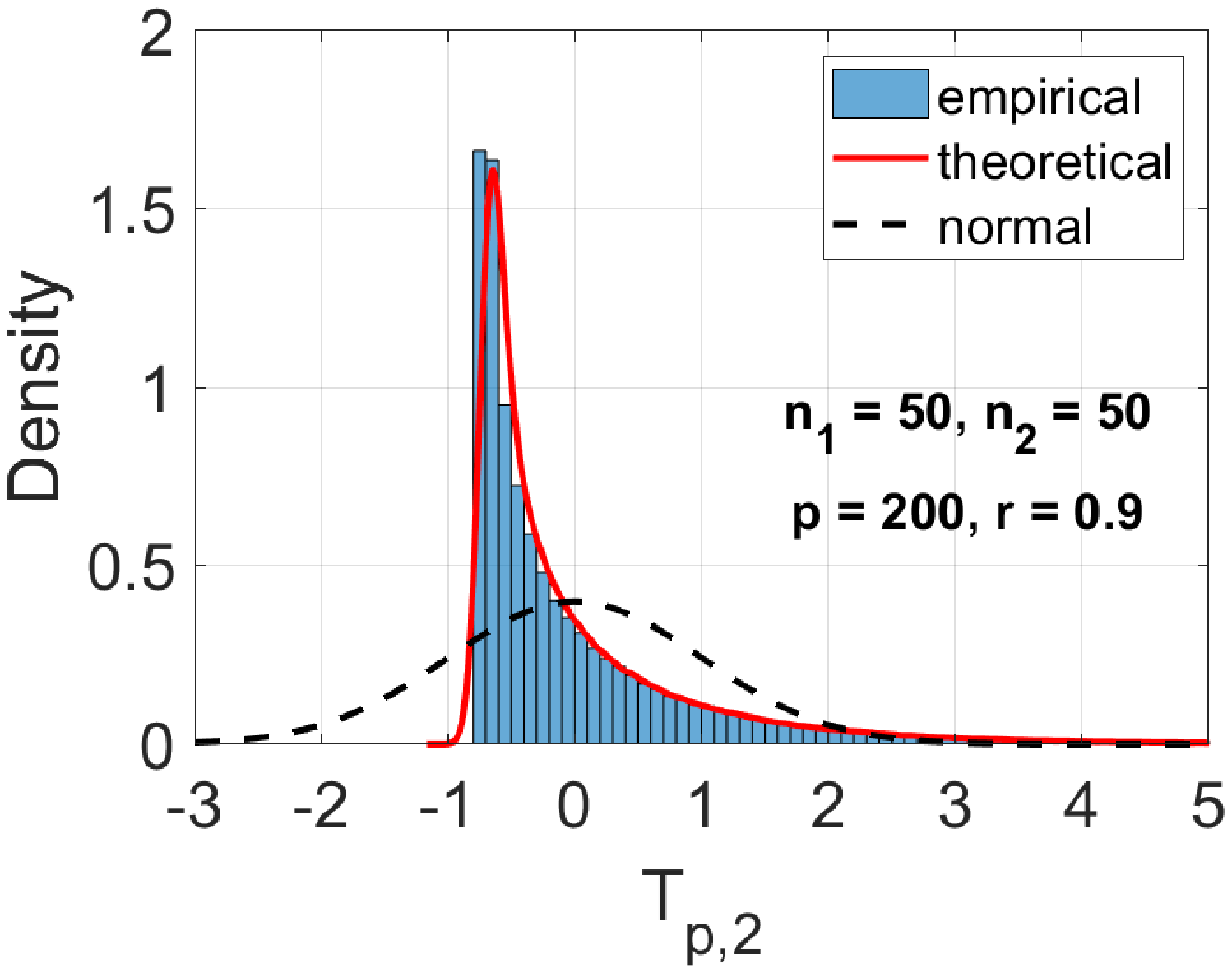}
}

\mbox{
\includegraphics[width=2.2in]{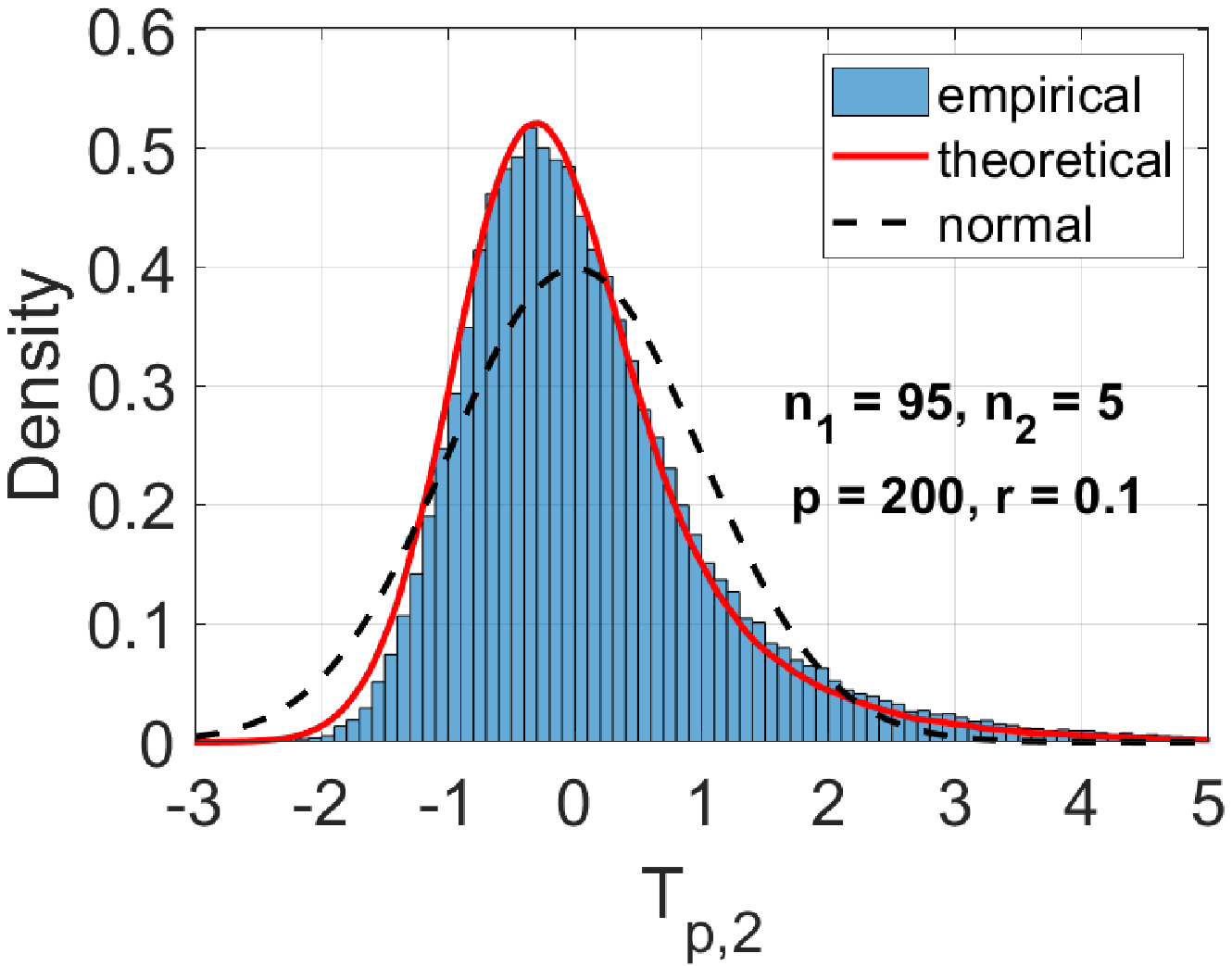}
\includegraphics[width=2.2in]{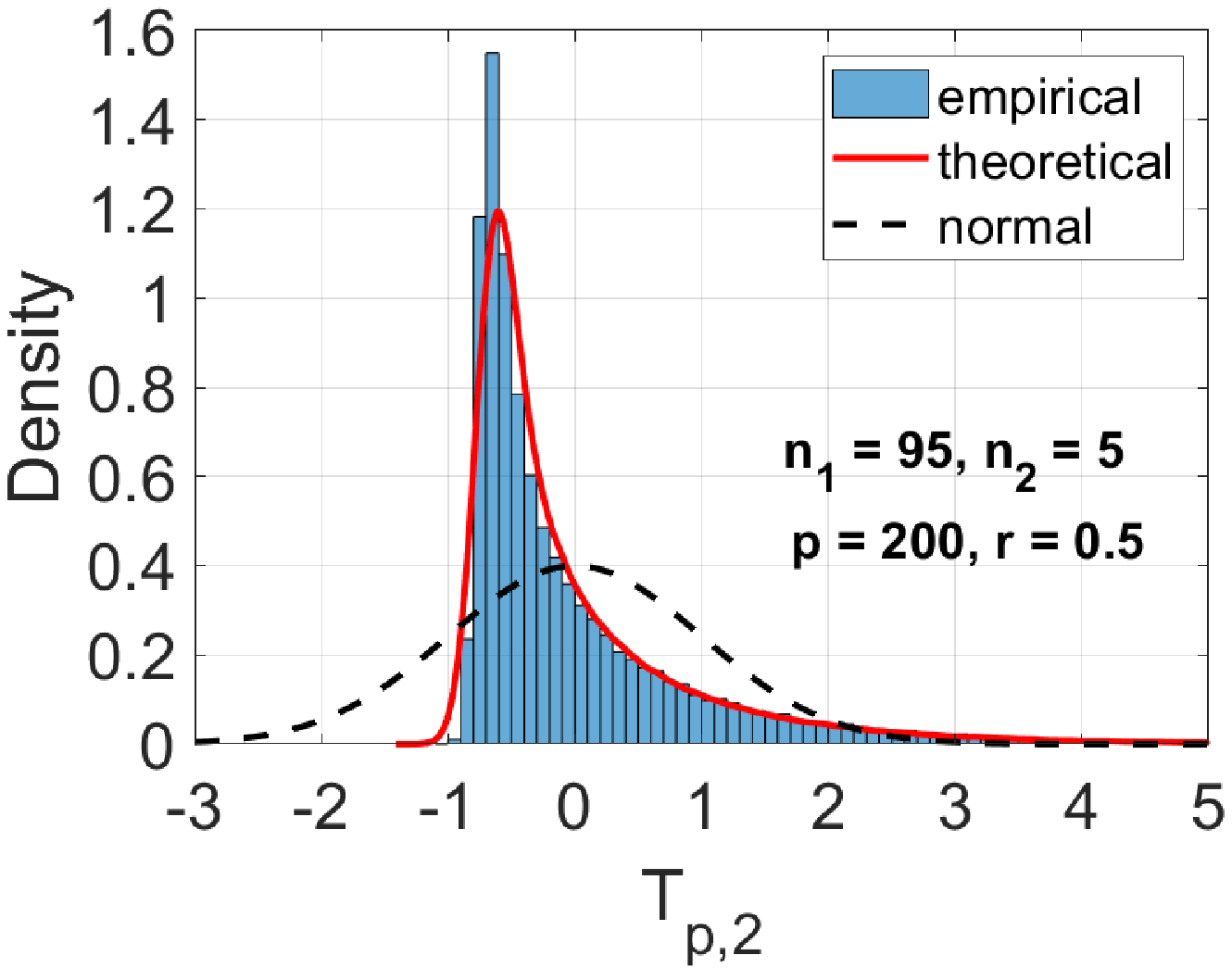}
\includegraphics[width=2.2in]{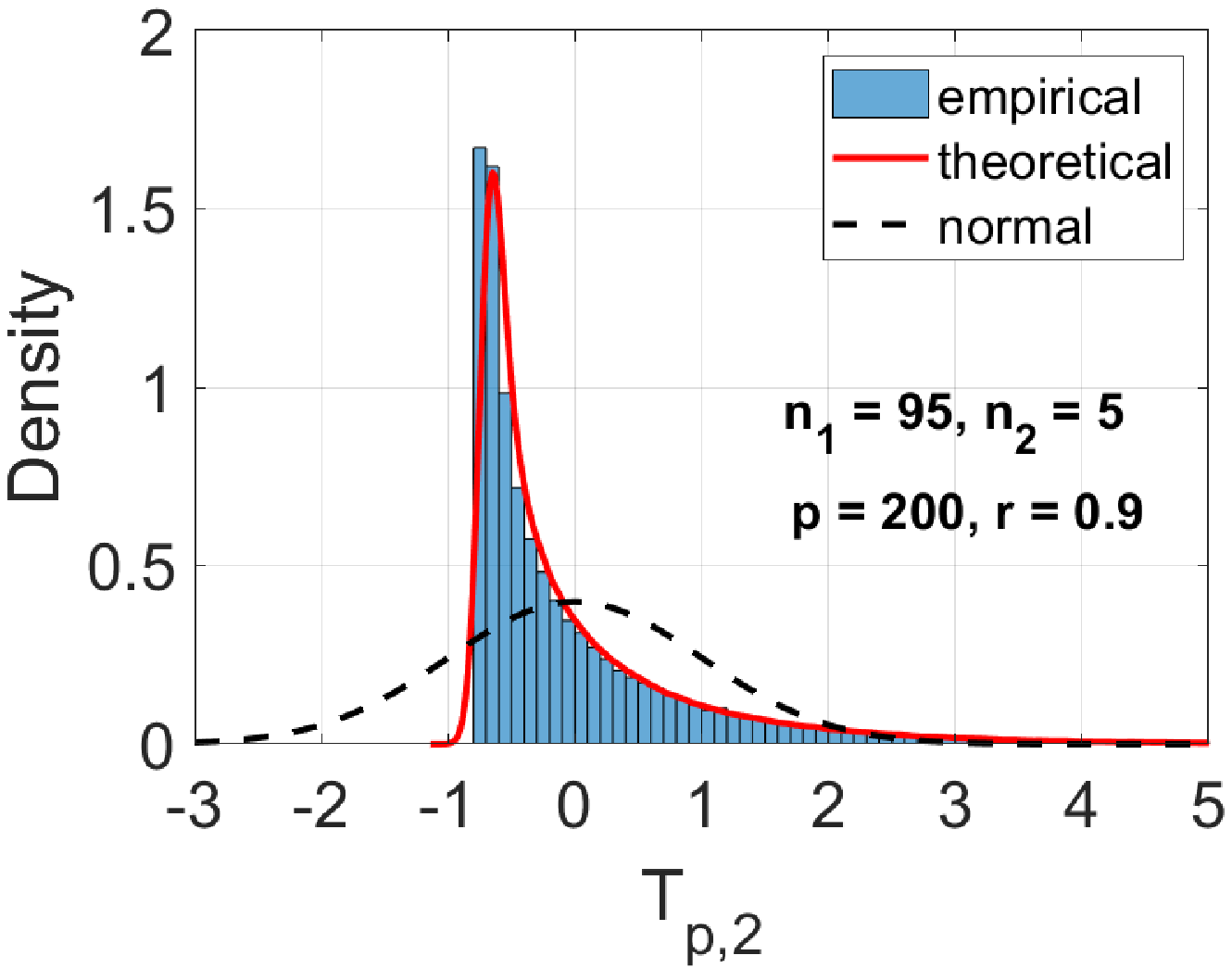}
}

    \caption{\small Example~\ref{Ex4}. Our theoretical curves always match the empirical ones well no matter the former is a normal-like or chi-square-like curve, even one of $n_1$ and $n_2$ being small. The normal curves by \cite{srivastava2008test} are farther from the empirical ones as correlation $r$ increases.}
    \label{fig:Ex2}
\end{figure}

\newpage

Similar to Example~\ref{Ex4}, modifications of Examples~\ref{Ex2} and~\ref{Ex3} can also be constructed for the two-sample case. It is straightforward, so we omit the detail. Given a data set, the spikes of the population correlation $\bd{R}$ have to be determined; see \cite{fan2020estimating} and \cite{morales2021asymptotics}.

Based on Example \ref{Ex4}, the statistic $T_{p,2}$ from \eqref{Statistics2} is simulated for $10^5$ times for each set of values of  $(n_1, n_2, p, r)$. The value of $r$ is chosen at $0,1, 0,5$ and $0.9$, designed for weakly dependence, dependence and strong dependence, respectively. Our theoretical curves are normal-like, mixing between $N(0,1)$ and chi-squared distribution and chi-squared-like curves. For weakly dependent case, by examining the four graphs in the first column, we see the normal approximation by \cite{srivastava2008test} becomes less and less accurate as all of the values of $n_1, n_2, p$ increase. However, our theoretical curve always match the empirical ones well. For dependent and very dependent cases appeared in the second and third columns, respectively, the normal curves are nowhere close to the empirical ones. It is good to see our theoretical curves are always close to their empirical counterparts.

\section{A Law of Large Numbers for Sample Correlation Matrices}\lbl{Sec_wlln}

Let $\X_1,\cdots,\X_n$ be a random sample from $N_p(\bmu, \bms)$ with correlation matrix $\bd{R}$. Let $\hat{\bd{R}}$ be the sample correlation matrix defined in \eqref{sample_corr_ma}. One of the crucial steps in proving  Theorems~\ref{Theorem1} and~\ref{Theorem2} is the use of an  asymptotically ratio-unbiased-estimator of $\mbox{tr}(\bd{R}^2)$. In their Lemma 3.2, \cite{srivastava2008test} state the following.

\begin{result}\lbl{AD2008Unbiased} If $n=O(p^{\zeta})$, $0<\zeta\leq 1$, under the condition \eqref{SD08a},
$(\mbox{tr}(\bd{R}^2)-p^2/n)/p$ converges to $\lim_{p\to\infty}\mbox{tr}(\bd{R}^2)/p$ in probability as $(n, p)\to\infty$ and thus can be considered as a consistent estimator of $\mbox{tr}(\bd{R}^2)/p$ as $(n, p)\to \infty$.
\end{result}
The proof of the above result is give on pages 400-402 from \cite{srivastava2008test}. It lacks a mathematical rigor.
By using a method from the random matrix theory, we rigorously obtain a more general result than the above in the following. Most efforts in the whole  proofs of this paper are devoted to this result, which is rather involved. The reason is that we assume no independence at all among the entries of the population distribution, i.e., $\bd{R}$ is arbitrary.

\begin{theorem}\lbl{Theorem3} Let $\X_1,\cdots,\X_n$ be a random sample from $N_p(\bmu, \bms)$ with correlation matrix $\bd{R}$. Let $\hat{\bd{R}}$ be the sample correlation matrix defined in \eqref{sample_corr_ma}. If $\lim_{p\to\infty}\frac{p}{n\|\bd{R}\|_F}= 0$ and $\lim_{p\to\infty}\frac{p}{n^{a}}=0$ for some constant $a>0$,  then
\beaa
\frac{1}{\mbox{tr}(\bd{R}^2)}\Big[\mbox{tr}(\hat{\bd{R}}^2)-\frac{p(p-1)}{n-1}
\Big] \to 1
\eeaa
in probability as $p\to\infty.$
\end{theorem}
Unlike $\mbox{tr}(\bd{R}^2)$, $\mbox{tr}(\bd{R}^3)$ and $\mbox{tr}(\bd{R}^4)$ needed in {\it Result~\ref{AD2008Unbiased}}, our conditions in the above theorem are imposed on  $\|\bd{R}\|_F=[\mbox{tr}(\bd{R}^2)]^{1/2}$ only.   Also, it is easy to see the difference between the assumptions on restrictions of $n$ and $p$:  {\it Result~\ref{AD2008Unbiased}} asks $n=O(p^{\zeta})$ for some $0<\zeta\leq 1$ and our Theorem~\ref{Theorem3} holds as long as $n$ and $p$ are in the order of a polynomial. We make some remarks on Theorem~\ref{Theorem3} next.

Naively,  an obvious ratio-unbiased-estimator of $\mbox{tr}(\bd{R}^2)$ is its sample version $\mbox{tr}(\hat{\bd{R}}^2)$. Theorem~\ref{Theorem3} indicates that it is not true for the high dimensional setting. Instead,  the modified version  ``$\mbox{tr}(\hat{\bd{R}}^2)-\frac{p(p-1)}{n-1}$" is a correct one. By using the fact $\mbox{tr}(\hat{\bd{R}}^2)\geq p$, we can see from Theorem~\ref{Theorem3} that $\mbox{tr}(\hat{\bd{R}}^2)$ is indeed a ratio-unbiased-estimator of $\mbox{tr}(\bd{R}^2)$ if $p$ is as small as  $p=o(n)$, that is, $\mbox{tr}(\hat{\bd{R}}^2)/\mbox{tr}(\bd{R}^2) \to 1$ in probability.

Sample correlation matrices are a special type of random matrices. The spectral distributions of  eigenvalue values are studied in  \cite{jiang2004limiting} and \cite{xiao2010almost}. The central limit theorems for  determinants under independent and correlated situations are obtained by \cite{jiang2013central} and \cite{jiang2009variance}, respectively. The Tracy-Widom law is derived by \cite{bao2012tracy}. The central limit theorem for linear statistics of eigenvalues is understood by \cite{gao2017high}.

An interesting fact is the derivation of Theorem~\ref{Theorem3}. Through the whole proof, contrary to standard techniques to handle randomness, we do not use/assume any independence from the population distribution $N_p(\bmu, \bms)$. In fact, in the most dependent case, that is, all of the $p$ entries of $N_p(\bmu, \bms)$ are identical, then $\hat{\bd{R}}=\bd{R}$ and all of their $p^2$ entries are equal to $1$, so  $\mbox{tr}(\hat{\bd{R}}^2)=\mbox{tr}(\bd{R}^2)=p^2$. One can see Theorem~\ref{Theorem3} trivially holds. On the other hand, for the most independent case, i.e., if $\bms=\bd{R}=\bd{I}_p$, then $\mbox{tr}(\bd{R}^2)=p.$ Assuming $p/n\to \rho\in (0,1]$, we have from Theorem~\ref{Theorem3} that $\mbox{tr}(\hat{\bd{R}}^2)/p\to 1+\rho.$ This actually can be  confirmed independently by a random matrix theory. In fact, let $\mu_p$ be the empirical distribution of the eigenvalues of  $\hat{\bd{R}}$. By Theorem 2 from \cite{jiang2004limiting}, $\mu_p$ converges weakly to the so-called Mar\v{c}henko-Pastur law $\mu$. The second moment of $\mu$ is equal to $1+\rho$ [Lemma 3.1 from \cite{bai2010spectral}]. Consequently, $\mbox{tr}(\hat{\bd{R}}^2)/p=\int x^2\,d\mu_p\to \int x^2\,d\mu=1+\rho.$

The proof of Theorem~\ref{Theorem3} is a bit technical and lengthy. We first give an accurate  estimate of $E\mbox{tr}(\hat{\bd{R}}^2)$ in Lemma~\ref{sdu329}. Then we only need to show that its variance go to zero. To do so, since $\hat{\bd{R}}^2=p+2\sum_{1\leq i< j \leq p}\hat{r}_{ij}^2$,  we need to understand the covariance between two sample correlations, say, $\hat{r}_{ij}^2$ and $\hat{r}_{kl}^2$ for any $1\leq i, j, k, l\leq p$. The precise result is given next.

\begin{theorem}\lbl{Theorem4} Assume $m\geq 5$. Let $\{(X_{1j}, X_{2j}, X_{3j}, X_{4j})^T \in \mathbb{R}^4;\, 1\leq j \leq m\}$ be i.i.d.  random vectors with distribution $N_4(\bd{0}, \bd{R})$, where $\bd{R}=(r_{ij})_{4\times 4}$ and $r_{ii}=1$ for each $i$. Set
\bea\lbl{sdi3208}
\hat{r}_{ij}=\frac{\sum_{k=1}^mX_{ik}X_{jk}}{(\sum_{k=1}^mX_{ik}^2)^{1/2}(\sum_{k=1}^mX_{jk}^2)^{1/2}}
\eea
for $1\leq i\leq j \leq 4.$  Then, for any $N\geq 1$,  $\mbox{Cov}(\hat{r}_{12}^2, \hat{r}_{34}^2)$ is equal to
\beaa
\varrho_{m,1}\cdot \sum_{1\leq i<j\leq 4}r_{ij}^2 + \varrho_{m,2}\cdot r_{12}^2r_{34}^2 &+& \varrho_{m,3}\cdot\big(r_{12}r_{23}r_{34}r_{41}+r_{12}r_{24}r_{43}r_{31}\big)\\
&+&\varrho_{m,4}\cdot\big(r_{13}r_{34}r_{41}+
r_{23}r_{34}r_{42}\big)r_{12}^2\\
&+ & \varrho_{m,5}\cdot\big(r_{12}r_{23}r_{31}+r_{12}r_{24}r_{41}\big)r_{34}^2\\
&+& \frac{\varrho_{m,6}}{m^{(N+1)/2}},
\eeaa
where $\{\varrho_{m,i};\, 3\leq i \leq 5\}$ are quantities not depending on $\bd{R}$,
\beaa
|\varrho_{m,1}| \leq Km^{-2}, \ \ |\varrho_{m,2}|\vee|\varrho_{m,3}|\vee|\varrho_{m,4}|\vee|\varrho_{m,5}| \leq K m^{-1},\ |\varrho_{m,6}| \leq K
\eeaa
and $K$ is a constant depending on $N$ but not on $m$ or $\bd{R}.$
\end{theorem}
The proof of Theorem~\ref{Theorem4} is rather involved. There are two reasons, one of them is that no assumption is imposed on $\bd{R}$, as a consequence, we are not able to use any techniques/theory related to independence. The second one is that sample correlation coefficients are more difficult to be handled than sample covariances. Look at \eqref{sdi3208} for $\hat{r}_{12}$ by taking $i=1$ and $j=2$. From the law of large numbers, $(1/m)\sum_{k=1}^mX_{1k}^2 \to 1$ and $(1/m)\sum_{k=1}^mX_{2k}^2 \to 1$. To understand $\hat{r}_{12}$, a naive idea is to replace the denominator in the expression of $\hat{r}_{12}$ by $m$. Interestingly, this works for the derivation of the Mar\v{c}henko-Pastur law (Jiang, 2004) when $\bd{R}=\bd{I}$. Now we elaborate this more for a further discussion. By the Taylor expansion, write
\bea\lbl{Taylor_ex}
\frac{1}{\sqrt{1+x}}=1-\frac{x}{2}+\frac{3}{8}x^2+\cdots + a_Nx^N+o(x^{N+1})
\eea
as $x$ is small. \cite{jiang2004limiting} uses the above expansion with $N=0$ and $x$ taking $(1/m)\sum_{k=1}^m(X_{1k}^2- 1)$ and $(1/m)\sum_{k=1}^m(X_{2k}^2-1)$, respectively. If $\bd{R}\ne\bd{I}$ but $\bd{R}=\bd{A}_p$ as in Example~\ref{Ex1}, \cite{fan2019largest} use the expansion \eqref{Taylor_ex} with $N=1$.  In the proof our Theorem~\ref{Theorem3}, since no structure of $\bd{R}$ is assumed, we have to use \eqref{Taylor_ex} for an arbitrary $N$. This is why the condition ``$\limsup_{p\to\infty}\frac{p}{n^{a}}=0$ for some constant $a>0$" is imposed in Theorems~\ref{Theorem1} and~\ref{Theorem3}. The new method here would be useful in the future for the study on sample correlation matrices. It is worthwhile to mention that handling sample correlation matrices is much more difficult than working on sample covariance matrices. The difference is obvious: the entries of a sample covariance matrix does not have the denominator as in the expression of $\hat{r}_{ij}$ from \eqref{sdi3208}. For instance, \cite{bai1996effect} also investigate a ratio-unbiased-estimator similar to Theorem~\ref{Theorem3} but for covariance matrices,  their argument is quick and short.

A byproduct of the proofs of Theorems~\ref{Theorem1} and~\ref{Theorem2} is the behavior of a quadratic form related to the diagonal matrix of $\bd{R}$. It would be useful for other research from this point forward. Assume $\X_1,\cdots,\X_n$ is a random sample from $N_p(\bmu, \bms)$ with correlation matrix $\bd{R}$. Review  $\D$ is the diagonal matrix of $\bms$ and  $\hat{\D}$ is the diagonal matrix of $\hat{\S}$ in \eqref{duck_yazi}.

\begin{lemma}\lbl{Jiaodong} Let $\bm{\eta}\sim N_p(\bd{0}, \bms)$ and $\bm{\eta}$ be independent of $\hat{\S}$. If $\lim_{p\to\infty}\frac{p}{n\|\bd{R}\|_F}=0$, then
\beaa
\frac{\bm{\eta}^T\hat{\D}^{-1}\bm{\eta}}{\sqrt{2\,\mbox{tr}(\bd{R}^2)}}
=\frac{\bm{\eta}^T\D^{-1}\bm{\eta}}{\sqrt{2\,\mbox{tr}(\bd{R}^2)}} + o_p(1)
\eeaa
as $p\to\infty.$
\end{lemma}
Since $\X_1,\cdots,\X_n$ is a random sample from a multivariate normal distribution,  $\bar{\X}$ is independent of $\hat{\S}$, and hence $\bar{\X}$ is independent of $\hat{\D}$. So the above conclusion holds if we take $\bm{\eta}=\sqrt{p}\bar{\X}$. This is actually the way we use this lemma in later proofs.

The limiting distributions in Theorems~\ref{Theorem1} and~\ref{Theorem2} essentially come from the following.

\begin{lemma}\lbl{vdu349j} For each $p\geq 1$, let $a_{p,1}\geq \cdots\geq a_{p, p}\geq 0$
 be constants with $a_{p,1}^2+ \cdots+ a_{p, p}^2=1$. Suppose $\lim_{p\to\infty}a_{p,i}=\rho_i\geq 0$ for each $i\geq 1$. Let $\xi_1, \xi_2, \cdots$ be i.i.d. with distribution $\chi^2(1)-1.$ Then $a_{p,1}\xi_1+\cdots + a_{p,p}\xi_p$ converges to  $[2(1-\sum_{i=1}^{\infty}\rho_i^2)]^{1/2}\eta+\sum_{i=1}^{\infty}\rho_i\xi_i$ in distribution as $p\to \infty$, where $\eta\sim N(0, 1)$ and $\eta$ is independent of $\{\xi_i;\, i\geq 1\}$.
\end{lemma}

The comment after Example~\ref{Ex3} is based on the following fact. It is interesting useful in its own right.

\begin{lemma}\lbl{dsih328} Let $\bd{M}$ be a $p\times p$ correlation matrix, namely, $\bd{M}$ is non-negative definite  and its diagonal entries are all equal to $1$. Suppose $\bd{M}$ has  eigenvalues $\lambda_1\geq \cdots\geq \lambda_p$. Then $\lambda_1+\cdots+ \lambda_k\geq k$ for $k=1, \cdots, p$ and $\lambda_1+ \cdots +\lambda_p=p$. Conversely, for any $\tau_1\geq \cdots \geq \tau_p\geq 0$ with $\tau_1+ \cdots+ \tau_k\geq k$ for each $k=1, \cdots, p$ and $\tau_1+ \cdots +\tau_p=p$, there always exists a correlation matrix with eigenvalues $\tau_1, \cdots, \tau_p$.
\end{lemma}

\section{Concluding Remarks and Discussion}\lbl{Sec_concluding}

1. In this paper we have studied one-sample and two-sample mean tests. For the multiple  population case, it becomes the classical MANOVA problem. In the ``large $p$, small $n$" situation, \cite{srivastava2006multivariate} consider the question by revising the classical $F$-test. They use functions of non-zero eigenvalues of pseudo-MANOVA random matrix as test statistics.  \cite{zhang2017linear} and \cite{chen2019two} study this case based on the $L_2$-type statistics. An extension of our work to the MANOVA case under ``large $p$, small $n$" situation is meaningful. One direction is to modify the classical tests, such as the Anderson test, the Pillal test, the Roy test, the Wilks test and the Olson test, by a method similar to those from \cite{srivastava2006multivariate}. The classical statistics take the correlations among population coordinates into account. It has kind of self-normalization, which is favorable.

2. All results in this paper are based on random samples from a multivariate Gaussian distribution. Is it possible to generalize this to non-Gaussian cases? The corresponding results are highly demanded simply because many data are not sampled from a Gaussian distribution. Also, conclusions on non-Gaussian scenario may help us understand robustness of our tests.

3. We study the two-sample mean test by assuming the two populations have the same covariance matrices. If the two covariance matrices are not identical, this is the multivariate Behrens-Fisher problem. Various methods are proposed, for example, by \cite{yao1965approximate}, \cite{johansen1980welch}, \cite{nel1986solution}, \cite{krishnamoorthy2004modified}, \cite{chen2010two} and \cite{chen2019two}. For one-dimensional case ($p=1$), the method initiated by \cite{welch1947generalization} is probably the most popular one. The author basically normalize the difference of sample means by its sample standard deviation. \cite{srivastava2013two} propose a statistic by replacing $\hat{\bd{D}}$ from \eqref{Statistics2} with a hybrid of two samples. Similar to the study in this paper, its properties are needed to be understood, too.

4. As far as proofs go, we spend most of our energy proving Theorem~\ref{Theorem4}. We do not use complicated technology. Instead, by we use brute force to compute mixed moments of multiple Gaussian random variables in combination with some random matrix theory. As a result, the argument is lengthy. It is possible to shorten the proof. To get the law of large numbers stated in Theorem~\ref{Theorem4}, one may like to try Gaussian concentration inequalities; see, for example, \cite{ledoux2001concentration} and \cite{boucheron2013concentration}. An alternative way is to get the joint density function of $\hat{r}_{12}$ and $\hat{r}_{34}$ similar to that of the marginal density of  $\hat{r}_{12}$ treated at \eqref{wifh} through an hypergeometric function.

5. Theorems~\ref{Theorem1} and~\ref{Theorem2} present the null distributions for the one-sample and two-sample tests, respectively. We actually have tried to derive the non-null limiting distribution to explore power functions of the tests. The argument is also very involved. The current paper is already very lengthy, so we postpone and leave it as a future work.

\section*{Acknowledgment}

The bulk of the work of Tiefeng Jiang was conducted during his visit at Baidu Research in 2019-2020. After completing an initial draft, we solicited comments from several colleagues and we appreciate their feedback.

\section{Proofs}\lbl{Sec_proofs}

One of the main steps of proving Theorems~\ref{Theorem1} and~\ref{Theorem2} is a weak law of large numbers for $\mbox{tr}(\hat{\bd{R}}^2)$ (Theorem~\ref{Theorem3}), where $\hat{\bd{R}}=\hat{\bd{R}}_p$ is the sample correlation matrix defined in \eqref{sample_corr_ma}. To derive this, we need to study the mean and variance of $\mbox{tr}(\hat{\bd{R}}^2)$. In Section~\ref{mean_study}, we will get an accurate estimate of $E[\mbox{tr}(\hat{\bd{R}}^2)]$ by an argument on hypergeometric functions. Then we will show $\mbox{Var}(\mbox{tr}(\hat{\bd{R}}^2))$ go to zero. The proof of Theorem~\ref{Theorem3} will be completed once this step is established. This is a rather involved step and it will be understood step by step in Sections~\ref{Variance_Study}-\ref{TPTRSCM}.

\subsection{Evaluation of the Mean of the Frobenius Norm of a Sample Correlation matrix}\lbl{mean_study}

In this section, we will work on the mean of $\mbox{tr}(\hat{\bd{R}}^2)$, where $\hat{\bd{R}}=\hat{\bd{R}}_p$ is the sample correlation matrix defined in \eqref{sample_corr_ma}. The critical tool is the hypergeometric function. We first need a preliminary result as follows.
\begin{lemma}\lbl{fwhuwefu} There exists a constant $C>0$ such that the following holds for all $k\geq 2$.
\beaa
\sum_{n=1}^{\infty}\binom{n+k}{k}^{-1}\leq Ck^{-1/4}.
\eeaa
\end{lemma}
\noindent\textbf{Proof of Lemma~\ref{fwhuwefu}}. For given $a>0$, set $f_a(x)=(1+ax^{-1})^{-x}$ for $x>0$. We claim that $f_a(x)$ is decreasing in $x\in(0, \infty).$  In fact, let $g(x)=\log f_a(x)$. Then
\beaa
g'(x)=\Big[-x\log \Big(1+\frac{a}{x}\Big)\Big]'=-\log \Big(1+\frac{a}{x}\Big)-x\frac{-\frac{a}{x^2}}{1+\frac{a}{x}}.
\eeaa
Then
\beaa
g'(x)=\log \Big(1-\frac{a}{x+a}\Big) +\frac{a}{x+a}\leq 0
\eeaa
for all $x>0$ since $\log (1+y)\leq y$ for all $y>-1$. Hence, $f_a(x)$ is decreasing in $x\in (0, \infty)$ for each $a>0.$ By the Stirling formula, $m!=\sqrt{2\pi m}m^me^{-m+\frac{\theta_m}{12m}}$ with $\theta_m\in (0,1)$ for each $m\geq 1$. Consequently,
\beaa
\sqrt{2\pi m}m^me^{-m}\leq m! \leq e\cdot \sqrt{2\pi m}m^me^{-m}
\eeaa
for all $m\geq 1.$ Therefore,
\beaa
\binom{n+k}{k}=\frac{(n+k)!}{n!k!}
&\geq & e^{-2}\cdot \frac{\sqrt{n+k}}{\sqrt{2\pi nk}}\cdot \frac{(n+k)^{n+k}e^{-(n+k)}}{k^kn^ne^{-(n+k)}}\\
&\geq&
 e^{-2}\frac{1}{\sqrt{2\pi k}}\cdot\frac{(n+k)^{n+k}}{k^kn^n}.
\eeaa
This implies that
\bea\lbl{kdify}
\frac{1}{\binom{n+k}{k}}\leq (2\pi e^2 k^{1/2})\Big(1+\frac{k}{n}\Big)^{-n}\Big(1+\frac{n}{k}\Big)^{-k}.
\eea
For fixed $k$, since $f_k(x)$ is decreasing in $x$, we have
\bea\lbl{fwefwu}
\sum_{n=1}^k\Big(1+\frac{k}{n}\Big)^{-n}\Big(1+\frac{n}{k}\Big)^{-k}
&\leq & \sum_{n=1}^k\Big(1+\frac{k}{n}\Big)^{-n}\nonumber\\
&\leq & \frac{1}{k+1}+\sum_{n=2}^k\int_{n-1}^{n}\Big(1+\frac{k}{x}\Big)^{-x}\,dx \nonumber\\
&=&\frac{1}{k+1}+\int_1^k\Big(1+\frac{k}{x}\Big)^{-x}\,dx.
\eea
Again, use the fact that $f_k(x)$ is decreasing to see
\beaa
\int_1^k\Big(1+\frac{k}{x}\Big)^{-x}\,dx
&=&\int_1^{k^{1/4}}\Big(1+\frac{k}{x}\Big)^{-x}\,dx
+\int_{k^{1/4}}^k\Big(1+\frac{k}{x}\Big)^{-x}\,dx\\
&\leq & \frac{k^{1/4}}{k+1}+ k\cdot (1+k^{3/4})^{-k^{1/4}}.
\eeaa
If $k\geq 81$, then $k^{1/4}\geq 3$. Consequently, $k\cdot (1+k^{3/4})^{-k^{1/4}}\leq k\cdot (k^{3/4})^{-3}=k^{-5/4}.$ From \eqref{fwefwu} we see
\bea\lbl{weuhefwo}
\sum_{n=1}^k\Big(1+\frac{k}{n}\Big)^{-n}\Big(1+\frac{n}{k}\Big)^{-k}\leq \frac{1}{k+1}+\frac{k^{1/4}}{k+1}+\frac{1}{k^{5/4}}\leq \frac{3}{k^{3/4}}.
\eea
Now, $(1+\frac{k}{n})^{-n}\leq 2^{-k}$ as $n\geq k$ by the fact $f_k(x)$ is decreasing. Also,  $(1+\frac{n}{k})^{-k}\leq (1+\frac{n}{k})^{-81}$ as $k\geq 81$. Hence
\bea\lbl{siu08}
\sum_{n=k}^{\infty}\Big(1+\frac{k}{n}\Big)^{-n}\Big(1+\frac{n}{k}\Big)^{-k}\leq 2^{-k}\sum_{n=k}^{\infty}\Big(1+\frac{n}{k}\Big)^{-81}\leq 2^{-k}k^{81}\sum_{n=1}^{\infty}\frac{1}{n^{81}}
\eea
as $k\geq 81.$ This joined with \eqref{kdify}, \eqref{weuhefwo} and \eqref{siu08} yields that
\beaa
\sum_{n=1}^{\infty}\frac{k^{1/4}}{\binom{n+k}{k}} \leq k^{1/4}\cdot (2\pi e^2k^{1/2})\cdot \big(3k^{-3/4}+ \zeta(81)2^{-k}k^{81} \big)\leq C'
\eeaa
for all $k\geq 81$, where $\zeta(81)=\sum_{n=1}^{\infty}\frac{1}{n^{81}}<\infty$ is the Riemann zeta function evaluated at $81$ and $C'$ is another numerical constant not depending on $k$. In summary,
\bea\lbl{huweio}
\sup_{k\geq 81}\sum_{n=1}^{\infty}k^{1/4}\binom{n+k}{k}^{-1}< \infty.
\eea
Note that
\beaa
\frac{\binom{n+i+1}{i+1}}{\binom{n+i}{i}}=\frac{n+i+1}{i+1}>1
\eeaa
for any $i\geq 1$ and $n\geq 1$. We know $\binom{n+i}{i}$ is increasing in $i$, and hence
\beaa
\sup_{2\leq k\leq 80}\sum_{n=1}^{\infty}k^{1/4}\binom{n+k}{k}^{-1}\leq 80^{1/4}\sum_{n=1}^{\infty}\binom{n+2}{2}^{-1}\leq 2\cdot 80^{1/4}\sum_{n=1}^{\infty}\frac{1}{n^2}<\infty.
\eeaa
This and \eqref{huweio} imply the conclusion. \hfill$\square$

\medskip

The following result quantifies the second moment of the sample correlation coefficient of a random sample of size $m$ up to an error $O(m^{-1/4}).$
\begin{lemma}\lbl{efwieii} Let $\{(X_i, Y_i)^T;\, 1\leq i \leq m\}$ be i.i.d. $2$-dimensional normal random vectors with $EX_1=EY_1=0$, $EX_1^2=EY_1^2=1$ and $\mbox{Cov}(X_1, Y_1)=r$. Set
\beaa
\hat{r}_m=\frac{\sum_{k=1}^mX_{i}Y_i}{(\sum_{k=1}^mX_i^2)^{1/2}(\sum_{k=1}^mY_i^2)^{1/2}}.
\eeaa
Write
\beaa
E\big(\hat{r}_m^2)=\frac{1}{m}+r^2+b_m(r)\cdot r^2
\eeaa
for $m\geq 4$. Then  $\sup_{m\geq 4,\, |r|\leq 1}|m^{1/4}b_m(r)|<\infty$.
\end{lemma}
\noindent\textbf{Proof of Lemma~\ref{efwieii}}. If $r=\pm 1$, since  $EX_1=EY_1=0$ and  $EX_1^2=EY_1^2=1$, we know $Y_i=\pm X_i$ for each $i$. By the definition of $\hat{r}_m$, trivially, $\hat{r}_m =\pm 1$. This implies that
\bea\lbl{wq09cn}
b_m(\pm 1)=-\frac{1}{m}
\eea
for each $m\geq 1.$ In the following we always assume $r^2<1.$
From \cite{ghosh1966asymptotic} or p. 156 in \cite{muirhead1982aspects},
\bea\lbl{wifh}
E\big(\hat{r}_m^2)=1-\frac{m-1}{m}(1-r^2)\cdot {}_2F_1\Big(1,1; \frac{1}{2}m+1; r^2\Big),
\eea
where ${}_2F_1(a, b; c; z)$ is an hypergeometric function defined by
\beaa
{}_2F_1(a, b; c; z)=\sum_{k=0}^{\infty}\frac{(a)_k(b)_k}{(c)_kk!}z^k,~~ |z|<1,
\eeaa
and where $(x)_0=1$ and $(x)_k=x(x+1)\cdots (x+k-1)$ for $k\geq 1$; see, for example, p. 20 from \cite{muirhead1982aspects}. Notice $(1)_k=k!$ and
\beaa
\Big(\frac{1}{2}m+1\Big)_k=\frac{m+2}{2}\cdot \frac{m+4}{2}\cdots \frac{m+2k}{2}.
\eeaa
Thus
\bea\lbl{feuofe45}
{}_2F_1\Big(1,1; \frac{1}{2}m+1; r^2\Big)
&=&1+\sum_{k=1}^{\infty}\frac{2^kk!}{(m+2)(m+4)\cdots (m+2k)}r^{2k} \nonumber\\
&=& 1+r^2\sum_{k=1}^{\infty}\frac{2^kk!}{(m+2)(m+4)\cdots (m+2k)}r^{2k-2}.
\eea
Evidently, the last sum is bounded by
\beaa
\sum_{k=1}^{\infty}\frac{2^kk!}{(m+2)(m+4)\cdots (m+2k)}.
\eeaa
Set $j=[\frac{m}{2}]$, where $[x]$ denotes the integer part of a real number $x\geq 0$. Then $j\leq \frac{m}{2}$, or equivalently,  $m\geq 2j$. It follows that
\beaa
\frac{2^kk!}{(m+2)(m+4)\cdots (m+2k)}\leq \frac{k!}{(j+1)(j+2)\cdots (j+k)}=\frac{j!k!}{(j+k)!}.
\eeaa
It follows that
\beaa
\sum_{k=1}^{\infty}\frac{2^kk!}{(m+2)(m+4)\cdots (m+2k)}\leq
\sum_{k=1}^{\infty}\frac{1}{\binom{j+k}{j}}\leq \frac{K_1}{j^{1/4}}\leq \frac{K_2}{m^{1/4}}
\eeaa
as $j\geq 2$, or equivalently, $m\geq 4$  by Lemma~\ref{fwhuwefu}, where $K_1$ and $K_2$ are constants not depending on $m$. From this and \eqref{feuofe45}, we are able to write
\beaa
{}_2F_1\Big(1,1; \frac{1}{2}m+1; r^2\Big)=1+a_mr^2,
\eeaa
where $0\leq a_m\leq K_2m^{-1/4}.$ Substitute this back to \eqref{wifh}, we get
\beaa
E\big(\hat{r}_m^2)
&=&1-\frac{m-1}{m}(1-r^2)\cdot(1+a_mr^2)\\
&= & 1-\frac{m-1}{m}(1-r^2)-\frac{m-1}{m}(1-r^2)\cdot a_mr^2\\
&=& \frac{1}{m}+r^2+b_m(r)r^2,
\eeaa
where
\beaa
b_m(r)=-\frac{1}{m}r^2-\frac{m-1}{m}(1-r^2)\cdot a_mr^2.
\eeaa
Obviously,
$\sup_{|r|\leq 1}|b_m(r)|\leq K_3m^{-1/4}$ for every $m\geq 4$, where $K_3$ is a constant not depending on $m$. This and \eqref{wq09cn} lead to the desired conclusion.  \hfill$\square$

\begin{lemma}\lbl{sdu329} Let $\X_1,\cdots,\X_n$ be a random sample from $N_p(\bmu, \bms)$ with correlation matrix $\bd{R}$. Let $\hat{\bd{R}}=\hat{\bd{R}}_p$ be the sample correlation matrix defined in \eqref{sample_corr_ma}. If $\lim_{p\to\infty}\frac{p}{n\|\bd{R}\|_F}= 0$ and $\limsup_{p\to\infty}\frac{p}{n^{a}}=0$ for some constant $a>0$,  then, as $p\to\infty$,
\beaa
E\,\mbox{tr}(\hat{\bd{R}}_p^2)
=\frac{p(p-1)}{n-1}+\mbox{tr}(\bd{R}_p^2)\cdot \big[1+O(m^{-1/4})\big].
\eeaa
\end{lemma}
\noindent\textbf{Proof of Lemma~\ref{sdu329}}. Set $m=n-1$. The notation $W_p(m, \bms)$ represents  the distribution of  the Wishart matrix $\bd{U}^T\bd{U}$, where $\bd{U}$ is an $m\times p$ matrix whose rows are i.i.d. with distribution $N_p(\bd{0}, \bms)$. Then, $n\hat{\S}$ has the Wishart distribution $W_p(m, \bms)$; see, for example, Theorem 3.1.2 from \cite{muirhead1982aspects}. That is, $n\hat{\S}\overset{d}{=} \bd{U}^T\bd{U}$. We claim
\bea\lbl{ifeie}
\hat{\bd{R}}_p=(\hat{r}_{ij})_{p\times p}\ \overset{d}{=}  \Big(\frac{\bd{v}_i^T\bd{v}_j}{\|\bd{v}_i^T\|\cdot\|\bd{v}_j\|}\Big)_{p\times p},
\eea
where the $m$ rows of $(\bd{v}_1, \cdots, \bd{v}_p)_{m\times p}$ are i.i.d. with distribution $N_p(\bd{0}, \bd{R})$.
In fact, write $\bd{U}=(u_{ij})=(\bd{u}_1, \cdots, \bd{u}_p)$ where $\bd{u}_i=(u_{1i}, \cdots, u_{mi})^T\in \mathbb{R}^m$ for each $i$. Then $\bd{u}_i^T\bd{u}_j$ is the $(i, j)$-entry of $\bd{U}^T\bd{U}$ and the diagonal entries are $\|\bd{u}_i\|^2$ for $1\leq i \leq p$. We know from~\eqref{sample_corr_ma} that
$\hat{\R}_p=\hat{\D}^{-1/2}\hat{\S}\hat{\D}^{-1/2}$,
where $\hat{\D}$ is the diagonal matrix of $\hat{\S}$. Then
\beaa
\hat{\bd{R}}_p=(\hat{r}_{ij})_{p\times p}\ \overset{d}{=}  \Big(\frac{\bd{u}_i^T\bd{u}_j}{\|\bd{u}_i^T\|\cdot\|\bd{u}_j\|}\Big)_{p\times p}.
\eeaa
Notice $(\bd{u}_i, \bd{u}_j)=(u_{ki}, u_{kj})_{1\leq k \leq m}$ for any $i<j$ and the $m$ rows are i.i.d. bivariate normal random variables with mean vector $\bd{0}$ and covariance matrix
\beaa
\begin{pmatrix}
\sigma_{ii}& \sigma_{ij}\\
\sigma_{ij} & \sigma_{jj}
\end{pmatrix}
\eeaa
where $\bms:=(\sigma_{ij})_{p\times p}$. Write $\bd{R}_p=(r_{ij})_{p\times p}$. Then  $r_{ij}=\sigma_{ij}(\sigma_{ii}\sigma_{jj})^{-1/2}$ by definition. Observe
\beaa
\frac{\bd{u}_i^T\bd{u}_j}{\|\bd{u}_i^T\|\cdot\|\bd{u}_j\|}=
\frac{(\sigma_{ii}^{-1/2}\bd{u}_i)^T(\sigma_{jj}^{-1/2}\bd{u}_j)}
{\|\sigma_{ii}^{-1/2}\bd{u}_i^T\|\cdot\|\sigma_{jj}^{-1/2}\bd{u}_j\|}.
\eeaa
Take $\bd{v}_i=\sigma_{ii}^{-1/2}\bd{u}_i$ to obtain \eqref{ifeie}.

By Lemma~\ref{efwieii},
\beaa
E\,\mbox{tr}(\hat{\bd{R}}_p^2)
&=&p+ 2\sum_{1\leq i<j \leq p}E\big(\hat{r}_{ij}^2\big)\\
&=&p+2\sum_{1\leq i<j \leq p}\Big[\frac{1}{m}+r_{ij}^2+b_m(r_{ij})r_{ij}^2\Big]\\
&=& \frac{p(p-1)}{m}+\mbox{tr}(\bd{R}_p^2) + 2\sum_{1\leq i<j \leq p}b_m(r_{ij})r_{ij}^2,
\eeaa
where  $\sup_{|r|\leq 1}|b_m(r)|\leq Km^{-1/4}$ and $K$ is a constant not depending on $m$ or $r$. Since
\beaa
2\Big|\sum_{1\leq i<j \leq p}b_m(r_{ij})r_{ij}^2\Big| \leq (2K)m^{-1/4}\cdot \mbox{tr}(\bd{R}_p^2),
\eeaa
the above two displays show that
\beaa
E\,\mbox{tr}(\hat{\bd{R}}_p^2)
=\frac{p(p-1)}{n-1}+\mbox{tr}(\bd{R}_p^2)\cdot \big[1+O(m^{-1/4})\big]
\eeaa
as $p\to \infty.$ The proof is completed.  \hfill$\square$

\subsection{The Proofs of Theorems~\ref{Theorem3} and~\ref{Theorem4}}\lbl{TRSCMTP}

\subsubsection{Mixing Moments of Gaussian Random Variables}\lbl{Variance_Study}

The following lemma is a very useful tool to compute the mean of the product of Gaussian random variables.
\begin{lemma}\lbl{Wick_formula}(Isserlis and Wick formula) Let $n\geq 2$ be an even integer and $(X_1, \cdots, X_n)^T\in \mathbb{R}^n$ follows a normal distribution with mean vector $\bd{0}$. Then
\beaa
E(X_1X_2\cdots X_n)=\sum_{p}\prod_{(i,j)\in p}E(X_iX_j),
\eeaa
where the sum runs over every pairing $p$ of $\{1, 2,\cdots, n\}$, that is, all distinct ways of partitioning $\{1,2,\cdots, n\}$ into pairs $\{i, j\}$, and the product is over the pairs contained in $p$. Sometimes we also use its equivalent form:
\beaa
E(X_1X_2\cdots X_n)=\sum_{i=2}^nE(X_1X_i)\cdot E\Big(\frac{X_2\cdots X_n}{X_i}\Big).
\eeaa
\end{lemma}

A seemingly more general, but actually an equivalent formula of Lemma~\ref{Wick_formula} is the following Lemma~\ref{hufweio}, which will be used only once  in a comment after the proof of Lemma~\ref{fourth_moment}.

\begin{lemma}\lbl{hufweio}\citep{guiard1986general}
Let $X=(X_1,\cdots,X_n)$ be n-dimensionally centralized normally distributed with $Cov(X_i,X_j)=\sigma_{ij}$. Let $\{\alpha_i;\, 1\leq i\leq n\}$ be positive integers.
 If $\sum_{i=1}^n\alpha_i$ is odd, then $\mathbb{E}(\prod_{i=1}^{n}X_i^{\alpha_{i}})=0.$
If $\sum_{i=1}^n\alpha_i$ is even, then
\beaa
    \mathbb{E}\left(\prod_{i=1}^{n}X_i^{\alpha_{i}}\right)=\sum_{2\beta_{ii}+\sum_{j:j\neq i}\beta_{ij}=\alpha_i,\forall i}\dfrac{\prod_{i=1}^n\alpha_i! }{\prod_{i=1}^n2^{\beta_{ii}} \prod_{1\leq i\leq j\leq n}\beta_{ij}!}\prod_{1\leq i\leq j\leq n}\sigma_{ij}^{\beta_{ij}},
\eeaa
where $\{\beta_{ij},1\leq i\leq j\leq n\}$ are non-negative integers.
\end{lemma}

In Lemmas~\ref{3rd_moment} and~\ref{fourth_moment} below, two identities on mixing moments of Gaussian random variables will be presented.

\newpage

\begin{lemma}\lbl{3rd_moment} Let $(X_1, X_2, X_3, X_4)^T \sim N_4(0, \bd{R})$. Assume $\bd{R}=(r_{ij})_{4\times 4}$ with $r_{ii}=1$ for each $i$. Then

(i) $E(X_1X_2X_3^2X_4^2)
= r_{12} + 2r_{12}r_{34}^2+2r_{13}r_{23} + 2r_{14}r_{24}+ 4r_{13}r_{24}r_{34}+ 4r_{14}r_{23}r_{34}$;

(ii) $E(X_{1}^2X_{2}^2X_{3}^2)= 1 +2r_{12}^2 + 2r_{13}^2+ 2r_{23}^2+ 8r_{12}r_{23}r_{31}.$

\end{lemma}
\noindent\textbf{Proof of Lemma~\ref{3rd_moment}}.
Let $(\xi_1, \cdots, \xi_6)^T$ be a multivariate normal with mean vector $\bd{0}$,  $\mbox{Cov}(\xi_i, \xi_j)=\sigma_{ij}$ and $\mbox{Var}(\xi_i)=1$ for each $i$. Then, by Lemma~\ref{Wick_formula},
\beaa
 E(\xi_1\xi_2 \xi_3\xi_4\xi_5\xi_6)
&=& \sigma_{12}\sigma_{34}\sigma_{56} + \sigma_{12}\sigma_{35}\sigma_{46} + \sigma_{12}\sigma_{36}\sigma_{45}+\\
&&  \sigma_{13}\sigma_{24}\sigma_{56} + \sigma_{13}\sigma_{25}\sigma_{46} + \sigma_{13}\sigma_{26}\sigma_{45}+\\
&& \sigma_{14}\sigma_{23}\sigma_{56} + \sigma_{14}\sigma_{25}\sigma_{36} + \sigma_{14}\sigma_{26}\sigma_{35} +\\
&& \sigma_{15}\sigma_{23}\sigma_{46} + \sigma_{15}\sigma_{24}\sigma_{36} + \sigma_{15}\sigma_{26}\sigma_{34} +\\
&& \sigma_{16}\sigma_{23}\sigma_{45} + \sigma_{16}\sigma_{24}\sigma_{35} + \sigma_{16}\sigma_{25}\sigma_{34}.
\eeaa
Now take $\xi_3=\xi_4$ and $\xi_5=\xi_6$ to see
\beaa
 E(\xi_1\xi_2 \xi_3^2\xi_5^2)
&=& \sigma_{12} + \sigma_{12}\sigma_{35}^2 + \sigma_{12}\sigma_{35}^2+\\
&&  \sigma_{13}\sigma_{23} + \sigma_{13}\sigma_{25}\sigma_{35} + \sigma_{13}\sigma_{25}\sigma_{35}+\\
&& \sigma_{13}\sigma_{23} + \sigma_{13}\sigma_{25}\sigma_{35} + \sigma_{13}\sigma_{25}\sigma_{35} +\\
&& \sigma_{15}\sigma_{23}\sigma_{35} + \sigma_{15}\sigma_{23}\sigma_{35} + \sigma_{15}\sigma_{25} +\\
&& \sigma_{15}\sigma_{23}\sigma_{35} + \sigma_{15}\sigma_{23}\sigma_{35} + \sigma_{15}\sigma_{25}\\
&= & \sigma_{12} + 2\sigma_{12}\sigma_{35}^2+2\sigma_{13}\sigma_{23} + 2\sigma_{15}\sigma_{25}+ 4\sigma_{13}\sigma_{25}\sigma_{35}+ 4\sigma_{15}\sigma_{23}\sigma_{35}.
\eeaa
This says that
\bea\lbl{long_4}
E(X_{1}X_{2}X_{3}^2X_{4}^2)=r_{12} + 2r_{12}r_{34}^2+2r_{13}r_{23} + 2r_{14}r_{24}+ 4r_{13}r_{24}r_{34}+ 4r_{14}r_{23}r_{34}.
\eea
We obtain (i). Now, take $X_2=X_1$ in \eqref{long_4}. By using the fact $r_{12}=1$ and by changing  ``$2$" to ``$1$" to the indices of $r$ we obtain
\beaa
E(X_{1}^2X_{3}^2X_{4}^2)
&=& 1 + 2r_{34}^2+2r_{13}^2 + 2r_{14}^2+ 4r_{13}r_{14}r_{34}+ 4r_{14}r_{13}r_{34}\\
& = & 1 +2r_{13}^2 + 2r_{14}^2+ 2r_{34}^2+ 8r_{13}r_{14}r_{34}.
\eeaa
In the above and change ``$3$" to ``$2$" and ``$4$" to ``$3$" to see
\beaa
E(X_{1}^2X_{2}^2X_{3}^2)= 1 +2r_{12}^2 + 2r_{13}^2+ 2r_{23}^2+ 8r_{12}r_{13}r_{23}.
\eeaa
The proof is completed. \hfill$\square$

\begin{lemma}\lbl{fourth_moment} Suppose the 4-dimensional random vector $(X_1, X_2, X_3, X_4)^T \sim N_4(0, \bd{R})$. Assume $\bd{R}=(r_{ij})_{4\times 4}$ with $r_{ii}=1$ for each $i$. Then
\beaa
 E(X_1^2X_2^2X_3^2X_4^2)
&=& 1+ 2\sum_{1\leq i<j\leq 4}r_{ij}^2 + 4(r_{12}^2r_{34}^2+ r_{13}^2r_{24}^2 + r_{14}^2r_{23}^2 ) \\
&&\ \ + 8(r_{12}r_{23}r_{31}+r_{12}r_{24}r_{41}+ r_{23}r_{34}r_{42}+r_{13}r_{34}r_{41})\\
&& \ \ + 16(r_{12}r_{23}r_{34}r_{41}+r_{12}r_{24}r_{43}r_{31}+  r_{13}r_{32}r_{24}r_{41}).
\eeaa
In particular, take $X_3=X_1$ and $X_4=X_2$ to see $E(X_1^4X_2^4)=9+72r_{12}^2+24r_{12}^4.$
\end{lemma}
\noindent\textbf{Proof of Lemma~\ref{fourth_moment}}. By Lemma~\ref{Wick_formula},
\beaa
E(X_1^2X_2^2X_3^2X_4^2)
&= & E(X_1^2)\cdot E(X_2^2X_3^2X_4^2) + 2E(X_1X_2)\cdot E(X_1X_2X_3^2X_4^2) \\
& & + 2E(X_1X_3)\cdot E(X_1X_3X_2^2X_4^2) + 2E(X_1X_4)\cdot E(X_1X_4X_2^2X_3^2).
\eeaa
By Lemma~\ref{3rd_moment}(i),
\bea\lbl{weijf2q0}
E(X_{1}X_{2}X_{3}^2X_{4}^2)=r_{12} + 2r_{12}r_{34}^2+2r_{13}r_{23} + 2r_{14}r_{24}+ 4r_{13}r_{24}r_{34}+ 4r_{14}r_{23}r_{34}.
\eea
Similarly, we obtain $E(X_1X_3X_2^2X_4^2)$ and $E(X_{1}X_4X_{2}^2X_{3}^2)$ by exchanging the roles of ``$X_2$" and ``$X_3$" and exchanging the roles of ``$X_2$" and ``$X_4$", respectively, from \eqref{weijf2q0}.  By Lemma~\ref{3rd_moment}(ii),
\beaa
E(X_{1}^2X_{2}^2X_{3}^2)= 1 +2r_{12}^2 + 2r_{13}^2+ 2r_{23}^2+ 8r_{12}r_{13}r_{23}.
\eeaa
Hence,
\beaa
E(X_1^2X_2^2X_3^2X_4^2)
&= &1 +2r_{23}^2 + 2r_{24}^2+ 2r_{34}^2+ 8r_{23}r_{24}r_{34}+\\
&& 2r_{12}\big(r_{12} + 2r_{12}r_{34}^2+2r_{13}r_{23} + 2r_{14}r_{24}+ 4r_{13}r_{24}r_{34}+ 4r_{14}r_{23}r_{34}\big)+\\
&& 2r_{13}\big(r_{13} + 2r_{13}r_{24}^2+2r_{12}r_{23} + 2r_{14}r_{34}+ 4r_{12}r_{34}r_{24}+ 4r_{14}r_{23}r_{24}\big)+\\
&& 2r_{14}\big(r_{14} + 2r_{14}r_{23}^2+2r_{13}r_{34} + 2r_{12}r_{24}+ 4r_{13}r_{24}r_{23}+ 4r_{12}r_{34}r_{23}\big).
\eeaa
Sorting them out, we have
\beaa
 E(X_1^2X_2^2X_3^2X_4^2)
&=& 1+ 2\sum_{1\leq i<j\leq 4}r_{ij}^2 + 4(r_{12}^2r_{34}^2+ r_{13}^2r_{24}^2 + r_{14}^2r_{23}^2 ) + \\
&&8(r_{12}r_{23}r_{31}+r_{12}r_{24}r_{41}+ r_{23}r_{34}r_{42}+r_{13}r_{34}r_{41})+\\
&& 16(r_{12}r_{23}r_{34}r_{41}+r_{12}r_{24}r_{43}r_{31}+ r_{13}r_{32}r_{24}r_{41}).
\eeaa
The proof is completed. \hfill$\square$

\medskip
Lemma~\ref{fourth_moment} studies $E(X_1^2X_2^2\cdots X_n^2)$ for $n=4$, which is sufficient for our purpose. It is interesting to see the formula for $n\geq 5$ by pure curiosity. If we argue the same way as in the proof of Lemma~\ref{fourth_moment} through Lemma~\ref{Wick_formula}, the sorting procedure would be messy. However, Lemma~\ref{hufweio} provides a way to do so by figuring out the non-negative integer solutions of the system equations  $2\beta_{ii}+\sum_{j:j\neq i}\beta_{ij}=2$ for $1\leq i \leq n.$

The following is a key step to study the covariance between two squared sample covariances stated in Lemma~\ref{sjwqi3}.

\begin{lemma}\lbl{shadage}  Let $\{(X_{1j}, X_{2j}, X_{3j}, X_{4j})^T \in \mathbb{R}^4;\, 1\leq j \leq m\}$ be i.i.d.  random vectors with distribution $N_4(\bd{0}, \bd{R})$, where $\bd{R}=(r_{ij})_{4\times 4}$ and $r_{ii}=1$ for each $i$. Set
\beaa
B_1=\frac{1}{m}\sum_{j=1}^mX_{1j}X_{2j} \ \ \ \ \mbox{and} \ \ \ \ B_2=\frac{1}{m}\sum_{j=1}^mX_{3j}X_{4j}.
\eeaa
Then,

(i) $\mbox{Var}(B_1)=\frac{1}{m}(1+r_{12}^2)$;

(ii) $\mbox{Cov}(B_1, B_2)=\frac{1}{m}(r_{13}r_{24}+r_{14}r_{23})$;

(iii) $E[(B_1-r_{12})(B_2-r_{34})^2]=\frac{2}{m^2}(r_{13}r_{23} + r_{14}r_{24}+ r_{13}r_{34}r_{42}+ r_{14}r_{43}r_{32})$;

(iv) $E[(B_1-r_{12})^2(B_2-r_{34})^2]$ is identical to
\beaa
&&\ \   \frac{1}{m^2}(1+r_{12}^2)(1+r_{34}^2) +
\frac{2}{m^2}\big[r_{13}^2r_{24}^2+r_{14}^2r_{23}^2+2r_{13}r_{32}r_{24}r_{41}\big] \\
& & +\frac{2}{m^3}\big[(r_{13}^2+r_{14}^2+r_{23}^2+r_{24}^2) + (r_{13}^2r_{24}^2 + r_{14}^2r_{23}^2)\\
&&+2(r_{12}r_{23}r_{31}+r_{12}r_{24}r_{41}+ r_{23}r_{34}r_{42}+r_{13}r_{34}r_{41})\\
&& +2r_{12}r_{23}r_{34}r_{41}+2r_{12}r_{24}r_{43}r_{31}+ 6 r_{13}r_{32}r_{24}r_{41}\big].
\eeaa
\end{lemma}

\noindent\textbf{Comment}. We now conduct an independent check of the accuracy of Lemma~\ref{shadage}(iv) for two special cases. In the above result, take $X_{1j}= X_{2j}= X_{3j}=X_{4j}$, then $r_{ij}=1$ for all $i,j$. Consequently, we get from Theorem~\ref{shadage}(iv) that
\beaa
E\Big[\frac{1}{m}\sum_{j=1}^m(X_{1j}^2-1)\Big]^4=E[(B_1-r_{12})^2(B_2-r_{34})^2]
=\frac{12}{m^2}+\frac{48}{m^3}.
\eeaa
Notice $\{X_{1j};\, 1\leq j \leq m\}$ are i.i.d. $N(0, 1)$. Then $\sum_{j=1}^m(X_{1j}^2-1)$ is a sum of i.i.d. random variables with $E(X_{11}^2-1)^2=2$ and $E(X_{11}^2-1)^4=60$. By a classical formula (see, for example, p. 69 from \cite{durrett2019probability}),
\beaa
E\Big[\frac{1}{m}\sum_{j=1}^m(X_{1j}^2-1)\Big]^4
&=&\frac{1}{m^4}\big\{m\cdot E(X_{11}^2-1)^4 +3m(m-1)\cdot \big[E(X_{11}^2-1)^2\big]^2\big\}\\
&=& \frac{12}{m^2}+\frac{48}{m^3}.
\eeaa
So Theorem~\ref{shadage}(iv) recovers the case for $r_{ij}=1$ for all $i,j$. On the other hand, assume the two $2$-dimensional random vectors $(X_{11}\, X_{21})^T$ and $(X_{31}\, X_{41})^T$ are independent, that is, $r_{13}=r_{14}=r_{23}=r_{24}=0$. By  Theorem~\ref{shadage}(iv),
\beaa
E[(B_1-r_{12})^2(B_2-r_{34})^2]=  \frac{1}{m^2}(1+r_{12}^2)(1+r_{34}^2).
\eeaa
On the other hand, by independence and Theorem~\ref{shadage}(i),
\beaa
E[(B_1-r_{12})^2(B_2-r_{34})^2]=E[(B_1-r_{12})^2]\cdot E[(B_2-r_{34})^2]=\frac{1}{m^2}(1+r_{12}^2)(1+r_{34}^2).
\eeaa
So this independent check indicates that Theorem~\ref{shadage}(iv) holds for the case $r_{13}=r_{14}=r_{23}=r_{24}=0$.

\noindent\textbf{Proof of Lemma~\ref{shadage}}. First, by Lemma~\ref{Wick_formula},
\bea
&& E(X_{11}^2X_{21}^2)=1+2r_{12}^2,\ \  E(X_{11}^3X_{21})=3r_{12}, \lbl{12340}\\
&& E(X_{11}^2X_{21}X_{31})=r_{23} + 2r_{12}r_{13},\ \ E(X_{11}X_{21}X_{31}X_{41})
=r_{12}r_{34}+r_{13}r_{24}+r_{14}r_{23}.\ \ \ \ \ \ \ \ \lbl{56780}
\eea
In fact, the two middle identities from the four in the above can be deduced immediately from the last one by taking $X_{11}=X_{21}=X_{31}$ and $X_{11}=X_{21}$, respectively.

(i) Notice $mB_1$ and $mB_2$ are sums of i.i.d. random variables with mean $r_{12}$ and $r_{34}$, respectively. Also, $X_{ij}\sim N(0, 1)$ for each $i,j.$ Then
$ E[(B_1-r_{12})^2]=\frac{1}{m}\mbox{Var}(X_{11}X_{21})$. Since $\mbox{Cov}(X_{11}, X_{21})=E(X_{11} X_{21})=r_{12}$, we have from \eqref{12340} that $\mbox{Var}(X_{11}X_{21})=1+r_{12}^2$.
So (i) follows.

(ii) By independence, $E[(B_1-r_{12})(B_2-r_{34})]=\frac{1}{m}\mbox{Cov}(X_{11}X_{21}, X_{31}X_{41}).$ From \eqref{56780},
\bea\lbl{catu1}
\mbox{Cov}(X_{11}X_{21}, X_{31}X_{41})=r_{12}r_{34}+r_{13}r_{24}+r_{14}r_{23}
-r_{12}r_{34}=r_{13}r_{24}+r_{14}r_{23}.
\eea

(iii) Write
\bea\lbl{chkas1}
(B_2-r_{34})^2=\frac{1}{m^2}\sum_{j=1}^m(X_{3j}X_{4j}-r_{34})^2 +\frac{2}{m^2}\sum_{1\leq k< l \leq m}(X_{3k}X_{4k}-r_{34})(X_{3l}X_{4l}-r_{34}).
\eea
By independence, the covariance between $B_1-r_{12}$ and any term from the last sum is zero. This implies
\bea\lbl{eklefw1}
E[(B_1-r_{12})(B_2-r_{34})^2]
&=&\frac{1}{m^2}E[(X_{11}X_{21}-r_{12})(X_{31}X_{41}-r_{34})^2].
\eea
Use expansion $(X_{31}X_{41}-r_{34})^2=X_{31}^2X_{41}^2-2r_{34}X_{31}X_{41}+r_{34}^2$ to see that
the last expectation in \eqref{eklefw1} is equal to $E[(X_{11}X_{21}-r_{12})(X_{31}^2X_{41}^2-2r_{34}X_{31}X_{41})]$, which is again equal to
\bea\lbl{hddu930}
&&E(X_{11}X_{21}X_{31}^2X_{41}^2)-2r_{34}E(X_{11}X_{21} X_{31}X_{41})-r_{12}E(X_{31}^2X_{41}^2)
+2r_{12}r_{34}E(X_{31}X_{41}) \nonumber\\
&= & E(X_{11}X_{21}X_{31}^2X_{41}^2)-2r_{34}(r_{12}r_{34}+r_{13}r_{24}+r_{14}r_{23})-
r_{12}(1+2r_{34}^2)+2r_{12}r_{34}^2 \nonumber\\
&= & E(X_{11}X_{21}X_{31}^2X_{41}^2)-2r_{34}(r_{12}r_{34}+r_{13}r_{24}+r_{14}r_{23})-
r_{12}
\eea
by \eqref{12340} and \eqref{56780}. From Lemma~\ref{3rd_moment}(i), we see that
\bea\lbl{eybc38}
E(X_{11}X_{21}X_{31}^2X_{41}^2)=r_{12} + 2r_{12}r_{34}^2+2r_{13}r_{23} + 2r_{14}r_{24}+ 4r_{13}r_{24}r_{34}+ 4r_{14}r_{23}r_{34}.
\eea
Plug this into the previous display,we arrive at
\beaa
E[(B_1-r_{12})(B_2-r_{34})^2]=\frac{1}{m^2}\big(2r_{13}r_{23} + 2r_{14}r_{24}+ 2r_{13}r_{24}r_{34}+ 2r_{14}r_{23}r_{34}\big).
\eeaa

(iv) By (i),
\bea\lbl{haoatx1}
&&E[(B_1-r_{12})^2(B_2-r_{34})^2]\nonumber\\
&=&\,\mbox{Cov}\big((B_1-r_{12})^2, (B_2-r_{34})^2\big) + E(B_1-r_{12})^2\cdot E(B_2-r_{34})^2 \nonumber\\
& = & \,\mbox{Cov}\big((B_1-r_{12})^2, (B_2-r_{34})^2\big) +\frac{1}{m^2}(1+r_{12}^2)(1+r_{34}^2).
\eea
Similar to \eqref{chkas1}, we have
\bea\lbl{woi85}
(B_1-r_{12})^2=\frac{1}{m^2}\sum_{j=1}^m(X_{1j}X_{2j}-r_{12})^2 +\frac{2}{m^2}\sum_{1\leq k< l \leq m}(X_{1k}X_{2k}-r_{12})(X_{1l}X_{2l}-r_{12}).
\eea
Recall the $m$ random variables  $\{(X_{1j}, X_{2j}, X_{3j}, X_{4j})^T;\, 1\leq j \leq m\}$ are i.i.d., thus each term from the double sums in \eqref{chkas1} and \eqref{woi85} is a product of two independent random variables. This implies that $E[(X_{3k}X_{4k}-r_{34})(X_{3l}X_{4l}-r_{34})]=0$ for any $k<l$, and the term $(X_{1j}X_{2j}-r_{12})^2$ is uncorrelated to $(X_{3k}X_{4k}-r_{34})(X_{3l}X_{4l}-r_{34})$ for any $j,k,l$ with $k<l.$ By the same spirit, it is easy to check that  $(X_{1k}X_{2k}-r_{12})(X_{1l}X_{2l}-r_{12})$ and $(X_{3k_1}X_{4k_1}-r_{34})(X_{3l_1}X_{4l_1}-r_{34})$ are uncorrelated  for any $k<l$ and $k_1<l_1$ as long as $(k,l)\ne (k_1, l_1)$. These yield
\bea
&&\mbox{Cov}\big((B_1-r_{12})^2, (B_2-r_{34})^2\big) \nonumber\\
&=& \frac{1}{m^3}
\mbox{Cov}\big((X_{11}X_{21}-r_{12})^2,(X_{31}X_{41}-r_{34})^2 \big)+\nonumber\\
&&\frac{4}{m^4}\sum_{1\leq k< l \leq m}\mbox{Cov}\big((X_{11}X_{21}-r_{12})(X_{12}X_{22}-r_{12}), (X_{31}X_{41}-r_{34})(X_{32}X_{42}-r_{34})\big)\nonumber\\
&= &  \frac{1}{m^3}\mbox{Cov}\big((X_{11}X_{21}-r_{12})^2,(X_{31}X_{41}-r_{34})^2 \big)+ \nonumber\\
&&\frac{2(m-1)}{m^3}\mbox{Cov}\big((X_{11}X_{21}-r_{12})(X_{12}X_{22}-r_{12}), (X_{31}X_{41}-r_{34})(X_{32}X_{42}-r_{34})\big). \lbl{cuybeuy1}
\eea
For brevity of notation, let $(X_1, X_2, X_3, X_4)\in \mathbb{R}^4$ and $(Y_1, Y_2, Y_3, Y_4)\in \mathbb{R}^4$ be i.i.d. random vectors with distribution $N_4(\bd{0}, \bd{R}).$ Then the last covariance in \eqref{cuybeuy1} is identical to
\bea\lbl{19tk}
&& E\big[(X_1X_2-r_{12})(Y_1Y_2-r_{12})(X_3X_4-r_{34})(Y_3Y_4-r_{34})\big] \nonumber\\
&=& \big\{E[(X_1X_2-r_{12})(X_3X_4-r_{34})]\big\}^2  \nonumber\\
& = & r_{13}^2r_{24}^2+r_{14}^2r_{23}^2+2r_{13}r_{32}r_{24}r_{41}
\eea
by independence and \eqref{catu1}. Now we calculate the first covariance in \eqref{cuybeuy1}. In  fact,
\bea\lbl{yibaka}
&& \mbox{Cov}\big((X_{11}X_{21}-r_{12})^2,(X_{31}X_{41}-r_{34})^2 \big)\nonumber\\
&=& \mbox{Cov}\big((X_{1}X_{2}-r_{12})^2,(X_{3}X_{4}-r_{34})^2 \big)\nonumber\\
&=& \mbox{Cov}\big(X_{1}^2X_{2}^2, X_{3}^2X_{4}^2)-2r_{12}\mbox{Cov}\big(X_{1}X_{2}, X_{3}^2X_{4}^2)-2r_{34}\mbox{Cov}\big(X_{1}^2X_{2}^2, X_{3}X_{4})+\nonumber\\
&&  4r_{12}r_{34}\mbox{Cov}\big(X_1X_2, X_3X_4\big)
\eea
since $(X_{1}X_{2}-r_{12})^2=X_{1}^2X_{2}^2-2r_{12}X_{1}X_{2}+r_{12}^2$ and $(X_{3}X_{4}-r_{34})^2=X_{3}^2X_{4}^2-2r_{34}X_{3}X_{4}+r_{34}^2$. From \eqref{12340} and Lemma~\ref{fourth_moment},
\beaa
\mbox{Cov}(X_1^2X_2^2, X_3^2X_4^2)
&=&E(X_1^2X_2^2X_3^2X_4^2)-E(X_1^2X_2^2)\cdot E(X_3^2X_4^2)\\
&=& 2(r_{13}^2+r_{14}^2+r_{23}^2+r_{24}^2) + 4(r_{13}^2r_{24}^2 + r_{14}^2r_{23}^2)\\
&&+8(r_{12}r_{23}r_{31}+r_{12}r_{24}r_{41}+ r_{23}r_{34}r_{42}+r_{13}r_{34}r_{41}) \\
&& +16(r_{12}r_{23}r_{34}r_{41}+r_{12}r_{24}r_{43}r_{31}+  r_{13}r_{32}r_{24}r_{41}).
\eeaa
Also, by \eqref{eybc38},
\beaa
\mbox{Cov}\,(X_1X_2, X_3^2X_4^2)
&=&E(X_1X_2X_3^3X_4^2)-r_{12}(1+2r_{34}^2)\\
&=& 2r_{13}r_{23} + 2r_{14}r_{24}+ 4r_{13}r_{24}r_{34}+ 4r_{14}r_{23}r_{34}.
\eeaa
By exchanging ``$1$" to ``$3$" and exchanging ``$2$" and ``$4$" in the above, we obtain
\beaa
\mbox{Cov}\,(X_1^2X_2^2, X_3X_4)=2r_{13}r_{14} + 2r_{23}r_{24}+ 4r_{13}r_{24}r_{12}+ 4r_{23}r_{14}r_{12}.
\eeaa
Combining the above computations with \eqref{catu1} and \eqref{yibaka}, we arrive at
\beaa
&&\mbox{Cov}\big((X_{11}X_{21}-r_{12})^2,(X_{31}X_{41}-r_{34})^2 \big)\\
&=& 2(r_{13}^2+r_{14}^2+r_{23}^2+r_{24}^2) + 4(r_{13}^2r_{24}^2 + r_{14}^2r_{23}^2)\\
&&+8(r_{12}r_{23}r_{31}+r_{12}r_{24}r_{41}+ r_{23}r_{34}r_{42}+r_{13}r_{34}r_{41}) \\
&& +16(r_{12}r_{23}r_{34}r_{41}+r_{12}r_{24}r_{43}r_{31}+  r_{13}r_{32}r_{24}r_{41})\\
&& -2r_{12}\big(2r_{13}r_{23} + 2r_{14}r_{24}+ 4r_{13}r_{24}r_{34}+ 4r_{14}r_{23}r_{34}\big)\\
&& -2r_{34}\big(2r_{13}r_{14} + 2r_{23}r_{24}+ 4r_{13}r_{24}r_{12}+ 4r_{23}r_{14}r_{12}\big)\\
&& + 4r_{12}r_{34}(r_{13}r_{24}+r_{14}r_{23}).
\eeaa
A careful cancellation leads to
\beaa
&& \mbox{Cov}\big((X_{11}X_{21}-r_{12})^2,(X_{31}X_{41}-r_{34})^2 \big)\\
&=& 2(r_{13}^2+r_{14}^2+r_{23}^2+r_{24}^2) + 4(r_{13}^2r_{24}^2 + r_{14}^2r_{23}^2)\\
&&+4(r_{12}r_{23}r_{31}+r_{12}r_{24}r_{41}+ r_{23}r_{34}r_{42}+r_{13}r_{34}r_{41})\\
&& +4r_{12}r_{23}r_{34}r_{41}+4r_{12}r_{24}r_{43}r_{31}+ 16 r_{13}r_{32}r_{24}r_{41}.
\eeaa
By rewriting $\frac{2(m-1)}{m^3}=\frac{2}{m^2}-\frac{2}{m^3}$, we see from the above, \eqref{cuybeuy1} and \eqref{19tk} that
\beaa
&&\mbox{Cov}\big((B_1-r_{12})^2, (B_2-r_{34})^2\big)\\
&= & \frac{1}{m^3}\big[2(r_{13}^2+r_{14}^2+r_{23}^2+r_{24}^2) + 4(r_{13}^2r_{24}^2 + r_{14}^2r_{23}^2)\\
&&+4(r_{12}r_{23}r_{31}+r_{12}r_{24}r_{41}+ r_{23}r_{34}r_{42}+r_{13}r_{34}r_{41})\\
&& +4r_{12}r_{23}r_{34}r_{41}+4r_{12}r_{24}r_{43}r_{31}+ 16 r_{13}r_{32}r_{24}r_{41}\big]\\
&&+\Big(\frac{2}{m^2}-\frac{2}{m^3}\Big)(r_{13}^2r_{24}^2+r_{14}^2r_{23}^2+2r_{13}r_{32}r_{24}r_{41}),
\eeaa
which is again equal to
\beaa
&&\frac{2}{m^2}\big[r_{13}^2r_{24}^2+r_{14}^2r_{23}^2+2r_{13}r_{32}r_{24}r_{41}\big]\\
&+ & \frac{1}{m^3}\big[2(r_{13}^2+r_{14}^2+r_{23}^2+r_{24}^2) + 2(r_{13}^2r_{24}^2 + r_{14}^2r_{23}^2)\\
&+&4(r_{12}r_{23}r_{31}+r_{12}r_{24}r_{41}+ r_{23}r_{34}r_{42}+r_{13}r_{34}r_{41})\\
&+& 4r_{12}r_{23}r_{34}r_{41}+4r_{12}r_{24}r_{43}r_{31}+ 12 r_{13}r_{32}r_{24}r_{41}\big].
\eeaa
This and \eqref{haoatx1} imply the desired conclusion. \hfill$\square$

Now we compute the covariance of two squared sample covariances.
\begin{lemma}\lbl{sjwqi3}  Let $\{(X_{1j}, X_{2j}, X_{3j}, X_{4j})^T \in \mathbb{R}^4;\, 1\leq j \leq m\}$ be i.i.d.  random vectors with distribution $N_4(\bd{0}, \bd{R})$, where $\bd{R}=(r_{ij})_{4\times 4}$ and $r_{ii}=1$ for each $i$. Let $B_1$ and $B_2$ be defined as in Lemma~\ref{shadage}. Then
\beaa
\mbox{Cov}(B_1^2, B_2^2)=\frac{1}{m}(r_{12}r_{23}r_{34}r_{41}+r_{12}r_{24}r_{43}r_{31}) +\frac{\delta_m}{m^2}\sum_{1\leq i< j \leq 4}r_{ij}^2,
\eeaa
where $|\delta_m|\leq \kappa$ and $\kappa$ is a numerical constant not depending on $m$ or $\bd{R}.$
\end{lemma}
\noindent\textbf{Proof of Lemma~\ref{sjwqi3}}. Write
\beaa
B_1^2=(B_1-r_{12})^2+2r_{12}(B_1-r_{12})+r_{12}^2 \ \ \mbox{and}\ \ B_2^2=(B_2-r_{34})^2+2r_{34}(B_2-r_{34})+r_{34}^2.
\eeaa
It follows that
\beaa
&& \mbox{Cov}(B_1^2, B_2^2)\\
&=&\mbox{Cov}\big((B_1-r_{12})^2, (B_2-r_{34})^2\big)+2r_{12}\cdot \mbox{Cov}\big((B_1-r_{12}), (B_2-r_{34})^2\big)\\
&& +2r_{34}\cdot \mbox{Cov}\big((B_1-r_{12})^2, (B_2-r_{34})\big)+ r_{12}r_{34}\cdot \mbox{Cov}\big(B_1, B_2).
\eeaa
By Lemma~\ref{shadage}(i) \& (iv),
\beaa
&& \mbox{Cov}\big((B_1-r_{12})^2, (B_2-r_{34})^2\big)\\
&=& E[(B_1-r_{12})^2(B_2-r_{34})^2]-E[(B_1-r_{12})^2]\cdot E[(B_2-r_{34})^2]\\
&=& \frac{2}{m^2}\big[r_{13}^2r_{24}^2+r_{14}^2r_{23}^2+2r_{13}r_{32}r_{24}r_{41}\big]  +\frac{2}{m^3}\big[(r_{13}^2+r_{14}^2+r_{23}^2+r_{24}^2) + (r_{13}^2r_{24}^2 + r_{14}^2r_{23}^2)\\
&&~~~~~~~~~~~~~~~~~~~~~~~~~~~~~~~~~~~~~~~~~~~~~+2(r_{12}r_{23}r_{31}+r_{12}r_{24}r_{41}+ r_{23}r_{34}r_{42}+r_{13}r_{34}r_{41})\\
&&~~~~~~~~~~~~~~~~~~~~~~~~~~~~~~~~~~~~~~~~~~~~~ +2r_{12}r_{23}r_{34}r_{41}+2r_{12}r_{24}r_{43}r_{31}+ 6 r_{13}r_{32}r_{24}r_{41}\big].
\eeaa
Use $x^2y^2\leq x^2$ and $|xyzu|\leq |xyz|\leq |xy|\leq x^2+y^2$ for any $x, y, z, u\in [-1, 1]$ to see that
\beaa
\big|\mbox{Cov}\big((B_1-r_{12})^2, (B_2-r_{34})^2\big)\big| \leq \frac{K_1}{m^2}\sum_{1\leq i< j \leq 4}r_{ij}^2,
\eeaa
where $K_1$ is a numerical constant not depending on $m$ or $\bd{R}.$ Same bounds hold with another numerical constant $K_2$ for $2r_{12}\cdot \mbox{Cov}\big((B_1-r_{12}), (B_2-r_{34})^2\big)$ and  $2r_{34}\cdot \mbox{Cov}\big((B_1-r_{12})^2, (B_2-r_{34})\big)$ by Lemma~\ref{shadage}(iii). Finally, by Lemma~\ref{shadage}(ii),
\beaa
r_{12}r_{34}\cdot \mbox{Cov}\big(B_1, B_2)=\frac{1}{m}(r_{12}r_{24}r_{43}r_{31}+r_{12}r_{23}r_{34}r_{41}).
\eeaa
The proof is finished by combining all of the estimates. \hfill$\square$

\medskip

Review sample correlation coefficient $\hat{r}_{ij}$ defined in Theorem~\ref{Theorem4}, the denominator is not easy to be handled. By the heuristic of the law of large numbers and the Taylor expansion, we are able to express $\hat{r}_{ij}^2$ as a polynomial of Gaussian random variables plus an error. The following result provides basic computations for this step. It will be used in the derivation of Lemma~\ref{youdaoa}.

\begin{lemma}\lbl{ufnuc} Let $\{(X_i, Y_i)^T;\, 1\leq i \leq m\}$ be i.i.d. $2$-dimensional normal random vectors with $EX_1=EY_1=0$, $EX_1^2=EY_1^2=1$ and $\mbox{Cov}(X_1, Y_1)=r$. Define
\beaa
 A=\frac{X_1^2+\cdots + X_m^2}{m}\ \ \ \ \mbox{and} \ \ \ \ B=\frac{Y_1^2+\cdots + Y_m^2}{m}.
 \eeaa
The following statements hold for all $m\geq 2$ with $|\delta_m|<\kappa$ where $\kappa$ is a constant not depending on $m$ or $r$.

(i) $E[(1-A)(1-B)]=\frac{2}{m}r^2$.

(ii) $E[(1-A)(1-B)^2]=-\frac{8}{m^2}r^2$.

(iii) $E[(1-A)(1-B)^3]=\frac{12r^2}{m^2} +\frac{\delta_m}{m^3}$.

(iv) $E[(1-A)^2(1-B)^2]=\frac{8r^4+4}{m^2}+ \frac{\delta_m}{m^3}$.

(v) $E[(1-A)^i(1-B)^j] =\frac{\delta_m}{m^3}$ for $(i, j)=(0,5), (1,4), (2,3).$
\end{lemma}
\noindent\textbf{Proof of Lemma~\ref{ufnuc}}. Write
\bea\lbl{cuybwr}
1-A=\frac{1}{m}\sum_{i=1}^m(1-X_i^2)\ \ \ \mbox{and}\ \ \ 1-B=\frac{1}{m}\sum_{i=1}^m(1-Y_i^2).
\eea

(i) Easily,
\bea
&& \mbox{Var}(X_1Y_1)=E\big[(X_1Y_1-r)^2\big]=E\big(X_1^2Y_1^2\big)-r^2=r^2+1;\lbl{ksaiub767}\\
&& E[(1-X_1^2)(1-Y_1^2)]=1-E(X_1^2)-E(Y_1^2)+E(X_1^2Y_1^2)=2r^2   \lbl{asji219}
\eea
by using \eqref{12340}. Due to independence,
\beaa
E[(1-A)(1-B)]=\frac{1}{m}E\big[(1-X_1^2)(1-Y_1^2)\big]=2r^2.
\eeaa

(ii)  Write
\bea\lbl{3874378}
m^2\cdot (1-B)^2=\sum_{i=1}^m(1-Y_i^2)^2+2\sum_{1\leq i<j\leq m}(1-Y_i^2)(1-Y_j^2).
\eea
Evidently, $1-X_k^2$ is uncorrelated to $(1-Y_i^2)(1-Y_j^2)$ for any $i, j, k$ with $i\ne j.$ Consequently,
\bea\lbl{0r32m}
E[(1-A)(1-B)^2]=\frac{1}{m^2}E\big[(1-X_1^2)(1-Y_1^2)^2\big].
\eea
Use the fact $(1-Y_1^2)^2=1-2Y_1^2+Y_1^4$ to see
\bea\lbl{ewn7r}
E\big[(1-X_1^2)(1-Y_1^2)^2\big]
&=& -2 E\big[(1-X_1^2)Y_1^2\big]+ E\big[(1-X_1^2)Y_1^4\big]\nonumber\\
&=& 2 E\big[(1-X_1^2)(1-Y_1^2)\big] + E\big(Y_1^4\big)-E\big(X_1^2Y_1^4\big).
\eea
Write $X_1=rY_1+\sqrt{1-r^2}Y_1' $, where $Y_1'\sim N(0, 1)$ and $Y_1'$ is independent of $X_1$. Then
\beaa
E(X_1^2Y_1^4)=r^2E(Y_1^6)+(1-r^2)E(Y_1'^2) E(Y_1^4)=15r^2+3(1-r^2)=3+12r^2.
\eeaa
Notice $E(Y_1^4)=3$ and $E(Y_1^6)=15$.
By \eqref{asji219} and \eqref{ewn7r},
\beaa
E\big[(1-X_1^2)(1-Y_1^2)^2\big]=2\cdot 2r^2+3-(3+12r^2)=-8r^2.
\eeaa
Then (ii) follows from \eqref{0r32m}.

(iii) Write
\bea\lbl{298ba}
m^3(1-B)^3
&=& \sum_{i=1}^m(1-Y_i^2)^3+3\sum_{1\leq i\ne j \leq m}(1-Y_i^2)^2(1-Y_j^2)\nonumber\nonumber\\
&&+ \sum_{1\leq i\ne j \ne k\leq m}(1-Y_i^2)(1-Y_j^2)(1-Y_k^2)\nonumber\\
&:=& U_1+U_2+U_3.
\eea
Notice $1-X_k^2$ is uncorrelated to $U_3$ and $1-X_i^2$ is uncorrelated to $(1-Y_i^2)^2(1-Y_j^2)$ if $k=i\ne j$. This gives
\beaa
&& m^4\cdot E[(1-A)(1-B)^3] \nonumber\\
&=&\sum_{i=1}^m E\big[(1-X_i^2)(1-Y_i^2)^3\big]
+3\sum_{1\leq i\ne j \leq m}E\big[(1-X_j^2)(1-Y_i^2)^2(1-Y_j^2)\big] \nonumber\\
& = & m E\big[(1-X_1^2)(1-Y_1^2)^3\big]+3m(m-1)E\big[(1-Y_1^2)^2\big]\cdot E\big[(1-X_2^2)(1-Y_2^2)\big].
\eeaa
The product of the last two expectations is equal to $2\cdot (2r^2)=4r^2$ by \eqref{asji219}. This implies
\beaa
E[(1-A)(1-B)^3]=\frac{12r^2}{m^2} +\frac{\delta_m}{m^3}
\eeaa
with $|\delta_m|<\kappa$ for all $m\geq 2$, where $\kappa$ here and later represents a constant not depending on $m$ or $r$, and can be different from line to line.

(iv) Similar to \eqref{3874378},
\bea\lbl{in3278}
m^2\cdot (1-A)^2
&=&\sum_{i=1}^m(1-X_i^2)^2+2\sum_{1\leq i<j\leq m}(1-X_i^2)(1-X_j^2) \nonumber\\
&:=& V_1+V_2.
\eea
Trivially, $(1-X_i^2)^2$ is uncorrelated to $(1-Y_j^2)(1-Y_k^2)$ for any $i,j,k$ with $j\ne k$. The same is true if $X$ and $Y$ are switched. Thus
\beaa
&& m^4\cdot \mbox{Cov}\big((1-A)^2, (1-B)^2\big)\\
&=& m\,\mbox{Cov}\big((1-X_1^2)^2, (1-Y_1^2)^2\big)+ 2m(m-1)\mbox{Cov}\big((1-X_1^2)(1-X_2^2), (1-Y_1^2)(1-Y_2^2)\big).
\eeaa
Easily, by independence and \eqref{asji219},
\beaa
 \mbox{Cov}\big((1-X_1^2)(1-X_2^2), (1-Y_1^2)(1-Y_2^2)\big)
=  \big[E(1-X_1^2)(1-Y_1^2)\big]^2= 4r^4.
\eeaa
Evidently, $E[(1-A)^2]=\frac{2}{m}$. It follows that
\beaa
E\big[(1-A)^2(1-B)^2\big]&=& \mbox{Cov}\big((1-A)^2, (1-B)^2\big) + E\big[(1-A)^2\big]\cdot E\big[(1-B)^2\big]\\
& = &  \frac{1}{m^4}\big[m\,\mbox{Cov}\big((1-X_1^2)^2, (1-Y_1^2)^2\big)+ 8 m(m-1)r^4\big]+\frac{4}{m^2}\\
& = & \frac{8r^4+4}{m^2}+ \frac{\delta_m}{m^3}.
\eeaa

(v) We will study the three cases one by one.

{\it Case 1: $(i, j)=(0,5)$}. Let $Z_1, \cdots, Z_m$ be i.i.d. random variables with mean zero. Write
\beaa
(Z_1+\cdots + Z_m)^5
&=& \sum_{i=1}^mZ_i^5+ c_1\sum_{i\ne j}Z_i^4Z_j+c_2\sum_{i\ne j}Z_i^3Z_j^2
+c_3\sum_{i\ne j\ne k}Z_i^3Z_jZ_k \\
&&+ c_4\sum_{i\ne j\ne k}Z_i^2Z_j^2Z_k +c_5\sum_{i\ne j\ne k}Z_i^2Z_jZ_kZ_l+c_6\sum_{a\ne i\ne j\ne k\ne l}Z_aZ_iZ_jZ_kZ_l,
\eeaa
where $c_1, \cdots, c_6$ are numerical coefficients not depending on $m$ or $r$. By independence, only the first and third sums  of the right hand side above are non-zero if the expectation is taken on both sides. The total number of the terms appeared in the two sums is $m+2\cdot \binom{m}{2}=m^2$. Take $Z_i=1-Y_i^2$ for each $i$ to have
\beaa
E[(1-B)^5]=\frac{\delta_m}{m^3}
\eeaa
for all $m\geq 2$.

{\it Case 2: $(i, j)=(1,4)$}. Similar to  {\it Case 1},
\bea\lbl{wuy9}
(Z_1+\cdots + Z_m)^4
&=& \sum_{i=1}^mZ_i^4+ c_1\sum_{i\ne j}Z_i^3Z_j+6\sum_{i< j}Z_i^2Z_j^2
+c_3\sum_{i\ne j\ne k}Z_i^2Z_jZ_k \nonumber\\
&&+ c_4\sum_{i\ne j\ne k\ne l}Z_iZ_jZ_kZ_l.
\eea
Take $Z_i=1-Y_i^2$ for each $i$. Recall \eqref{cuybwr}, each term $1-X_i^2$ is uncorrelated to any term in the last two sums. By the same argument as in {\it Case 1}, we know $E[(1-A)(1-B)^4]=\frac{1}{m^3}\delta_m$ for all $m\geq 2.$

{\it Case 3: $(i, j)=(2,3)$}. Review \eqref{298ba} and \eqref{in3278}. Easily, $U_3$ is uncorrelated to both $V_1$ and $V_2$;   $U_1$ is uncorrelated to $V_2$. It follows that
\beaa
m^5\cdot \mbox{Cov}\big((1-A)^2, (1-B)^3\big)=\mbox{Cov}\big(U_1, V_1\big) + \mbox{Cov}\big(U_2, V_1\big) + \mbox{Cov}\big(U_2, V_2\big).
\eeaa
By independence and the same argument as before, it is easy to see that $\mbox{Cov}\big(U_1, V_1\big)=m\delta_m$,  $\mbox{Cov}\big(U_2, V_1\big)=m^2\delta_m'$ and $\mbox{Cov}\big(U_2, V_2\big)=m^2\delta_m''$, where $|\delta_m|+|\delta_m'|+|\delta_m''|\leq \kappa$ and $\kappa$ is a constant not depending on $m$ or $r.$ Consequently,
\bea\lbl{cq827}
\mbox{Cov}\big((1-A)^2, (1-B)^3\big)=\frac{\delta_m}{m^3}.
\eea
From \eqref{in3278} and \eqref{298ba}, it is readily seen $E[(1-A)^2]=\frac{2}{m}$ and $E[(1-B)^3]=\frac{1}{m^2}\cdot E(1-Y_1^2)^3$. The conclusion then follows from these facts, \eqref{cq827} and the formula
\beaa
E\big[(1-A)^2(1-B)^3\big]=\mbox{Cov}\big((1-A)^2, (1-B)^3\big) + E\big[(1-A)^2\big]\cdot E\big[(1-B)^3].
\eeaa
The proof is finished. \hfill$\square$

\subsubsection{Evaluations of Covariances between Monomials of Gaussian Random Variables }\lbl{Monomials}

Let $(X_1, X_2, X_3, X_4)^T$ be a $4$-dimensional random vector with distribution $N_4(\bd{0}, \bd{R})$ where $\bd{R}=(r_{ij})_{4\times 4}$ and $r_{ii}=1$ for each $i$. One of major tasks in this section is showing $|\mbox{Cov}(X_1^{\alpha}X_2^{\beta}, X_3^{\gamma}X_4^{\delta})|\leq C\sum_{1\leq i <j \leq 4}r_{ij}^2$ for some constant $C$ depending only on non-negative integers $\alpha, \beta, \gamma, \delta$. Review the $t$-norm $\|U\|_t=[E(|U|^t)]^{1/t}$ for any random variable $U$ and $t\geq 1$. In the following lemma, by convention we regard $\prod_{i=1}^0a_i=1$ for any $a_i$.

\begin{lemma}\lbl{garbage} Let $\{X_1, \cdots, X_k, Y_1, \cdots, Y_l\}$ be random variables with $k\geq 0$, $l\geq 0$ and $k+l\geq 1$. Let $p_1, \cdots, p_k, q_1, \cdots, q_l$ be positive integers. Assume, for each $i$,
(a) $X_i=\frac{1}{\sqrt{p_i}}\sum_{j=1}^{p_i}\xi_{ij}$ where $\{\xi_{i1}, \cdots, \xi_{ip_i}\}$ are i.i.d. with mean $0$, and (b) $Y_i=\frac{1}{q_i}\chi^2(q_i)$. Then
\beaa
E\Big[\prod_{i=1}^k|X_i|^{\alpha_i}\cdot \prod_{j=1}^{l}Y_j^{\beta_j}\Big] <C\cdot \prod_{i=1}^k\|\,|\xi_{i1}|^{\alpha_i}\|_{k+l}
\eeaa
for any $\alpha_i\geq 2/(k+l)$,  $\beta_j>(2-q_j)/[2(k+l)]$, $1\leq i \leq k$ and $1\leq j \leq l$, where $C$ is a constant depending on $k$, $l$, $\alpha_i$ and $\beta_j$ but not depending on $p_i$ or $q_j$.
\end{lemma}
\noindent\textbf{Proof of Lemma~\ref{garbage}}. First we assume $k\geq 1$ and $l\geq 1$. By applying the H\"{o}lder inequality to the product of $k+l$ terms, we see
\bea\lbl{8n509}
E\Big[\prod_{i=1}^k|X_i|^{\alpha_i}\cdot \prod_{j=1}^{l}Y_j^{\beta_j}\Big]\leq
\prod_{i=1}^k\big(E|X_i|^{\alpha_i(k+l)}\big)^{1/(k+l)}\cdot \prod_{j=1}^{l}\big(E Y_j^{\beta_j(k+l)}\big)^{1/(k+l)}.
\eea
By assumption, $\alpha_i(k+l) \geq 2$ for each $i$. According to the Marcinkiewicz-Zygmund inequality [see, for example,  the proof of Corollary 2 on p. 387 from \cite{chow1997probability}],
\bea\lbl{sagwn}
E|X_i|^{\alpha_i(k+l)} \leq C(\alpha_i, k,l)\cdot E\big(|\xi_{i1}|^{\alpha_i(k+l)}\big),
\eea
where $C(\alpha_i, k,l)$ is a constant depending on $\alpha_i, k,l$ only. As a result,
\bea\lbl{ackutkru}
\prod_{i=1}^k\big(E|X_i|^{\alpha_i(k+l)}\big)^{1/(k+l)}\leq
C(\alpha_1,\cdots, \alpha_k, k, l)\cdot \prod_{i=1}^k\||\xi_{i1}|^{\alpha_i}\|_{k+l}.
\eea
Second, by recalling the density of $\chi^2(m)$ is $f(x):=\frac{1}{\Gamma(\frac{m}{2})2^{m/2}}x^{(m/2)-1}e^{-x/2}$ for $x>0$, we have
\bea
E[\chi^2(m)^{\beta}]
&=&\frac{1}{\Gamma(\frac{m}{2})2^{m/2}}\int_0^{\infty}x^{(m+2\beta)/2-1}e^{-x/2}\,dx \nonumber\\
& = & \frac{\Gamma(\frac{m}{2}+\beta)}{\Gamma(\frac{m}{2})2^{-\beta}} \lbl{consist_with}
\eea
for any real number $\beta > -m/2$. From the fact that $\lim_{x\to \infty}\frac{\Gamma(x+a)}{\Gamma(x)x^a}=1$ for any number $a\in \mathbb{R}$, it is seen that   $E[\chi^2(m)/m]^\beta\to 1$ as $m\to \infty.$ This implies that
\bea\lbl{cdh8293f}
\sup E\Big[\frac{\chi^2(m)}{m}\Big]^\beta= C_{\beta}<\infty,
\eea
where the supremum is taken over all $m$ satisfying $m\geq 1$ and $m\geq 2(1-\beta)$ and where $C_{\beta}$ depends on $\beta$ only. The reason we choose $m$  such that $\frac{m}{2}+\beta\geq 1$ is because we need $\frac{m}{2}+\beta$ to stay away from the only singular point $0$ of  $\Gamma(x)$ defined on $[0, \infty).$  It follows that
\beaa
E\big(|Y_j|^{\beta_j(k+l)}\big)=E\Big[\Big(\frac{\chi^2(q_j)}{q_j}\Big)^{\beta_j(k+l)}\Big]
\leq C_{\beta_j(k+l)}<\infty
\eeaa
provided $q_j\geq 2[1-\beta_j(k+l)]$, or equivalently, $\beta_j\geq (2-q_j)/[2(k+l)]$. This asserts that
\bea\lbl{c29}
\prod_{j=1}^{l}\big(E Y_j^{\beta_j(k+l)}\big)^{1/(k+l)}=C(k,l,\beta_1, \cdots, \beta_l)<\infty
\eea
under the assumption  $\beta_j\geq (2-q_j)/[2(k+l)]$ for each $1\leq j \leq l$.
 This joined with \eqref{8n509} and \eqref{ackutkru} leads to the desired inequality.

If $k\geq 1$ and $l=0$, then \eqref{8n509}-\eqref{ackutkru} still hold. If $k=0$ and $l\geq 1$, then \eqref{8n509} and \eqref{c29} are also true. The desired statements are then derived. \hfill$\square$

Recall sample correlation coefficient $\hat{r}_{ij}$ defined in Theorem~\ref{Theorem4} from Section~\ref{TRSCMTP}, the next result provides an estimate for the dependency between two denominators.

\begin{lemma}\lbl{youdaoa} Let $\{(X_i, Y_i)^T;\, 1\leq i \leq m\}$ be i.i.d. $2$-dimensional normal random vectors with $EX_1=EY_1=0$, $EX_1^2=EY_1^2=1$ and $\mbox{Cov}(X_1, Y_1)=r$. Write
\bea\lbl{yilong}
E\Big(\frac{m}{X_1^2+\cdots + X_m^2}\cdot \frac{m}{Y_1^2+\cdots + Y_m^2}\Big)=1+ \frac{4+2r^2}{m}+\frac{12+8r^2+8r^4}{m^2}+\frac{\delta_m}{m^3}.
\eea
Then $|\delta_m|\leq \kappa$ for all $m\geq 11$, where $\kappa$ is a constant not depending on $m$ or $r$.
\end{lemma}

Before we present the proof let us have a quick check for two cases with $r=0$ and $r=1$. Recall
\eqref{consist_with}. We have $E\frac{1}{\chi^2(m)}=\frac{1}{m-2}$ and $E\frac{1}{\chi^2(m)^2}=\frac{1}{(m-2)(m-4)}$ for any $m\geq 5$. If $r=0$, then $\{X_i, Y_i;\, 1\leq i \leq m\}$ are i.i.d. $N(0, 1)$. Consequently the left hand side of \eqref{yilong} is identical to
\beaa
m^2\Big[E\frac{1}{\chi^2(m)}\Big]^2=\frac{m^2}{(m-2)^2}=1+\frac{4}{m}+\frac{12}{m^2} +\frac{32m-48}{m^2(m-2)^2},
\eeaa
which corresponds to the right hand side of \eqref{yilong} with $r=0$. If $r=1$, then $X_i=Y_i$ for each~$i$. Thus the left hand side of \eqref{yilong} becomes
\beaa
m^2\cdot E\frac{1}{\chi^2(m)^2}=\frac{m^2}{(m-2)(m-4)}=1+\frac{6}{m} + \frac{28}{m^2}+\frac{120m-224}{m^2(m-2)(m-4)},
\eeaa
which is equal to the right hand side of \eqref{yilong} with $r=1$.
\medskip

\noindent\textbf{Proof of Lemma~\ref{youdaoa}}. Define
\beaa
A=\frac{X_1^2+\cdots + X_m^2}{m}\ \ \mbox{and}\ \ B=\frac{Y_1^2+\cdots + Y_m^2}{m}.
\eeaa
Set  $A_1=1-A$ and $B_1=1-B$. By the formula $x^{-1}=1+(1-x)+(1-x)^2+(1-x)^3+(1-x)^4+(1-x)^5+x^{-1}(1-x)^6$, we write
\beaa
\frac{1}{AB}=
\Big(\frac{A_1^6}{A}+\sum_{i=0}^5A_1^i\Big)\Big( \frac{B_1^6}{B}+\sum_{i=0}^5B_1^i\Big).
\eeaa
We will expand $(\sum_{i=0}^5A_1^i)(\sum_{i=0}^5B_1^i)$ and write it as the sum of $A_1^aB_1^b$, and furthermore break the sum into two sums, the first of which is for the terms with $a+b\leq 5$ and the second of which is for those with  $a+b\geq 6$. Therefore,
\bea\lbl{ewuy3r8}
\frac{1}{AB}
&=& 1+ (A_1+B_1)+(A_1^2+A_1B_1+A_2^2)  + \big(A_1^3+A_1^2B_1+A_1B_1^2+B_1^3\big)\nonumber\\
&& +\big(A_1^4+A_1^3B_1+A_1^2B_1^2+A_1B_1^3+B_1^4\big)\nonumber\\
&& + \big(A_1^5+A_1^4B_1+A_1^3B_1^2+A_1^2B_1^3+A_1B_1^4 +B_1^5\big)
+\epsilon_m,
\eea
where
\beaa
\epsilon_m:=\frac{A_1^6}{A}\sum_{i=0}^5B_1^i+ \frac{B_1^6}{B}\sum_{i=0}^5A_1^i +\frac{A_1^6B_1^6}{AB}+p_m
\eeaa
and where $p_m=\sum A_1^aB_1^b$ with the sum running over all $a\geq 1$ and $b\geq 1$ satisfying $6\leq a+b\leq 10$. Obviously, by counting the number of the terms on the right hand side of \eqref{ewuy3r8}, we know the total number of the terms in the sum of $p_m$ is  $36-(1+2+3+4+5+6)=15$. Notice
\beaa
A\sim \frac{\chi^2(m)}{m},\ \ \ B\sim \frac{\chi^2(m)}{m},\ \ \ \sqrt{m}A_1=\frac{1}{\sqrt{m}}\sum_{i=1}^m(1-X_i^2),
\ \ \ \sqrt{m}B_1=\frac{1}{\sqrt{m}}\sum_{i=1}^m(1-Y_i^2).
\eeaa
Review the notation in Lemma~\ref{garbage}. We first consider the term $E(A_1^6B_1^{i}A^{-1})$. If $i=0$, take $k=l=1$, $\alpha_1=6$, $\beta_1=-1$ and $q_1=m$. Then $\alpha_1\geq 2/(k+l)=1$ and $\beta_1>(2-q_1)/[2(k+l)]=(2-m)/4$ as  $m\geq 7.$ Similarly, if $i\geq 1$, take $k=2$, $l=1$, $\alpha_1=6$, $\alpha_2=i$, $\beta_1=-1$ and $q_1=m$. It is always true that $\alpha_1\wedge\alpha_2\geq 2/(k+l)=2/3.$ Also, $\beta_1> (2-q_1)/[2(k+l)]=(2-m)/6$ if $m\geq 9.$ A similar but easier check can be done for $A_1^aB_1^b$ with $a\geq 1$, $b\geq 1$ and $a+b\geq 6$. It then follows from Lemma~\ref{garbage} that
\bea\lbl{shwno87}
E\Big(\frac{A_1^6}{A}B_1^i\Big)=\frac{\delta(m,i)}{m^3}\ \ \mbox{and}\ \ E(A_1^aB_1^b)=\frac{\delta(m,a, b)'}{m^3}
\eea
as $m\geq 9$, where  $|\delta(m,i)|+|\delta(m,a,b)'|\leq \kappa$ for all $1\leq i \leq 5$ and  integers $a\geq 1$ and $b\geq 1$ with $a+b\geq 6$. Here and later $\kappa$  represents a constant not depending on $m$ or $r$ and can be different from line to line.

Now we turn to look at the term $A_1^6B_1^6A^{-1}B^{-1}$. Take $k=l=2$, $\alpha_1=\alpha_2=6$, $\beta_1=\beta_2=-1$ and $q_1=q_2=m$ to see that, for each $i=1,2$, we have $\alpha_i\geq 2/(k+l)=1/2$ and $\beta_i> (2-q_i)/[2(k+l)]=(2-m)/8$ provided $m\geq 11.$ From Lemma~\ref{garbage} again we obtain
\bea\lbl{dc239}
E\frac{A_1^6B_1^6}{AB}=\frac{\delta_m}{m^6}
\eea
with $|\delta(m)|\leq \kappa$ as $m\geq 11.$ Combining \eqref{shwno87} and \eqref{dc239} and using the symmetry of $A$ and $B$ as well as that of $A_i$ and $B_i$, we see that
\bea\lbl{wueowe}
E\epsilon_m=\frac{\delta_m}{m^3}
\eea
with $|\delta_m|\leq \kappa$ for all $m\geq 11$.

Back to \eqref{ewuy3r8}, let us examine the expectation of each monomial of $A_1$ and $B_1$ on the right hand side. By independence,
\bea\lbl{c348}
E\big(A_1^2\big)=\frac{1}{m^2}\sum_{i=1}^mE\big[(1-X_i^2)^2\big]=\frac{2}{m}
\eea
since $\mbox{Var}(X_1^2)=2.$ Recall
\beaa
(t_1+\cdots +t_m)^3=\sum_{i=1}^mt_i^3+3\sum_{1\leq i\ne j \leq m}t_i^2t_j+ \sum_{1\leq i\ne j \ne k\leq m}t_it_jt_k
\eeaa
for any real numbers $t_1, \cdots, t_m$. Take $t_i=1-X_i^2$  and use independence to see
\beaa
E\big(A_1^3\big)=\frac{1}{m^2}E\big[(1-X_1^2)^3\big]
=\frac{1}{m^2}\big[1-3E(X_1^2)+3E(X_1^4)-E(X_1^6)\big]=-\frac{8}{m^2}
\eeaa
since $E(X_1^2)=1$, $E(X_1^4)=3$ and $E(X_1^6)=15$. Review
\beaa
(t_1+\cdots + t_m)^4
&=& \sum_{i=1}^mt_i^4+ c_1\sum_{i\ne j}t_i^3t_j+6\sum_{i< j}t_i^2t_j^2
+c_3\sum_{i\ne j\ne k}t_i^2t_jt_k \\
&&+ c_4\sum_{i\ne j\ne k\ne l}t_it_jt_kt_l
\eeaa
for any $t_1, \cdots, t_m$. By taking $t_i=1-X_i^2$ and using independence, we get
\beaa
E\big(A_1^4\big)
&=& \frac{1}{m^3}\cdot E\big[(1-X_1^2)^4\big] + \frac{6}{m^4}\cdot
\binom{m}{2}\cdot E\big[(1-X_1^2)^2(1-X_2^2)^2\big]\\
&= & \frac{12}{m^2}+\frac{\delta_m}{m^3}
\eeaa
with $|\delta_m|\leq \kappa$ for all $m\geq 1$ since $E\big[(1-X_1^2)^2(1-X_2^2)^2\big]=[E(1-X_1^2)^2]^2=4.$ Now, by \eqref{ewuy3r8},
\beaa
E\Big(\frac{1}{AB}-\epsilon_m\Big)
&=& 1+E\big(A_1^2+A_1B_1+B_1^2\big) + E\big(A_1^3+A_1^2B_1+A_1B_1^2+B_1^3\big)\\
&& +E\big(A_1^4+A_1^3B_1+A_1^2B_1^2+A_1B_1^3+B_1^4\big)\\
&& + E\big(A_1^5+A_1^4B_1+A_1^3B_1^2+A_1^2B_1^3+A_1B_1^4 +B_1^5\big).
\eeaa
From Lemma~\ref{ufnuc}, the expressions of $EA_1^i$ for $i=2,3,4$ above as well as the symmetry of $A$ and $B$, we obtain
\beaa
E\Big(\frac{1}{AB}-\epsilon_m\Big)
&=& 1+\Big(\frac{4}{m}+\frac{2}{m}r^2\Big)+ \Big(-\frac{16}{m^2}-\frac{16}{m^2}r^2\Big)\\
&& \ \ \,+ \Big(\frac{24}{m^2}+\frac{24r^2}{m^2}+\frac{8r^4+4}{m^2}\Big)+\frac{\delta_m}{m^3}\\
& = & 1+ \frac{4+2r^2}{m}+\frac{12+8r^2+8r^4}{m^2}+\frac{\delta_m}{m^3}
\eeaa
with $|\delta_m|\leq \kappa$ for all $m\geq 2$. This and \eqref{wueowe} conclude
\beaa
E\frac{1}{AB}=1+  \frac{4+2r^2}{m}+\frac{12+8r^2+8r^4}{m^2}+\frac{\delta_m}{m^3}
\eeaa
with $|\delta_m|\leq \kappa$ for all $m\geq 11$. \hfill$\square$

\medskip

Recall the notation $(2i-1)!!=1\cdot 3\cdots(2i-1)$ for any integer $i\geq 1$. We set  $(-1)!!=1$ by convention.
\begin{lemma}\lbl{wdh120} Let $X$ and $Y$ be $N(0, 1)$.  Assume they are jointly normal with covariance $r$. Then the following hold.

(i) For any integers $i\geq 1$ and $j\geq 1$, there exists a polynomial $f(x)$ depending on $i$ and $j$ such that $f(0)=1$ and $E(X^{2i-1}Y^{2j-1})=(2i-1)!!\cdot (2j-1)!!\cdot r f(r^2).$

(ii) For any integers $i\geq 0$ and $j\geq 0$, there exists a polynomial $g(x)$ depending on $i$ and $j$ such that $E(X^{2i}Y^{2j})=(2i-1)!!\cdot (2j-1)!!+r^2g(r^2)$.
\end{lemma}
\noindent\textbf{Proof of Lemma~\ref{wdh120}}. Let $Z$ be a random variable with distribution  $N(0, 1)$ and $Z$ be independent of $X$ and $Y$. Then we may write $Y=rX+r'Z$, where $r'=\sqrt{1-r^2}.$

(i) Easily,
\beaa
Y^{2j-1}=\sum_{k=0}^{2j-1}\binom{2j-1}{k}r^kr'^{2j-k-1}X^kZ^{2j-k-1}.
\eeaa
By independence and the fact $E[N(0,1)^n]=0$ for any odd number $n\geq 1$, we have
\bea
E(X^{2i-1}Y^{2j-1})
&=&\sum_{k=0}^{2j-1}\binom{2j-1}{k}r^kr'^{2j-k-1}E\big(X^{2i+k-1}Z^{2j-k-1}\big)\nonumber\\
&=& \sum_{l=1}^{j}\binom{2j-1}{2l-1}r^{2l-1}r'^{2(j-l)}E\big(X^{2(i+l-1)}\big)\cdot E\big(Z^{2(j-l)}\big)\nonumber\\
&=& r\sum_{l=1}^{j}\binom{2j-1}{2l-1}(r^2)^{l-1}(1-r^2)^{j-l}\cdot (2(i+l)-3)!!\cdot(2(j-l)-1)!!\nonumber\\
&:=& rh(r^2) \lbl{1290},
\eea
where we set $k=2l-1$ in the second identity and the fact $E[N(0,1)^n]=(n-1)!!$ for even number $n\geq 1$ is used in the last display. Obviously the last term in \eqref{1290} is obviously a polynomial. Then
\beaa
h(0)
&=& \sum_{l=1}^{j}\binom{2j-1}{2l-1}(r^2)^{l-1}(1-r^2)^{j-l}\cdot (2(i+l)-3)!!\cdot(2(j-l)-1)!!\Big|_{r=0}\\
&= & (2j-1)\cdot (2(i+1)-3)!!\cdot(2(j-1)-1)!!\\
&=& (2i-1)!!\cdot (2j-1)!!.
\eeaa
Set $f(x)=g(x)[(2i-1)!!\cdot (2j-1)!!]^{-1}$. The desired conclusion follows.

(ii) Recall $E(X^{2i})=(2i-1)!!$ and $E(Y^{2j})=(2j-1)!!$. The conclusion is obviously true if $(i,j)=(0,0)$, $(i,j)=(0,1)$ and $(i,j)=(1,0)$, in which cases $g(x)=0$ for each $x\in \mathbb{R}$.  We next assume $i\geq 1$ and $j\geq 1$. By the same argument as in (i),
\beaa
Y^{2j}=(rX+r'Z)^{2j}=\sum_{k=0}^{2j}\binom{2j}{k}r^kr'^{2j-k}X^kZ^{2j-k}.
\eeaa
Then $E(X^{2i}Y^{2j})$ is equal to
\beaa
&&\sum_{k=0}^{2j}\binom{2j}{k}r^kr'^{2j-k}E\big(X^{2i+k}Z^{2j-k}\big)\\
&=& \sum_{l=0}^{j}\binom{2j}{2l}r^{2l}r'^{2(j-l)}E\big(X^{2(i+l)}\big)\cdot E\big(Z^{2(j-l)}\big)\\
&=& (1-r^2)^jE\big(X^{2i}\big)\cdot E\big(Z^{2j}\big)+r^2\sum_{l=1}^{j}\binom{2j}{2l}(r^2)^{l-1}(1-r^2)^{j-l}E\big(X^{2(i+l)}\big)\cdot E\big(Z^{2(j-l)}\big),
\eeaa
where we set $k=2l$ in the first identity and  single out the term with $l=0$ in the last step. Obviously the last sum is a polynomial of $r^2$, say, $h_1(r^2).$
Write  $(1-r^2)^j=1+r^2h_2(r^2)$ with the function $h_2(x)$ being a polynomial. Consequently,
\beaa
E(X^{2i}Y^{2j})
&=& [1+r^2h_2(r^2)]E(X^{2i})\cdot E(Z^{2j}) +r^2h_1(r^2)\\
&=&(2i-1)!!\cdot (2j-1)!!+r^2g(r^2),
\eeaa
where $g(x):=(2i-1)!!\cdot (2j-1)!!h_2(x)+h_1(x)$ is a polynomial. The proof is completed. \hfill$\square$

\medskip

With the help of Lemma~\ref{wdh120}, we will obtain Lemmas~\ref{dw0ms}-\ref{29385} in the following. They  estimate the size of $|\mbox{Cov}(X_1^{d_1}X_2^{d_2}, X_3^{d_3}X_4^{d_4})|$ for non-negative integers $d_1, d_2, d_3, d_4$ with special requirements.

\begin{lemma}\lbl{dw0ms} Let $(X_1, X_2, X_3, X_4)^T$ be a $4$-dimensional random vector with distribution $N_4(\bd{0}, \bd{R})$ where $\bd{R}=(r_{ij})_{4\times 4}$ and $r_{ii}=1$ for each $i$. Then, there exists a constant $C>0$ depending on $i, j,k, l$ but not on $\bd{R}$ such that
$|\mbox{Cov}(X_1^{2i}X_2^{2j}, X_3^{2k}X_4^{2l})|\leq C\sum_{1\leq i<j\leq 4}r_{ij}^2$ for all non-negative integers $i,j,k$ and $l$.
\end{lemma}
\noindent\textbf{Proof of Lemma~\ref{dw0ms}}.  If $i+j=0$, then $i=j=0$, and hence $X_1^{2i}X_2^{2j}=1$. The conclusion obviously holds. The same is true if $k+l=0$. So we assume
$i+j\geq 1$ and $k+l\geq 1$ next. Write
\bea\lbl{9147}
\mbox{Cov}(X_1^{2i}X_2^{2j}, X_3^{2k}X_4^{2l})
=E\big(X_1^{2i}X_2^{2j} X_3^{2k}X_4^{2l}\big)-E\big(X_1^{2i}X_2^{2j}\big)\cdot E\big(X_3^{2k}X_4^{2l}\big).
\eea
By Lemma~\ref{wdh120}(ii), there exists polynomials $g_1(x)$ and $g_2(x)$ for $x\in \mathbb{R}$ such that
\bea
&& E\big(X_1^{2i}X_2^{2j}\big)=(2i-1)!!\cdot (2j-1)!!+r_{12}^2\cdot g_1(r_{12}^2);\lbl{fwq0}\\
&& E\big(X_3^{2k}X_4^{2l}\big)=(2k-1)!!\cdot (2l-1)!!+r_{34}^2\cdot g_2(r_{34}^2).\lbl{aiqw09}
\eea
Set $r=(r_{12}, r_{13}, r_{14}, r_{23}, r_{24}, r_{34})^T\in \mathbb{R}^6$. Then, by Lemma~\ref{Wick_formula}, $E(X_1^{2i}X_2^{2j} X_3^{2k}X_4^{2l})$ is a multivariate polynomial of $r_{ij}, 1\leq i<j \leq 4$. So we are able to write
\bea\lbl{a7280}
E\big(X_1^{2i}X_2^{2j} X_3^{2k}X_4^{2l}\big)=C_0+\sum_{1\leq i< j \leq 4}C_{ij}r_{ij}+ c(r),
\eea
where $C_0$ and $C_{ij}$ are constants and  $c(r)$ is a linear combination  of $\prod_{1\leq i<j\leq 4}r_{ij}^{\alpha_{ij}}$ with $2\leq \sum_{1\leq i<j \leq 4}\alpha_{ij}\leq i+j+k+l.$ Next  we will use the fact that \eqref{a7280} holds for every $r_{ij}\in [-1, 1]$ and $1\leq i<j \leq 4$ to identify the values of $C_0$ and every $C_{ij}$.

Set $r_{ij}=0$ for all $i\ne j$. Then $X_1, X_2, X_3, X_4$ are i.i.d. $N(0,1)$. It follows from \eqref{a7280} that
\beaa
C_0
&=& E(X_1^{2i})\cdot E(X_2^{2j})\cdot E(X_3^{2k})\cdot E(X_4^{2l})\\
&=& (2i-1)!!\cdot (2j-1)!!\cdot (2k-1)!!\cdot (2l-1)!!.
\eeaa
Review \eqref{a7280}. We claim that $C_{ij}=0$ for all $1\leq i< j \leq 4$. Take $r_{ij}=0$ for all $1\leq i<j\leq 4$ except $(i, j)=(1,2)$. Then the three random variables $(X_1, X_2)^T\in \mathbb{R}^2,\, X_3$ and $X_4$ are  independent. By~\eqref{a7280},
\bea\lbl{wdqh09}
&&E\big(X_1^{2i}X_2^{2j}\big)\cdot E(X_3^{2k})\cdot E(X_4^{2l})\nonumber\\
&=& (2i-1)!!\cdot (2j-1)!!\cdot (2k-1)!!\cdot (2l-1)!!+C_{12}r_{12}+c_1(r_{12}),
\eea
where $c_1(x)=c_2x^2+c_3x^3+\cdots + c_px^p$ with $2\leq p \leq i+j+k+l$. Use the fact $E(X_3^{2k})=(2k-1)!!$, $E(X_4^{2l})=(2l-1)!!$ and \eqref{fwq0} to see
\beaa
&& E\big(X_1^{2i}X_2^{2j}\big)\cdot E(X_3^{2k})\cdot E(X_4^{2l})\\
&=& \big[(2i-1)!!\cdot (2j-1)!!+r_{12}^2\cdot g_1(r_{12}^2)\big]\cdot (2k-1)!!\cdot (2l-1)!!.
\eeaa
Observe there are no linear terms of $r_{12}$ on the right hand side. Compare this with \eqref{wdqh09}, we see $C_{12}=0$. By symmetry, $C_{ij}=0$ for all  $1\leq i< j \leq 4$. Combining this, \eqref{9147}-\eqref{a7280}, we get $\mbox{Cov}(X_1^{2i}X_2^{2j}, X_3^{2k}X_4^{2l})$ is identical to
\bea
&& (2i-1)!!\cdot (2j-1)!!\cdot (2k-1)!!\cdot (2l-1)!!+ c(r) \nonumber\\
 &-&\big[(2i-1)!!\cdot (2j-1)!!+r_{12}^2\cdot g_1(r_{12}^2)\big]\cdot \big[(2k-1)!!\cdot (2l-1)!!+r_{34}^2\cdot g_2(r_{34}^2)\big] \nonumber\\
&:=& c_2(r).\lbl{yw8946}
\eea
Observe that $\mbox{Cov}(X_1^{2i}X_2^{2j}, X_3^{2k}X_4^{2l})=c_2(r)$ is a linear combination  of $\prod_{1\leq i<j\leq 4}r_{ij}^{\alpha_{ij}}$ with $2\leq \sum_{1\leq i<j \leq 4}\alpha_{ij}\leq i+j+k+l.$ Notice $|r_{ij}|\leq 1$, and hence $|r_{ij}^{\alpha_{ij}}|\leq r_{ij}^2$ for all $\alpha_{ij}\geq 2$ and $\prod_{1\leq i<j\leq 4}r_{ij}^{\alpha_{ij}}\leq |r_{kl}|^{\alpha_{kl}}|r_{uv}|^{\alpha_{uv}}$ for any $1\leq k<l\leq 4$, $1\leq u< v\leq 4$, $(k, l) \ne (u, v)$ with $\alpha_{kl}+\alpha_{uv}\geq 2$. Use the formula $|xy|\leq \frac{1}{2}(x^2+y^2)$ to get
\beaa
|r_{kl}|^{\alpha_{kl}}|r_{uv}|^{\alpha_{uv}}\leq r_{kl}^2+r_{uv}^2\leq \sum_{1\leq i<j\leq 4}r_{ij}^2.
\eeaa
This says that there exists a constant $C>0$ depending on $i, j,k, l$ but not on $\bd{R}$ such that
\bea\lbl{qi105}
|\mbox{Cov}(X_1^{2i}X_2^{2j}, X_3^{2k}X_4^{2l})|=|c_2(r)|\leq C\sum_{1\leq i<j\leq 4}r_{ij}^2.
\eea
The proof is completed. \hfill$\square$

\begin{lemma}\lbl{dhwnal} Let $(X_1, X_2, X_3, X_4)^T$ be the same as in Lemma~\ref{dw0ms}.  Then, there exists a constant $C>0$ depending on $i, j,k, l$ but not on $\bd{R}$ such that  $|\mbox{Cov}(X_1^{2i-1}X_2^{2j-1}, X_3^{2k}X_4^{2l})|\leq C\sum_{1\leq i<j\leq 4}r_{ij}^2$ for all non-negative integers $i,j,k,l$ with $i\geq 1$ and $j\geq 1$.
\end{lemma}
\noindent\textbf{Proof of Lemma~\ref{dhwnal}}. If $k+l=0$, then $k=l=0$, and hence $X_3^{2k}X_4^{2l}=1$. The conclusion trivially holds. So we assume   and $k+l\geq 1$ next.

By Lemma~\ref{wdh120}(i),  there exists a polynomial $f(x)$ for $x\in \mathbb{R}$ with $f(0)=1$ and
\bea\lbl{fiure0}
E(X_1^{2i-1}X_2^{2j-1})=(2i-1)!!\cdot (2j-1)!!\cdot r_{12} f(r_{12}^2)
\eea
for all integers $i\geq 1$ and $j\geq 1$. Write
\bea\lbl{vee958}
&&\mbox{Cov}(X_1^{2i-1}X_2^{2j-1}, X_3^{2k}X_4^{2l}) \nonumber\\
&=& E\big(X_1^{2i-1}X_2^{2j-1} X_3^{2k}X_4^{2l}\big)-E\big(X_1^{2i-1}X_2^{2j-1}\big)\cdot E\big(X_3^{2k}X_4^{2l}\big).
\eea
 Then, by Lemma~\ref{Wick_formula},
\bea\lbl{cl120}
E\big(X_1^{2i-1}X_2^{2j-1} X_3^{2k}X_4^{2l}\big)=D_0+\sum_{1\leq i< j \leq 4}D_{ij}r_{ij}+ d(r),
\eea
where $D_0$ and $D_{ij}$'s are constants and  $d(r)$ is a linear combination  of $\prod_{1\leq i<j\leq 4}r_{ij}^{\alpha_{ij}}$ with $2\leq \sum_{1\leq i<j \leq 4}\alpha_{ij}\leq i+j+k+l-1.$ We now evaluate the values of $D_0$ and $D_{ij}$'s.

Take $r_{ij}=0$ for all $i\ne j$. Then $X_1, X_2, X_3, X_4$ are i.i.d. with distribution $N(0, 1)$. From \eqref{cl120}, we see
\beaa
D_0=E(X_1^{2i-1})\cdot E(X_2^{2j-1})\cdot E( X_3^{2k})\cdot E(X_4^{2l})=0
\eeaa
since $E(X_1^{2i-1})=0.$ Then \eqref{cl120} becomes
\bea\lbl{u2ei0}
E\big(X_1^{2i-1}X_2^{2j-1} X_3^{2k}X_4^{2l}\big)=\sum_{1\leq i< j \leq 4}D_{ij}r_{ij}+ d(r).
\eea
Take $r_{ij}=0$ for all $1\leq i<j\leq 4$ except $(i, j)=(1,2)$. Then the three random variables $(X_1, X_2)^T\in \mathbb{R}^2,\, X_3$ and $X_4$ are  independent. It follows from \eqref{fiure0} and \eqref{u2ei0} that
\beaa
(2i-1)!!\cdot (2j-1)!!\cdot  (2k-1)!!\cdot (2l-1)!!\cdot  f(r_{12}^2)r_{12}=D_{12}r_{12}+d_1(r_{12})
\eeaa
where $d_1(x)=c_2x^2+c_3x^3+\cdots + c_px^p$ with $2\leq p \leq i+j+k+l-1$. Set $d_2(x)=c_2x+c_3x^2+\cdots + c_px^{p-1}$. The above implies
\beaa
(2i-1)!!\cdot (2j-1)!!\cdot  (2k-1)!!\cdot (2l-1)!!\cdot  f(r_{12}^2)=D_{12}+d_2(r_{12}).
\eeaa
Take $r_{12}=0$. Recall $f(0)=1$. Then
\beaa
D_{12}=(2i-1)!!\cdot (2j-1)!!\cdot  (2k-1)!!\cdot (2l-1)!!.
\eeaa
Now we claim $D_{ij}=0$ for all $(i, j)\ne (1, 2).$ In fact, set $r_{ij}=0$ for all $1\leq i < j \leq 4$ except $(i, j)=(1, 3)$. Then $X_2$ is independent of $(X_1, X_3, X_4)^T$. Use the fact $E(X_2^{2j-1})=0$ and independence to see  from \eqref{u2ei0} that
\bea\lbl{u2e1093}
0=D_{13}r_{13}+d_3(r_{13}),
\eea
where $d_3(x)=a_2x^2+\cdots + a_px^{p}$ with $2\leq p \leq i+j+k+l-1$. Divide both sides of \eqref{u2e1093} by $r_{13}$ and then set $r_{13}=0$, we get $D_{13}=0$. Similarly, $D_{ij}=0$ for all $(i, j)\ne (1, 2).$ Therefore, by \eqref{u2ei0},
\bea\lbl{hnjun}
E\big(X_1^{2i-1}X_2^{2j-1} X_3^{2k}X_4^{2l}\big)=(2i-1)!!\cdot (2j-1)!!\cdot  (2k-1)!!\cdot (2l-1)!!\cdot r_{12}+ d(r).
\eea
Since $f(x)$ in \eqref{fiure0} is a polynomial with $f(0)=1$, we are able to write  $f(x)=1+xf_1(x)$ for each $x\in \mathbb{R}$, where $f_1(x)$ ia a polynomial. Hence,
\beaa
&&E(X_1^{2i-1}X_2^{2j-1})\\
&=& (2i-1)!!\cdot (2j-1)!!\cdot r_{12}+ (2i-1)!!\cdot (2j-1)!!\cdot r_{12}^3f_1(r_{12}^2).
\eeaa
Joining this by \eqref{aiqw09}, \eqref{vee958} and \eqref{hnjun}, we have that $\mbox{Cov}(X_1^{2i-1}X_2^{2j-1}, X_3^{2k}X_4^{2l})$ is equal to
\beaa
&&(2i-1)!!\cdot (2j-1)!!\cdot  (2k-1)!!\cdot (2l-1)!!\cdot r_{12}+ d(r)\\
&-& \big[(2i-1)!!\cdot (2j-1)!!\cdot r_{12}+ (2i-1)!!\cdot (2j-1)!!\cdot r_{12}^3f_1(r_{12}^2)\big]\\
& & \cdot \big[(2k-1)!!\cdot (2l-1)!!+r_{34}^2\cdot g_2(r_{34}^2)\big]\\
&:=& d_4(r),
\eeaa
where $g_2(x)$ is a polynomial. Observe there are no linear terms of $r_{ij},\, 1\leq i<j \leq 4$ in the expression of $d_4(r)$. Consequently, $d_4(r)$ is a linear combination  of $\prod_{1\leq i<j\leq 4}r_{ij}^{\alpha_{ij}}$ with $2\leq \sum_{1\leq i<j \leq 4}\alpha_{ij}\leq i+j+k+l-1.$ The conclusion then follows from the same argument between \eqref{yw8946} and \eqref{qi105}. \hfill$\square$

\begin{lemma}\lbl{29385} Let $(X_1, X_2, X_3, X_4)^T$ be the same as in Lemma~\ref{dw0ms}.  Then, there exists a constant $C>0$ depending on $i, j,k, l$ but not on $\bd{R}$ such that $|\mbox{Cov}(X_1^{2i-1}X_2^{2j-1}, X_3^{2k-1}X_4^{2l-1})|\leq C\sum_{1\leq i<j\leq 4}r_{ij}^2$ for all positive integers $i,j,k,l$.
\end{lemma}
\noindent\textbf{Proof of Lemma~\ref{29385}}. Write
\bea\lbl{vesho8}
&&\mbox{Cov}(X_1^{2i-1}X_2^{2j-1}, X_3^{2k-1}X_4^{2l-1}) \nonumber\\
&=& E\big(X_1^{2i-1}X_2^{2j-1} X_3^{2k-1}X_4^{2l-1}\big)-E\big(X_1^{2i-1}X_2^{2j-1}\big)\cdot E\big(X_3^{2k-1}X_4^{2l-1}\big).
\eea
By Lemma~\ref{wdh120}(i),
\beaa
E(X_3^{2k-1}X_4^{2l-1})=(2k-1)!!\cdot (2l-1)!!\cdot r_{34} g(r_{34}^2)
\eeaa
for some polynomial $g(x)$. This together with \eqref{fiure0} implies that
\bea\lbl{uhqwq03}
\big|E\big(X_1^{2i-1}X_2^{2j-1}\big)\cdot E\big(X_3^{2k-1}X_4^{2l-1}\big)\big|\leq C_1|r_{12}r_{34}|\leq C\sum_{1\leq i<j \leq 4}r_{ij}^2,
\eea
where $C_1$ is a constant depending on $i, j,k, l$ but not on $\bd{R}$. By Lemma~\ref{Wick_formula},
\bea\lbl{yfoewf8e}
E\big(X_1^{2i-1}X_2^{2j-1} X_3^{2k-1}X_4^{2l-1}\big)=B_0+\sum_{1\leq i< j \leq 4}B_{ij}r_{ij}+ b(r),
\eea
where $B_0$ and $B_{ij}$'s are constants and  $b(r)$ is a linear combination  of $\prod_{1\leq i<j\leq 4}r_{ij}^{\alpha_{ij}}$ with $2\leq \sum_{1\leq i<j \leq 4}\alpha_{ij}\leq i+j+k+l-2.$ We claim that
\bea\lbl{jh1206}
B_0=B_{ij}=0
\eea
for all $1\leq i<j\leq 4$. In fact, take $r_{ij}=0$ for all $1\leq i<j\leq 4$, then $X_1, X_2, X_3, X_4$ are i.i.d. $N(0, 1)$-distributed random variables. The left hand side of \eqref{yfoewf8e} is zero because of independence. The right hand side of \eqref{yfoewf8e} is equal to $B_0$. This concludes $B_0=0$. Then \eqref{yfoewf8e} becomes
\bea\lbl{s218}
E\big(X_1^{2i-1}X_2^{2j-1} X_3^{2k-1}X_4^{2l-1}\big)=\sum_{1\leq i< j \leq 4}B_{ij}r_{ij}+ b(r).
\eea
For this identity we take $r_{ij}=0$ for all $1\leq i<j\leq 4$ except $(i, j)=(1,2)$. Then the three random variables $(X_1, X_2)^T\in \mathbb{R}^2,\, X_3$ and $X_4$ are  independent. By independence and the fact $E(X_3^{2k-1})=0$, the left hand side of \eqref{s218} is zero. Hence, \eqref{s218} is reduced to
\beaa
0=B_{12}r_{12}+b_1(r_{12}),
\eeaa
where $b_1(x)=c_2x^2+c_3x^3+\cdots + c_px^p$ with $2\leq p \leq i+j+k+l-2$. Set $b_2(x)=c_2x+c_3x^2+\cdots + c_px^{p-1}$. Then $0=B_{12}+b_2(r_{12}).$ Take $r_{12}=0$ to get $B_{12}=0$. By symmetry, we know $B_{ij}=0$ for all $1\leq i<j\leq 4$. So \eqref{jh1206} has been verified. It follows that \eqref{yfoewf8e} is reduced to
$E\big(X_1^{2i-1}X_2^{2j-1} X_3^{2k-1}X_4^{2l-1}\big)= b(r)$. By employing the same argument between \eqref{yw8946} and \eqref{qi105}, we get
\beaa
\big|E\big(X_1^{2i-1}X_2^{2j-1} X_3^{2k-1}X_4^{2l-1}\big)\big| \leq C_2\sum_{1\leq i<j\leq 4}r_{ij}^2,
\eeaa
where $C_2$ is a constant depending on $i, j,k, l$ but not on $\bd{R}$. This joined~\eqref{vesho8} and \eqref{uhqwq03} yields the desired inequality. \hfill$\square$

\medskip

The following fact supplies a convenient tool to handle the covariance between products of independent random variables. Recall $\|\xi\|_q=(E|\xi|^q)^{1/q}$ for any random variable $\xi$ and $q\geq 1$.
\begin{lemma}\lbl{xintong} For $k\geq 2$, let $\{(U_i, V_i)\in \mathbb{R}^2;\, 1\leq i \leq k\}$ be independent random vectors. Assume $C:=3\prod_{i=1}^k(1+\|U_i\|_k)(1+\|V_i\|_k)<\infty$ and $K=C^k$. Then
\beaa
\Big|\mbox{Cov}\Big(\prod_{i=1}^kU_i, \prod_{i=1}^kV_i\Big)\Big|\leq K\cdot \sum_{i=1}^k|\mbox{Cov}(U_i, V_i)|.
\eeaa
\end{lemma}
\noindent\textbf{Proof of Lemma~\ref{xintong}}. The inequality holds obviously for $k=1$. Assume now $k\geq 2$. By definition,
\bea\lbl{dh205}
\mbox{Cov}(U_1U_2, V_1V_2)
&=& E(U_1U_2 V_1V_2)-E(U_1U_2)\cdot E(V_1V_2)\nonumber\\
& = & E(U_1V_1)\cdot E(U_2V_2)-EU_1\cdot EU_2\cdot EV_1\cdot EV_2.
\eea
Write $E(U_1V_1)=\mbox{Cov}(U_1, V_1)+EU_1\cdot EV_1$ and $E(U_2V_2)=\mbox{Cov}(U_2, V_2)+EU_2\cdot EV_2$. Plug the two identities into~\eqref{dh205} to see
\bea\lbl{ajs6}
\mbox{Cov}(U_1U_2, V_1V_2)
&=& \mbox{Cov}(U_1, V_1)\cdot \mbox{Cov}(U_2, V_2)+EU_2\cdot EV_2\cdot \mbox{Cov}(U_1, V_1)\nonumber\\
&&  ~~~~~~~~~~~~~~~~~~~~~~~~~~~~~~~~+EU_1\cdot EV_1\cdot \mbox{Cov}(U_2, V_2).
\eea
Since $k\geq 2$, we have from the H\"{o}lder inequality that $E|\xi|\leq \|\xi\|_F\leq \|\xi\|_k$ for any random variable $\xi$. By the Cauchy-Schwartz inequality, $E|U_1V_1|\leq \|U_1\|_2\cdot \|V_1\|_2\leq \|U_1\|_k\cdot \|V_1\|_k$. As a result, $|\mbox{Cov}(U_1, V_1)|$ is bounded by
\bea\lbl{efjwl3}
 E|U_1V_1|+E|U_1|\cdot E|V_1| \leq 2\|U_1\|_k \|V_1\|_k.
\eea
Set $\tau_j=3\prod_{i=1}^j(1+\|U_i\|_k)(1+\|V_i\|_k)$ for $j=1, \cdots, k$. We claim that
\bea\lbl{euoewf9}
\Big|\mbox{Cov}\Big(\prod_{i=1}^kU_i, \prod_{i=1}^kV_i\Big)\Big|\leq (\tau_1\cdots \tau_k)\cdot \sum_{i=1}^k|\mbox{Cov}(U_i, V_i)|.
\eea
In fact, by applying \eqref{efjwl3} to \eqref{ajs6}, we see
\bea
&&|\mbox{Cov}(U_1U_2, V_1V_2)|\nonumber \\
&\leq & 2\|U_1\|_k \|V_1\|_k \cdot|\mbox{Cov}(U_2, V_2)|+\|U_2\|_k \|V_2\|_k\cdot |\mbox{Cov}(U_1, V_1)|+ \|U_1\|_k\|V_1\|_k\cdot |\mbox{Cov}(U_2, V_2)| \nonumber\\
&\leq & \Big[3\prod_{i=1}^2\big(1+\|U_i\|_k\big)\big(1+\|V_i\|_k\big)\Big]\cdot \big[ |\mbox{Cov}(U_1, V_1)|+ |\mbox{Cov}(U_2, V_2)|\big].\lbl{209r5}
\eea
So claim \eqref{euoewf9} holds for $k=2$ due to the fact $\tau_j\geq 1$ for each $j$. Now we assume $k\geq 3$ and use induction to complete the proof. Assume
\bea\lbl{weuo32r09}
\Big|\mbox{Cov}\Big(\prod_{i=1}^jU_i, \prod_{i=1}^jV_i\Big)\Big|\leq (\tau_1\cdots \tau_j)\cdot \sum_{i=1}^j|\mbox{Cov}(U_i, V_i)|
\eea
for some $2\leq j < k.$ By assumption, $(U_{j+1}, V_{j+1})$ and $(\prod_{i=1}^jU_i, \prod_{i=1}^jV_i)$ are independent. We obtain from~\eqref{209r5} that
\bea
&& \Big|\mbox{Cov}\Big(\prod_{i=1}^{j+1}U_i, \prod_{i=1}^{j+1}V_i\Big)\Big|\nonumber\\
& \leq & \Big[3\big(1+\|U_{j+1}\|_k\big)\big(1+\|V_{j+1}\|_k\big)\big(1+\|U_1\cdots U_j\|_k\big)\big(1+\|V_1\cdots V_j\|_k\big)\Big]\cdot \nonumber\\
&&~~~~~~~~~~~~~~~~~~~~~~~~~~~~~~~\Big[ |\mbox{Cov}(U_{j+1}, V_{j+1})|+\Big|\mbox{Cov}\Big(\prod_{i=1}^jU_i, \prod_{i=1}^jV_i\Big)\Big|\Big].\lbl{wfqiu92}
\eea
Since $2\leq j < k$, then we have from the H\"{o}lder inequality that
\beaa
\|U_1\cdots U_j\|_k=\big[E\big(U_1\cdots U_j\cdot \underbrace{1\cdots 1}_{k-j}\big)\big]^{1/k}
\leq \prod_{i=1}^j\|U_i\|_k.
\eeaa
The above also holds if the symbol ``$U$" is replaced by ``$V$". Thus,
\beaa
3\big(1+\|U_1\cdots U_j\|_k\big)\big(1+\|V_1\cdots V_j\|_k\big)\big(1+\|U_{j+1}\|_k\big)\big(1+\|V_{j+1}\|_k\big)\leq \tau_{j+1}.
\eeaa
By \eqref{weuo32r09} and \eqref{wfqiu92},
\beaa
\Big|\mbox{Cov}\Big(\prod_{i=1}^{j+1}U_i, \prod_{i=1}^{j+1}V_i\Big)\Big|
&\leq & \tau_{j+1}\cdot \Big[|\mbox{Cov}(U_{j+1}, V_{j+1})|+(\tau_1\cdots \tau_j)\cdot \sum_{i=1}^j|\mbox{Cov}(U_i, V_i)|\Big]\\
& \leq & (\tau_1\cdots \tau_{j+1})\cdot \sum_{i=1}^{j+1}|\mbox{Cov}(U_i, V_i)|.
\eeaa
This confirms \eqref{euoewf9}. The proof is completed by taking $C=\tau_k$ and $K=C^k.$  \hfill$\square$

\subsubsection{Combinatorics}\lbl{Combinatorics111}

In this section we will work on some combinatorics problems. They will be used to evaluate    covariances between squared sample correlations coefficients in Section~\ref{Squared_SCC}. We always assume $\alpha_1,\cdots,\alpha_m, \beta_1,\cdots, \beta_m, \gamma_1,\cdots, \gamma_m, \delta_1, \cdots, \delta_m$ are non-negative integers.
Set $\bm{\alpha}=(\alpha_1,\cdots,\alpha_m)$, $\bm{\beta}= (\beta_1,\cdots, \beta_m)$, $\bm{\gamma}=(\gamma_1,\cdots, \gamma_m)$ and $\bm{\delta}=(\delta_1, \cdots, \delta_m)$.

\begin{lemma}\lbl{wqjiefg1k} Let $m\geq 2$, $a\geq 0$ and $b\geq 1$  be  integers. Then the following hold with constant $K$ depending on $a$ and $b$ but not $m$.

(i) Let $N_1$ be the total number of non-negative integer solutions $(x_1, \cdots, x_m)$ of $x_1+\cdots+x_m=a$, then $N_1\leq Km^a$.

(ii) Let $N_2$ be the total number of non-negative integer solutions $(x_1, \cdots, x_m)$ of $x_1+\cdots+x_m=a$ with $x_1\geq b$. Then $N_2\leq Km^{a-b}$.

(iii) Given $1\leq n<m$ and $c_1\geq 1, \cdots, c_n\geq 1$ with $c_1+\cdots +c_n\leq a$,
let $N_3$ the total number of non-negative integer solutions $(x_1, \cdots, x_m)$ of $x_1+\cdots+x_m=a$ with $x_i\geq c_i$ for  $1\leq i\leq n$.   Then $N_3\leq K m^{a-c_1-\cdots - c_n}$.

\end{lemma}
A quick comment is that (ii) is a special case of (iii). We single it out because $N_2$ has a much neater statement and it will be used very frequently.

\medskip

\noindent\textbf{Proof of Lemma~\ref{wqjiefg1k}}. (i) If $a=0$, the only non-negative integer solution of $x_1+\cdots+x_m=a$ is $(0,\cdots, 0)$. Then $N_1=1$ and the conclusion follows with any constant $K\geq 1.$ We assume next that $a\geq 1$. It is well-known that
\bea\lbl{gyg20}
N_1=\binom{m+a-1}{a}\leq (m+a-1)^a\leq (1+a)^am^a.
\eea

(iii)
Set $y_i=x_i-c_i$ for $i=1,\cdots, n$ and $y_i=x_i$ for $n+1\leq i\leq m$. Then $N_3$ is equal to the total number of non-negative integer solutions $(y_1, \cdots, y_m)$ of $y_1+\cdots + y_m=a-c_1-\cdots - c_n.$ From (i) we see $N_3\leq K m^{a-c_1-\cdots - c_n}$. The proof is completed.

The statement (ii) follows because it is a special case of (iii). \hfill$\square$

\medskip

In the following when we say a non-negative integer solution $(\bm{\alpha}, \bm{\beta}, \bm{\gamma}, \bm{\delta})$ of a certain equation, we mean $(\alpha_1, \cdots, \alpha_m, \beta_1, \cdots, \beta_m, \gamma_1, \cdots, \gamma_m, \delta_1, \cdots, \delta_m)$ satisfies that equation with each of $\{\alpha_i, \beta_i, \gamma_i, \delta_i;\, 1\leq i \leq m\}$ being a non-negative integer.

\begin{lemma}\lbl{wqjiefg} Let $m\geq 4$ and $\alpha, \beta, \gamma, \delta$  be non-negative integers. Let $N_1$ be the total number of non-negative integer solutions of $(\bm{\alpha}, \bm{\beta}, \bm{\gamma}, \bm{\delta})$ satisfying
\bea\lbl{husoq922}
\sum_{i=1}^m\alpha_i=\alpha,\ \ \sum_{i=1}^m\beta_i=\beta,\ \ \sum_{i=1}^m\gamma_i=\gamma,\ \  \sum_{i=1}^m\delta_i=\delta.
\eea
Set $I_1:= \{1\}\cup\{2\leq i\leq m;\, (\bm{\alpha}, \bm{\beta}, \bm{\gamma}, \bm{\delta})\ \mbox{satisfies}\ \eqref{husoq922}\ \mbox{and}\ \alpha_i+\beta_i\geq 1\}$ and
\beaa
I_2:= \{2, 3\}\cup\big\{i\in\{1,4,5,\cdots, m\};\, (\bm{\alpha}, \bm{\beta}, \bm{\gamma}, \bm{\delta})\ \mbox{satisfies}\ \eqref{husoq922}\ \mbox{and}\ \gamma_i+\delta_i\geq 1\big\}.
\eeaa
Let $N_2$ be the total number of non-negative integer solutions of $(\bm{\alpha}, \bm{\beta}, \bm{\gamma}, \bm{\delta})$ satisfying \eqref{husoq922} and $I_1\cap I_2\ne \emptyset$.
Then, there exists a constant $K$ depending on $\alpha, \beta, \gamma, \delta$ but not  $m$ such that (i) $N_1\leq K\cdot m^{\alpha+\beta+\gamma+\delta}$; (ii)  $N_2\leq K\cdot m^{\alpha+\beta+\gamma+\delta-1}.$
\end{lemma}
\noindent\textbf{Proof of Lemma~\ref{wqjiefg}}. First, recall the fact that the total number of non-negative integer solutions of $x_1+\cdots + x_m=k$ for any non-negative integer $k$ is $\binom{m+k-1}{k}.$ Therefore, considering the four equations in \eqref{husoq922} separately, the total numbers of non-negative integer solutions are
\bea\lbl{ehuweui}
\binom{m+\alpha-1}{\alpha},\
\binom{m+\beta-1}{\beta},\ \binom{m+\gamma-1}{\gamma}\ \ \mbox{and}\ \ \binom{m+\delta-1}{\delta},
\eea
respectively.

(i) By \eqref{gyg20} and \eqref{ehuweui},
\bea
N_1 &\leq & \binom{m+\alpha-1}{\alpha}\binom{m+\beta-1}{\beta}
\binom{m+\gamma-1}{\gamma}\binom{m+\delta-1}{\delta} \nonumber\\
&\leq & (1+\alpha)^{\alpha}(1+\beta)^{\beta}(1+\gamma)^{\gamma}(1+\delta)^{\delta}\cdot m^{\alpha+\beta+\gamma+\delta}.
\eea
The conclusion follows by taking $K=(1+\alpha)^{\alpha}(1+\beta)^{\beta}(1+\gamma)^{\gamma}(1+\delta)^{\delta}$.

(ii) If $\alpha+\beta=0$ and $\gamma+\delta=0$, then $I_1=\{1\}$ and $I_2= \{2, 3\}$, and hence  $I_1\cap I_2= \emptyset$. Thus $N=0$. The conclusion holds.  So we assume next that either $\alpha+\beta\geq 1$ or $\gamma+\delta\geq 1$.

Notice $I_1\cap I_2=A_1 \cup A_2 \cup A_3$, where
\beaa
A_1&=&\{i\in\{1\};\, (\bm{\alpha}, \bm{\beta}, \bm{\gamma}, \bm{\delta})\ \mbox{satisfies}\ \eqref{husoq922}\ \mbox{and}\ \gamma_i+\delta_i\ne 0\big\};\\
A_2 &=&  \{\{i\in \{2, 3\};\, (\bm{\alpha}, \bm{\beta}, \bm{\gamma}, \bm{\delta})\ \mbox{satisfies}\ \eqref{husoq922}\ \mbox{and}\ \alpha_i+\beta_i\ne 0\};\\
A_3 &=& \big\{i\in\{4,5,\cdots, m\};\, (\bm{\alpha}, \bm{\beta}, \bm{\gamma}, \bm{\delta})\ \mbox{satisfies}\ \eqref{husoq922}\ \mbox{and}\ \alpha_i+\beta_i\ne 0\ \mbox{and}\ \gamma_i+\delta_i\ne 0\big\}.
\eeaa
Hence, if $I_1\cap I_2 \ne \emptyset$, then either $A_1\ne\emptyset$, $A_2\ne\emptyset$ or $A_3\ne\emptyset$. Let us consider the three scenarios one by one next.

\textbf{Scenario 1: $A_1\ne\emptyset$}. In this situation, $\gamma_1+\delta_1 \geq 1$. Consequently, either $\gamma_1\geq 1$ or $\delta_1\geq 1$. Thus, taking $b=1$ in Lemma~\ref{wqjiefg1k}(i) and (ii), we know the total number of non-negative integer solutions of $(\bm{\alpha}, \bm{\beta}, \bm{\gamma}, \bm{\delta})$ satisfying \eqref{husoq922} and $A_1\ne\emptyset$ is bounded by
\beaa
&& K_1m^{\alpha}\cdot K_1m^{\beta} \cdot K_1m^{\gamma-1}\cdot K_1m^{\delta}+ K_1m^{\alpha}\cdot K_1m^{\beta} \cdot K_1m^{\gamma}\cdot K_1m^{\delta-1}\\
&=& 2(K_1)^4\cdot m^{\delta+\beta+\gamma+\delta-1}
\eeaa
where $K_1$ here and below is a constant depending on $\alpha, \beta, \gamma, \delta$ but not $m$.

\textbf{Scenario 2: $A_2\ne\emptyset$}. In this situation, either $\alpha_2+\beta_2\geq 1$ or $\alpha_3+\beta_3\geq 1$. Similar to the first case, the total number of non-negative integer solutions of $(\bm{\alpha}, \bm{\beta}, \bm{\gamma}, \bm{\delta})$ satisfying \eqref{husoq922} and $A_2\ne\emptyset$ is bounded by $2(K_1)^4\cdot m^{\delta+\beta+\gamma+\delta-1} + 2(K_1)^4\cdot m^{\delta+\beta+\gamma+\delta-1}=4(K_1)^4\cdot m^{\delta+\beta+\gamma+\delta-1}$.

\textbf{Scenario 3: $A_3\ne\emptyset$}. In this situation,  there exists $i\in\{4,5,\cdots, m\}$ such that $\alpha_i+\beta_i\geq 1$ and $\gamma_i+\delta_i\geq 1$. For fixed $i$, if $\alpha_i+\beta_i\geq 1$ and $\gamma_i+\delta_i\geq 1$, then one of the four cases must be true: (a) $\alpha_i\geq 1$ and $\gamma_i\geq 1$; (b) $\alpha_i\geq 1$ and $\delta_i\geq 1$; (c) $\beta_i\geq 1$ and $\gamma_i\geq 1$; (d) $\beta_i\geq 1$ and $\delta_i\geq 1$. From Lemma~\ref{wqjiefg1k}(i) and (ii) again, we have  that the total number of non-negative integer solutions of $(\bm{\alpha}, \bm{\beta}, \bm{\gamma}, \bm{\delta})$ satisfying \eqref{husoq922} and (a) is dominated by
$K_1m^{\alpha-1}\cdot K_1m^{\beta} \cdot K_1m^{\gamma-1}\cdot K_1m^{\delta}=
(K_1)^4\cdot m^{\delta+\beta+\gamma+\delta-2}$. By symmetry, the same inequality holds if ``(a)" is replaced by (b), (c) and (d), respectively. In conclusion, for fixed $i$, the total number of non-negative integer solutions of $(\bm{\alpha}, \bm{\beta}, \bm{\gamma}, \bm{\delta})$ satisfying \eqref{husoq922} and $\alpha_i+\beta_i\geq 1$ and $\gamma_i+\delta_i\geq 1$ is controlled by $4(K_1)^4\cdot m^{\delta+\beta+\gamma+\delta-2}$. Now, $i\in\{4,5,\cdots, m\}$ has at most $m$ choices. Then the total number of non-negative integer solutions of $(\bm{\alpha}, \bm{\beta}, \bm{\gamma}, \bm{\delta})$ satisfying \eqref{husoq922} and $A_3\ne\emptyset$ is bounded by $m\cdot 4(K_1)^4\cdot m^{\delta+\beta+\gamma+\delta-2}=4(K_1)^4\cdot m^{\delta+\beta+\gamma+\delta-1}.$

Finally, add the bounds up in the above three scenarios, we get $N_2\leq 10(K_1)^4\cdot m^{\delta+\beta+\gamma+\delta-1}$. The proof is completed by taking $K=10(K_1)^4$. \hfill$\square$

\begin{lemma}\lbl{cjhew0} Assume $m\geq 5$ and $\alpha, \beta, \gamma, \delta$ are non-negative integers. Define
\beaa
S = \{5\leq i \leq m; (\bm{\alpha}, \bm{\beta}, \bm{\gamma}, \bm{\delta})\ \mbox{satisfies}\ \eqref{husoq922},\,  \alpha_i+\beta_i\geq 1\ \mbox{and}\  \gamma_i+\delta_i\geq 1\}.
\eeaa
 Then the following statements hold with a constant $K$ depending on $\alpha, \beta, \gamma, \delta$ but not $m$.

(i) The total number of solutions of  $(\bm{\alpha}, \bm{\beta}, \bm{\gamma}, \bm{\delta})$ satisfying \eqref{husoq922} and $S\ne \emptyset$ is bounded by $Km^{\alpha+\beta+\gamma+\delta-1}.$

(ii) The total number of solutions of $(\bm{\alpha}, \bm{\beta}, \bm{\gamma}, \bm{\delta})$ satisfying~\eqref{husoq922} with $\gamma_1+\delta_1 \geq 2$ is bounded by
$K\cdot m^{\alpha+\beta+\gamma+\delta-2}$.
\end{lemma}
\noindent\textbf{Proof of Lemma~\ref{cjhew0}}. (i) Since $S\ne \emptyset$, then there exists some  $5\leq i \leq m$ such that $\alpha_i+\beta_i\geq 1$ and  $\gamma_i+\delta_i\geq 1$.
According to Scenario 3 in the proof of Lemma~\ref{wqjiefg}, the total number of solutions of  $(\bm{\alpha}, \bm{\beta}, \bm{\gamma}, \bm{\delta})$ satisfying \eqref{husoq922}, $\alpha_i+\beta_i\geq 1$ and  $\gamma_i+\delta_i\geq 1$ is dominated by $K\cdot m^{\delta+\beta+\gamma+\delta-2}$, where $K$ is a constant depending on $\alpha, \beta, \gamma, \delta$ but not $m$. Noticing $5\leq i \leq m$, then the desired number is bounded by $(m-4)\cdot Km^{\delta+\beta+\gamma+\delta-2}\leq Km^{\delta+\beta+\gamma+\delta-1}$.

(ii) Let $K$ be a constant in Lemma~\ref{wqjiefg1k}(i) with $a=\alpha$ or $\beta$. Also, the $K$ satisfies Lemma~\ref{wqjiefg1k}(ii) with $a\in \{\gamma, \delta\}$ and $b\in \{1,2\}$. Since $\gamma_1+\delta_1 \geq 2$, then one of the three cases must be true: (a) $\gamma_1\geq 2$, (b) $\delta_1\geq 2$ or  (c) $\gamma_1\geq 1$ and $\delta_1\geq 1$ simultaneously. By Lemma~\ref{wqjiefg1k}(ii), the total number of solutions  $(\bm{\gamma}, \bm{\delta})$ of the last two equations from \eqref{husoq922} with $\gamma_1\geq 2$ is no more than $Km^{\gamma+\delta-2}$. The same holds if ``$\gamma_1\geq 2$" is replaced by ``$\delta_1\geq 2$". Similarly, by Lemma~\ref{wqjiefg1k}(ii) again, the total number of solutions  $(\bm{\gamma}, \bm{\delta})$ of the last two equations from \eqref{husoq922} satisfying $\gamma_1\geq 1$ and $\delta_1\geq 1$ is bounded by $Km^{\gamma-1}\cdot Km^{\delta-1}=K^2m^{\gamma+\delta-2}$. Consequently, the total number of solutions of $(\bm{\alpha}, \bm{\beta}, \bm{\gamma}, \bm{\delta})$ satisfying~\eqref{husoq922} with $\gamma_1+\delta_1 \geq 2$ is bounded by
\beaa
Km^{\alpha}\cdot Km^{\beta}\cdot (Km^{\gamma+\delta-2}+Km^{\gamma+\delta-2}+K^2m^{\gamma+\delta-2})
=(K^4+2K^2)m^{\alpha+\beta+\gamma+\delta-2}.
\eeaa
Therefore, the desired conclusion follows by regarding $K^4+2K^2$ as new constant $K$. \hfill$\square$

\begin{lemma}\lbl{wejheif4} Assume $m\geq 5$ and $\alpha, \beta, \gamma, \delta$ are non-negative integers.    Define
\beaa
S &=& \{i\in \{3, 4\}; (\bm{\alpha}, \bm{\beta}, \bm{\gamma}, \bm{\delta})\ \mbox{satisfies}\ \eqref{husoq922}\ \mbox{and}\ \alpha_i+\beta_i\geq 1\}\cup\\
&& \{j\in \{1, 2\}; (\bm{\alpha}, \bm{\beta}, \bm{\gamma}, \bm{\delta})\ \mbox{satisfies}\ \eqref{husoq922}\ \mbox{and}\ \gamma_j+\delta_j\geq 1\}.
\eeaa
 Let $T_{m,1}$ be the set of $(\bm{\alpha}, \bm{\beta}, \bm{\gamma}, \bm{\delta})$ satisfying \eqref{husoq922} and $|S|=1$. Let  $T_{m,2}$ be the set of $(\bm{\alpha}, \bm{\beta}, \bm{\gamma}, \bm{\delta})$ satisfying \eqref{husoq922} and $|S|\geq 2$. Let $T_{m, 3}$ be the set of $(\bm{\alpha}, \bm{\beta}, \bm{\gamma}, \bm{\delta})$ satisfying \eqref{husoq922}, $|S|=1$ and one of the following:

(1) $\alpha_i+\beta_i\geq 1$ for some $i\in \{1,2 \}$;

(2) $\gamma_j+\delta_j\geq 1$ for some $j \in \{3, 4\}$;

(3) $\alpha_k+\beta_k\geq 2$ for some $5\leq k\leq m$;

(4) $\gamma_l+\delta_l\geq 2$ for some $5\leq l\leq m$;

(5) $\alpha_t+\beta_t=1$\ \mbox{and}\ $\gamma_t+\delta_t=1$ simultaneously for some $5\leq t\leq m$. Then, there exists a constant $K$ depending on $\alpha, \beta, \gamma, \delta$ but not $m$ such that $|T_{m, 1}|\leq Km^{\alpha+\beta+\gamma+\delta-1}$, $|T_{m, 2}|\leq Km^{\alpha+\beta+\gamma+\delta-2}$ and $|T_{m, 3}|\leq K\cdot m^{\alpha+\beta+\gamma+\delta-2}.$
\end{lemma}
\noindent\textbf{Proof of Lemma~\ref{wejheif4}}. If $\alpha=\beta=\gamma=\delta=0$, then $T_{m,1}=T_{m,2}=T_{m,3}=\emptyset$, the conclusion obviously holds.  So we assume next that at least one of the four numbers is positive. Note that the bounds in the conclusions are $Km^{\alpha+\beta+\gamma+\delta-1}$ and $Km^{\alpha+\beta+\gamma+\delta-2}$. So, in case one of $\{\alpha, \beta, \gamma, \delta\}$ is zero, say, $\alpha=0$, any discussions below related to $\alpha$  will disappear by convention. In the following we will prove the three conclusions one by one.

{\it The bound for $T_{m,1}$}.  If $|S|=1$, then one of the following four situations must occur:  (a) $\alpha_i\geq 1$ for some $i\in \{3, 4\}$;  (b) $\beta_i\geq 1$ for some $i\in \{3, 4\}$;  (c) $\gamma_j\geq 1$ for some $j\in \{1, 2\}$; (d) $\delta_j\geq 1$ for some $j\in \{1, 2\}$. If $\alpha_i\geq 1$, by Lemma~\ref{wqjiefg1k}(ii), the total number of non-negative integer solutions $\bm{\alpha}$ of $\alpha_1+\cdots + \alpha_m=\alpha$ is no more than $K_1m^{\alpha-1}$. Here and later $K_1$ represents a constant depending on $\alpha, \beta, \gamma, \delta$ but not  $m$, and could be different from line to line. By Lemma~\ref{wqjiefg1k}(i), the total number of points $(\bm{\alpha}, \bm{\beta}, \bm{\gamma}, \bm{\delta})$ satisfying \eqref{husoq922} and (a) is controlled by
\beaa
 K_1m^{\alpha-1}\cdot K_1m^{\beta}\cdot K_1m^{\gamma}\cdot K_1m^{\delta}=(K_1)^4\cdot m^{\alpha+\beta+\gamma+\delta-1}.
\eeaa
Likewise the same bound holds if ``(a)" is replaced by ``(b)", ``(c)" or ``(d)".  This implies $|T_{m, 1}|$ is dominated by the sum of the four bounds, that is, $4(K_1)^4\cdot m^{\alpha+\beta+\gamma+\delta-1}.$

{\it The bound for $T_{m, 2}$}. The assumption $|S|\geq 2$ implies one of the following three statements must be true: (e) $\alpha_3+\beta_3\geq 1$ and $\alpha_4+\beta_4\geq 1$; (f) $\gamma_1+\delta_1\geq 1$ and $\gamma_2+\delta_2\geq 1$; (g) $\alpha_i+\beta_i\geq 1$ for some $i\in \{3, 4\}$ and $\gamma_j+\delta_j$ for some $j\in \{1, 2\}$.

Under (e), one of the next four  cases has to be true: (e1) $\alpha_3\geq 1$ and $\alpha_4\geq 1$; (e2) $\alpha_3\geq 1$ and $\beta_4\geq 1$; (e3) $\alpha_4\geq 1$ and  $\beta_3\geq 1$; (e4) $\beta_3\geq 1$ and $\beta_4\geq 1$. By Lemma~\ref{wqjiefg1k}(i) and (ii), the total number of points $(\bm{\alpha}, \bm{\beta}, \bm{\gamma}, \bm{\delta})$ satisfying \eqref{husoq922} and (e1)  is no more than
\beaa
 K_1m^{\alpha-2}\cdot K_1m^{\beta}\cdot K_1m^{\gamma}\cdot K_1m^{\delta}=(K_1)^4\cdot m^{\alpha+\beta+\gamma+\delta-2}.
\eeaa
By the same spirit, the total number of points $(\bm{\alpha}, \bm{\beta}, \bm{\gamma}, \bm{\delta})$ satisfying \eqref{husoq922} and (e2) is controlled by
\beaa
 K_1m^{\alpha-1}\cdot K_1m^{\beta-1}\cdot K_1m^{\gamma}\cdot K_1m^{\delta}=(K_1)^4\cdot m^{\alpha+\beta+\gamma+\delta-2}.
\eeaa
By similar discussions, the same conclusion above also holds if ``(e2)" is replaced by ``(e3)" and ``(e4)", respectively, and ``$(K_1)^4$ is replaced by another polynomial of $K_1$. In conclusion, by summing the four bounds corresponding to $(e1)-(e4)$, we see that the total number of points $(\bm{\alpha}, \bm{\beta}, \bm{\gamma}, \bm{\delta})$ satisfying \eqref{husoq922} and (e) is no more than $K_1\cdot m^{\alpha+\beta+\gamma+\delta-2}.$ Similarly, the same conclusion holds if ``(e)" is replaced by ``(f)" and ``(g)", respectively. The desired conclusion is then yielded by adding up the three bounds corresponding to (e), (f) and (g).

{\it The bound for $T_{m, 3}$}. Let $A_1$ be  the set of $(\bm{\alpha}, \bm{\beta}, \bm{\gamma}, \bm{\delta})$ satisfying \eqref{husoq922}, $|S|=1$  and (1). Similarly, we define $A_2, A_3, A_4, A_5$ with ``(1)" is replaced by ``(2)", ``(3)", ``(4)", ``(5)", respectively. It suffices to show
\bea\lbl{ewyu348}
|A_i|\leq  C_i\cdot m^{\alpha+\beta+\gamma+\delta-2}
\eea
for $i=1,2,3,4,5$, where $C_i$ is a constant depending on $\alpha, \beta, \gamma, \delta$ but not  $m$.

We first look into $A_1$ and $A_2$. Assuming $|S|=1$  and (1), then there are two possibilities: $\alpha_i+\beta_i\geq 1$ and $\alpha_j+\beta_j\geq 1$ for a pair $(i, j)$ with $i, j\in \{1,2,3,4\}$ and $i\ne j$; $\alpha_i+\beta_i\geq 1$ and $\gamma_j+\delta_j\geq 1$ for a pair $(i, j)$ with $i, j\in \{1,2\}$. Review the analysis of case (e) and (g) above and  the conclusion that the total number of points $(\bm{\alpha}, \bm{\beta}, \bm{\gamma}, \bm{\delta})$ satisfying \eqref{husoq922} and (e) is no more than $K_1\cdot m^{\alpha+\beta+\gamma+\delta-2}.$ We know \eqref{ewyu348} is true for $i=1$. By symmetry, \eqref{ewyu348} is also true for $i=2$.

Now we work on $A_3$ and $A_4$. For fixed $k\in \{5, \cdots,m\}$, the assumption   $\alpha_k+\beta_k\geq 2$ hints that either $\alpha_k\geq 2$, $\beta_k\geq 2$ or the third possibility that $\alpha_k\geq 1$ and $\beta_k\geq 1$. On the other hand, the condition $|S|=1$ implies that
either $\alpha_i+\beta_i\geq 1$ for some $i \in \{3, 4\}$ or $\gamma_j+\delta_j\geq 1$ for some $j\in \{1, 2\}$. In total we see $3\times 2=6$ scenarios. The only scenario we have not encountered so far comes from the combination $\alpha_k\geq 2$ and $\alpha_i+\beta_i\geq 1$ for some $i \in \{3, 4\}$. In this case, either $\alpha_k\geq 2$ and $\alpha_i\geq 1$ for some $i \in \{3, 4\}$ or the second possibility $\alpha_k\geq 2$ and $\beta_i\geq 1$ for some $i \in \{3, 4\}$. By Lemma~\ref{wqjiefg1k}(ii) and (iii), the total number of points $(\bm{\alpha}, \bm{\beta})$ satisfying \eqref{husoq922} and this  combination is bounded by $K_1(2m^{\alpha-3}\cdot m^{\beta}+ m^{\alpha-2}\cdot (2m^{\beta-1})=(4K_1)m^{\alpha+\beta-3}.$ By using this and earlier argument, we have the same bound for any of the six scenarios. Adding them up and noting $k$ has $m-4$ choices, we obtain \eqref{ewyu348} for $i=3.$ Similarly, \eqref{ewyu348} also holds for $i=4.$

Finally we study $A_5$. Fix $5\leq t\leq m$. Then the assumptions that $\alpha_t+\beta_t=1$\ \mbox{and}\ $\gamma_t+\delta_t=1$ have four possibilities: $\alpha_t= 1$ and $\gamma_t=1$; $\alpha_t= 1$ and $\delta_t=1$; $\beta_t= 1$ and $\gamma_t=1$; $\beta_t= 1$ and $\delta_t=1$. As aforementioned,  the condition $|S|=1$ implies that
either $\alpha_i+\beta_i\geq 1$ for some $i \in \{3, 4\}$ or $\gamma_j+\delta_j\geq 1$ for some $j\in \{1, 2\}$. So there are eight combinations with a common feature that the values of three different members of $\{\alpha_i, \beta_i, \gamma_i, \delta_i;\, 1\leq i \leq m\}$ are required to be at least $1$. By Lemma~\ref{wqjiefg1k} and the assumption that $t$ has no more than $m$ choices, we know \eqref{ewyu348} is true for $i=5.$

After the verification of \eqref{ewyu348} for $i=1,2,3,4,5$, we obtain the bound for $T_{m,3}.$ Observe the three upper bounds for $T_{m,1}$, $T_{m,2}$ and $T_{m,3}$ are involved with polynomials of $K_1$. We choose $K$ to be the maximum of the three polynomials. The whole proof is completed. \hfill$\square$

\newpage

\begin{lemma}\lbl{feufeuiqfP} Assume $m\geq 5$ and $\alpha, \beta, \gamma, \delta$ are non-negative integers. Let $S$ be the set of $(\bm{\alpha}, \bm{\beta}, \bm{\gamma}, \bm{\delta})$ satisfying \eqref{husoq922} and one of the following holds:

(i) $\alpha_i+\beta_i\geq 1$ for some $i\in \{1, 2, 3\};$

(ii) $\gamma_i+\delta_i\geq 1$ for some $i\in \{1, 2, 3\};$

(iii) $\alpha_i+\beta_i\geq 2$ or $\gamma_i+\delta_i\geq 2$ for some $4\leq i \leq m$;

(iv) $\alpha_i+\beta_i\geq 1$ and $\gamma_i+\delta_i\geq 1$ simultaneously for some $1\leq i \leq m.$\\
Then $|S|\leq Km^{\alpha+\beta+\gamma+\delta-1}$ for some constant $K$ depending on $\alpha, \beta, \gamma, \delta$ but not  $m$.
\end{lemma}
\noindent\textbf{Proof of Lemma~\ref{feufeuiqfP}}. The proof is very similar to that of Lemma~\ref{wejheif4} and is even easier. We omit the details.  \hfill$\square$

\subsubsection{Evaluation of Covariances between Polynomials of Gaussian Random Variables
}\lbl{Squared_SCC}

With the previous preparation, we are now ready to study covariances between polynomials of Gaussian random variables. The basic setting is that
\bea\lbl{dida58}
&& \mbox{Let}\ m\geq 5\ \mbox{and}\ \{(X_{1j}, X_{2j}, X_{3j}, X_{4j})^T \in \mathbb{R}^4;\, 1\leq j \leq m\}\ \mbox{be i.i.d.  random vectors} \nonumber\\
&&\mbox{with distribution}\ N_4(\bd{0}, \bd{R}),\ \mbox{where}\ \bd{R}=(r_{ij})_{4\times 4}\ \mbox{and}\ r_{ii}=1\ \mbox{for each}\ i.\ \ \
\eea
In this section, $K$ and $K_1$ always represent constants depending on $\alpha, \beta, \gamma, \delta$ but not $m$ or $\bd{R}$, and can be different from line to line. Review the notation $\bm{\alpha}, \bm{\beta}, \bm{\gamma}, \bm{\delta}$ before the statement of Lemma~\ref{wqjiefg1k}.

\begin{lemma}\lbl{fjii09f} Assume \eqref{dida58} holds. Let $\{a, b, c, d, a_i, b_i, c_i, d_i;\, 1\leq i \leq m\}$ be non-negative integers with $a=\sum_{i=1}^ma_i$, $b=\sum_{i=1}^mb_i$, $c=\sum_{i=1}^mc_i$, $d=\sum_{i=1}^md_i$. Define
$U_i=X_{i1}^{a_i}X_{i2}^{b_i}$ and $V_i=X_{i3}^{c_i}X_{i4}^{d_i}$ for $1\leq i \leq m$.
If $a_i+b_i$ and $c_i+d_i$ are both even for each $1\leq i\leq m$, then
\beaa
\Big|\mbox{Cov}\Big(\prod_{i=1}^mU_i, \prod_{i=1}^mV_i\Big)\Big| \leq K\sum_{1\leq i<j \leq 4}r_{ij}^2.
\eeaa
\end{lemma}
\noindent\textbf{Proof of Lemma~\ref{fjii09f}}. If $a+b=0$, then $U_i=1$ for each $i$, and the conclusion trivially holds. If $c+d=0$, then $V_i=1$ for each $i$, and the conclusion is still valid. Now we assume that both $a+ b\geq 1$ and $c+d \geq 1$. For any random variable $\xi$, its $L_s$-norm $\|\xi\|_s=(E|\xi|^s)^{1/s}$ is non-decreasing in $s\geq 1$ by H\"{o}lder's inequality. Furthermore, by the same inequality, since $X_{ij}\sim N(0,1)$ for each $i, j$,  we have
$E(|X_{11}|^{2a_is})\leq E(|X_{11}|^{(2a+1)s})^{a_i/(2a+1)}\leq 1+E(|X_{11}|^{(2a+1)s}).$ A similar conclusion holds for $E(|X_{11}|^{2b_is})$. Consequently,
\beaa
\|U_i\|_s^s=E\big(|X_{i1}|^{sa_i}\cdot |X_{i2}^{sb_i}|\big)
&\leq & \big[E\big(|X_{11}|^{2a_i s}\big)\big]^{1/2}\cdot \big[E\big(|X_{11}|^{2b_is}\big)\big]^{1/2}\\
&\leq & \big[1+E(|X_{11}|^{(2a+1)s})\big]\cdot \big[1+E(|X_{11}|^{(2b+1)s})\big].
\eeaa
Hence,
\bea\lbl{wufeowh}
\max\{\|U_i\|_s, \|V_i\|_s;\, 1\leq i \leq m\} \leq K_1,
\eea
where $K_1$ depends on $a, b, c, d$ and $s\geq 1$. By definition, $a=\sum_{i=1}^ma_i$, hence $|\{1\leq i \leq m;\, a_i\geq 1\}| \leq a$. The same is also true for the analogue of $b$, $c$ and $d$, respectively. Set  $\Psi=\{1\leq i \leq m;\, a_i+b_i\geq 1\, \mbox{or}\ c_i+d_i\geq 1\}$. For any $i\in \Psi$, either $a_i\geq 1$, $b_i\geq 1$, $c_i\geq 1$ or $d_i\geq 1$, it follows that $|\Psi| \leq a+b+c+d$. On the contrary, if $i\notin \Psi$ then $U_i=V_i=1$, therefore
\beaa
\mbox{Cov}\Big(\prod_{i=1}^mU_i, \prod_{i=1}^mV_i\Big) = \mbox{Cov}\Big(\prod_{i\in \Psi}U_i, \prod_{i\in \Psi}V_i\Big).
\eeaa
Set $k=|\Psi|\geq 1.$ Then $C(k):=3\prod_{i=1}^k(1+\|U_i\|_k)(1+\|V_i\|_k)\leq C(l)$ with $l=a+b+c+d$ since $\|\cdot\|_s$ is non-decreasing in $s$.  By Lemma~\ref{xintong}, there exists a constant $K>0$ depending on $a, b, c ,d $ but not $m$ such that
\beaa
\Big|\mbox{Cov}\Big(\prod_{i\in \Psi}U_i, \prod_{i\in \Psi}V_i\Big)\Big|
\leq K\cdot \sum_{i\in \Psi}\Big|\mbox{Cov}(U_i, V_i)\big|\leq K\cdot |\Psi|\cdot \max_{1\leq i\leq m}|\mbox{Cov}(U_i, V_i)|.
\eeaa
Use the fact  $|\Psi| \leq a+b+c+d$ to see
\beaa
\Big|\mbox{Cov}\Big(\prod_{i=1}^mU_i, \prod_{i=1}^mV_i\Big)\Big|
\leq (a+b+c+d)K\cdot \max_{1\leq i\leq m}|\mbox{Cov}(U_i, V_i)|.
\eeaa
For non-negative integers $x$ and $y$ with $x+y$ being even, we know that $x$ and $y$ have to be  both even or both odd. The conclusion then follows from Lemmas~\ref{dw0ms}-\ref{29385}. \hfill$\square$

\begin{lemma}\lbl{uwee9} Assume the setting in \eqref{dida58}. Define
\bea\lbl{hjwq38}
&& A_i=\frac{1}{m}\sum_{j=1}^mX_{ij}^2,\ i=1,2,3,4.
\eea
For given non-negative integers $\alpha, \beta, \gamma, \delta$ and $q\in \{1,2\}$, we have
\bea\lbl{dsy238}
 \big|\mbox{Cov}\big((X_{11}X_{21})^2A_1^{\alpha}A_2^{\beta}, (X_{3q}X_{4q})^{2}A_3^{\gamma}A_4^{\delta}\big)\big|
\leq  K \sum_{1\leq i<j \leq 4}r_{ij}^2.
\eea
\end{lemma}
\noindent\textbf{Proof of Lemma~\ref{uwee9}}. Write
\bea
&&(mA_1)^{\alpha}=\Big(\sum_{j=1}^mX_{1j}^2\Big)^{\alpha}=\sum\frac{\alpha!}{\alpha_1!\cdots \alpha_m!}X_{11}^{2\alpha_1}\cdots X_{1m}^{2\alpha_m};\lbl{eiy1}\\
&&(mA_2)^{\beta}=\Big(\sum_{j=1}^mX_{2j}^2\Big)^{\beta}=\sum\frac{\beta!}{\beta_1!\cdots \beta_m!}X_{21}^{2\beta_1}\cdots X_{2m}^{2\beta_m};\lbl{eiy2}\\
&&(mA_3)^{\gamma}=\Big(\sum_{j=1}^mX_{3j}^2\Big)^{\gamma}=\sum\frac{\gamma!}{\gamma_1!\cdots \gamma_m!}X_{31}^{2\gamma_1}\cdots X_{3m}^{2\gamma_m};\lbl{eiy3}\\
&&(mA_4)^{\delta}=\Big(\sum_{j=1}^mX_{4j}^2\Big)^{\delta}=\sum\frac{\delta!}{\delta_1!\cdots \delta_m!}X_{41}^{2\delta_1}\cdots X_{4m}^{2\delta_m},\lbl{eiy4}
\eea
where $\alpha_i$, $\beta_i$, $\gamma_i$ and $\delta_i$ are non-negative integers for each $i$ satisfying
\beaa
\sum_{i=1}^m\alpha_i=\alpha,\ \ \sum_{i=1}^m\beta_i=\beta,\ \ \sum_{i=1}^m\gamma_i=\gamma,\ \  \sum_{i=1}^m\delta_i=\delta,
\eeaa
respectively. This restriction is exactly the same as \eqref{husoq922}. To avoid repetition in the future, once this restriction is used, we will always quote \eqref{husoq922}.

First, we consider the case  $q=1$. Notice
\bea\lbl{stake8}
m^{\alpha+\beta+\gamma+\delta}\cdot \mbox{Cov}\big((X_{11}X_{21})^2A_1^{\alpha}A_2^{\beta}, (X_{3q}X_{4q})^{2}A_3^{\gamma}A_4^{\delta}\big)
\eea
is a linear combination of $N_1$ terms of the form
\beaa
\mbox{Cov}\Big(\prod_{i=1}^mU_i, \prod_{i=1}^mV_i\Big)
\eeaa
with positive coefficients no more than $\alpha!\beta!\gamma!\delta!$, where
\bea
U_1&=&X_{11}^{2\alpha_1+2}X_{21}^{2\beta_1+2}\ \ \mbox{and}\ \ U_i=X_{1i}^{2\alpha_i}X_{2i}^{2\beta_i}; \nonumber\\
V_1&=&X_{31}^{2\gamma_1+2}X_{41}^{2\delta_1+2}\ \ \mbox{and}\ \ V_i=X_{3i}^{2\gamma_i}X_{4i}^{2\delta_i} \lbl{wdqu}
\eea
for $2\leq i \leq m$. Here $N_1$ is the total number of non-negative integer solutions of $(\bm{\alpha}, \bm{\beta}, \bm{\gamma}, \bm{\delta})$ satisfying the set of equations from~\eqref{husoq922}. By Lemma~\ref{wqjiefg}(i),
\bea\lbl{gyq7}
N_1\leq K_1\cdot m^{\alpha+\beta+\gamma+\delta}.
\eea

\bea\lbl{sAHWDYT}
&&m^{\alpha+\beta+\gamma+\delta}\cdot \big|\mbox{Cov}\big((X_{11}X_{21})^2A_1^{\alpha}A_2^{\beta}, (X_{3q}X_{4q})^{2}A_3^{\gamma}A_4^{\delta}\big)\big| \nonumber\\
&\leq & (K_1\alpha!\beta!\gamma!\delta!)\cdot m^{\alpha+\beta+\gamma+\delta}\cdot \max \Big|\mbox{Cov}\Big(\prod_{i=1}^mU_i, \prod_{i=1}^mV_i\Big)\Big|,
\eea
where the maximum is taken over all $(\bm{\alpha}, \bm{\beta}, \bm{\gamma}, \bm{\delta})$ satisfying \eqref{husoq922}. Note that
\bea\lbl{gi3298gq3}
&&(2\alpha_1+2) +\sum_{i=2}^m2\alpha_i=2\alpha+2;~~ (2\beta_1+2) +\sum_{i=2}^m2\beta_i=2\beta+2; \nonumber\\
&&(2\gamma_1+2) +\sum_{i=2}^m2\gamma_i=2\gamma+2;~~ (2\delta_1+2) +\sum_{i=2}^m2\delta_i=2\delta+2.
\eea
By Lemma~\ref{fjii09f},
\bea\lbl{uh32843}
\Big|\mbox{Cov}\Big(\prod_{i=1}^mU_i, \prod_{i=1}^mV_i\Big)\Big|\leq K_1\sum_{1\leq i<j \leq 4}r_{ij}^2.
\eea
Combining this with \eqref{sAHWDYT}, we arrive at
\bea\lbl{cr3289}
&& m^{\alpha+\beta+\gamma+\delta}\cdot \big|\mbox{Cov}\big((X_{11}X_{21})^2A_1^{\alpha}A_2^{\beta}, (X_{3q}X_{4q})^{2}A_3^{\gamma}A_4^{\delta}\big)\big| \nonumber\\
&\leq & K\cdot m^{\alpha+\beta+\gamma+\delta}\cdot \sum_{1\leq i<j \leq 4}r_{ij}^2,
\eea
where $K=K_1^2\alpha!\beta!\gamma!\delta!$. So \eqref{dsy238} follows for the case  $q=1$.

For the case $q=2$, we keep $U_i$  in \eqref{wdqu} unchanged but  modify $V_i$ such that
$V_2=X_{32}^{2\gamma_2+2}X_{42}^{2\delta_2+2}$ and  $V_i=X_{3i}^{2\gamma_i}X_{4i}^{2\delta_i}$ for all $i=1,3,\cdots, m.$ By Lemma~\ref{fjii09f}, \eqref{uh32843} still holds.  From \eqref{sAHWDYT} we then get \eqref{dsy238} for the case  $q=2$. The proof is completed. \hfill$\square$

\newpage

\begin{lemma}\lbl{jeu8} Assume the setting in \eqref{dida58}.
Let $A_i$ be defined as in \eqref{hjwq38}. Given non-negative integers $\alpha, \beta, \gamma, \delta$, set
\beaa
I_m(a, b)=\mbox{Cov}\big((X_{11}X_{21})^2A_1^{\alpha}A_2^{\beta}, (X_{3a}X_{4a})(X_{3b}X_{4b})A_3^{\gamma}A_4^{\delta}\big)
\eeaa
for integers $a\geq 1$ and $b\geq 1$. Then
\beaa
|I_m(a, b)|\leq
\begin{cases}
K\sum_{1\leq i<j \leq 4}r_{ij}^2 & \text{if $(a, b)=(1, 2)$};\\
\frac{K}{m}\sum_{1\leq i<j \leq 4}r_{ij}^2 & \text{if $(a, b)=(2, 3)$}.
\end{cases}
\eeaa
\end{lemma}
\noindent\textbf{Proof of Lemma~\ref{jeu8}}. We will use the same notation as in the proof of Lemma~\ref{uwee9}. Review \eqref{husoq922} and \eqref{eiy1}. We will consider the two cases for $(a, b)$ separately, that is, $(a, b)=(1, 2)$ or $(a, b)=(2, 3)$.

{\it Case 1: $(a, b)=(1, 2)$.} Set
\bea
U_1&=&X_{11}^{2(\alpha_1+1)}X_{21}^{2(\beta_1+1)}\ \ \mbox{and}\ \ U_i=X_{1i}^{2\alpha_i}X_{2i}^{2\beta_i}\ \mbox{for}\ 2\leq i \leq m; \lbl{fuyw9}\\
V_1&=&X_{31}^{2\gamma_1+1}X_{41}^{2\delta_1+1},\ \ V_2=X_{32}^{2\gamma_2+1}X_{42}^{2\delta_2+1}\ \ \mbox{and}\ \ V_i=X_{3i}^{2\gamma_i}X_{4i}^{2\delta_i} \lbl{wdqu1}
\eea
for $3\leq i \leq m$. As before, let $N_1$ be the total number of solutions $(\bm{\alpha}, \bm{\beta}, \bm{\gamma}, \bm{\delta})$ of \eqref{husoq922} with  a bound provided in \eqref{gyq7}. Then
\beaa
m^{\alpha+\beta+\gamma+\delta}\cdot \mbox{Cov}\big((X_{11}X_{21})^2A_1^{\alpha}A_2^{\beta}, (X_{31}X_{41})(X_{32}X_{42})A_3^{\gamma}A_4^{\delta}\big)
\eeaa
is a linear combination of $N_1$ terms of the form $\mbox{Cov}(\prod_{i=1}^mU_i, \prod_{i=1}^mV_i)$
with positive coefficients no more than $\alpha!\beta!\gamma!\delta!$. From the restriction in \eqref{husoq922}, we know $\alpha_i\in \{0, \cdots, \alpha\}$, $\beta_i\in \{0, \cdots, \beta\}$, $\gamma_i\in \{0, \cdots, \gamma\}$ and $\delta_i\in \{0, \cdots, \delta\}$ for each $i$. By Lemma~\ref{fjii09f} and a discussion similar to \eqref{gi3298gq3}, we obtain
\beaa
\max\Big|\mbox{Cov}(\prod_{i=1}^mU_i, \prod_{i=1}^mV_i)\Big|\leq K_1\sum_{1\leq i<j \leq 4}r_{ij}^2,
\eeaa
where the maximum is taken over all $(\bm{\alpha}, \bm{\beta}, \bm{\gamma}, \bm{\delta})$ satisfying \eqref{husoq922}. By \eqref{sAHWDYT}, we obtain the bound for $I_m(1,2).$

{\it Case 2: $(a, b)=(2, 3)$.} Again,
\bea\lbl{sqwio}
m^{\alpha+\beta+\gamma+\delta}\cdot \mbox{Cov}\big((X_{11}X_{21})^2A_1^{\alpha}A_2^{\beta}, (X_{32}X_{42})(X_{33}X_{43})A_3^{\gamma}A_4^{\delta}\big)
\eea
is a linear combination of $N_1$ terms of the form $\mbox{Cov}(\prod_{i=1}^mU_i, \prod_{i=1}^mV_i)$
with positive coefficients no more than $\alpha!\beta!\gamma!\delta!$, where $U_i$ is as in \eqref{fuyw9} and
\bea
V_1'=X_{31}^{2\gamma_1}X_{41}^{2\delta_1},\ \ V_2'=X_{32}^{2\gamma_2+1}X_{42}^{2\delta_2+1},\ \ V_3'=X_{33}^{2\gamma_3+1}X_{43}^{2\delta_3+1}\ \mbox{and}\ \ V_i'=X_{3i}^{2\gamma_i}X_{4i}^{2\delta_i}\ \ \ \ \ \   \lbl{wdhqu1}
\eea
for $4\leq i \leq m$.
Set
\beaa
I_1&=& \{1\}\cup\{2\leq i \leq m;\, (\bm{\alpha}, \bm{\beta}, \bm{\gamma}, \bm{\delta})\ \mbox{satisfies}\ \eqref{husoq922}\ \mbox{and}\ \alpha_i+\beta_i\geq 1\};\\
I_2&=& \{2, 3\}\cup\big\{i\in\{1,4,5,\cdots, m\};\, (\bm{\alpha}, \bm{\beta}, \bm{\gamma}, \bm{\delta})\ \mbox{satisfies}\ \eqref{husoq922}\ \mbox{and}\ \gamma_i+\delta_i\geq  1\big\}.
\eeaa
Recalling \eqref{dida58},  $\{(U_i, V_i)^T;\, 1\leq i\leq m\}$ are independent and $m\geq 5$. Reviewing the form of $U_i$ from \eqref{fuyw9} and $V_i'$ from \eqref{wdhqu1}, we see $U_i=1$ if $\alpha_i+\beta_i=0$ for $2\leq i\leq m$ and $V_i'=1$ if $\gamma_i+\delta_i=0$ for  $i\in \{1,4,5, \cdots, m\}.$ Consequently, if  $I_1\cap I_2=\emptyset$,  Then $\alpha_2=\beta_2=\alpha_3=\beta_3=0$ and  $\gamma_1=\delta_1=0$. This says that $U_2=U_3=V_1'=1$. Also, for each $4\leq i \leq m$, the following have to be true: $\alpha_i+\beta_i=0$ if $\gamma_i+\delta_i\geq 1$ and $\gamma_i+\delta_i=0$ if $\alpha_i+\beta_i\geq 1$. These imply $\{U_i, V_i;\, 1\leq i \leq m\}$ are independent,  and hence $\mbox{Cov}(\prod_{i=1}^mU_i, \prod_{i=1}^mV_i')=0$. Thus, we only need to study the situation $I_1\cap I_2 \ne \emptyset$.
Let $N_2$ be defined as in Lemma~\ref{wqjiefg}. Thus, the quantity from \eqref{sqwio} is a linear combination of $N_2$ terms of the form $\mbox{Cov}(\prod_{i=1}^mU_i, \prod_{i=1}^mV_i)$. From Lemma~\ref{wqjiefg},  $N_2\leq K_1\cdot m^{\alpha+\beta+\gamma+\delta-1}.$
By Lemma~\ref{fjii09f} and applying the same argument of~\eqref{gi3298gq3} to \eqref{wdhqu1}, we have
\beaa
\max\Big|\mbox{Cov}(\prod_{i=1}^mU_i, \prod_{i=1}^mV_i)\Big|\leq K_1\sum_{1\leq i<j \leq 4}r_{ij}^2
\eeaa
where the maximum is taken over all $(\bm{\alpha}, \bm{\beta}, \bm{\gamma}, \bm{\delta})$ satisfying \eqref{husoq922}. Combining all of these we get
\beaa
&&m^{\alpha+\beta+\gamma+\delta}\cdot \big|\mbox{Cov}\big((X_{11}X_{21})^2A_1^{\alpha}A_2^{\beta}, (X_{32}X_{42})(X_{33}X_{43})A_3^{\gamma}A_4^{\delta}\big)\big|\\
&\leq & K_1^2(1+\alpha+\beta)\alpha!\beta!\gamma!\delta!\cdot m^{\alpha+\beta+\gamma+\delta-1}\cdot \sum_{1\leq i<j \leq 4}r_{ij}^2.
\eeaa
This gives the bound for $I_m(2, 3)$. \hfill$\square$

\begin{lemma}\lbl{dhjwe54} Assume the setting in \eqref{dida58}. Let $A_i$ be defined as in \eqref{hjwq38}. Given non-negative integers $\alpha, \beta, \gamma, \delta$, set
\bea\lbl{qwhwqoig}
J_m(a, b)=\mbox{Cov}\big((X_{11}X_{21})(X_{12}X_{22})A_1^{\alpha}A_2^{\beta}, (X_{3a}X_{4a})(X_{3b}X_{4b})A_3^{\gamma}A_4^{\delta}\big)
\eea
for integers $a\geq 1$ and $b\geq 1$. Then
\beaa
|J_m(1, 2)| \leq K\sum_{1\leq i<j \leq 4}r_{ij}^2.
\eeaa
\end{lemma}
\noindent\textbf{Proof of Lemma~\ref{dhjwe54}}. Set
\bea
U_1&=&X_{11}^{2\alpha_1+1}X_{21}^{2\beta_1+1},\ \ U_2=X_{12}^{2\alpha_2+1}X_{22}^{2\beta_2+1}\ \mbox{and}\ \ U_i=X_{1i}^{2\alpha_i}X_{2i}^{2\beta_i}; \nonumber\\
V_1&=&X_{31}^{2\gamma_1+1}X_{41}^{2\delta_1+1},\ \ V_2=X_{32}^{2\gamma_2+1}X_{42}^{2\delta_2+1}\ \mbox{and}\ \ V_i=X_{3i}^{2\gamma_i}X_{4i}^{2\delta_i} \lbl{tn98}
\eea
for $3\leq i \leq m$. Let $N_1$ be the total number of solutions $(\bm{\alpha}, \bm{\beta}, \bm{\gamma}, \bm{\delta})$ satisfying \eqref{husoq922}. From \eqref{gyq7}, we see $N_1\leq K_1\cdot m^{\alpha+\beta+\gamma+\delta}$. Review the formulas between \eqref{eiy1} and \eqref{eiy4}. We have
\bea\lbl{dsjso}
m^{\alpha+\beta+\gamma+\delta}\cdot \mbox{Cov}\big((X_{11}X_{21})(X_{12}X_{22})A_1^{\alpha}A_2^{\beta}, (X_{31}X_{41})(X_{32}X_{42})A_3^{\gamma}A_4^{\delta}\big)
\eea
is a linear combination of $N_1$ terms of the form $\mbox{Cov}(\prod_{i=1}^mU_i, \prod_{i=1}^mV_i)$ with positive coefficients no more than $\alpha!\beta!\gamma!\delta!$.  Recall the notation $J_m(a, b)$ and \eqref{dsjso}. We then have
\beaa
m^{\alpha+\beta+\gamma+\delta}\cdot \big|J_m(1, 2)\big|
\leq  N_1\cdot (\alpha!\beta!\gamma!\delta!)\cdot \max \Big|\mbox{Cov}\Big(\prod_{i=1}^mU_i, \prod_{i=1}^mV_i\Big)\Big|,
\eeaa
where the maximum is taken over all $\{\alpha_i, \beta_i, \gamma_i, \delta_i;\, 1\leq i \leq m\}$ satisfying \eqref{husoq922}. By \eqref{gyq7}, we have
\bea\lbl{feefii}
\big|J_m(1, 2)\big|
\leq  K_1\cdot \max \Big|\mbox{Cov}\Big(\prod_{i=1}^mU_i, \prod_{i=1}^mV_i\Big)\Big|,
\eea
where the maximum is taken over all $(\bm{\alpha}, \bm{\beta}, \bm{\gamma}, \bm{\delta})$ satisfying \eqref{husoq922}. By Lemma~\ref{fjii09f} and applying the same argument of~\eqref{gi3298gq3} to \eqref{tn98}, we have
\bea\lbl{jhfewu}
\max|\mbox{Cov}(U_i, V_i)| \leq K_1\sum_{1\leq i<j \leq 4}r_{ij}^2,
\eea
where the maximum is taken over all $\{\alpha_i, \beta_i, \gamma_i, \delta_i;\, 1\leq i \leq m\}$ satisfying \eqref{husoq922}. This and \eqref{feefii} conclude
\beaa
\big|J_m(1, 2)\big|\leq K_1^2\sum_{1\leq i<j \leq 4}r_{ij}^2.
\eeaa
This proves the inequality for $(a, b)=(1, 2).$ \hfill $\square$

\begin{lemma}\lbl{dduh7} Assume the setting in \eqref{dida58}. Give non-negative integers $\alpha_i, \beta_i, \gamma_j, \delta_j$ for $i=1, 2$ and $j=3, 4$, set
\bea
U_1&=&X_{11}X_{21},\ \ U_2=X_{12}X_{22}\ \mbox{and}\ \ U_i=X_{1i}^{2\alpha_i}X_{2i}^{2\beta_i},\ i \in \{3, 4\}; \nonumber\\
V_3&=&X_{33}X_{43},\ \ V_4=X_{34}X_{44}\ \mbox{and}\ \ V_i=X_{3i}^{2\gamma_i}X_{4i}^{2\delta_i},\ i\in \{1, 2\}.
\eea
Define
\beaa
S &=& \{i\in \{3, 4\}; (\bm{\alpha}, \bm{\beta}, \bm{\gamma}, \bm{\delta})\ \mbox{satisfies}\ \eqref{husoq922}\ \mbox{and}\ \alpha_i+\beta_i\geq 1\}\cup\\
&& \{i\in \{1, 2\}; (\bm{\alpha}, \bm{\beta}, \bm{\gamma}, \bm{\delta})\ \mbox{satisfies}\ \eqref{husoq922}\ \mbox{and}\ \gamma_i+\delta_i\geq 1\}.
\eeaa
Then the following hold.

(i) If $S=\{1\}$ and $\gamma_1+\delta_1=1$, then
\beaa
\mbox{Cov}\Big(\prod_{i=1}^4U_i, \prod_{i=1}^4V_i\Big)
=
\begin{cases}
2(r_{12}r_{23}r_{31})r_{34}^2, & \text{if $\gamma_1=1$ and $\delta_1=0$};\\
2(r_{12}r_{24}r_{41})r_{34}^2, &\text{if $\gamma_1=0$ and $\delta_1=1$}.
\end{cases}
\eeaa

(ii) If $S=\{2\}$ and $\gamma_2+\delta_2=1$, then
\beaa
\mbox{Cov}\Big(\prod_{i=1}^4U_i, \prod_{i=1}^4V_i\Big)
=
\begin{cases}
2(r_{12}r_{23}r_{31})r_{34}^2, & \text{if $\gamma_2=1$ and $\delta_2=0$};\\
2(r_{12}r_{24}r_{41})r_{34}^2, &\text{if $\gamma_2=0$ and $\delta_2=1$}.
\end{cases}
\eeaa

(iii) If $S=\{3\}$ and $\alpha_3+\beta_3=1$, then
\beaa
\mbox{Cov}\Big(\prod_{i=1}^4U_i, \prod_{i=1}^4V_i\Big)
=
\begin{cases}
2(r_{13}r_{34}r_{41})r_{12}^2, & \text{if $\alpha_3=1$ and $\beta_3=0$};\\
2(r_{23}r_{34}r_{42})r_{12}^2, &\text{if $\alpha_3=0$ and $\beta_3=1$}.
\end{cases}
\eeaa

(iv) If $S=\{4\}$ and $\alpha_4+\beta_4=1$, then
\beaa
\mbox{Cov}\Big(\prod_{i=1}^4U_i, \prod_{i=1}^4V_i\Big)
=
\begin{cases}
2(r_{13}r_{34}r_{41})r_{12}^2, & \text{if $\alpha_4=1$ and $\beta_4=0$};\\
2(r_{23}r_{34}r_{42})r_{12}^2, &\text{if $\alpha_4=0$ and $\beta_4=1$}.
\end{cases}
\eeaa
\end{lemma}
\noindent\textbf{Proof of Lemma~\ref{dduh7}}. By assumption, $\{(X_{1j}, X_{2j}, X_{3j}, X_{4j})^T \in \mathbb{R}^4;\, 1\leq j \leq m\}$ are i.i.d.  random vectors with distribution $N_4(\bd{0}, \bd{R})$ where $\bd{R}=(r_{ij})_{4\times 4}$ and $r_{ii}=1$ for each $i.$ Thus,
\bea\lbl{euoqe863}
EU_1=EU_2=r_{12}\ \ \mbox{and}\ \ EV_3=EV_4=r_{34}.
\eea

(i) Under the case $S=\{1\}$ and $\gamma_1+\delta_1=1$, we know that $\alpha_3=\beta_3=\alpha_4=\beta_4=\gamma_2=\delta_2=0$ and that  $(\gamma_1, \delta_1)$ is equal to $(1, 0)$ or $(0, 1)$. Hence
\bea\lbl{ddjdwi}
&&U_1=X_{11}X_{21},\ U_2=X_{12}X_{22},\ U_3=1,\ U_4=1;\nonumber\\
&&V_1=X_{31}^2\, \mbox{or}\, X_{41}^2 ,\ V_2=1,\ V_3=X_{33}X_{43},\ V_4=X_{34}X_{44}.
\eea
This implies that $\{(U_1, V_1)^T, U_i, V_i;\, 2\leq i \leq 4\}$ are independent. By \eqref{euoqe863} and by Lemma~\ref{Wick_formula},
\beaa
&& \mbox{Cov}(U_1, X_{31}^2)=E(X_{11}X_{21}X_{31}^2)-r_{12}=2r_{13}r_{23};\\
&& \mbox{Cov}(U_1, X_{41}^2)=E(X_{11}X_{21}X_{41}^2)-r_{12}=2r_{14}r_{24}.
\eeaa
Notice
\bea\lbl{wefhe24}
\mbox{Cov}(\xi_1\eta_1, \xi_2\eta_2)=E\xi_1\cdot E\xi_2\cdot \mbox{Cov}(\eta_1, \eta_2)
\eea
if $\xi_1$ and $\xi_2$ are independent and $\{\xi_1, \xi_2\}$ are independent of $\{\eta_1, \eta_2\}$. Note $V_1=X_{31}^2$ if $(\gamma_1, \delta_1)=(1, 0)$ and $V_1=X_{41}^2$ if $(\gamma_1, \delta_1)=(0, 1)$. Then
\beaa
\mbox{Cov}\Big(\prod_{i=1}^4U_i, \prod_{i=1}^4V_i\Big)
&=&\mbox{Cov}(U_1, V_1)\cdot  \prod_{i=2, 3, 4}\big(EU_i\cdot EV_i\big)\\
&=&
\begin{cases}
2(r_{12}r_{23}r_{31})r_{34}^2, & \text{if $\gamma_1=1$ and $\delta_1=0$};\\
2(r_{12}r_{24}r_{41})r_{34}^2, &\text{if $\gamma_1=0$ and $\delta_1=1$}.
\end{cases}
\eeaa

(ii) Under the case $S=\{2\}$ and $\gamma_2+\delta_2=1$, we know that $\alpha_3=\beta_3=\alpha_4=\beta_4=\gamma_1=\delta_1=0$ and that  $(\gamma_2, \delta_2)$ is equal to $(1, 0)$ or $(0, 1)$. Hence
\beaa
&&U_1=X_{11}X_{21},\ U_2=X_{12}X_{22},\ U_3=1,\ U_4=1;\nonumber\\
&&V_1=1,\ V_2=X_{32}^2\, \mbox{or}\, X_{42}^2 ,\ V_3=X_{33}X_{43},\ V_4=X_{34}X_{44}.
\eeaa
Then, $U_1$, $(U_2, V_2)^T$, $V_3$ and $V_4$ are independent.  By \eqref{euoqe863},
\beaa
&& \mbox{Cov}(U_2, X_{32}^2)=E(X_{12}X_{22}X_{32}^2)-r_{12}=2r_{13}r_{23};\\
&& \mbox{Cov}(U_2, X_{42}^2)=E(X_{12}X_{22}X_{42}^2)-r_{12}=2r_{14}r_{24}.
\eeaa
By \eqref{wefhe24},
\beaa
\mbox{Cov}\Big(\prod_{i=1}^4U_i, \prod_{i=1}^4V_i\Big)
&=&\mbox{Cov}(U_2, V_2)\cdot  \prod_{i=1,3,4}^m\big(EU_i\cdot EV_i\big)\\
&=&
\begin{cases}
2(r_{12}r_{23}r_{31})r_{34}^2, & \text{if $\gamma_2=1$ and $\delta_2=0$};\\
2(r_{12}r_{24}r_{41})r_{34}^2, &\text{if $\gamma_2=0$ and $\delta_2=1$}.
\end{cases}
\eeaa

(iii) Under the case $S=\{3\}$ and $\alpha_3+\beta_3=1$, we know that $\alpha_4=\beta_4=\gamma_1=\delta_1=\gamma_2=\delta_2=0$ and that  $(\alpha_3, \beta_3)$ is equal to $(1, 0)$ or $(0, 1)$. Hence
\beaa
&&U_1=X_{11}X_{21},\ U_2=X_{12}X_{22},\ U_3=X_{13}^2\ \mbox{or}\ X_{23}^2,\ U_4=1;\nonumber\\
&&V_1=1,\ V_2=1 ,\ V_3=X_{33}X_{43},\ V_4=X_{34}X_{44}.
\eeaa
Then, $U_1$, $U_2$, $(U_3, V_3)^T$ and $V_4$ are independent. By \eqref{euoqe863},
\beaa
&& \mbox{Cov}(X_{13}^2, V_3)=E(X_{33}X_{43}X_{13}^2)-r_{34}=2r_{13}r_{14};\\
&& \mbox{Cov}(X_{23}^2, V_3)=E(X_{33}X_{43}X_{23}^2)-r_{34}=2r_{23}r_{24}.
\eeaa
By \eqref{wefhe24},
\beaa
\mbox{Cov}\Big(\prod_{i=1}^4U_i, \prod_{i=1}^4V_i\Big)
&=&\mbox{Cov}(U_3, V_3)\cdot  \prod_{i=1,2,4}^m\big(EU_i\cdot EV_i\big)\\
&=&
\begin{cases}
2(r_{13}r_{34}r_{41})r_{12}^2, & \text{if $\alpha_3=1$ and $\beta_3=0$};\\
2(r_{23}r_{34}r_{42})r_{12}^2, &\text{if $\alpha_3=0$ and $\beta_3=1$}.
\end{cases}
\eeaa

(iv) Under the case $S=\{4\}$ and $\alpha_4+\beta_4=1$, we know that $\alpha_3=\beta_3=\gamma_1=\delta_1=\gamma_2=\delta_2=0$ and that  $(\alpha_4, \beta_4)$ is equal to $(1, 0)$ or $(0, 1)$. Hence
\beaa
&&U_1=X_{11}X_{21},\ U_2=X_{12}X_{22},\ U_3=1,\ U_4=X_{14}^2\ \mbox{or}\ X_{24}^2;\nonumber\\
&&V_1=1,\ V_2=1 ,\ V_3=X_{33}X_{43},\ V_4=X_{34}X_{44}.
\eeaa
Then, $U_1$, $U_2$, $V_3$, $(U_4, V_4)^T$ are independent. By \eqref{euoqe863},
\beaa
&& \mbox{Cov}(X_{14}^2, V_4)=E(X_{34}X_{44}X_{14}^2)-r_{34}=2r_{13}r_{14};\\
&& \mbox{Cov}(X_{24}^2, V_4)=E(X_{34}X_{44}X_{24}^2)-r_{34}=2r_{23}r_{24}.
\eeaa
By \eqref{wefhe24},
\beaa
\mbox{Cov}\Big(\prod_{i=1}^4U_i, \prod_{i=1}^4V_i\Big)
&=&\mbox{Cov}(U_4, V_4)\cdot  \prod_{i=1,2,3}^m\big(EU_i\cdot EV_i\big)\\
&=&
\begin{cases}
2(r_{13}r_{34}r_{41})r_{12}^2, & \text{if $\alpha_4=1$ and $\beta_4=0$};\\
2(r_{23}r_{34}r_{42})r_{12}^2, &\text{if $\alpha_4=0$ and $\beta_4=1$}.
\end{cases}
\eeaa
The verification is finished. \hfill$\square$

Let $\alpha_1,\cdots,\alpha_m, \beta_1,\cdots, \beta_m, \gamma_1,\cdots, \gamma_m, \delta_1, \cdots, \delta_m$ be non-negative integers, review the notation
$\bm{\alpha}=(\alpha_1,\cdots,\alpha_m)$, $\bm{\beta}= (\beta_1,\cdots, \beta_m)$, $\bm{\gamma}=(\gamma_1,\cdots, \gamma_m)$ and $\bm{\delta}=(\delta_1, \cdots, \delta_m)$. For non-negative integers $\alpha, \beta, \gamma, \delta$ and  $(\bm{\alpha}, \bm{\beta}, \bm{\gamma}, \bm{\delta})$ satisfying \eqref{husoq922}, define
\bea\lbl{f8439}
C(\bm{\alpha}, \bm{\beta}, \bm{\gamma}, \bm{\delta})=\frac{\alpha!}{\alpha_1!\cdots \alpha_m!}\cdot
\frac{\beta!}{\beta_1!\cdots \beta_m!}\cdot \frac{\gamma!}{\gamma_1!\cdots \gamma_m!}\cdot \frac{\delta!}{\delta_1!\cdots \delta_m!}.
\eea

\begin{lemma}\lbl{shj2389} Assume the setting in \eqref{dida58}. Define $\{U_i, V_i;\, 1\leq i \leq m\}$ such that
\beaa
U_1&=&X_{11}^{2\alpha_1+1}X_{21}^{2\beta_1+1},\ \ U_2=X_{12}^{2\alpha_2+1}X_{22}^{2\beta_2+1}\ \mbox{and}\ \ U_i=X_{1i}^{2\alpha_i}X_{2i}^{2\beta_i};\nonumber\\
V_3&=&X_{33}^{2\gamma_3+1}X_{43}^{2\delta_3+1},\ \ V_4=X_{34}^{2\gamma_4+1}X_{44}^{2\delta_4+1}\ \mbox{and}\ \ V_j=X_{3j}^{2\gamma_j}X_{4j}^{2\delta_j}
\eeaa
for $3\leq i \leq m$ and  $j\in \{1,2,\cdots, m\}\backslash \{3, 4\}$, where $\alpha_i, \beta_i, \gamma_i, \delta_i$ satisfies~\eqref{husoq922}. Define
\beaa
S&= & \{i\in \{3, 4\}; (\bm{\alpha}, \bm{\beta}, \bm{\gamma}, \bm{\delta})\ \mbox{satisfies}\ \eqref{husoq922}\ \mbox{with}\ \alpha_i+\beta_i\geq 1\}\cup\\
  && \{j\in \{1, 2\}; (\bm{\alpha}, \bm{\beta}, \bm{\gamma}, \bm{\delta})\ \mbox{satisfies}\ \eqref{husoq922}\ \mbox{with}\ \gamma_j+\delta_j\geq 1\}.
\eeaa
Obviously, $S\subset \{1, 2, 3, 4\}$. Then
\beaa
&&\sum_{(\bm{\alpha}, \bm{\beta}, \bm{\gamma}, \bm{\delta}): |S|=1}C(\bm{\alpha}, \bm{\beta}, \bm{\gamma}, \bm{\delta})\cdot \mbox{Cov}\Big(\prod_{i=1}^mU_i, \prod_{i=1}^mV_i\Big)\\
& = & \rho_{m, 1} \sum_{1\leq i<j\leq 4}r_{ij}^2
 +\rho_{m, 2} \cdot\big[(r_{12}r_{23}r_{31})r_{34}^2+(r_{12}r_{24}r_{41})r_{34}^2\big] \\
&& ~~~~~~~~~~~~~~~~~+\rho_{m, 3} \cdot\big[(r_{13}r_{34}r_{41})r_{12}^2+
(r_{23}r_{34}r_{42})r_{12}^2\big]
\eeaa
where $\max\{m^2|\rho_{m, 1}|, m|\rho_{m, 2}|, m|\rho_{m, 3}|\}\leq K$ and $\rho_{m, 2}$ and $\rho_{m, 3}$ do not depend on $\bd{R}.$
\end{lemma}
\noindent\textbf{Proof of Lemma~\ref{shj2389}}. First, $|S|=1$ implies that $S=\{1\}$, $S=\{2\}$, $S=\{3\}$ or $S=\{4\}.$  We will first examine the case $S=\{1\}$ next.

Assume now $S=\{1\}$. Then $\gamma_1+\delta_1\geq 1$ and $\alpha_3=\beta_3=\alpha_4=\beta_4=\gamma_2=\delta_2=0$. Hence
\bea\lbl{fhuwqu}
&& U_3=U_4=V_2=1\ \mbox{and}\ U_2, V_3 ,V_4\ \mbox{are independent themselves and they are}\nonumber\\
&& \mbox{also independent of}\ \{(U_i, V_i)^T; i=1, 5, 6, \cdots, m\}
\eea
by the fact $\{(U_i, V_i)^T;\, 1\leq i \leq m\}$ are independent aforementioned. By Lemma~\ref{fjii09f} and applying the same argument of~\eqref{gi3298gq3} to $\{U_i, V_i,\, 1\leq i \leq m\}$, we obtain
\bea\lbl{dhh1du}
\max|\mbox{Cov}(U_i, V_i)| \leq K_1\sum_{1\leq i<j \leq 4}r_{ij}^2,
\eea
where the maximum is taken over all $\{\alpha_i, \beta_i, \gamma_i, \delta_i;\, 1\leq i \leq m\}$ satisfying \eqref{husoq922}.
Next we bound
\bea\lbl{jifewieo}
\sum\Big|\mbox{Cov}\Big(\prod_{i=1}^mU_i, \prod_{i=1}^mV_i\Big)\Big|,
\eea
where
the sum runs over $(\bm{\alpha}, \bm{\beta}, \bm{\gamma}, \bm{\delta})$ satisfing\eqref{husoq922} and $S=\{1\}$. When  $S=\{1\}$, we know $\gamma_1+\delta_1 \geq 1$. We will distinguish two cases: $\gamma_1+\delta_1 \geq  2$ and $\gamma_1+\delta_1 =1$. Recall the definition of $T_{m, 3}$ from Lemma~\ref{wejheif4}. For the case $\gamma_1+\delta_1 =1$, we will divide it into another two case: $T_{m, 3}$ and  $T_{m, 3}^c$. The derivation of bounds for \eqref{jifewieo} under $\gamma_1+\delta_1 \geq  2$ and $T_{m, 3}^c$ are easier than that under  $T_{m, 3}$. We will take two steps next two handle the two cases.

{\it Step 1}. By Lemma~\ref{cjhew0}(ii),  the total number of solutions $(\bm{\alpha}, \bm{\beta}, \bm{\gamma}, \bm{\delta})$ of \eqref{husoq922} with $\gamma_1+\delta_1 \geq 2$ is bounded by
$K_1\cdot m^{\alpha+\beta+\gamma+\delta-2}$.
This joined with \eqref{dhh1du} implies that
\bea\lbl{wehjefwif}
\sum\Big|\mbox{Cov}\Big(\prod_{i=1}^mU_i, \prod_{i=1}^mV_i\Big)\Big| \leq K_1\cdot m^{\alpha+\beta+\gamma+\delta-2}\sum_{1\leq i<j\leq 4}r_{ij}^2,
\eea
where the sum runs over all $(\bm{\alpha}, \bm{\beta}, \bm{\gamma}, \bm{\delta})$  satisfying \eqref{husoq922},  $S=\{1\}$ and $\gamma_1+\delta_1 \geq 2$.

Review the definition of $T_{m, 3}$ in Lemma~\ref{wejheif4}. We have $|T_{m, 3}|\leq K_1\cdot m^{\alpha+\beta+\gamma+\delta-2}.$  This together with \eqref{dhh1du} yields
\bea\lbl{dsg14}
\sum\Big|\mbox{Cov}\Big(\prod_{i=1}^mU_i, \prod_{i=1}^mV_i\Big)\Big| \leq K_1\cdot m^{\alpha+\beta+\gamma+\delta-2}\sum_{1\leq i<j\leq 4}r_{ij}^2,
\eea
where the sum runs over all $(\bm{\alpha}, \bm{\beta}, \bm{\gamma}, \bm{\delta})$  satisfying \eqref{husoq922},  $S=\{1\}$ and $(\bm{\alpha}, \bm{\beta}, \bm{\gamma}, \bm{\delta})\in T_{m, 3}.$

{\it Step 2}. We now estimate \eqref{jifewieo} as the index  $(\bm{\alpha}, \bm{\beta}, \bm{\gamma}, \bm{\delta})$  satisfies \eqref{husoq922}, the event $S=\{1\}$ holds and $(\bm{\alpha}, \bm{\beta}, \bm{\gamma}, \bm{\delta})\notin T_{m, 3}$.  Review the definition of $T_{m, 3}$ and expressions of $U_i$ and $V_i$, under the new conditions, $U_i$ and $V_i$ take much simpler form:
\bea\lbl{duwe8}
&&U_1=X_{11}X_{21},\ U_2=X_{12}X_{22},\ U_3=1,\ U_4=1;\nonumber\\
&&V_1=X_{31}^2\, \mbox{or}\, X_{41}^2 ,\ V_2=1,\ V_3=X_{33}X_{43},\ V_4=X_{34}X_{44}.
\eea
Furthermore, if $(\bm{\alpha}, \bm{\beta}, \bm{\gamma}, \bm{\delta})\notin T_{m, 3}$, then
$\alpha_k+\beta_k\leq  1$ for all $5\leq k\leq m$, $\gamma_l+\delta_l\leq 1$ for all $5\leq l\leq m$ and the two identities $\alpha_t+\beta_t=1$\ \mbox{and}\ $\gamma_t+\delta_t=1$ cannot occur at the same time for any $5\leq t\leq m$. The key observation is that, if  $\alpha_t+\beta_t=1$ then $U_t \sim \chi^2(1)$ and $V_t=1$. Similarly, if  $\gamma_t+\delta_t=1$ then $U_t=1$ and $V_t\sim \chi^2(1)$.
 Therefore, the $2m-8$ random variables in
 $\{U_i, V_i;\, 5\leq i \leq m\}$ are  independent random variables, each has mean $1$.
As used earlier, $\{(U_i, V_i)^T;\, 1\leq i \leq m\}$ are independent.
This and the special structures in \eqref{duwe8} imply that the $2m-1$ random quantities in $\{(U_1, V_1)^T, U_i, V_i;\, 2\leq i \leq m\}$ are independent. By Lemma~\ref{dduh7}(i) and \eqref{wefhe24},
\beaa
\mbox{Cov}\Big(\prod_{i=1}^mU_i, \prod_{i=1}^mV_i\Big)
&=& \mbox{Cov}\Big(\prod_{i=1}^4U_i, \prod_{i=1}^4V_i\Big)\\
&=&
\begin{cases}
2(r_{12}r_{23}r_{31})r_{34}^2, & \text{if $\gamma_1=1$ and $\delta_1=0$};\\
2(r_{12}r_{24}r_{41})r_{34}^2, &\text{if $\gamma_1=0$ and $\delta_1=1$}.
\end{cases}
\eeaa
 This says that
\bea\lbl{sdhi280}
\sum  \mbox{Cov}\Big(\prod_{i=1}^mU_i, \prod_{i=1}^mV_i\Big) = L_{m,1} \cdot\big[(r_{12}r_{23}r_{31})r_{34}^2+(r_{12}r_{24}r_{41})r_{34}^2\big],
\eea
where the sum runs over $\Gamma$, defined to be the set of $(\bm{\alpha}, \bm{\beta}, \bm{\gamma}, \bm{\delta})$ satisfying \eqref{husoq922}, $|S|=1$  with $\gamma_1+\delta_1=1$, and $(\bm{\alpha}, \bm{\beta}, \bm{\gamma}, \bm{\delta}) \notin T_{m, 3}$; $L_{m,1}:=2|\Gamma|$. Obviously, $L_{m,1}$ does not depend on the matrix $\bd{R}=(r_{ij})$. By the bound on $T_{m, 1}$ in Lemma~\ref{wejheif4}, we have $L_{m,1}\leq K_1m^{\alpha+\beta+\gamma+\delta-1}.$ Notice, if $(\bm{\alpha}, \bm{\beta}, \bm{\gamma}, \bm{\delta})$  satisfies \eqref{husoq922}, the event $S=\{1\}$ holds and $(\bm{\alpha}, \bm{\beta}, \bm{\gamma}, \bm{\delta})\notin T_{m, 3}$, then any one from $\{\alpha_i, \beta_i, \gamma_i, \delta_i;\, 1\leq i \leq m\}$ is either $1$ or $0$. According to \eqref{f8439}, $C(\bm{\alpha}, \bm{\beta}, \bm{\gamma}, \bm{\delta}) = \alpha!\beta!\gamma!\delta!$. Then \eqref{sdhi280} becomes
\beaa
\sum C(\bm{\alpha}, \bm{\beta}, \bm{\gamma}, \bm{\delta}) \mbox{Cov}\Big(\prod_{i=1}^mU_i, \prod_{i=1}^mV_i\Big) =\alpha!\beta!\gamma!\delta!\cdot L_{m,1} \cdot\big[(r_{12}r_{23}r_{31})r_{34}^2+(r_{12}r_{24}r_{41})r_{34}^2\big].
\eeaa
Thus, combining this with \eqref{wehjefwif} and \eqref{dsg14} and using the trivial fact that $C(\bm{\alpha}, \bm{\beta}, \bm{\gamma}, \bm{\delta}) \leq \alpha!\beta!\gamma!\delta!$, we arrive at
\bea\lbl{comay273}
\sum C(\bm{\alpha}, \bm{\beta}, \bm{\gamma}, \bm{\delta})\cdot \mbox{Cov}\Big(\prod_{i=1}^mU_i, \prod_{i=1}^mV_i\Big) = \tau_{m,1} \cdot\big[(r_{12}r_{23}r_{31})r_{34}^2+(r_{12}r_{24}r_{41})r_{34}^2\big] +\tau_{m, 1}',
\eea
where the sum runs over the set of $(\bm{\alpha}, \bm{\beta}, \bm{\gamma}, \bm{\delta})$ satisfying \eqref{husoq922} and  $S=\{1\}$, and where $\tau_{m,1}$ does not depend on  $r_{ij}$ and
\bea\lbl{wq9rjvns}
|\tau_{m,1}|\leq K_1m^{\alpha+\beta+\gamma+\delta-1} \ \ \mbox{and}\ \ |\tau_{m, 1}'|\leq K_1\cdot m^{\alpha+\beta+\gamma+\delta-2}\sum_{1\leq i<j\leq 4}r_{ij}^2.
\eea

By applying the same argument as the derivation of \eqref{comay273} to  $S=\{2\}$ and using Lemma~\ref{dduh7}(ii), we get an analogue of~\eqref{comay273} as the sum runs over the set of $(\bm{\alpha}, \bm{\beta}, \bm{\gamma}, \bm{\delta})$ satisfying \eqref{husoq922} and $S=\{2\}$, where $\tau_{m,1}$ and $\tau_{m1}'$ will be replaced by two corresponding symbols but still satisfy \eqref{wq9rjvns}.

By applying the same argument as the derivation of \eqref{comay273} to $S=\{3\}$ and using Lemma~\ref{dduh7}(iii), we get
\bea\lbl{sy37k}
\sum\mbox{Cov}\Big(\prod_{i=1}^mU_i, \prod_{i=1}^mV_i\Big) = \tilde{\tau}_{m,1} \cdot\big[
(r_{13}r_{34}r_{41}+r_{23}r_{34}r_{42})r_{12}^2\big] + \tilde{\tau}'_{m,1},
\eea
where the sum runs over the set of $(\bm{\alpha}, \bm{\beta}, \bm{\gamma}, \bm{\delta})$ satisfying \eqref{husoq922} and $S=\{3\}$, and the inequalities from \eqref{wq9rjvns} still hold as  ``$(\tau_{m,1}, \tau_{m,1}')$" is replaced by ``$(\tilde{\tau}_{m,1}, \tilde{\tau}'_{m,1})$".

By applying the same argument as the derivation of \eqref{comay273} to  $S=\{4\}$ and using Lemma~\ref{dduh7}(iv), we get an analogue of~\eqref{sy37k} as the sum runs over the set of $(\bm{\alpha}, \bm{\beta}, \bm{\gamma}, \bm{\delta})$ satisfying \eqref{husoq922} and $S=\{4\}$, the quantities ``$(\tilde{\tau}_{m,1}, \tilde{\tau}_{m,1}')$" is replaced by ``$(T_{m,1}, T_{m1}'')$", and $T_{m,1}$ does not depend on $r_{ij}$ and
\beaa
|T_{m,1}|\leq K_1m^{\alpha+\beta+\gamma+\delta-1} \ \ \mbox{and}\ \ |T_{m, 1}'|\leq K_1\cdot m^{\alpha+\beta+\gamma+\delta-2}\sum_{1\leq i<j\leq 4}r_{ij}^2.
\eeaa
The proof is completed by summing the above four upper bounds corresponding to $S=\{1\}, \{2\}, \{3\}$ and $\{4\}$.   \hfill$\square$

\begin{lemma}\lbl{ghuw2} Assume the setting in \eqref{dida58}. Let $J_m(a, b)$ be defined as in Lemma~\ref{dhjwe54}. Then
\beaa
J_m(3, 4)
&=&\tau_{m,1}r_{12}^2r_{34}^2 + \tau_{m, 2}\sum_{1\leq i<j\leq 4}r_{ij}^2
 +\tau_{m,3} \cdot\big[(r_{12}r_{23}r_{31})r_{34}^2+(r_{12}r_{24}r_{41})r_{34}^2\big] \\
&& ~~~~~~~~~~~~~~~~~~~~~~~~~~~~~~~~+\tau_{m, 4}\cdot\big[(r_{13}r_{34}r_{41})r_{12}^2+
(r_{23}r_{34}r_{42})r_{12}^2\big],
\eeaa
where $\max\{m|\tau_{m,1}|, m^2|\tau_{m, 2}|, m|\tau_{m, 3}|, m|\tau_{m, 4}|\}\leq K$ and $\tau_{m, 3}$ and $\tau_{m, 4}$ do not depend on $\bd{R}.$
\end{lemma}
\noindent\textbf{Proof of Lemma~\ref{ghuw2}}. Set
\bea
U_1&=&X_{11}^{2\alpha_1+1}X_{21}^{2\beta_1+1},\ \ U_2=X_{12}^{2\alpha_2+1}X_{22}^{2\beta_2+1}\ \mbox{and}\ \ U_i=X_{1i}^{2\alpha_i}X_{2i}^{2\beta_i};\nonumber\\
V_3&=&X_{33}^{2\gamma_3+1}X_{43}^{2\delta_3+1},\ \ V_4=X_{34}^{2\gamma_4+1}X_{44}^{2\delta_4+1}\ \mbox{and}\ \ V_j=X_{3j}^{2\gamma_j}X_{4j}^{2\delta_j} \lbl{huhu}
\eea
for $3\leq i \leq m$ and  $j\in \{1,2,,\cdots, m\}\backslash \{3, 4\}$. Let
$N_1$ be the total number of solutions  $(\bm{\alpha}, \bm{\beta}, \bm{\gamma}, \bm{\delta})$ for~\eqref{husoq922} with a bound provided in \eqref{gyq7}. Review the discussions between \eqref{eiy1} and \eqref{wdqu}. We know that $m^{\alpha+\beta+\gamma+\delta} J_m(3,4)$
is a linear combination of $N_1$ terms of the form
$\mbox{Cov}(\prod_{i=1}^mU_i, \prod_{i=1}^mV_i)$
with positive coefficients no more than $\alpha!\beta!\gamma!\delta!$.
 Define
\beaa
 S &= & \{i\in \{3, 4\}; (\bm{\alpha}, \bm{\beta}, \bm{\gamma}, \bm{\delta})\ \mbox{satisfies}\ \eqref{husoq922}\ \mbox{with}\ \alpha_i+\beta_i\geq 1\}\cup\\
  && \{j\in \{1, 2\}; (\bm{\alpha}, \bm{\beta}, \bm{\gamma}, \bm{\delta})\ \mbox{satisfies}\ \eqref{husoq922}\ \mbox{with}\ \gamma_j+\delta_j\geq 1\};\\
   S_1 &=& \{5\leq i \leq m;\, (\bm{\alpha}, \bm{\beta}, \bm{\gamma}, \bm{\delta})\ \mbox{satisfies}\ \eqref{husoq922}\ \mbox{with either}\ \alpha_i+\beta_i\geq 1\ \mbox{or}\  \delta_i+\gamma_i\geq 1\}.
\eeaa
Similar to the proof of Lemma~\ref{fjii09f}, we have
\bea\lbl{ewuwq9}
|S_1|
& \leq & \alpha+\beta +\gamma +\delta.
\eea
We now  estimate $\mbox{Cov}(\prod_{i=1}^mU_i, \prod_{i=1}^mV_i)$ by differentiating three cases: $|S|=0$, $|S|=1$ and $|S| \geq 2$. Quickly, for the case $|S|=1$, by reviewing  $C(\bm{\alpha}, \bm{\beta}, \bm{\gamma},  \bm{\delta})$ from \eqref{f8439}, we have from  Lemma~\ref{shj2389} that
\bea\lbl{dhiewfuo}
&& \sum_{(\bm{\alpha}, \bm{\beta}, \bm{\gamma}, \bm{\delta}): |S|=1}C(\bm{\alpha}, \bm{\beta}, \bm{\gamma},  \bm{\delta})\cdot\mbox{Cov}\Big(\prod_{i=1}^mU_i, \prod_{i=1}^mV_i\Big)\nonumber\\
& = & \rho_{m, 1} \sum_{1\leq i<j\leq 4}r_{ij}^2
 +\rho_{m, 2} \cdot\big[(r_{12}r_{23}r_{31})r_{34}^2+(r_{12}r_{24}r_{41})r_{34}^2\big] \nonumber \\
&& ~~~~~~~~~~~~~~~~~+\rho_{m, 3} \cdot\big[(r_{13}r_{34}r_{41})r_{12}^2+
(r_{23}r_{34}r_{42})r_{12}^2\big],
\eea
where $\rho_{m, 1} \leq Km^{-2}$ and $\rho_{m, 2} \vee\rho_{m, 3}  \leq m^{-1}K$, and where $\rho_{m, 2}$ and $\rho_{m, 3}$ do not depend on $\bd{R}.$ To finish the proof, it remains to study the cases ``$|S|=0$" and ``$|S| \geq 2$". This will be worked out in two steps.

{\textbf{Step 1}:
First, the condition $|S|=0$ implies $\alpha_3=\alpha_4=\beta_3=\beta_4=\gamma_1=\gamma_2=\delta_1=\delta_2=0$, and hence we have from \eqref{huhu} that
\beaa
&& U_1=X_{11}^{2\alpha_1+1}X_{21}^{2\beta_1+1},\ \ U_2=X_{12}^{2\alpha_2+1}X_{22}^{2\beta_2+1},\ U_3=1, U_4=1\ \mbox{and}\ \ U_i=X_{1i}^{2\alpha_i}X_{2i}^{2\beta_i};\\
&& V_1= 1, V_2=1, V_3=X_{33}^{2\gamma_3+1}X_{43}^{2\delta_3+1},\ \ V_4=X_{34}^{2\gamma_4+1}X_{44}^{2\delta_4+1}\ \mbox{and}\ \ V_i=X_{3i}^{2\gamma_i}X_{4i}^{2\delta_i}
\eeaa
for $i=5, \cdots, m$. By assumption \eqref{dida58}, $\{(X_{1i}, X_{2i}, X_{3i}, X_{4i})^T \in \mathbb{R}^4;\, 1\leq i \leq m\}$ are  i.i.d.  random vectors with distribution $N_4(\bd{0}, \bd{R})$, where $\bd{R}=(r_{ij})_{4\times 4}\ \mbox{and}\ r_{ii}=1$ for each $i$. In particular,  $\{(U_i, V_i)^T;\, 1\leq i \leq m\}$ are independent. As a consequence, $U_1, U_2, V_3, V_4$ are themselves independent, and furthermore $\{U_1, U_2, V_3, V_4\}$ are also independent of $\{U_i, V_i;\,  5\leq i \leq m\}$. By \eqref{wefhe24},
\bea\lbl{wqw3r09}
\mbox{Cov}\Big(\prod_{i=1}^mU_i, \prod_{i=1}^mV_i\Big)
&=&EU_1EU_2EV_3EV_4\cdot \mbox{Cov}\Big(\prod_{i=5}^mU_i, \prod_{i=5}^mV_i\Big)\nonumber\\
&=& EU_1EU_2EV_3EV_4\cdot \mbox{Cov}\Big(\prod_{i\in S_1}U_i, \prod_{i\in S_1}V_i\Big).
\eea
By definition of $S_1$, we see that $\sum_{i\in S_1}\alpha_i=\alpha$, $\sum_{i\in S_1}\beta_i=\beta$, $\sum_{i\in S_1}\gamma_i=\gamma$ and $\sum_{i\in S_1}\delta_i=\delta$. By Lemma~\ref{fjii09f},
\beaa
\big|\mbox{Cov}\Big(\prod_{i\in S_1}U_i, \prod_{i\in S_1}V_i\Big)\big|
\leq K_1\sum_{1\leq i<j \leq 4}r_{ij}^2\leq 6K_1.
\eeaa
Bounds for $EU_1, EU_2, EV_3, EV_4$ are given in
Lemma~\ref{wdh120}. By the lemma, we see
\bea\lbl{wh367d}
\Big|\mbox{Cov}\Big(\prod_{i=1}^mU_i, \prod_{i=1}^mV_i\Big)\Big|
&\leq &  K_1r_{12}^2r_{34}^2\cdot \Big|\mbox{Cov}\Big(\prod_{i=5}^mU_i, \prod_{i=5}^mV_i\Big)\Big|\nonumber\\
& \leq & K_1\cdot r_{12}^2r_{34}^2.
\eea
Let $S_2$ be the set of solutions $(\bm{\alpha}, \bm{\beta}, \bm{\gamma}, \bm{\delta})$ satisfying~\eqref{husoq922} with  $\alpha_i+\beta_i\geq 1$ and $\gamma_i+\delta_i\geq 1$ simultaneously for some $5\leq i \leq m$. By Lemma~\ref{feufeuiqfP}, we have $|S_2| \leq Km^{\alpha+\beta+\gamma+\delta-1}$. Recall $U_i=X_{1i}^{2\alpha_i}X_{2i}^{2\beta_i}$ and $V_i=X_{3i}^{2\gamma_i}X_{4i}^{2\delta_i}$ for $5\leq i \leq m.$ If $\alpha_i=\beta_i=0$ then $U_i=1$. Likewise, $V_i=1$ if $\gamma_i=\delta_i=0.$ This together with the fact $\{(U_i, V_i);\, 1\leq i \leq m\}$ are independent implies $\prod_{i=5}^mU_i$ and $\prod_{i=5}^mV_i$ are independent (hence their covariance is zero) if there is no $i\in \{5,\cdots, m\}$ such that $\alpha_i+\beta_i\geq 1$ and $\gamma_i+\delta_i\geq 1$ at the same time. Therefore, by~\eqref{wh367d},
\bea\lbl{uqwu394}
\sum\Big|\mbox{Cov}\Big(\prod_{i=1}^mU_i, \prod_{i=1}^mV_i\Big)\Big| \leq K_1\cdot m^{\alpha+\beta+\gamma+\delta-1}r_{12}^2r_{34}^2
\eea
where the sum runs over all $(\bm{\alpha}, \bm{\beta}, \bm{\gamma}, \bm{\delta})$ satisfying \eqref{husoq922} and $|S|=0$.

{\textbf{Step 2}: Assume the index $(\bm{\alpha}, \bm{\beta}, \bm{\gamma}, \bm{\delta})$ satisfies that $|S|\geq 2$. Review the structures of $U_i$ and $V_i$ from \eqref{huhu}, we have from Lemma~\ref{fjii09f}  that
\beaa
\max\Big|\mbox{Cov}\Big(\prod_{i=1}^mU_i, \prod_{i=1}^mV_i\Big)\Big| \leq K_1\cdot \sum_{1\leq i<j \leq 4}r_{12}^2,
\eeaa
 where the sum runs over all $(\bm{\alpha}, \bm{\beta}, \bm{\gamma}, \bm{\delta})$ satisfying \eqref{husoq922}. Review the definition of $T_{m, 2}$ from Lemma \eqref{wejheif4}, we have from the lemma that $|T_{m, 2}|\leq K_{13}m^{\alpha+\beta+\gamma+\delta-2}$. It follows that
 \bea\lbl{sdhwq8}
&&\sum\Big|\mbox{Cov}\Big(\prod_{i=1}^mU_i, \prod_{i=1}^mV_i\Big)\Big|\leq K_{1}m^{\alpha+\beta+\gamma+\delta-2}\sum_{1\leq i<j \leq 4}r_{12}^2,
\eea
where the sum runs over all $(\bm{\alpha}, \bm{\beta}, \bm{\gamma}, \bm{\delta})$  satisfying \eqref{husoq922} and $|S|\geq 2$.

Finally, we add up the three bounds from \eqref{dhiewfuo}, \eqref{uqwu394} and \eqref{sdhwq8}. By changing ``$K_1$" to ``$\tau_{m,1}$", ``$\rho_{m, 1}+K_1$" to ``$\tau_{m, 2}$", ``$\rho_{m, 2}$" to ``$\tau_{m,3}$" and ``$\rho_{m, 3}$" to ``$\tau_{m,4}$", we complete the proof.

\hfill$\square$

\begin{lemma}\lbl{whwqo} Assume the setting in \eqref{dida58}. Let all notation be the same as those in Lemma~\ref{dhjwe54}.
Then
\beaa
J_m(2, 3)&=&\tau_{m, 1}\sum_{1\leq i<j\leq 4}r_{ij}^2
 +\tau_{m,2} \cdot\big(r_{12}r_{23}r_{34}r_{41}+r_{12}r_{24}r_{43}r_{31}\big),
\eeaa
where $\max\{m|\tau_{m,1}|, |\tau_{m, 2}|\}\leq K$ and $\tau_{m, 2}$ does not depend on  $\bd{R}$.
\end{lemma}
\noindent\textbf{Proof of Lemma~\ref{whwqo}}. Set
\bea
U_1&=&X_{11}^{2\alpha_1+1}X_{21}^{2\beta_1+1},\ \ U_2=X_{12}^{2\alpha_2+1}X_{22}^{2\beta_2+1}\ \mbox{and}\ \ U_i=X_{1i}^{2\alpha_i}X_{2i}^{2\beta_i}; \nonumber\\
V_2&=&X_{32}^{2\gamma_2+1}X_{42}^{2\delta_2+1},\ \ V_3=X_{33}^{2\gamma_3+1}X_{43}^{2\delta_3+1}\ \mbox{and}\ \ V_j=X_{3j}^{2\gamma_j}X_{4j}^{2\delta_j} \lbl{huhu1}
\eea
for $3\leq i \leq m$ and $j\in \{1,2,\cdots, m\}\backslash \{2,3\}.$ Review $C(\bm{\alpha}, \bm{\beta}, \bm{\gamma}, \bm{\delta})$ in \eqref{f8439}. By \eqref{eiy1}-\eqref{eiy4} and~\eqref{qwhwqoig},
\bea\lbl{fwek3}
J_m(2, 3)
&=&\mbox{Cov}\big((X_{11}X_{21})(X_{12}X_{22})A_1^{\alpha}A_2^{\beta}, (X_{32}X_{42})(X_{33}X_{43})A_3^{\gamma}A_4^{\delta}\big)\nonumber\\
&=& \frac{1}{m^{\alpha+\beta+\gamma+\delta}}\sum C(\bm{\alpha}, \bm{\beta}, \bm{\gamma}, \bm{\delta})\cdot \mbox{Cov}\Big(\prod_{i=1}^mU_i, \prod_{i=1}^mV_i\Big),
\eea
where the sum runs over all $(\bm{\alpha}, \bm{\beta}, \bm{\gamma}, \bm{\delta})$  satisfying \eqref{husoq922}.
Let $S$ be defined as in Lemma~\ref{feufeuiqfP}. By the lemma, $|S|\leq K_1m^{\alpha+\beta+\gamma+\delta-1}$. By Lemma~\ref{fjii09f} and structures of $U_i$ and $V_i$ in \eqref{huhu1},
\beaa
\max\Big|\mbox{Cov}\Big(\prod_{i=1}^mU_i, \prod_{i=1}^mV_i\Big)\Big| \leq K_{1}\cdot \sum_{1\leq i<j \leq 4}r_{12}^2,
\eeaa
 where the maximum is taken over all $(\bm{\alpha}, \bm{\beta}, \bm{\gamma}, \bm{\delta})$ satisfying \eqref{husoq922}. Combining the two facts, we obtain that
 \beaa
 \sum_{(\bm{\alpha}, \bm{\beta}, \bm{\gamma}, \bm{\delta})\in S}\Big|\mbox{Cov}\Big(\prod_{i=1}^mU_i, \prod_{i=1}^mV_i\Big)\Big|\leq K_{1}m^{\alpha+\beta+\gamma+\delta-1}\sum_{1\leq i<j \leq 4}r_{12}^2.
 \eeaa
As used before, $1\leq C(\bm{\alpha}, \bm{\beta}, \bm{\gamma}, \bm{\delta}) \leq \alpha!\beta!\gamma!$ for any  $(\bm{\alpha}, \bm{\beta}, \bm{\gamma}, \bm{\delta})$  satisfies~\eqref{husoq922}. We rewrite the above as
\bea\lbl{viuew9}
\frac{1}{m^{\alpha+\beta+\gamma+\delta}}\sum_{(\bm{\alpha}, \bm{\beta}, \bm{\gamma}, \bm{\delta})\in S} C(\bm{\alpha}, \bm{\beta}, \bm{\gamma}, \bm{\delta})\cdot \mbox{Cov}\Big(\prod_{i=1}^mU_i, \prod_{i=1}^mV_i\Big)
=\tau_{m, 1}\sum_{1\leq i<j \leq 4}r_{12}^2,
\eea
where $|\tau_{m,1}| \leq K_1m^{-1}.$

Now, if $(\bm{\alpha}, \bm{\beta}, \bm{\gamma}, \bm{\delta})$ satisfies \eqref{husoq922} but not in $S$, then (i) $\alpha_i+\beta_i=0$ for each $i\in \{1, 2, 3\}$; (ii) $\gamma_i+\delta_i=0$ for each  $i\in \{1, 2, 3\}$; (iii) $\alpha_i+\beta_i\leq 1$ and  $\gamma_i+\delta_i \leq 1$ for each $4\leq i \leq m$; (iv) for each $4\leq i \leq m$, if $\alpha_i+\beta_i=1$ then $\gamma_i+\delta_i=0$, and if $\gamma_i+\delta_i=1$ then $\alpha_i+\beta_i=0$.  This implies that $\{U_i, V_i;\, 4\leq i \leq m\}$ are independent random variables and each of them is either $1$ or $\chi^2(1).$ Keep in mind that $\{(U_i, V_i);\, 1\leq i\leq m\}$ are independent random variables and $E(\chi^2(1))=1$. Furthermore, it is readily seen from \eqref{huhu1} that
\beaa
U_1&=&X_{11}X_{21},\ \ U_2=X_{12}X_{22}\ \mbox{and}\ \ U_3=1; \nonumber\\
V_1&=&1,\ \  V_2=X_{32}X_{42},\ \ V_3=X_{33}X_{43}.
\eeaa
Since  $\{(U_i, V_i)^T;\, 1\leq i \leq m\}$ are independent by the assumption from \eqref{dida58},  the three random quantities $\{U_1, (U_2, V_2)^T, V_3\}$ are independent and they are independent of $\{(U_i, V_i)^T;\, 4\leq i \leq m\}$. Thus, it follows from \eqref{wefhe24} that
\beaa
\mbox{Cov}\Big(\prod_{i=1}^mU_i, \prod_{i=1}^mV_i\Big)
&=& E\Big(\prod_{i=4}^mU_i\Big)\cdot E\Big( \prod_{i=4}^mV_i\Big)\cdot\mbox{Cov}\Big(\prod_{i=1}^3U_i, \prod_{i=1}^3V_i\Big) \\
&=& EU_1\cdot EV_3\cdot \mbox{Cov}(U_2, V_2).
\eeaa
Recall $\{(X_{1j}, X_{2j}, X_{3j}, X_{4j})^T \in \mathbb{R}^4;\, 1\leq j \leq m\}$ are i.i.d.  random vectors with distribution $N_4(\bd{0}, \bd{R})$, where $\bd{R}=(r_{ij})_{4\times 4}$ and $r_{ii}=1$ for each $i$. Then $EU_1=r_{12}$, $EV_3=r_{34}$ and
\beaa
\mbox{Cov}(U_2, V_2)=E\big(X_{12}X_{22}X_{32}X_{42}\big)-r_{12}r_{34}=r_{13}r_{24}+r_{14}r_{23}
\eeaa
by Lemma~\ref{Wick_formula}.
Thus,
\beaa
\mbox{Cov}\Big(\prod_{i=1}^mU_i, \prod_{i=1}^mV_i\Big)=r_{12}r_{24}r_{43}r_{31}+r_{12}r_{23}r_{34}r_{41}.
\eeaa
Recall \eqref{f8439}, $C(\bm{\alpha}, \bm{\beta}, \bm{\gamma}, \bm{\delta})=\alpha!\beta!\gamma!\delta!$ in this case, that is, $(\bm{\alpha}, \bm{\beta}, \bm{\gamma}, \bm{\delta})$ satisfies \eqref{husoq922} but not in $S$. Let $N_1$ be the total number of solutions $(\bm{\alpha}, \bm{\beta}, \bm{\gamma}, \bm{\delta})$ satisfying \eqref{husoq922}. From \eqref{gyq7}, we see $N_1\leq K_1\cdot m^{\alpha+\beta+\gamma+\delta}$. Therefore, there exists a constant $\tau_{m, 2}$ not depending $r_{ij}$ and $|\tau_{m, 2}|\leq K_1$ such that
\beaa
 \frac{1}{m^{\alpha+\beta+\gamma+\delta}}\sum C(\bm{\alpha}, \bm{\beta}, \bm{\gamma}, \bm{\delta})\cdot\mbox{Cov}\Big(\prod_{i=1}^mU_i, \prod_{i=1}^mV_i\Big) =\tau_{m, 2}\big(r_{12}r_{23}r_{34}r_{41}+r_{12}r_{24}r_{43}r_{31}\big)
 \eeaa
where the sum is taken over every $(\bm{\alpha}, \bm{\beta}, \bm{\gamma}, \bm{\delta})$ satisfying \eqref{husoq922} and the restriction in $S^c$. By connecting this fact  to \eqref{fwek3}  and \eqref{viuew9}, we get desired the conclusion.  \hfill$\square$

\subsubsection{A Study on Correlation Matrices}\lbl{Fescm}

As needed in the proof of Lemma~\ref{main_an} later, we have to handle certain functions of the entries of sample correlation matrices. They are interesting on their own. Through the whole section, we assume $\bd{R}=(r_{ij})_{p\times p}$ is a non-negative definite matrix with $r_{ii}=1$ for $1\leq i \leq p$. Review the Frobenius norm $\|\bd{R}\|_F=[\mbox{tr}(\bd{R}^2)]^{1/2}=(\sum_{1\leq i, j\leq p}r_{ij}^2)^{1/2}$.

\begin{lemma}\lbl{dhdqwi} Assume $\{m_p;\, p\geq 1\}$ are positive constants with $\lim_{p\to\infty}m_p= \infty$. Define
\beaa
W_1=\sum_{1\leq i, j, k\leq p}r_{ij}r_{jk}r_{ki}\ \ \mbox{and}\ \ \
W_2=\sum_{1\leq i, j, k, l\leq p}r_{ij}r_{jk}r_{kl}r_{li}.
\eeaa
Then $\lim_{p\to\infty}\frac{W_i}{m\|\bd{R}\|_F^4}= 0$ for $i=1, 2.$
\end{lemma}
\noindent\textbf{Proof of Lemma~\ref{dhdqwi}}. Let $\lambda_1\geq 0, \cdots, \lambda_p \geq 0$ be the eigenvalues of $\bd{R}$. Write $\bd{R}=\bd{O}^T\mbox{diag}(\lambda_1, \cdots, \lambda_p)\bd{O}$, where $\bd{O}$ is a $p\times p$ orthogonal matrix. Recall the fact
\bea\lbl{qwdygw8}
\mbox{tr}\big(\bd{R}^s\big)=\sum_{1\leq i_1, i_2, \cdots, i_s\leq p}r_{i_1i_2}r_{i_2i_3}r_{i_3i_4}\cdots r_{i_si_1}
\eea
for any integer $s\geq 2$.   Easily,
\beaa
W_1=\,\mbox{tr}(\bd{R}^3)=\lambda_1^3+\cdots +\lambda_p^3\leq  \big(\lambda_1^2+\cdots +\lambda_p^2\big)^{3/2}.
\eeaa
In addition, $\|\bd{R}\|_F^2\geq \sum_{i=1}^pr_{ii}^2=p$. Therefore,
\beaa
\frac{|W_1|}{[\mbox{tr}(\bd{R}^2)]^2} \leq \frac{[\mbox{tr}(\bd{R}^2)]^{3/2}}{[\mbox{tr}(\bd{R}^2)]^2}=\frac{1}{[\mbox{tr}(\bd{R}^2)]^{1/2}}\leq \frac{1}{\sqrt{p}}.
\eeaa
It follows that $\lim_{p\to\infty}\frac{W_1}{m\|\bd{R}\|_F^4}= 0$. Now prove the second conclusion. Note
\beaa
W_2=\,\mbox{tr}(\bd{R}^4)=\lambda_1^4+\cdots + \lambda_p^4\leq (\lambda_1^2+\cdots +\lambda_p^2)^2.
\eeaa
Consequently,
\beaa
\frac{|W_2|}{m[\mbox{tr}(\bd{R}^2)]^2} \leq \frac{ [\mbox{tr}(\bd{R}^2)]^2}{m[\mbox{tr}(\bd{R}^2)]^2}=\frac{1}{m}\to 0
\eeaa
as $p\to\infty.$ The proof is completed. \hfill$\square$

\begin{lemma}\lbl{wdhdwqo} Given  integer $\alpha\geq 0$, define
\beaa
S_m=\sum_{1\leq i,j,k,l\leq p}r_{ij}r_{jk}r_{ki}r_{kl}^{\alpha}.
\eeaa
If $\lim_{p\to\infty}\frac{p}{m\|\bd{R}\|_F}= 0$, then $\lim_{p\to\infty}\frac{S_m}{m\|\bd{R}\|_F^4}= 0$.
\end{lemma}
\noindent\textbf{Proof of Lemma~\ref{wdhdwqo}}. Set $a_k=\sum_{l=1}^pr_{kl}^{\alpha}$ for $k=1,2,\cdots, p$. Then $|a_k|\leq p$ for each $k$.  Define a $p\times p$ matrix  $\bd{D}=\mbox{diag}(a_1, \cdots, a_p).$ Then the $(k, i)$-entry of $\bd{R}\bd{D}$ is $r_{ki}a_i$. It follows that
\beaa
S_m=\sum_{1\leq i,j,k\leq p}r_{ij}r_{jk}r_{ki}a_k=\sum_{1\leq i,j,k\leq p}r_{ij}r_{jk}\big(\bd{R}\bd{D}\big)_{ki}.
\eeaa
Therefore, $S_m=\mbox{tr}\,(\bd{R}\bd{R}(\bd{R}\bd{D}))=\mbox{tr}\,(\bd{R}^3\bd{D})$. Set $(b_{ij})_{p\times p}=\bd{B}= \bd{R}^3$. Then $\bd{B}$ is a non-negative definite matrix due to the fact that $\bd{R}$ is non-negative definite. Hence, $b_{ii}\geq 0$ for each $i$ and
\beaa
|\mbox{tr}\,(\bd{R}^3\bd{D})|=|\mbox{tr}\,(\bd{B}\bd{D})|=\big|\sum_{i=1}^pb_{ii}a_i\big|\leq p\sum_{i=1}^pb_{ii}.
\eeaa
This shows that $|S_m| \leq p\cdot \mbox{tr}\,(\bd{R}^3)$. Let $\lambda_1\geq 0, \cdots, \lambda_p \geq 0$ be the eigenvalues of $\bd{R}$. Easily,
\beaa
\mbox{tr}\,(\bd{R}^3)=\lambda_1^3+\cdots +\lambda_p^3\leq  \big(\lambda_1^2+\cdots +\lambda_p^2\big)^{3/2}.
\eeaa
Therefore $|S_m|\leq p\cdot [\mbox{tr}(\bd{R}^2)]^{3/2}.$
Consequently,
\beaa
\frac{|S_m|}{m\,[\mbox{tr}(\bd{R}^2)]^2} \leq \frac{p}{m\,[\mbox{tr}(\bd{R}^2)]^{1/2}}\to 0
\eeaa
by assumption. The proof is finished. \hfill$\square$

\begin{lemma}\lbl{ewf394} Define $\Lambda_p:=\big\{(i, j, k, l); \, 1\leq i\ne j \leq p\ \mbox{and}\ 1\leq k\ne l\leq p\big\}$,
\beaa
V_{p,1}&=&\frac{1}{m\|\bd{R}\|_F^4}\cdot \sum_{(i,j,k,l)\in \Lambda_p}r_{ij}r_{jk}r_{kl}r_{li}, \nonumber\\
V_{p,2}&=&\frac{1}{m\|\bd{R}\|_F^4}\cdot\sum_{(i,j,k,l)\in \Lambda_p}(r_{ik}r_{kl}r_{li})r_{ij}^2.
\eeaa
If $\lim_{p\to\infty}\frac{p}{m\|\bd{R}\|_F}= 0$, then
$\lim_{p\to\infty}V_{p,i} = 0$ for $i=1, 2.$
\end{lemma}
\noindent\textbf{Proof of Lemma~\ref{ewf394}}. By assumption, $r_{ii}=1$ for any $1\leq i \leq p$. Then
\beaa
\sum_{1\leq i,j,k,l\leq m}r_{ij}r_{jk}r_{kl}r_{li}
&=&\sum_{1\leq i=j,k,l\leq p}r_{ij}r_{jk}r_{kl}r_{li}+\sum_{1\leq i\ne j,k,l\leq p}r_{ij}r_{jk}r_{kl}r_{li}\\
& = & \sum_{1\leq i,k,l\leq p}r_{ik}r_{kl}r_{li}+ \sum_{1\leq i\ne j,k=l\leq p}r_{ij}r_{jk}r_{kl}r_{li}+\sum_{(i,j,k,l)\in \Lambda_p}r_{ij}r_{jk}r_{kl}r_{li}.
\eeaa
Now
\beaa
\sum_{1\leq i\ne j,k=l\leq m}r_{ij}r_{jk}r_{kl}r_{li}=\sum_{1\leq i\ne j,k\leq m}r_{ij}r_{jk}r_{ki}=\sum_{1\leq i, j,k\leq m}r_{ij}r_{jk}r_{ki}-\sum_{1\leq j,k\leq m}r_{jk}^2.
\eeaa
Recall \eqref{qwdygw8}. The above two identities imply that
\beaa
\sum_{(i,j,k,l)\in \Lambda_p}r_{ij}r_{jk}r_{kl}r_{li}=\, \mbox{tr}\big(\bd{R}^4\big)-2\,\mbox{tr}\big(\bd{R}^3\big)+\mbox{tr}\big(\bd{R}^2\big).
\eeaa
Then $V_{p,1}\to 0$ by using the fact $\|\bd{R}\|_F^2\geq p$ and the conclusions for $W_1$ and $W_2$ in Lemma~\ref{dhdqwi}.

Now we prove $V_{p,2}\to 0$. Note that
\bea\lbl{Hulatang}
\sum_{1\leq i,j,k,l\leq p}(r_{ik}r_{kl}r_{li})r_{ij}^2
&=&\sum_{1\leq j,k,l\leq p}r_{jk}r_{kl}r_{lj}
+ \sum_{1\leq i\ne j,k,l\leq p}(r_{ik}r_{kl}r_{li})r_{ij}^2 \nonumber\\
& = & \mbox{tr}\big(\bd{R}^3\big)+\sum_{1\leq i\ne j,k\leq p}r_{ik}^2r_{ij}^2+\sum_{(i,j,k,l)\in \Lambda_p}(r_{ik}r_{kl}r_{li})r_{ij}^2.
\eea
Now
\beaa
0\leq \sum_{1\leq i\ne j,k\leq p}r_{ik}^2r_{ij}^2\leq \sum_{i=1}^p\Big(\sum_{j=1}^pr_{ij}^2\Big)^2
\leq \Big(\sum_{i=1}^p\sum_{j=1}^pr_{ij}^2\Big)^2=\, \big[\mbox{tr}\big(\bd{R}^2\big)\big]^2.
\eeaa
Hence,
\beaa
\frac{1}{m\,[\mbox{tr}(\bd{R}^2)]^2}\sum_{1\leq i\ne j,k\leq p}r_{ik}^2r_{ij}^2\leq \frac{1}{m}\to 0.
\eeaa
Also, $\frac{1}{m[\mbox{tr}(\bd{R}^2)]^2}\mbox{tr}(\bd{R}^3)\to 0$ by using the conclusion for $W_1$ in Lemma~\ref{dhdqwi}. Then the conclusion $V_{p,2}\to 0$ follows from~\eqref{Hulatang} and Lemma~\ref{wdhdwqo} with $\alpha=2.$ \hfill$\square$

\newpage

\noindent\textbf{Proof of Lemma~\ref{dsih328}}. Write $\bd{M}=(m_{ij})$. Then $m_{ii}=1$ for each $i$. By Theorem 4.3.26 from \cite{horn2012matrix}, $\{\lambda_1,\cdots, \lambda_p\}$ majorizes $\{m_{11}, \cdots, m_{pp}\}$, the diagonal matrix of $\bd{M}$. By the definition of majorization, $\lambda_1+\cdots + \lambda_k\geq m_{11}+\cdots+m_{kk}=k$ for each $1\leq k \leq p$ and  $\lambda_1+\cdots + \lambda_p= m_{11}+\cdots+m_{pp}=p$.

On the other hand, by the definition of majorization, the $p$ numbers $\{1, \cdots, 1\}$ majorizes $\{\tau_1, \cdots, \tau_p\}$. By Theorem 4.3.32 from \cite{horn2012matrix}, there exists a symmetric matrix $\bd{B}=(b_{ij})_{p\times p}$ such that $b_{ii}=1$ for each $i$ and that $\bd{B}$ has non-negative eigenvalues $\tau_1, \cdots, \tau_p$. By definition, $\bd{B}$ is a correlation matrix. \hfill$\square$

\subsubsection{The Proofs of  Theorems~\ref{Theorem3} and~\ref{Theorem4}}\lbl{TPTRSCM}

In this part, by using the preliminary results developed in Sections~\ref{Variance_Study}-\ref{Fescm}, we are now ready to prove the two main results Theorems~\ref{Theorem4} and ~\ref{Theorem3} stated in Section~\ref{TRSCMTP}.

\begin{lemma}\lbl{main_an} Assume the setting in \eqref{dida58} with $\bd{R}=(r_{ij})_{4\times4}$. Recall $A_i$ in \eqref{hjwq38}. Define
\beaa
B_1=\frac{1}{m}\sum_{j=1}^mX_{1j}X_{2j}\ \ \mbox{and}\ \ B_2=\frac{1}{m}\sum_{j=1}^mX_{3j}X_{4j}.
\eeaa
Given integer $N\geq 1$, the covariance between
\beaa
\sum_{0\leq j, k\leq N} (1-A_1)^j(1-A_2)^kB_1^2\ \ \mbox{and}\ \
\sum_{0\leq j', k'\leq N}(1-A_3)^{j'}(1-A_4)^{k'}B_2^2
\eeaa
is equal to
\beaa
\varrho_{m,1}\cdot \sum_{1\leq i<j\leq 4}r_{ij}^2 + \varrho_{m,2}\cdot r_{12}^2r_{34}^2 &+& \varrho_{m,3}\cdot\big(r_{12}r_{23}r_{34}r_{41}+r_{12}r_{24}r_{43}r_{31}\big)\\
&+&\varrho_{m,4}\cdot\big[(r_{13}r_{34}r_{41})r_{12}^2+
(r_{23}r_{34}r_{42})r_{12}^2\big]\\
&+ & \varrho_{m,5}\cdot\big[(r_{12}r_{23}r_{31})r_{34}^2+(r_{12}r_{24}r_{41})r_{34}^2\big],
\eeaa
where  $\{\varrho_{m,i};\, 3\leq i \leq 5\}$ do not depend on $\bd{R}$,
\beaa
|\varrho_{m,1}| \leq Km^{-2},  \ \ |\varrho_{m,2}|\vee|\varrho_{m,3}|\vee|\varrho_{m,4}|\vee|\varrho_{m,5}| \leq K m^{-1}
\eeaa
and $K$ is a constant depending on $N$ but not on $m$ or $\bd{R}.$
\end{lemma}
\noindent\textbf{Proof of Lemma~\ref{main_an}}. For convenience, we use $\Delta_m$ to denote the covariance between
\beaa
\sum_{0\leq j, k\leq N} (1-A_1)^j(1-A_2)^kB_1^2\ \ \mbox{and}\ \
\sum_{0\leq j', k'\leq N}(1-A_3)^{j'}(1-A_4)^{k'}B_2^2.
\eeaa
Then
\bea
\Delta_m =
 \sum\mbox{Cov}\big((1-A_1)^j(1-A_2)^kB_1^2, (1-A_3)^{j'}(1-A_4)^{k'}B_2^2\big),\ \ \ \ \lbl{efh19}
\eea
where the sum runs over all non-negative integers $j, j', k, k'$ such that $0\leq j, k\leq N$ and  $0\leq j', k'\leq N$.  For each $i=1, 2, 3, 4$, write
\beaa
(1-A_i)^l=1+\sum_{\alpha=1}^l(-1)^{\alpha}\binom{l}{\alpha}A_i^{\alpha}
\eeaa
for any $l\geq 1$. Trivially,  $\mbox{Cov}(U_1+h_1, U_2+h_2)=\mbox{Cov}(U_1, U_2)$ for any random variables $U_1$ and $U_2$ and constants $h_1$ and $h_2.$ Then the last covariance from \eqref{efh19} is
\beaa
&& \mbox{Cov}\big((1-A_1)^j(1-A_2)^kB_1^2, (1-A_3)^{j'}(1-A_4)^{k'}B_2^2\big)\\
& = &~ \mbox{a finite linear combination of}~ H~ \mbox{terms of the form }
\mbox{Cov}\big(B_1^2A_1^{\alpha}A_2^{\beta}, B_2^2A_3^{\gamma}A_4^{\delta}\big)
\eeaa
where the coefficients in the linear combination depend on $\alpha, \beta, \gamma, \delta$ but not $m$ or $\bd{R}$, and $H:=(j+1)(j'+1)(k+1)(k'+1)$. This and \eqref{efh19} imply that
\bea
 \Delta_m= \mbox{a finite linear combination of}~ N'~ \mbox{terms of the form } \mbox{Cov}\big(B_1^2A_1^{\alpha}A_2^{\beta}, B_2^2A_3^{\gamma}A_4^{\delta}\big)~~~~~ \lbl{dg328}
\eea
where $0\leq \alpha+\beta\leq N$ and $0\leq \gamma+\delta\leq N$ and the coefficients in the linear combination depend on $N$ but not on $m$ or $\bd{R}$,  and $N'$ is bounded by
\beaa
\sum_{j, k\leq N, j', k'\leq N}(j+1)(j'+1)(k+1)(k'+1)
&\leq & (N+1)^4\sum_{0\leq j, j', k, k'\leq N}1\\
&\leq & (N+1)^8.
\eeaa
As $\alpha+\beta=0$ and $\gamma+\delta=0$, the covariance becomes
\bea\lbl{ewje9}
\mbox{Cov}(B_1^2, B_2^2)=\frac{1}{m}(r_{12}r_{23}r_{34}r_{41}+r_{12}r_{24}r_{43}r_{31}) +\frac{\delta_m}{m^2}\sum_{1\leq i< j \leq 4}r_{ij}^2
\eea
by Lemma~\ref{sjwqi3}, where $|\delta_m|\leq \kappa$ and $\kappa$ is a numerical constant not depending on $m$, $\bd{R}$ or $N$. So we next only need to study $\mbox{Cov}\big(B_1^2A_1^{\alpha}A_2^{\beta}, B_2^2A_3^{\gamma}A_4^{\delta}\big)$ from \eqref{dg328} with an extra assumption that either $\alpha+\beta\geq 1$ or $\gamma+\delta\geq 1.$

Write
\bea
&& m^2\cdot B_1^2=\sum_{j=1}^m(X_{1j}X_{2j})^2 +2\sum_{1\leq k< l \leq m}(X_{1k}X_{2k})(X_{1l}X_{2l}); \lbl{zenmeshde1}\\
&& m^2\cdot B_2^2=\sum_{q=1}^m(X_{3q}X_{4q})^2 +2\sum_{1\leq a< b \leq m}(X_{3a}X_{4a})(X_{3b}X_{4b}).\lbl{zenmeshde2}
\eea
Then
\bea\lbl{we6ft3}
m^4\cdot \mbox{Cov}\big(B_1^2A_1^{\alpha}A_2^{\beta}, B_2^2A_3^{\gamma}A_4^{\delta}\big)
= D_1+2D_2+2D_3+4D_4,
\eea
where
\beaa
D_1 &=& \sum_{1\leq j, q\leq m}\,\mbox{Cov}\big((X_{1j}X_{2j})^2A_1^{\alpha}A_2^{\beta}, (X_{3q}X_{4q})^2A_3^{\gamma}A_4^{\delta}\big),\\
D_2 &=& \sum_{1\leq j\leq m, 1\leq a< b \leq m}\,\mbox{Cov}\big((X_{1j}X_{2j})^2A_1^{\alpha}A_2^{\beta}, (X_{3a}X_{4a})(X_{3b}X_{4b})A_3^{\gamma}A_4^{\delta}\big),\\
D_3 &= & \sum_{1\leq q\leq m, 1\leq k< l \leq m}\,\mbox{Cov}\big((X_{1k}X_{2k})(X_{1l}X_{2l})A_1^{\alpha}A_2^{\beta}, (X_{3q}X_{4q})^2A_3^{\gamma}A_4^{\delta}\big)
\eeaa
and
\beaa
D_4=\sum_{1\leq k< l \leq m, 1\leq a<b\leq m} \mbox{Cov}\big((X_{1k}X_{2k})(X_{1l}X_{2l})A_1^{\alpha}A_2^{\beta}, (X_{3a}X_{4a})(X_{3b}X_{4b})A_3^{\gamma}A_4^{\delta}\big).
\eeaa
We next study the four terms in steps

{\it Step 1: the estimate of $D_1$}. Write
\beaa
&&\sum_{1\leq j, q\leq m}\,\mbox{Cov}\big((X_{1j}X_{2j})^2A_1^{\alpha}A_2^{\beta}, (X_{3q}X_{4q})^2A_3^{\gamma}A_4^{\delta}\big)\\
&=& \sum_{j=1}^ m\,\mbox{Cov}\big((X_{1j}X_{2j})^2A_1^{\alpha}A_2^{\beta}, (X_{3j}X_{4j})^2A_3^{\gamma}A_4^{\delta}\big)\\
&&+\sum_{1\leq j\ne q\leq m}\,\mbox{Cov}\big((X_{1j}X_{2j})^2A_1^{\alpha}A_2^{\beta}, (X_{3q}X_{4q})^2A_3^{\gamma}A_4^{\delta}\big).
\eeaa
Review $\{(X_{1j}, X_{2j}, X_{3j}, X_{4j})^T \in \mathbb{R}^4;\, 1\leq j \leq m\}$ are i.i.d.  random vectors. It follows that
\beaa
D_1 &=& m\cdot \mbox{Cov}\big((X_{11}X_{21})^2A_1^{\alpha}A_2^{\beta}, (X_{31}X_{41})^2A_3^{\gamma}A_4^{\delta}\big)\\
&& +m(m-1)\cdot\mbox{Cov}\big((X_{11}X_{21})^2A_1^{\alpha}A_2^{\beta}, (X_{32}X_{42})^2A_3^{\gamma}A_4^{\delta}\big).
\eeaa
By Lemma~\ref{uwee9},
\bea\lbl{fewjieowfj}
|D_1|
\leq  K_1m^2 \sum_{1\leq i<j \leq 4}r_{ij}^2,
\eea
where $K_1$ here and  later denotes a constant depending on $\alpha, \beta, \gamma, \delta$ but not $m$ or $\bd{R}$, and can be different from line to line.

{\it Step 2: the estimate of $D_2$}. Write
\beaa
&&D_2\\
 &=& \sum_{1\leq a< b \leq m}\Big(\sum_{j\in \{a, b\}} +\sum_{1\leq j \leq m, j\notin \{a, b\}}\Big) \,\mbox{Cov}\big((X_{1j}X_{2j})^2A_1^{\alpha}A_2^{\beta}, (X_{3a}X_{4a})(X_{3b}X_{4b})A_3^{\gamma}A_4^{\delta}\big)\\
 &= & 2\cdot \binom{m}{2}\cdot I_m(1, 2) + (m-2)\binom{m}{2}\cdot I_m(2, 3),
\eeaa
where
\beaa
I_m(a, b):=\mbox{Cov}\big((X_{11}X_{21})^2A_1^{\alpha}A_2^{\beta}, (X_{3a}X_{4a})(X_{3b}X_{4b})A_3^{\gamma}A_4^{\delta}\big).
\eeaa
By Lemma~\ref{jeu8},
\bea\lbl{sahuiu}
|D_2| \leq K_1m^2\sum_{1\leq i<j \leq 4}r_{ij}^2.
\eea

{\it Step 3: the estimate of $D_3$}. By switching the roles of ``$(X_{1j}, X_{2j}, A_1, A_2, \alpha, \beta)$'' and ``$(X_{3j}, X_{4j}, A_3, A_4, \gamma, \delta)$'' in {\it Step 2}, and using \eqref{sahuiu}, we obtain
\bea\lbl{jweoiW}
|D_3| \leq K_1m^2\sum_{1\leq i<j \leq 4}r_{ij}^2.
\eea

{\it Step 4: the estimate of $D_4$}. Rewrite
\beaa
D_4=\sum_{1\leq a<b\leq m}\Big(\sum_{\Gamma_1}+\sum_{\Gamma_2}+ \sum_{\Gamma_3}\Big) \mbox{Cov}\big((X_{1k}X_{2k})(X_{1l}X_{2l})A_1^{\alpha}A_2^{\beta}, (X_{3a}X_{4a})(X_{3b}X_{4b})A_3^{\gamma}A_4^{\delta}\big),
\eeaa
where
\beaa
&& \Gamma_1=\{(k, l): (k, l)=(a,b)\},\ \ \Gamma_2=\{(k, l):\, 1\leq k< l \leq m,\, |\{k, l\}\cap \{a,b\}|=1\},\\
&&  \Gamma_3=\{(k, l): 1\leq k< l \leq m,\, \{k, l\}\cap \{a,b\}=\emptyset\}.
\eeaa
Given $1\leq a<b\leq m$, it is easy to see $|\Gamma_1|=1$, $|\Gamma_2| \leq 2m$ and $\Gamma_3\leq m^2.$  Set
\beaa
J_m(a, b)=\mbox{Cov}\big((X_{11}X_{21})(X_{12}X_{22})A_1^{\alpha}A_2^{\beta}, (X_{3a}X_{4a})(X_{3b}X_{4b})A_3^{\gamma}A_4^{\delta}\big)
\eeaa
for $a\geq 1$ and $b\geq 1$. Then
\beaa
D_4=\frac{1}{2}m(m-1)\big[J_m(1, 2)+ |\Gamma_2|\cdot J_m(2,3) + |\Gamma_3|\cdot J_m(3, 4)\big].
\eeaa
By Lemma~\ref{dhjwe54},
\beaa
\big|J_m(1, 2)\big| \leq K_1\sum_{1\leq i<j \leq 4}r_{ij}^2.
\eeaa
By Lemma~\ref{whwqo},
\beaa
J_m(2, 3)&=&\tau_{m, 1}\sum_{1\leq i<j\leq 4}r_{ij}^2
 +\tau_{m,2} \cdot\big(r_{12}r_{23}r_{34}r_{41}+r_{12}r_{24}r_{43}r_{31}\big)
\eeaa
where $|\tau_{m,1}|\leq K_1m^{-1}$,  $|\tau_{m, 2}| \leq K_1$ and $\tau_{m, 2}$ does not depend on $\bd{R}$.  By Lemma~\ref{ghuw2},
\beaa
J_m(3, 4)
&=&\tau_{m,1}'r_{12}^2r_{34}^2 + \tau_{m, 2}'\sum_{1\leq i<j\leq 4}r_{ij}^2
 +\tau_{m,3} \cdot\big[(r_{12}r_{23}r_{31})r_{34}^2+(r_{12}r_{24}r_{41})r_{34}^2\big] \\
&& ~~~~~~~~~~~~~~~~~~~~~~~~~~~~~~~~+\tau_{m, 4}\cdot\big[(r_{13}r_{34}r_{41})r_{12}^2+
(r_{23}r_{34}r_{42})r_{12}^2\big],
\eeaa
where $|\tau_{m,1}'|\leq K_1m^{-1}$, $|\tau_{m, 2}'| \leq K_1m^{-2}$, $|\tau_{m, 3}|\vee|\tau_{m, 4}| \leq K_1m^{-1}$, and $\tau_{m, 3}$ and $\tau_{m, 4}$ do not depend on $\bd{R}$. Combining all of the above we get
\beaa
D_4=\rho_{m,1}\cdot \sum_{1\leq i<j\leq 4}r_{ij}^2 + \rho_{m,2}\cdot r_{12}^2r_{34}^2
&+&\rho_{m,3}\cdot\big(r_{12}r_{23}r_{34}r_{41}+r_{12}r_{24}r_{43}r_{31}\big)\\
&+&\rho_{m,4}\cdot\big[(r_{13}r_{34}r_{41})r_{12}^2+
(r_{23}r_{34}r_{42})r_{12}^2\big]\\
&+ & \rho_{m,5}\cdot\big[(r_{12}r_{23}r_{31})r_{34}^2+(r_{12}r_{24}r_{41})r_{34}^2\big],
\eeaa
where $\{\rho_{m,i},\, 1\leq i \leq 5\}$ satisfy that
\beaa
|\rho_{m,1}| \leq K_1m^2,\ \
\max\big\{|\rho_{m,2}|, |\rho_{m,3}|, |\rho_{m,4}|, |\rho_{m,5}|\big\} \leq K_1m^3
\eeaa
and $\rho_{m,3}$, $\rho_{m,4}$ and $\rho_{m,5}$ do not depend on $\bd{R}$.
Through combining the estimates of $D_1, D_2, D_3, D_4$ and \eqref{we6ft3}, we see $\mbox{Cov}\big(B_1^2A_1^{\alpha}A_2^{\beta}, B_2^2A_3^{\gamma}A_4^{\delta}\big)$ is equal to
\beaa
 \rho_{m,1}'\cdot \sum_{1\leq i<j\leq 4}r_{ij}^2 + \rho_{m,2}'\cdot r_{12}^2r_{34}^2 &+& \rho_{m,3}'\cdot\big(r_{12}r_{23}r_{34}r_{41}+r_{12}r_{24}r_{43}r_{31}\big)\\
&+&\rho_{m,4}'\cdot\big[(r_{13}r_{34}r_{41})r_{12}^2+
(r_{23}r_{34}r_{42})r_{12}^2\big]\\
&+ & \rho_{m,5}'\cdot\big[(r_{12}r_{23}r_{31})r_{34}^2+(r_{12}r_{24}r_{41})r_{34}^2\big],
\eeaa
where $\{\rho_{m,i}';\, 1\leq i \leq 5\}$ satisfy that
\beaa
|\rho_{m,1}'| \leq K_1m^{-2}, \ \ |\rho_{m,2}'|\vee|\rho_{m,3}'|\vee|\rho_{m,4}'|\vee|\rho_{m,5}'| \leq K_1m^{-1},
\eeaa
and $\rho_{m,3}', \rho_{m,4}'$ and $\rho_{m,5}'$ do not depend on $\bd{R}$.
Recalling~\eqref{dg328}, we arrive at that $\Delta_m$ is equal to
\beaa
\varrho_{m,1}\cdot \sum_{1\leq i<j\leq 4}r_{ij}^2 + \varrho_{m,2}\cdot r_{12}^2r_{34}^2 &+& \varrho_{m,3}\cdot\big(r_{12}r_{23}r_{34}r_{41}+r_{12}r_{24}r_{43}r_{31}\big)\\
&+&\varrho_{m,4}\cdot\big[(r_{13}r_{34}r_{41})r_{12}^2+
(r_{23}r_{34}r_{42})r_{12}^2\big]\\
&+ & \varrho_{m,5}\cdot\big[(r_{12}r_{23}r_{31})r_{34}^2+(r_{12}r_{24}r_{41})r_{34}^2\big],
\eeaa
where
\beaa
|\varrho_{m,1}| \leq K_1m^{-2}, \ \ |\varrho_{m,2}|\vee|\varrho_{m,3}|\vee|\varrho_{m,4}|\vee|\varrho_{m,5}| \leq K_1m^{-1},
\eeaa
and $\varrho_{m,3}, \varrho_{m,4}$ and $\varrho_{m,5}$ do not depend on $\bd{R}$. The proof is completed. \hfill$\square$

\medskip
We will first prove  Theorem~\ref{Theorem4} and then  prove~\ref{Theorem3}.
\medskip

\noindent\textbf{Proof of Theorem~\ref{Theorem4}}. Recall the earlier notation that
\beaa
&& A_i=\frac{1}{m}\sum_{j=1}^mX_{ij}^2,\ \ \ \ \ \ B_1=\frac{1}{m}\sum_{j=1}^mX_{1j}X_{2j}, \ \ \  \ \ \ \ B_2=\frac{1}{m}\sum_{j=1}^mX_{3j}X_{4j}
\eeaa
for $i=1,2,3,4.$ Then
\bea\lbl{iewfiu}
\hat{r}_{12}=\frac{B_1}{\sqrt{A_1A_2}}\ \ \ \mbox{and}\ \ \ \ \hat{r}_{34}=\frac{B_2}{\sqrt{A_3A_4}}.
\eea
Given $N\geq 1$, write
\beaa
\frac{1}{x}=1+(1-x)+\cdots +(1-x)^N+\frac{1}{x}(1-x)^{N+1}
\eeaa
for $x\ne 0$. Thus
\beaa
\frac{1}{A_1A_2}
&=&\Big[\frac{(1-A_1)^{N+1}}{A_1}+\sum_{i=0}^N(1-A_1)^i\Big]\cdot \Big[\frac{(1-A_2)^{N+1}}{A_2}+\sum_{j=0}^N(1-A_2)^j\Big]\\
&= & \epsilon_{m, 1}+ \sum_{0\leq i, j \leq N}(1-A_1)^i(1-A_2)^j,
\eeaa
where
\beaa
\epsilon_{m, 1}=\frac{(1-A_1)^{N+1}(1-A_2)^{N+1}}{A_1A_2}
&+&
\sum_{j=0}^N \frac{(1-A_1)^{N+1}(1-A_2)^j}{A_1}\\
& + & \sum_{i=0}^N \frac{(1-A_1)^i(1-A_2)^{N+1}}{A_2}.
\eeaa
Similarly,
\beaa
\frac{1}{A_3A_4}= \epsilon_{m, 2}+ \sum_{0\leq i, j \leq N}(1-A_3)^i(1-A_4)^j
\eeaa
where
\beaa
\epsilon_{m, 2}=\frac{(1-A_3)^{N+1}(1-A_4)^{N+1}}{A_3A_4}
&+&
\sum_{j=0}^N \frac{(1-A_3)^{N+1}(1-A_4)^j}{A_3}\\
& + & \sum_{i=0}^N \frac{(1-A_3)^i(1-A_4)^{N+1}}{A_4}.
\eeaa
By \eqref{iewfiu},
\bea
\mbox{Cov}(\hat{r}_{12}^2, \hat{r}_{34}^2)
&=&\,\mbox{Cov}\Big(\sum_{0\leq i, j \leq N}(1-A_1)^i(1-A_2)^jB_1^2,\, \sum_{0\leq i, j \leq N}(1-A_3)^i(1-A_4)^jB_2^2\Big) \nonumber\\
&&+ \,\mbox{Cov}\Big(\sum_{0\leq i, j \leq N}(1-A_1)^i(1-A_2)^jB_1^2,\, \epsilon_{m, 2}B_2^2\Big) \nonumber\\
&&+\,\mbox{Cov}\Big(\epsilon_{m, 1}B_1^2,\, \sum_{0\leq i, j \leq N}(1-A_3)^i(1-A_4)^jB_2^2\Big) \nonumber\\
& & +\, \mbox{Cov}\big(\epsilon_{m, 1}B_1^2,\, \epsilon_{m, 2}B_2^2\big).\lbl{ewy129}
\eea
We claim that
\bea\lbl{djcjwe}
\mbox{the absolute value of each of the last three covariances in \eqref{ewy129}}\ \leq \frac{K_1}{m^{(N+1)/2}}
\eea
where $K_1$ is a constant depending on $N$ but not  $m$ or $\bd{R}.$ In fact, by writing $\bar{B}_1=B_1-r_{12}$ and $\bar{B}_2=B_2-r_{34}$, then $B_1^2=\bar{B}_1^2+2r_{12}\bar{B}_1+r_{12}^2$ and $B_2^2=\bar{B}_2^2+2r_{34}\bar{B}_2+r_{34}^2$, and hence by linearity of the  covariance, each of the last three covariances from~\eqref{ewy129} is a linear combination of $N'$ terms of the form
\bea\lbl{ewuewu}
r_{12}^ar_{34}^b\cdot E\frac{\bar{B}_1^{t_1}\bar{B}_2^{t_2}(1-A_1)^{n_1}(1-A_2)^{n_2}(1-A_3)^{n_3}
(1-A_4)^{n_4}}{A_1^{s_1}A_2^{s_2}A_3^{s_3}A_4^{s_4}}
\eea
where $N'$ depends only on $N$ and all powers are non-negative integers with
\beaa
&& a, b, t_i\in \{0, 1, 2\} ~~~ \mbox{and}\ ~~~ s_1+s_2+s_3+s_4\geq 1;\\
&& n_i\leq N~~~~~~~~~~~~~~~ \mbox{and}\ ~~~  n_1+n_2+n_3+n_4\geq N+1
\eeaa
for each possible $i$. The crucial observation is that $n_1+n_2+n_3+n_4\geq N+1$. If some of $\{t_1, t_2, n_1, n_2, n_3, n_4\}$ are zero, the corresponding terms simply disappear. Set
\beaa
&& k=~\mbox{the count of positive values from}~\{t_1, t_2, n_1, n_2, n_3, n_4\};\\
&& l=~\mbox{the count of positive values from}~\{s_1, s_2, s_3, s_4\}.
\eeaa
Then $k\geq 1$ and $l\geq 1$.  Take
\beaa
&&X_i=
\begin{cases}
\sqrt{m}\bar{B}_i, & \text{$i=1, 2$};\\
\sqrt{m}(1-A_{i-2}),& \text{$i=3, 4, 5, 6$}
\end{cases}
\ \ \ \ \
\mbox{and}\ \ \ \
\alpha_i=
\begin{cases}
t_i, & \text{$i=1, 2$};\\
n_{i-2},& \text{$i=3, 4, 5, 6$}
\end{cases}
\eeaa
and $p_i=m$ for $i=1, \cdots, 6.$ Furthermore, take
$Y_j=A_j, q_j=m$ and $\beta_j=s_j$ for $j=1, 2,3, 4$. Easily, $2/(k+l)\leq 1$ and $(2-q_j)/[2(k+l)]<0$ as $m\geq 5.$ Observe
\beaa
&& X_1=\frac{1}{\sqrt{m}}\sum_{j=1}^m (X_{1j}X_{2j}-r_{12}),\ \ X_2=\frac{1}{\sqrt{m}}\sum_{j=1}^m(X_{3j}X_{4j}-r_{34}),\\
&& X_{i+2}=\frac{1}{\sqrt{m}}\sum_{j=1}^m(1-X_{ij}^2),\ \  \ \ \ \ \ \ \ \ Y_j\sim \frac{1}{m}\chi^2(m)
\eeaa
for $i=1, 2, 3, 4$ and $j=1,2,3, 4.$ Set $\xi_{1j}=X_{1j}X_{2j}-r_{12}$,  $\xi_{2j}=X_{3j}X_{4j}-r_{34}$ and $\xi_{i+2\, j}=1-X_{ij}^2$ for $i=1, 2, 3, 4$. Notice $X_{ij}\sim N(0, 1)$ and $|r_{ij}|\leq 1$ for each $i,j$. By Lemma~\ref{garbage},
\beaa
\Big|r_{12}^ar_{34}^b\cdot E\frac{\bar{B}_1^{t_1}\bar{B}_2^{t_2}(1-A_1)^{n_1}(1-A_2)^{n_2}(1-A_3)^{n_3}
(1-A_4)^{n_4}}{A_1^{s_1}A_2^{s_2}A_3^{s_3}A_4^{s_4}}\Big|\leq \frac{K_1}{m^{(N+1)/2}},
\eeaa
where $K_1$ is a constant depending on $N$ but not on $m$ or $\bd{R}.$ This confirms claim \eqref{djcjwe}.  The first covariance on the right hand side of \eqref{ewy129} is studied in Lemma ~\ref{main_an}. Combining this lemma and \eqref{djcjwe}, we finish the proof. \hfill$\square$

\begin{lemma}\lbl{main_theorem1} Let $\X_1,\cdots,\X_n$ be a random sample from $N_p(\bmu, \bms)$ with correlation matrix $\bd{R}$. Let $\hat{\bd{R}}$ be defined in \eqref{sample_corr_ma}. Assume, for some $a>0$, $p\leq n^{a}$ for each $p\geq 1$. If\,
  $\lim_{p\to\infty}\frac{p}{n\|\bd{R}\|_F}= 0$, then  $\mbox{Var}(\mbox{tr}(\hat{\bd{R}}^2))\cdot \|\bd{R}\|_F^{-2}$ goes to zero as $p\to \infty$.
\end{lemma}
\noindent\textbf{Proof of Lemma~\ref{main_theorem1}}.  Set $m=n-1$. By \eqref{ifeie} from the proof of Lemma~\ref{sdu329},
\beaa
\hat{\bd{R}}=\hat{\bd{R}}_p=(\hat{r}_{ij})_{p\times p}\ \overset{d}{=}  \Big(\frac{\bd{v}_i^T\bd{v}_j}{\|\bd{v}_i^T\|\cdot\|\bd{v}_j\|}\Big)_{p\times p},
\eeaa
where the $m$ rows of $(\bd{v}_1, \cdots, \bd{v}_p)_{m\times p}$ are i.i.d. with distribution $N_p(\bd{0}, \bd{R})$. Write
$\mbox{tr}(\hat{\bd{R}}^2)=p+\sum_{1\leq i\ne j \leq m}\hat{r}_{ij}^2.$
Then,
\bea\lbl{gywul}
\mbox{Var}\big(\mbox{tr}(\hat{\bd{R}}^2)\big)
=\mbox{Cov}\Big(\sum_{1\leq i\ne j \leq p}\hat{r}_{ij}^2, \sum_{1\leq k\ne l \leq p}\hat{r}_{kl}^2\Big)=\sum\,\mbox{Cov}\big(\hat{r}_{ij}^2, \hat{r}_{kl}^2\big),
\eea
where the last sum runs over all $(i, j, k, l)\in \Lambda_p$, where
\bea\lbl{simple1}
\Lambda_p:=\big\{(i, j, k, l); \, 1\leq i\ne j \leq p\ \mbox{and}\ 1\leq k\ne l\leq p\big\}.
\eea
Review Theorem~\ref{Theorem4}. We never impose any condition on the $4\times 4$ correlation matrix $\bd{R}_{4\times4}$ (not confuse the $p\times p$ correlation matrix $\bd{R}$ here) in the proposition. For example, if all of the entries of $\bd{R}_{4\times4}$ are equal to $1$, then the four random variables are actually equal. Keeping this understanding in mind,
by changing ``$(1,2,3,4)$" in Theorem~\ref{Theorem4} to ``$(i, j, k, l)$"  and taking $N=7$ in the proposition, we see $\mbox{Cov}(\hat{r}_{ij}^2, \hat{r}_{kl}^2)$ is equal to
\bea
\varrho_{p,1}\cdot \sum_{u, v\in \{i, j, k, l\}, u\ne v}r_{uv}^2 + \varrho_{p,2}\cdot r_{ij}^2r_{kl}^2 &+& \varrho_{p,3}\cdot\big(r_{ij}r_{jk}r_{kl}r_{li}+r_{ij}r_{jl}r_{lk}r_{ki}\big) \nonumber\\
&+&\varrho_{p,4}\cdot\big[(r_{ik}r_{kl}r_{li})r_{ij}^2+
(r_{jk}r_{kl}r_{lj})r_{ij}^2\big]\nonumber\\
&+ & \varrho_{p,5}\cdot\big[(r_{ij}r_{jk}r_{ki})r_{kl}^2+(r_{ij}r_{jl}r_{li})r_{kl}^2\big] \nonumber\\
&+& \frac{\varrho_{p,6}}{m^{4}}, \lbl{hwqiu}
\eea
where the sum runs over the six pairs from $\{i, j, k, l\}$,
\bea\lbl{wvhwq76}
|\varrho_{p,1}| \leq Km^{-2}; \ \ |\varrho_{p,2}|\vee|\varrho_{p,3}|\vee|\varrho_{p,4}|\vee|\varrho_{p,5}| \leq K m^{-1},\ |\varrho_{p,6}| \leq K,
\eea
$\{\varrho_{p,i};\, 3\leq i \leq 5\}$ do not depend on $\bd{R}$, and $K$ is a constant not depending on $m$ or $\bd{R}.$ Notice
\beaa
\sum_{(i,j,k,l)\in \Lambda_p}\sum_{u, v\in \{i, j, k, l\}, u\ne v}r_{uv}^2\leq 6p^2\sum_{1\leq i, j\leq p}r_{ij}^2=6p^2\cdot \mbox{tr}(\bd{R}^2).
\eeaa
Also,
\beaa
\sum_{(i,j,k,l)\in \Lambda_p}r_{ij}^2r_{kl}^2\leq \sum_{1\leq i,j\leq p}r_{ij}^2\cdot \sum_{1\leq k,l\leq p}r_{kl}^2=\, \big[\mbox{tr}(\bd{R}^2)\big]^2.
\eeaa
Recall $\|\bd{R}\|_F^4=[\mbox{tr}(\bd{R}^2)]^2$. The two facts together with \eqref{gywul} and \eqref{hwqiu} imply that
\beaa
&&\frac{1}{\|\bd{R}\|_F^4}\cdot \mbox{Var}\big(\mbox{tr}(\hat{\bd{R}}^2)\big)\\
& \leq & (6K)\cdot \frac{p^2}{m^2\cdot\mbox{tr}(\bd{R}^2)} + \frac{K}{m}  + \Big(\sum_{i=3}^5|m\varrho_{p,i}|\cdot |Q_{p,i}|\Big)+  \frac{Kp^4}{[\mbox{tr}(\bd{R}^2)]^2\cdot m^{4}}
\eeaa
where
\beaa
Q_{p,3}&=&\frac{1}{m[\mbox{tr}(\bd{R}^2)]^2}\cdot \sum_{(i,j,k,l)\in \Lambda_p}\big(r_{ij}r_{jk}r_{kl}r_{li}+r_{ij}r_{jl}r_{lk}r_{ki}\big), \nonumber\\
Q_{p,4}&=&\frac{1}{m[\mbox{tr}(\bd{R}^2)]^2}\cdot\sum_{(i,j,k,l)\in \Lambda_p}\big[(r_{ik}r_{kl}r_{li})r_{ij}^2+
(r_{jk}r_{kl}r_{lj})r_{ij}^2\big],\nonumber\\
Q_{p,5}&=& \frac{1}{m[\mbox{tr}(\bd{R}^2)]^2}\cdot\sum_{(i,j,k,l)\in \Lambda_p}\big[(r_{ij}r_{jk}r_{ki})r_{kl}^2+(r_{ij}r_{jl}r_{li})r_{kl}^2\big].
\eeaa
Now, by assumption, $p=o(m\|\bd{R}\|_F)$, hence
\beaa
\frac{p^2}{m^2\cdot\mbox{tr}(\bd{R}^2)}\to 0\ \ \mbox{and}\ \ \frac{Kp^4}{[\mbox{tr}(\bd{R}^2)]^2\cdot m^{4}}\to 0.
\eeaa
Because of \eqref{wvhwq76}, to prove the conclusion, it suffices to show $\lim_{p\to\infty}Q_{p,i}=0$ for $i=3, 4, 5.$ Recall \eqref{simple1}. By switching ``$k$" and ``$l$" in $r_{ij}r_{jl}r_{lk}r_{ki}$,  switching ``$k$" and ``$l$" in $(r_{jk}r_{kl}r_{lj})r_{ij}^2$ and switching ``$k$" and ``$l$" in $(r_{ij}r_{jl}r_{li})r_{kl}^2$, respectively, we obtain
\beaa
Q_{p,3}&=&\frac{2}{m[\mbox{tr}(\bd{R}^2)]^2}\cdot \sum_{(i,j,k,l)\in \Lambda_p}r_{ij}r_{jk}r_{kl}r_{li}, \nonumber\\
Q_{p,4}&=&\frac{2}{m[\mbox{tr}(\bd{R}^2)]^2}\cdot\sum_{(i,j,k,l)\in \Lambda_p}(r_{ik}r_{kl}r_{li})r_{ij}^2,\nonumber\\
Q_{p,5}&=& \frac{2}{m[\mbox{tr}(\bd{R}^2)]^2}\cdot\sum_{(i,j,k,l)\in \Lambda_p}(r_{ij}r_{jk}r_{ki})r_{kl}^2.
\eeaa
By interchanging ``$(i, j)$" with ``$(k, l)$" in the last sum, we see $Q_{p,4}=Q_{p,5}$. Finally, we see from Lemma~\ref{ewf394} that $\lim_{p\to\infty}Q_{p,i}=0$ for $i=3, 4, 5.$ \hfill$\square$

\noindent\textbf{Proof of Theorem~\ref{Theorem3}}. Recall $\hat{\bd{R}}=\hat{\bd{R}}_p$. By Lemma~\ref{main_theorem1},
\bea\lbl{ffwqi2}
\frac{\mbox{tr}(\hat{\bd{R}}_p^2)-E\,\mbox{tr}(\hat{\bd{R}}_p^2)}{\mbox{tr}(\bd{R}_p^2)} \to 0
\eea
in probability as $p\to\infty.$ By Lemma~\ref{sdu329}, under the assumption $\limsup_{p\to\infty}\frac{p}{n^{a}}=0$ for some constant $a>0$, we have
\beaa
E\,\mbox{tr}(\hat{\bd{R}}_p^2)
=\frac{p(p-1)}{n-1}+\mbox{tr}(\bd{R}_p^2)\cdot \big[1+O(m^{-1/4})\big].
\eeaa
This implies that
\beaa
\frac{1}{\mbox{tr}(\bd{R}_p^2)}\cdot \Big[E\,\mbox{tr}(\hat{\bd{R}}_p^2)-\frac{p(p-1)}{n-1}-\mbox{tr}(\bd{R}_p^2)\Big]\to 0.
\eeaa
The proof is completed by adding this and that from \eqref{ffwqi2}. \hfill$\square$

\newpage

\subsection{The Proofs of Theorems~\ref{Theorem1} and~\ref{Theorem2} and Proposition~\ref{Lemma_Remark_1}}\lbl{wefuewf}

Let  $\bm{\xi}_1, \cdots, \bm{\xi}_n$ be i.i.d. $p$-dimensional random vectors with distribution $N_p(\bm{\mu}, \bms)$. Let $\D$ be the diagonal matrix of $\bms$. The $p\times p$ population correlation matrix is defined by $\R=\D^{-1/2}\bms\D^{-1/2}$. The sample mean is  $\bar{\bm{\xi}}=\frac{1}{n}(\bm{\xi}_1+\cdots \bm{\xi}_n)$ and the sample covariance matrix is defined by
\bea\lbl{daloe}
\hat{\S}=\frac{1}{n}\sum_{i=1}^n (\bm{\xi}_i-\bar{\bm{\xi}})(\bm{\xi}_i-\bar{\bm{\xi}})^T.
\eea
Review $W_p(m, \bms)$ stands for the distribution of  the Wishart matrix $\bd{U}^T\bd{U}$ for any $m\geq 1$, where $\bd{U}$ is an $m\times p$ matrix whose rows are i.i.d. with distribution $N_p(\bd{0}, \bms)$. Then, $n\hat{\S}$ has the Wishart distribution $W_p(n-1, \bms)$; see, for example, Theorem 3.1.2 from \cite{muirhead1982aspects}. Let $\hat{\D}$ be the diagonal matrix of $\hat{\S}$. Then  $\hat{\R}:=\hat{\D}^{-1/2}\hat{\S}\hat{\D}^{-1/2}$ is the
 sample correlation matrix generated by $\bm{\xi}_1, \cdots, \bm{\xi}_n$. Before proving Theorems~\ref{Theorem1} and~\ref{Theorem2}, we first will reduce the test statistic appearing in Theorem~\ref{Theorem1} to a simple form. Recall we assume $n$ depends on $p$ and sometimes write $n_p$ if there is any possible confusion. Also, the Frobenius norm $\|\bd{R}\|_F=[\mbox{tr}(\bd{R}^2)]^{1/2}$ and the notation $o_p(1)$ representing a random variable converging to $0$ in probability.

\noindent\textbf{Proof of Lemma~\ref{Jiaodong}}. We need to show
\bea\lbl{wgyfu}
\frac{\bm{\eta}^T\hat{\D}^{-1}\bm{\eta}-\bm{\eta}^T\D^{-1}\bm{\eta}}{\sqrt{\,\mbox{tr}(\bd{R}^2)}}
\to 0
\eea
in probability as $p\to \infty.$ To do so, it suffices to prove that both its mean and variance converging to $0$.

{\it Step 1: the mean of random variable from \eqref{wgyfu}}. Write $\bm{\eta}=\bms^{1/2}\bm{\theta}$ where  $\bm{\theta}\sim N_p(\bd{0}, \bm{I}_p)$ and $\bms^{1/2}$ is a non-negative definite matrix satisfying $\bms^{1/2}\cdot \bms^{1/2}=\bms$.   By assumption, $\bm{\theta}$ is independent of $\hat{\S}$. In particular, $\bm{\theta}$ is independent of $\hat{\D}$, the diagonal matrix of $\hat{\S}$. Notice
\bea\lbl{dshewiu}
\bm{\eta}^T\hat{\D}^{-1}\bm{\eta}=\bm{\theta}^T\big(\bms^{1/2}\hat{\D}^{-1}\bms^{1/2}\big)\bm{\theta}
\ \ \mbox{and}\ \ \bm{\eta}^T\D^{-1}\bm{\eta}=\bm{\theta}^T\big(\bms^{1/2}\D^{-1}\bms^{1/2}\big)\bm{\theta}.
\eea
For any $p\times p$ symmetric matrix $\bd{A}$ with eigenvalues $\lambda_1, \cdots, \lambda_p$, by the orthogonal invariance of $N_p(\bd{0}, \bm{I}_p)$, we know $\bm{\theta}^T\bd{A}\bm{\theta}$ and $\lambda_1\theta_1^2+\cdots +\lambda_p\theta_p^2$ have the same distribution, where $\theta_1, \cdots, \theta_p$ are i.i.d. $N(0, 1)$-distributed random variables. Consequently,
\bea\lbl{uh2380}
E(\bm{\theta}^T\bd{A}\bm{\theta})=\, \mbox{tr}(\bd{A})\ \ \mbox{and}\ \ \mbox{Var}(\bm{\theta}^T\bd{A}\bm{\theta})=2\,\mbox{tr}(\bd{A}^2).
\eea
It follows from independence and conditioning on $\hat{\D}$ that
\bea\lbl{zhemekuai}
E(\bm{\eta}^T\hat{\D}^{-1}\bm{\eta})=E\,\mbox{tr}\big(\bms^{1/2}\hat{\D}^{-1}\bms^{1/2}\big)
=\mbox{tr}\big[\bms^{1/2}E\big(\hat{\D}^{-1}\big)\bms^{1/2}\big]
\eea
by linearity of expectations and traces, where $E\big(\hat{\D}^{-1}\big)$ is the entry-wise expectation of the diagonal matrix $\hat{\D}^{-1}$.  Set $\bms=(\sigma_{ij})_{p\times p}$. Then $\bd{D}=\mbox{diag}(\sigma_{11}, \cdots, \sigma_{pp})$. Set $m=n-1$. It is known
\bea\lbl{hongji}
n\hat{\S} \overset{d}{=}\sum_{j=1}^{m} \hat{\bm{\xi}}_j\hat{\bm{\xi}}_j^T
\eea
and $\hat{\S}$ is independent of $\bar{\bm{\xi}}=\frac{1}{n}(\bm{\xi}_1+\cdots \bm{\xi}_n)$, where $\hat{\bm{\xi}}_1, \cdots, \hat{\bm{\xi}}_m$ are i.i.d. $N_p(\bd{0}, \bm{\Sigma})$-distributed random vectors; see, for example, Theorem 3.1.2 from \cite{muirhead1982aspects}. Write $\hat{\bm{\xi}}_j=(\xi_{1j}, \cdots, \xi_{pj})^T$ for each $j.$ Then the $(i, i)$-entry of $\hat{\bm{\xi}}_j\hat{\bm{\xi}}_j^T$ is equal to $\xi_{ij}^2$. As a result,
\bea\lbl{evu2389}
\mbox{the}\ \mbox{$(i, i)$-entry of}\ n\hat{\S}\ \mbox{is}\ \sum_{j=1}^m\xi_{ij}^2 \sim \sigma_{ii}^2\cdot \chi^2(m)
\eea
for each $1\leq i\leq p$. Since $\hat{\D}=\mbox{diag}(s_{11}, \cdots, s_{pp})$ is the diagonal matrix of $\hat{\S}:=(s_{ij})_{p\times p}$, we know $n s_{ii}/\sigma_{ii}\sim \chi^2(m)$ for each $i$.
It is known that
\bea\lbl{kewu}
E\frac{1}{\chi^2(k)}=\frac{1}{k-2}\ \ \mbox{and}\ \ \mbox{Var}\Big(\frac{1}{\chi^2(k)}\Big)=\frac{2}{(k-2)^2(k-4)}
\eea
for any integer $k\geq 3$. Therefore,
\bea\lbl{pagezi}
E\frac{1}{s_{ii}}=\frac{n}{(m-2)\sigma_{ii}}\ \ \mbox{and}\ \ \ \mbox{Var}\Big(\frac{1}{s_{ii}}\Big)=\frac{2n^2}{(m-2)^2(m-4)\sigma_{ii}^2}.
\eea
It follows that $E(\hat{\D}^{-1})=\frac{n}{m-2}\bd{D}^{-1}.$ Observe   $\mbox{tr}(\bms^{1/2}\bd{D}^{-1}\bms^{1/2})=
\mbox{tr}(\bd{D}^{-1/2}\bms\bd{D}^{-1/2})=\mbox{tr}(\bd{\bd{R}})=p$. From \eqref{zhemekuai}, we have
\bea\lbl{xiwanga}
E(\bm{\eta}^T\hat{\D}^{-1}\bm{\eta})=\frac{n}{m-2}\cdot \mbox{tr}\big(\bms^{1/2}\bd{D}^{-1}\bms^{1/2}\big)=\frac{np}{n-3}.
\eea
Similarly, we have from \eqref{dshewiu} that
\beaa
E(\bm{\eta}^T\D^{-1}\bm{\eta})=E(\bm{\theta}^T\bms^{1/2}\D^{-1}\bms^{1/2}\bm{\theta})
=\,\mbox{tr}(\bms^{1/2}\D^{-1}\bms^{1/2})=p.
\eeaa
Therefore,
\beaa
E(\bm{\eta}^T\hat{\D}^{-1}\bm{\eta})-E(\bm{\eta}^T\D^{-1}\bm{\eta})=\frac{np}{m-2}-p=\frac{3p}{m-2}.
\eeaa
It follows that
\bea\lbl{ehuowe0}
\frac{1}{\sqrt{2\,\mbox{tr}(\bd{R}^2)}}
\big[E(\bm{\eta}^T\hat{\D}^{-1}\bm{\eta})-E(\bm{\eta}^T\D^{-1}\bm{\eta})\big] =\frac{\sqrt{4.5}\, p}{(m-2)\sqrt{\mbox{tr}(\bd{R}^2)}}\to 0
\eea
by the assumption $\lim_{p\to\infty}\frac{p}{m\|\bd{R}\|_F}=0$.

{\it Step 2: the variance of random variable from \eqref{wgyfu}}. Set $\B=\bms^{1/2}(\hat{\D}^{-1}-\D^{-1})\bms^{1/2}$. It is seen from \eqref{dshewiu} that
\bea\lbl{uhwiu}
\bm{\theta}^T\B\bm{\theta}=\bm{\eta}^T\hat{\D}^{-1}\bm{\eta}-\bm{\eta}^T\D^{-1}\bm{\eta}.
\eea
Recall the formula $\mbox{Var}(v)=E\mbox{Var}(v|\B)+ \mbox{Var}(E(v|\B))$ for any random variable $v$. Then, by the independence between $\bm{\theta}$ and $\B$ as well as \eqref{uh2380},
\bea\lbl{jinxindongpo}
\mbox{Var}(\bm{\theta}^T\B\bm{\theta})=2E\,\mbox{tr}(\B^2) + \mbox{Var}(\mbox{tr}(\B)).
\eea
Our focus next will be the evaluation of the two terms.

{\it Step 3: the evaluation of $E\mbox{tr}(\B^2)$ from \eqref{jinxindongpo}}.
Let us consider the last two terms one by one. First,
\beaa
\mbox{tr}(\B^2)&=&\mbox{tr}\big(\bms^{1/2}(\hat{\D}^{-1}-\D^{-1})
\bms(\hat{\D}^{-1}-\D^{-1})\bms^{1/2}\big)\\
& = & \mbox{tr}\big((\hat{\D}^{-1}-\D^{-1})
\bms(\hat{\D}^{-1}-\D^{-1})\bms\big).
\eeaa
Let $\bd{Q}(i, j)$ denote the $(i,j)$-entry of a matrix $\bd{Q}$. For any matrices $\bd{Q}_1, \bd{Q}_2$, $\bd{Q}_3, \bd{Q}_4$, we have $\mbox{tr}(\bd{Q}_1\bd{Q}_2\bd{Q}_3\bd{Q}_4)=\sum \bd{Q}_1(i, j)\bd{Q}_2(j, k)\bd{Q}_3(k, l)\bd{Q}_4(l, i)$, where the sum runs over all possible indices $i, j, k, l$. It follows that
\bea\lbl{vi54}
\mbox{tr}(\B^2)
&=&\sum_{1\leq i, j \leq p}\sigma_{ij}^2\Big(\frac{1}{s_{ii}}-\frac{1}{\sigma_{ii}}\Big)\Big(\frac{1}{s_{jj}}
-\frac{1}{\sigma_{jj}}\Big) \nonumber\\
&=& \sum_{1\leq i, j \leq p}r_{ij}^2\Big(\frac{\sigma_{ii}}{s_{ii}}-1\Big)\Big(\frac{\sigma_{jj}}{s_{jj}}
-1\Big)
\eea
since $r_{ij}=\sigma_{ij}(\sigma_{ii}\sigma_{jj})^{-1/2}$. By \eqref{pagezi},
\beaa
\frac{\sigma_{ii}}{s_{ii}}-1=\frac{\sigma_{ii}}{s_{ii}}-E\frac{\sigma_{ii}}{s_{ii}} + \frac{3}{m-2}.
\eeaa
Therefore,
\bea\lbl{ciub}
E\Big(\frac{\sigma_{ii}}{s_{ii}}-1\Big)\Big(\frac{\sigma_{jj}}{s_{jj}}
-1\Big)=\mbox{Cov}\Big(\frac{\sigma_{ii}}{s_{ii}}, \frac{\sigma_{jj}}{s_{jj}}\Big) + \frac{9}{(m-2)^2}.
\eea
The fact from \eqref{evu2389} implies that $ns_{ii}=X_1^2+\cdots + X_{m}^2$ and $ns_{jj}=Y_1^2+\cdots + Y_{m}^2$, where $(X_1, Y_1)^T, \cdots, (X_{m}, Y_{m})^T$ are i.i.d. $2$-dimensional normal random vectors with $EX_1=EY_1=0$, $EX_1^2=\sigma_{ii}$ and $EY_1^2=\sigma_{jj}$ and $\mbox{Cov}(X_1, Y_1)=\sigma_{ij}$. Recall $\R=(r_{ij})_{p\times p}$ with $r_{ij}=\sigma_{ij}/\sqrt{\sigma_{ii}\sigma_{jj}}.$ Then $\mbox{Cov}(X_1/\sqrt{\sigma_{ii}}, Y_1/\sqrt{\sigma_{jj}})=r_{ij}$. By Lemma~\ref{youdaoa}, we have
\beaa
\frac{m^2}{n^2}\cdot E\Big(\frac{\sigma_{ii}}{s_{ii}}\cdot\frac{\sigma_{ii}}{s_{jj}}\Big)
= 1+ \frac{4+2r_{ij}^2}{m}+\frac{12+8r_{ij}^2+8r_{ij}^4}{m^2}+\frac{\delta_m(i, j)}{m^3},
\eeaa
where  $\max_{1\leq i, j\leq p}|\delta_m(i, j)|\leq C$ for all $m\geq 11$, where $C$ is a constant not depending on $m$ or $\bd{R}=(r_{ij})$. This and \eqref{pagezi} conclude that $\mbox{Cov}(\frac{\sigma_{ii}}{s_{ii}},\, \frac{\sigma_{jj}}{s_{jj}})$ is identical to
\bea\lbl{maoju}
&& \Big[1+ \frac{4+2r_{ij}^2}{m}+\frac{12+8r_{ij}^2+8r_{ij}^4}{m^2}+\frac{\delta_m(i, j)}{m^3}\Big]\cdot \Big(\frac{m+1}{m}\Big)^2
-\frac{m+1}{m-2}\cdot \frac{m+1}{m-2} \nonumber\\
&=& \Big[1+ \frac{4+2r_{ij}^2}{m}+\frac{\delta_m(i, j)'}{m^2}\Big]\cdot \Big(1+\frac{2}{m}+\frac{1}{m^2}\Big)-\Big[1+\frac{6}{m}+\frac{21m-24}{m(m-2)^2}\Big]\nonumber\\
&=& \frac{2r_{ij}^2}{m}+\frac{\delta_m(i, j)''}{m^2},
\eea
where
\bea\lbl{vop485}
\max_{1\leq i,j\leq p}\big\{|\delta_m(i,j)'|, |\delta_m(i,j)''|\big\}\leq K_1
\eea
for all $m\geq 11$ and $K_1$ here and later represents a constant not depending on $m$ or $r_{ij}$, and can be different from line to line.   This,~\eqref{vi54} and \eqref{ciub} conclude
\bea\lbl{eutvyqbi}
E\,\mbox{tr}(\B^2)
&=&\frac{2}{m}\Big(\sum_{1\leq i, j \leq p}r_{ij}^4\Big) + \sum_{1\leq i, j \leq p}r_{ij}^2\Big[\frac{\delta_m(i,j)'}{m^2}+ \frac{9}{(m-2)^2}\Big] \nonumber\\
& \leq & \frac{2}{m}\,\mbox{tr}(\R^2) + \frac{K_1}{m^2}\,\mbox{tr}(\R^2) \nonumber\\
& \leq & \frac{3}{m}\,\mbox{tr}(\R^2)
\eea
as $m$ is sufficiently large,  which is guaranteed as $p\to\infty$ since $\lim_{p\to\infty}n_p=\infty$ and $m=n_p-1$. In the second step above we use the fact   $\sum_{1\leq i, j \leq p}r_{ij}^4\leq \sum_{1\leq i, j \leq p}r_{ij}^2=\mbox{tr}(\R^2)$.

{\it Step 4: the evaluation of $\mbox{Var}(\mbox{tr}(\B))$ from \eqref{jinxindongpo}}.  Note
\beaa
\mbox{tr}(\B)=\mbox{tr}\big(\bms^{1/2}(\hat{\D}^{-1}-\D^{-1})\bms^{1/2}\big)
&= &
\mbox{tr}\big((\hat{\D}^{-1}-\D^{-1})\bms\big)\\
& = & \sum_{i=1}^m\Big(\frac{1}{s_{ii}}-\frac{1}{\sigma_{ii}}\Big)\sigma_{ii}.
\eeaa
Recall $\frac{ns_{ii}}{\sigma_{ii}} \sim \chi^2(m)$ for each $i$.
It then follows from \eqref{pagezi} and \eqref{vop485} that
\beaa
\mbox{Var}(\mbox{tr}(\B))
&=&\mbox{Var}\Big(\sum_{i=1}^m\frac{\sigma_{ii}}{s_{ii}}\Big)\\
&=& \sum_{i=1}^m\mbox{Var}\Big(\frac{\sigma_{ii}}{s_{ii}}\Big) + 2\sum_{1\leq i< j \leq p}\mbox{Cov}\Big(\frac{\sigma_{ii}}{s_{ii}}, \frac{\sigma_{jj}}{s_{jj}}\Big)\\
& \leq & m\cdot \frac{2n^2}{(m-2)^2(m-4)} + 2\sum_{1\leq i< j \leq p}\Big[\frac{2r_{ij}^2}{m} + \frac{\delta_m(i,j)''}{m^2}\Big].
\eeaa
Thus,
\beaa
\mbox{Var}(\mbox{tr}(\B)) \leq 3 + \frac{2}{m}\,\mbox{tr}(\R^2) + \frac{K_1p^2}{m^2}
\eeaa
as $m$ is sufficiently large.

Finally, combining the analysis of the two terms from \eqref{jinxindongpo} in {\it Step 3} and {\it Step 4}, we eventually obtain
\beaa
\mbox{Var}(\bm{\theta}^T\B\bm{\theta})
&\leq & \frac{6}{m}\,\mbox{tr}(\R^2) +3 + \frac{2}{m}\,\mbox{tr}(\R^2) + K_1\frac{p^2}{m^2}\\
& = & 3+\frac{8}{m}\,\mbox{tr}(\R^2) + K_1\frac{p^2}{m^2}
\eeaa
as $p$ is sufficiently large. Easily, $\mbox{tr}(\R^2)\geq p$. It follows that
\beaa
\mbox{Var}\Big(\frac{\bm{\theta}^T\B\bm{\theta}}{\sqrt{2\,\mbox{tr}(\bd{R}^2)}}\Big)
=\frac{\mbox{Var}(\bm{\theta}^T\B\bm{\theta}) }{2\,\mbox{tr}(\bd{R}^2)} \leq \frac{3}{p} + \frac{4}{m}+K_1\frac{p^2}{m^2\,\mbox{tr}(\bd{R}^2)}\to 0
\eeaa
since $\lim_{p\to\infty}\frac{p}{m\|\bd{R}\|_F}=0$ and $\|\bd{R}\|_F^2=\,\mbox{tr}(\bd{R}^2)$. This joined \eqref{ehuowe0} concludes \eqref{wgyfu}. \hfill$\square$

\newpage

\noindent\textbf{Proof of Lemma~\ref{vdu349j}}. First, by the monotone property of $a_{p,i}$, we obtain $\rho_1\geq \rho_2\geq \cdots$. Moreover, $1=a_{p,1}^2+ \cdots+ a_{p, p}^2\geq a_{p,1}^2+ \cdots+ a_{p, i}^2\geq ia_{p, i}^2$ for any $1\leq i \leq p.$ This implies that $0\leq a_{p,i}\leq i^{-1/2}$ for each $1\leq i \leq p.$ Take $p\to\infty$ to obtain $0\leq \rho_i\leq i^{-1/2}$ for each $i \geq 1.$ Also, by using the fact $1\geq a_{p,1}^2+ \cdots+ a_{p, i}^2$ for any $1\leq i \leq p$, and letting $p\to\infty$ first and then $i\to\infty$, we get $\sum_{i=1}^{\infty}\rho_i^2\leq 1.$

We first handle a trivial case: $\rho_1=0$. By monotonicity, $\rho_i=0$ for each $i\geq 1$. Notice $E\xi_1=0$ and $\mbox{Var}(\xi_1)=2$. Thus, $s_n^2:=\mbox{Var}(a_{p,1}\xi_1+\cdots + a_{p,p}\xi_p)=2(a_{p,1}^2+\cdots + a_{p,p}^2)=2.$ Easily,
\beaa
\frac{1}{s_n^3}\sum_{i=1}^pE(|a_{p,i}\xi_i|^3)
=\frac{E(|\xi_1|^3)}{2\sqrt{2}}\sum_{i=1}^pa_{p,i}^3\leq \frac{E(|\xi_1|^3)}{2\sqrt{2}}\cdot a_{p,1}\cdot \sum_{i=1}^pa_{p,i}^2,
\eeaa
which goes to zero by the assumption $\lim_{p\to\infty}a_{p,1}=\rho_1= 0$. The desired result follows from the Lyapunov central limit theorem. From now on, we assume $\rho_1>0$.

In the following a useful fact will be derived first. For each $p\geq 1$, let $b_{p,1}\geq b_{p,2} \geq \cdots \geq 0$ be constants satisfying  $\sum_{i=1}^{\infty}b_{p,i}^2\leq 1$. We claim that
\bea\lbl{edoue90}
\prod_{i=m}^{\infty}\Big[e^{-tb_{p,i}}\big(1-2tb_{p,i}\big)^{-1/2}\Big]
=e^{\gamma_{p,m}}\cdot\exp\Big(t^2\sum_{i=m}^{\infty}b_{p,i}^2\Big)
\eea
for all $m\geq 16$ and $|t|<1$, where $\sup_{p\geq 1}|\gamma_{p,m}|\leq \frac{8}{\sqrt{m}}$. In fact,  write $\log (1-x)=-\sum_{i=1}^{\infty}\frac{1}{i}x^i:=-x-\frac{1}{2}x^2-B(x)$ for $|x|<1$. Then
\bea\lbl{dvhu380}
|B(x)|\leq \sum_{i=3}^{\infty}\frac{1}{i}|x|^i\leq \sum_{i=3}^{\infty}|x|^i\leq \frac{|x|^3}{1-|x|}\leq 2|x|^3
\eea
if $|x|\leq \frac{1}{2}.$ By the same argument as that in the beginning, we know  $0\leq b_{p,i}\leq \frac{1}{\sqrt{i}}$ for each $i\geq 1$ and $p\geq 1$. Observe
\beaa
\prod_{i=m}^{\infty}\Big[e^{-tb_{p,i}}\big(1-2tb_{p,i}\big)^{-1/2}\Big]
&=&\prod_{i=m}^{\infty}\exp\Big[-tb_{p,i}-\frac{1}{2}\log \big(1-2tb_{p,i}\big)\Big] \nonumber\\
&=& \prod_{i=m}^{\infty}\exp\Big[t^2b_{p,i}^2+\frac{1}{2}B(2tb_{p,i})\Big].\lbl{3vohq}
\eeaa
By the monotone property, $\max_{i\geq m}|2tb_{p,i}|=2|t|b_{p,m}\leq \frac{2|t|}{\sqrt{m}}\leq \frac{1}{2}|t|$ for $m\geq 16$. This and \eqref{dvhu380} say that
\beaa
\sum_{i=m}^{\infty}\frac{1}{2}|B(2tb_{p,i})| \leq 8|t|^3\sum_{i=m}^{\infty}b_{p,i}^3\leq \frac{ 8|t|^3}{\sqrt{m}}\sum_{i=m}^{\infty}b_{p,i}^2\leq \frac{8}{\sqrt{m}}
\eeaa
for any $t$ with $|t|<1$ and $p\geq 1$. These lead to \eqref{edoue90}.  In two steps next we will apply \eqref{edoue90}  to $a_{p,1}\xi_1+\cdots + a_{p,p}\xi_p$ and its limit stated in the lemma, respectively. The limit case goes first.

{\it Step 1}. Set $b=[2(1-\sum_{i=1}^{\infty}\rho_i^2)]^{1/2}$ and $X=b\eta+ \sum_{i=1}^{\infty}\rho_i\xi_i$, where $\eta\sim N(0, 1)$ and $\eta$ is independent of $\{\xi_i;\, i\geq 1\}$. Then, by independence and the fact $E\exp(t\chi^2(1))=(1-2t)^{-1/2}$ for $t<\frac{1}{2}$, we see
\bea\lbl{dvpij90}
Ee^{tX}=\Big(\prod_{i=1}^{\infty}Ee^{t\rho_i\xi_1}\Big)\cdot Ee^{tb\eta}=e^{b^2t^2/2}\cdot \prod_{i=1}^{\infty}\Big[e^{-t\rho_i}\big(1-2t\rho_i\big)^{-1/2}\Big]
\eea
for $t$ with $|t\rho_i|<\frac{1}{2}$ for each $i\geq 1$, which holds as $|t|<\frac{1}{2\rho_1}.$
Take $b_{p,i}=\rho_i$ for all $i\geq 1$ and $p\geq 1$ in \eqref{edoue90} to see
\bea\lbl{yidw73b}
\prod_{i=m}^{\infty}\Big[e^{-t\rho_i}\big(1-2t\rho_i\big)^{-1/2}\Big]
=e^{\gamma_{m}}\cdot\exp\Big(t^2\sum_{i=m}^{\infty}\rho_i^2\Big)
\eea
for $m\geq 16$ and $|t|<1$, where $|\gamma_{m}|\leq \frac{8}{\sqrt{m}}$. This and \eqref{dvpij90} especially indicate $Ee^{tX}<\infty$ for every $|t|<\frac{1}{2}.$ Recall $\sum_{i=1}^{\infty}\rho_i^2\leq 1$, by sending $m\to\infty$ we see the left hand side of \eqref{yidw73b} goes to $1$. Therefore,
\bea\lbl{dj38da}
\prod_{i=1}^{m-1}\Big[e^{-t\rho_i}\big(1-2t\rho_i\big)^{-1/2}\Big]\to \prod_{i=1}^{\infty}\Big[e^{-t\rho_i}\big(1-2t\rho_i\big)^{-1/2}\Big]
\eea
as $m\to\infty$ for every $|t|<\frac{1}{2}$.

{\it Step 2}. Evidently,
\bea\lbl{duoweu}
E^{t(a_{p,1}\xi_1+\cdots + a_{p,p}\xi_p)}
=\prod_{i=1}^{p}Ee^{ta_{p,i}\xi_1}
=\prod_{i=1}^{p}\Big[e^{-ta_{p,i}}\big(1-2ta_{p,i}\big)^{-1/2}\Big]
\eea
provided $|t|<\frac{1}{2a_{p,1}}$. In particular, this holds if $|t|<\frac{1}{2}$.
Now, by taking $b_{p,i}=a_{p,i}$ for $1\leq i \leq p$ and  $b_{p,i}=0$ for $i>p$ from \eqref{edoue90}, we obtain
\beaa
E^{t(a_{p,1}\xi_1+\cdots + a_{p,p}\xi_p)}
=e^{\gamma_{p,m}}\cdot \exp\Big(t^2\sum_{i=m}^pa_{p,i}^2\Big)\cdot \prod_{i=1}^{m-1}\Big[e^{-ta_{p,i}}\big(1-2ta_{p,i}\big)^{-1/2}\Big]
\eeaa
for any $m$ with $16\leq m\leq p$ and $|t|<\frac{1}{2}$, where $\sup_{p\geq 1}|\gamma_{p,m}|\leq \frac{8}{\sqrt{m}}$. Consequently, if $16\leq m\leq p$ and $|t|<\frac{1}{2}$ then
\bea\lbl{guangbangzi1}
E^{t(a_{p,1}\xi_1+\cdots + a_{p,p}\xi_p)}\leq e^{8/\sqrt{m}}\cdot \exp\Big[t^2\Big(1-\sum_{i=1}^{m-1}a_{p,i}^2\Big)\Big]\cdot \prod_{i=1}^{m-1}\Big[e^{-ta_{p,i}}\big(1-2ta_{p,i}\big)^{-1/2}\Big]
\eea
by the assumption $a_{p, 1}^2+ \cdots+ a_{p, p}^2=1$, and
\bea\lbl{guangbangzi2}
E^{t(a_{p,1}\xi_1+\cdots + a_{p,p}\xi_p)}\geq e^{-8/\sqrt{m}}\cdot \exp\Big[t^2\Big(1-\sum_{i=1}^{m-1}a_{p,i}^2\Big)\Big]\cdot \prod_{i=1}^{m-1}\Big[e^{-ta_{p,i}}\big(1-2ta_{p,i}\big)^{-1/2}\Big].
\eea

With the two steps established above, we are now ready to complete the proof.
Recall the assumption $\lim_{p\to\infty}a_{p,i}=\rho_i$ for each $i\geq 1$.  For fixed $m\geq 16$ we send $p\to\infty$ and then send $m\to\infty$ in \eqref{guangbangzi1} and \eqref{guangbangzi2}, we have from \eqref{dj38da} and then \eqref{dvpij90} that
\beaa
E^{t(a_{p,1}\xi_1+\cdots + a_{p,p}\xi_p)}\to \exp\Big[t^2\Big(1-\sum_{i=1}^{\infty}\rho_{i}^2\Big)\Big]\cdot \prod_{i=1}^{\infty}\Big[e^{-t\rho_{i}}\big(1-2t\rho_{i}\big)^{-1/2}\Big]=Ee^{tX}
\eeaa
as $p\to\infty$ for $|t|<\frac{1}{2}$. The desired conclusion then follows from the uniqueness of the moment generating function. \hfill$\square$

\medskip

Recall $F(1, m)$ stands for the $F$-distribution with degrees of freedoms $1$ and $m$.

\begin{lemma}\lbl{vhufew08} Let $m=m_p\to \infty$ as $p\to\infty$.  For each $p\geq 1$, let $X_{p,1}, \cdots, X_{p,p}$ be i.i.d. with distribution $F(1, m)$. Then $(2p)^{-1/2}[X_{p,1}+ \cdots+ X_{p,p}-mp(m-2)^{-1}]\to N(0, 1)$
in distribution as $p\to\infty.$
\end{lemma}
\noindent\textbf{Proof of Lemma~\ref{vhufew08}}. First, by the property of $F$-distribution,
\beaa
EX_{p,1}=\frac{m}{m-2} ~~~ \mbox{and} ~~~\mbox{Var}(X_{p,1})=\frac{2m^2(m-1)}{(m-2)^2(m-4)}
\eeaa
for $m\geq 5$. By definition, we write $X_{p,1}=\frac{m\xi_0^2}{\xi_1^2+\cdots + \xi_m^2}$, where $\xi_0, \xi_1, \cdots, \xi_m$ are i.i.d. $N(0, 1)$. Then
\beaa
E\big(X_{p,1}-EX_{p,1}\big)^4
&= &E\Big[\Big(\frac{m}{\xi_1^2+\cdots + \xi_m^2}-1\Big)\xi_0^2+\xi_0^2-\frac{m}{m-2}\Big]^4\\
& \leq & 3^3E\Big[\Big(\frac{m}{\xi_1^2+\cdots + \xi_m^2}-1\Big)^4\xi_0^8\Big] + 3^3E(\xi_0^8)
+3^3\Big(\frac{m}{m-2}\Big)^4.
\eeaa
By using the Cauchy-Schwartz inequality twice,
\beaa
E\Big[\Big(\frac{m}{\xi_1^2+\cdots + \xi_m^2}-1\Big)^4\xi_0^8\Big]
& \leq & \Big\{E\Big[\Big(\frac{\xi_1^2+\cdots + \xi_m^2-m}{\xi_1^2+\cdots + \xi_m^2}\Big)^8\Big]\Big\}^{1/2}\cdot \big(E\xi_0^{16}\big)^{1/2}\\
& \leq & K\cdot \Big[E\big(\xi_1^2+\cdots + \xi_m^2-m\big)^{16}\Big]^{1/4}\cdot \Big[E\frac{1}{(\xi_1^2+\cdots + \xi_m^2)^{16}}\Big]^{1/4},
\eeaa
where $K$ here and later is a constant free of $m$ and $p$, and can be different from line to line.
By using the Marcinkiewicz-Zygmund inequality [see, for example,  the proof of Corollary 2 on p. 387 from \cite{chow1997probability}], $E\big(\xi_1^2+\cdots + \xi_m^2-m\big)^{16}\leq Km^8$.  Furthermore, take $\beta=-16$ in \eqref{cdh8293f} to see
 $E[(\xi_1^2+\cdots + \xi_m^2)^{-16}]\leq Km^{-16}$ for all $m\geq 34$. Combining all of the above calculation, we see $E\big(X_{p,1}-EX_{p,1}\big)^4\leq K$ as $m\geq 34$. Notice $\mbox{Var}(X_{p,1}) \to 2$ as $p\to \infty$. Then
 \beaa
 \frac{1}{(p\mbox{Var}(X_{p,1}))^2}\sum_{i=1}^pE(X_{p,i}-EX_{p,i})^4 =O\Big(\frac{1}{p}\Big) \to 0
 \eeaa
as $p\to\infty$. By the Lyapunov CLT, we obtain the desired result. \hfill$\square$

\medskip

\noindent\textbf{Proof of Theorem~\ref{Theorem1}}. First, by Theorem~\ref{Theorem3},
\bea\lbl{csiy289}
\frac{1}{\mbox{tr}(\bd{R}^2)}\Big[\mbox{tr}(\hat{\bd{R}}^2)-\frac{p(p-1)}{n-1}
\Big] \to 1
\eea
in probability as $p\to\infty$. In the following we will use this fact twice to show $T_{SD}$ and $T_{p,1}$ are equivalent. First, it follows from the assumption $\lim_{p\to\infty}\frac{p}{n\|\bd{R}\|_F}= 0$ that
\beaa
\frac{\tr(\hat{\R}^2)-p^2(n-1)^{-1}}{\tr(\R^2)}=
\frac{\tr(\hat{\R}^2)-p(p-1)(n-1)^{-1}}{\tr(\R^2)}-
\frac{p(n-1)^{-1}}{\tr(\R^2)}=1+o_p(1).
\eeaa
As a consequence,
\beaa
H_{p}:=\Big[\frac{\tr(\hat{\R}^2)-p^2(n-1)^{-1}}{\tr(\R^2)}\Big]^{-1/2}=1+o_p(1).
\eeaa
Review \eqref{jin_wuzu}. We have
\bea\lbl{duowe8}
T_{SD}
&=&\frac{[n\bar{\X}^{T} \hat{\D}^{-1}\bar{\X}-pn(n-3)^{-1}]+p(n-3)^{-1}}{\sqrt{2\tr(\R^2)}}
\cdot H_{p} \nonumber\\
&=& \frac{[n\bar{\X}^{T} \hat{\D}^{-1}\bar{\X}-pn(n-3)^{-1}]}{\sqrt{2\tr(\R^2)}}\cdot[1+ o_p(1)] + o_p(1)
\eea
by  the assumption $\frac{p}{n\|\bd{R}\|_F}\to 0$ and the notation  $\|\bd{R}\|_F^2=\mbox{tr}(\bd{R}^2)$. By \eqref{Statistics1},
\beaa
T_{p,1}=\frac{n\bar{\X}^{T}\hat{\D}^{-1}\bar{\X}-pn(n-3)^{-1}}{\sqrt{ 2 \big|\tr(\hat{\R}^2)-p(p-1)(n-1)^{-1}\big|}}.
\eeaa
It follows from \eqref{csiy289} that
\beaa
T_{p,1}=\frac{[n\bar{\X}^{T} \hat{\D}^{-1}\bar{\X}-pn(n-3)^{-1}]}{\sqrt{2\tr(\R^2)}}\cdot[1+ o_p(1)].
\eeaa
Comparing this with \eqref{duowe8}, we obtain
\bea\lbl{fehi32r80}
T_{SD}=T_{p,1}\cdot[1+ o_p(1)] + o_p(1).
\eea
So to finish the proof, by using the Slutsky lemma, it suffices to prove that $T_{p,1}\to  (1-\sum_{i=1}^{\infty}\rho_i^2)^{1/2}\xi_0+\frac{1}{\sqrt{2}}\sum_{i=1}^{\infty}\rho_i(\xi_i^2-1)$ in distribution.

Set $m=n-1$. Then
\bea\lbl{wdquigr9}
&& \sqrt{n}\bar{\X} \sim N_p(\bd{0}, \bms),\   n\hat{\S} \sim W_p(m, \bms),\ \mbox{and}\ \bar{\X}\ \mbox{and}\ \hat{\S}\ \mbox{are independent},
\eea
where $W_p(m, \bms)$ is the Wishart distribution defined after \eqref{daloe}; see, for example, Theorem 3.1.2 from \cite{muirhead1982aspects}. This implies
\bea\lbl{diuewi}
\hat{\D}\,\overset{d}{=}\, \mbox{the diagonal matrix of}\ \frac{1}{n} W_p(m, \bms).
\eea
In particular, $\bar{\X}$ is independent of $\hat{\D}$. By Lemma~\ref{Jiaodong} and assumption  $\lim_{p\to\infty}\frac{p}{m\|\bd{R}\|_F}=0$,
\beaa
\frac{n\bar{\X}^{T} \hat{\D}^{-1}\bar{\X}}{\sqrt{2\,\mbox{tr}(\bd{R}^2)}}
=\frac{n\bar{\X}^{T} \D^{-1}\bar{\X}}{\sqrt{2\,\mbox{tr}(\bd{R}^2)}} + o_p(1)
\eeaa
as $p\to\infty$, where $\D$ is the diagonal matrix of the population covariance matrix $\bms$.  Now we will analyze the behavior of $n\bar{\X}^{T} \D^{-1}\bar{\X}$   in two steps. Once they are established, the limiting distribution of $T_{p,1}$ will be identified quickly.

{\it Step 1. The exact distribution of $n\bar{\X}^{T} \D^{-1}\bar{\X}$}. By assumption, $\lambda_1, \cdots, \lambda_p$ are the eigenvalues of $\bd{R}=\D^{-1/2}\bms\D^{-1/2}$. We claim that
\bea\lbl{h28p0}
n\bar{\X}^{T} \D^{-1}\bar{\X} \overset{d}{=}
\sum_{i=1}^p\lambda_i\xi_i^2,
\eea
where  $\xi_1, \cdots, \xi_p$ are i.i.d. $N(0, 1)$ and where ``$\overset{d}{=}$" means both sides of ``$=$" have the same distribution. We  show \eqref{h28p0} next.

Since $\sqrt{n}\bar{\X} \sim N_p(\bd{0}, \bms)$, we are able to write $\sqrt{n}\bar{\X}=\bms^{1/2}\bm{\xi}$, where $\bms^{1/2}$ is a non-negative definite matrix satisfying $\bms^{1/2}\cdot \bms^{1/2}=\bms$ and $\bm{\xi}\sim N_p(\bd{0}, \bd{I}_p).$ It follows that
\bea\lbl{daqiao}
n\bar{\X}^{T} \D^{-1}\bar{\X}=\bm{\xi}^T\big(\bms^{1/2}\D^{-1}\bms^{1/2}\big)\bm{\xi}.
\eea
 Since $\bd{A}\bd{B}$ and $\bd{B}\bd{A}$ have the same eigenvalues for any square matrix $\bd{A}$ and $\bd{B}$, it is easy to see that $\lambda_1, \cdots, \lambda_p$ are also the eigenvalues of $\bms^{1/2}\D^{-1}\bms^{1/2}$. Write $\bms^{1/2}\D^{-1}\bms^{1/2}=\bd{O}^T\mbox{diag}(\lambda_1, \cdots, \lambda_p)\bd{O}$ for some orthogonal matrix $\bd{O}$. Then by the orthogonal invariance property, we know $\bd{O}\bm{\xi}\overset{d}{=}\bm{\xi}$. Consequently,
\bea\lbl{huwqiwq9}
\bm{\xi}^T\big(\bms^{1/2}\D^{-1}\bms^{1/2}\big)\bm{\xi}=(\bd{O}\bm{\xi})^T\mbox{diag}(\lambda_1, \cdots, \lambda_p)(\bd{O}\bm{\xi})\overset{d}{=}\sum_{i=1}^p\lambda_i\xi_i^2
\eea
where $\xi_1, \cdots, \xi_p$ are i.i.d. $N(0, 1)$. We then get \eqref{h28p0}.

{\it Step 2: the limiting distribution of $n\bar{\X}^{T} \D^{-1}\bar{\X}$}. By \eqref{h28p0} and the fact $\mbox{tr}(\bd{R})=\lambda_1+\cdots +\lambda_p=p$, we know
\beaa
\frac{n\bar{\X}^{T} \D^{-1}\bar{\X}-p}{\sqrt{2\,\mbox{tr}(\bd{R}^2)}}\overset{d}{=}
\frac{1}{\sqrt{2}}\sum_{i=1}^p\frac{\lambda_i}{\|\bd{R}\|_F}(\xi_i^2-1).
\eeaa
Set $a_{p, i}=\frac{\lambda_i}{\|\bd{R}\|_F}$ for $1\leq i \leq p$. Then $a_{p, 1}\geq \cdots \geq a_{p,p}\geq 0$ and $a_{p,1}^2+\cdots + a_{p,p}^2=1$. By assumption  (a), $\lim_{p\to\infty}a_{p,i}=\rho_i \geq 0$ for each $i\geq 1$. From Lemma~\ref{vdu349j},
\beaa
\frac{1}{\sqrt{2}}\sum_{i=1}^p\frac{\lambda_i}{\|\bd{R}\|_F}(\xi_i^2-1) \to \Big(1-\sum_{i=1}^{\infty}\rho_i^2\Big)^{1/2}\xi_0 + \frac{1}{\sqrt{2}}\sum_{i=1}^{\infty}\rho_i(\xi_i^2-1)
\eeaa
in distribution as $p\to\infty$, where $\xi_0\sim N(0,1)$ and $\xi_0$ is independent of $\{\xi_i;\, i\geq 1\}$. Therefore,
\bea\lbl{dvuh32r8}
\frac{n\bar{\X}^{T} \D^{-1}\bar{\X}-p}{\sqrt{2\,\mbox{tr}(\bd{R}^2)}} \to \Big(1-\sum_{i=1}^{\infty}\rho_i^2\Big)^{1/2}\xi_0 + \frac{1}{\sqrt{2}}\sum_{i=1}^{\infty}\rho_i(\xi_i^2-1)
\eea
in distribution.

With {\it Step 1} and {\it Step 2} completed, let us now proceed to finish the proof. In fact, by
assumption (b),
\beaa
\frac{1}{\sqrt{2}\|\bd{R}\|_F}\big[p-pn(n-3)^{-1}\big]=-\frac{3p}{\sqrt{2}(n-3)\|\bd{R}\|_F}\to 0.
\eeaa
Summing this and \eqref{dvuh32r8}, we see from the Slutsky lemma that
\beaa
\frac{n\bar{\X}^{T} \D^{-1}\bar{\X}-pn(n-3)^{-1}}{\sqrt{2}\|\bd{R}\|_F}\to
\Big(1-\sum_{i=1}^{\infty}\rho_i^2\Big)^{1/2}\xi_0 + \frac{1}{\sqrt{2}}\sum_{i=1}^{\infty}\rho_i(\xi_i^2-1)
\eeaa
in distribution. Combine this with \eqref{csiy289} to see
\beaa
T_{p,1},
\to \Big(1-\sum_{i=1}^{\infty}\rho_i^2\Big)^{1/2}\xi_0 + \frac{1}{\sqrt{2}}\sum_{i=1}^{\infty}\rho_i(\xi_i^2-1)
\eeaa
in distribution as $p\to \infty.$ The proof is completed.\hfill$\square$

\medskip

Now we begin to prove  Proposition~\ref{Lemma_Remark_1}. As far as the proofs go, parts (i) and (ii) from the lemma  have different natures, to make the presentation clear, we will handle the two parts separately.

\noindent\textbf{Proof of Proposition~\ref{Lemma_Remark_1}(i)}. Assume $\bm{\mu}=0$. Then $\bd{X}_1, \cdots, \bd{X}_n$ are i.i.d. with distribution $N_p(\bd{0}, \bms)$, where all of the $p^2$ entries of $\bms$ are identical to $\sigma^2$. For this reason, we write
$\X_i=\xi_i(1, \cdots, 1)^T\in \mathbb{R}^p$ for $i=1,\cdots, n$, where $\xi_1, \cdots, \xi_n$ are i.i.d. $N(0, \sigma^2)$. By \eqref{wuoi0} and \eqref{sample_corr_ma}, $\bar{\X}=\frac{1}{n}\sum_{i=1}^n\X_i$,
\beaa
\hat{\S}=\frac{1}{n}\sum_{i=1}^n (\X_i-\bar{\X})(\X_i-\bar{\X})^T ~~~~~ \mbox{and} ~~~~~
\hat{\R}=\hat{\D}^{-1/2}\hat{\S}\hat{\D}^{-1/2},
\eeaa
where $\hat{\D}$ is the diagonal matrix of $\hat{\S}$. Set $\bar{\xi}=\frac{1}{n}\sum_{i=1}^n\xi_i$ and $W=\frac{1}{n}\sum_{i=1}^n(\xi_i-\bar{\xi})^2$. Then $\sqrt{n}\bar{\xi}/\sigma\sim N(0, 1)$,  $nW/\sigma^2\sim \chi^2(n-1)$, and $\bar{\xi}$ and $W$ are independent. Since $\bar{\X}=\bar{\xi}(1, \cdots, 1)^T\in \mathbb{R}^p$, we see  $\X_i-\bar{\X}=(\xi_i-\bar{\xi})^T(1, \cdots, 1)^T\in \mathbb{R}^p$. Hence, $\hat{\S}=W\cdot \bd{J}$, where $\bd{J}$ is a $p\times p$ matrix whose entries are all equal to $1$. It follows $\hat{\D}=W\cdot \bd{I}_p$ and $\hat{\R}=\bd{J}$. In particular, $\tr(\hat{\R}^2)=p^2$. Consequently,
\bea\lbl{dvsipewp}
n\bar{\X}^{T} \hat{\D}^{-1}\bar{\X}=\frac{np}{n-1}\cdot \frac{(\sqrt{n}\bar{\xi}/\sigma)^2}{nW/[(n-1)\sigma^2]} \overset{d}{=} \frac{np}{n-1}\cdot F_{1, n-1}.
\eea
By \eqref{jin_wuzu},
\bea\lbl{spartan}
T_{SD}=\frac{n\bar{\X}^{T} \hat{\D}^{-1}\bar{\X}-p(n-1)(n-3)^{-1}}{\sqrt{2[\tr(\hat{\R}^2)-p^2(n-1)^{-1}]}}
&\overset{d}{=} &
\frac{\frac{n}{n-1}\cdot F_{1, n-1}-\frac{n-1}{n-3}}{\sqrt{2[1-(n-1)^{-1}]}}.
\eea
Notice $F_{1, n-1}\to \chi^2(1)$ as $n\to\infty$. We see from the Slutsky lemma that $T_{SD}\to \frac{1}{\sqrt{2}}\cdot [\chi^2(1)-1]$ as $p\to\infty$. Finally, by \eqref{Statistics1} and \eqref{dvsipewp},
\beaa
T_{p,1}=\frac{n\bar{\X}^{T} \hat{\D}^{-1}\bar{\X}-pn(n-3)^{-1}}{\sqrt{ 2 \big|\tr(\hat{\R}^2)-p(p-1)(n-1)^{-1}\big|}}\overset{d}{=}\frac{\frac{np}{n-1}\cdot F_{1, n-1}-pn(n-3)^{-1}}{\sqrt{ 2 \big|p^2-p(p-1)(n-1)^{-1}\big|}}.
\eeaa
As a consequence,
\beaa
T_{p,1}
&\overset{d}{=} &
\frac{\frac{n}{n-1}\cdot F_{1, n-1}-n(n-3)^{-1}}{\sqrt{2\{1-(p-1)/[(n-1)p]\}}}\to \frac{1}{\sqrt{2}}\cdot \big[\chi^2(1)-1\big].
\eeaa
The proof is completed. \hfill$\square$

\medskip

\noindent\textbf{Proof of Proposition~\ref{Lemma_Remark_1}(ii)}.   Recall
\beaa
T_{SD}=\frac{n\bar{\X}^{T} \hat{\D}^{-1}\bar{\X}-p(n-1)(n-3)^{-1}}{\sqrt{2[\tr(\hat{\R}^2)-p^2(n-1)^{-1}]}}.
\eeaa
Since $T_{SD}$ is scale-invariant, without loss of generality, we assume $\bms=\bd{I}_p$. By \eqref{wdquigr9},
\beaa
&& \sqrt{n}\bar{\X} \sim N_p(\bd{0}, \bd{I}_p),\   n\hat{\S} \sim W_p(m, \bd{I}_p),\ \mbox{and}\ \bar{\X}\ \mbox{and}\ \hat{\S}\ \mbox{are independent},
\eeaa
where $m:=n-1$ and $W_p(m, \bms)$ is the Wishart distribution defined after \eqref{daloe}. Also, by \eqref{diuewi}, $\hat{\D}$ and  the diagonal matrix of $\frac{1}{n} W_p(m, \bd{I}_p)$ have the same distribution. By definition of $W_p(m, \bd{I}_p)$, its $p$ diagonal entries are i.i.d. $\chi^2(m)$. Therefore, $n\bar{\X}^{T} \hat{\D}^{-1}\bar{\X}$ is a sum of $p$ i.i.d. random variables with distribution $n\cdot \frac{N(0,1)^2}{\chi^2(m)}$, where the numerator and denominator are independent. Thus, we are able to write
\bea\lbl{dog_bite}
n\bar{\X}^{T} \hat{\D}^{-1}\bar{\X}=\frac{n}{m}\sum_{i=1}^pX_{p,i},
\eea
where $X_{p,1}, \cdots, X_{p,p}$ are i.i.d. with distribution $F(1, m)$. By Lemma~\ref{vhufew08},
\bea\lbl{GMBueffet}
U_p:=\frac{1}{\sqrt{2p}}\Big(\sum_{i=1}^pX_{p,i}-\frac{mp}{m-2}\Big)\to N(0, 1)
\eea
as $p\to\infty$, regardless of the speeds of $n$ and $p$ going to infinity. Solve $\sum_{i=1}^pX_{p,i}$ in terms of $U_p$ to see $ \sum_{i=1}^pX_{p,i}=\sqrt{2p}U_p+\frac{mp}{m-2}$. By plugging this into \eqref{dog_bite}, we obtain
\bea\lbl{fuhw8}
n\bar{\X}^{T} \hat{\D}^{-1}\bar{\X}=\frac{n}{m}\sqrt{2p}U_p+\frac{np}{m-2}.
\eea
Consequently,
\beaa
n\bar{\X}^{T} \hat{\D}^{-1}\bar{\X}-p(n-1)(n-3)^{-1}
&=&\frac{n}{m}\sqrt{2p}U_p+\frac{np}{m-2}-\frac{mp}{m-2}\\
&=& \sqrt{2p}U_p + \frac{\sqrt{2p} }{m}U_p+ \frac{p}{m-2}.
\eeaa
It follows that
\bea\lbl{fewuoefw90}
\frac{1}{\sqrt{2p}}\big[n\bar{\X}^{T} \hat{\D}^{-1}\bar{\X}-p(n-1)(n-3)^{-1}\big]=U_p+\frac{1 }{m}U_p+\frac{\sqrt{p}}{\sqrt{2}(m-2)}.
\eea

On the other hand,  $\hat{\R}=\hat{\D}^{-1/2}\hat{\S}\hat{\D}^{-1/2}$ by \eqref{sample_corr_ma}. Write $n\hat{\S}=(\bm{\xi}_1, \cdots, \bm{\xi}_p)^T(\bm{\xi}_1, \cdots, \bm{\xi}_p)$, where $\bm{\xi}_1, \cdots, \bm{\xi}_p$ are i.i.d. $N_m(0, \bd{I}_m)$. Set
$\bd{e}_i=\frac{\bm{\xi}_i}{\|\bm{\xi}_1\|}$ for $i=1,\cdots, p$. Then $\{\bd{e}_1, \cdots, \bd{e}_p\}$ are i.i.d. uniformly distributed on the $m$-dimensional sphere; see, for example, p. 38 from \cite{muirhead1982aspects}. Also, $\hat{\R}=(\bd{e}_i^T\bd{e}_j)$. Set
$V_p:=\frac{1}{p}\sum_{1\leq i< j \leq p}\big[m(\bd{e}_i'\bd{e}_j)^2-1\big].$
We claim
\bea\lbl{fhir0}
\sup_{p\geq 2}\mbox{Var}(V_p)\leq \frac{3}{2}.
\eea
In fact, by Lemma 11 from \cite{cai2013distributions} or (22) in Lemma 4.1 from \cite{cai2012phase}, $\{\bd{e}_i^T\bd{e}_j;\, 1\leq i< j \leq p\}$ are pairwise i.i.d. Thus, $\mbox{Var}(V_p)=\frac{1}{p^2}\cdot \frac{1}{2}p(p-1)\cdot m^2\cdot\mbox{Var}((\bd{e}_1'\bd{e}_2)^2).$ Since $\bd{e}_1$ and $\bd{e}_2$ are independent, by Theorems 1.5.6 and 1.5.7 from \cite{muirhead1982aspects}, we know $\bd{e}_1^T\bd{e}_2 \overset{d}{=}\xi_1(\xi_1^2+\cdots +\xi_m^2)^{-1/2}$, where $\xi_1, \cdots, \xi_m$ are i.i.d. $N(0, 1)$. In particular, if $n=2$ then $m=1$, $\bd{e}_1^T\bd{e}_2$ is a symmetric Bernoulli random variable, and hence $(\bd{e}_1^T\bd{e}_2)^2=1$ and $\mbox{Var}((\bd{e}_1^T\bd{e}_2))=0$. Now, for $n\geq 3$,
\beaa
\mbox{Var}((\bd{e}_1'\bd{e}_2)^2)=\mbox{Var}\Big(\frac{\xi_1^2}{\xi_1^2+\cdots +\xi_m^2}\Big)
\leq E\frac{\xi_1^4}{(\xi_1^2+\cdots +\xi_m^2)^2}=\frac{3}{m(m+2)}
\eeaa
by taking $a_1=2$ and other $a_i=0$ in Lemma 2.4 from Jiang (2012). The inequality is  true for all $m\geq 1$. Combining this and the earlier expression of $\mbox{Var}(V_p)$, we obtain \eqref{fhir0}. In particular, by the Chebyshev inequality, \eqref{fhir0} indicates
\bea\lbl{dskhwq8}
\frac{V_p}{m}\to 0
\eea
as $p\to\infty$ regardless of the speeds of $n$ and $p$ going to infinity. Write
\beaa
\sum_{1\leq i< j \leq p}(\bd{e}_i'\bd{e}_j)^2=\frac{p}{m}V_p+\frac{p(p-1)}{2m}.
\eeaa
Since $\hat{\R}$ is a symmetric matrix whose diagonal entries are all equal to $1$. The above implies
\bea\lbl{dvhw80}
\tr(\hat{\R}^2)=p+\frac{2p}{m}V_p+\frac{p(p-1)}{m}.
\eea
Consequently, we have from \eqref{dskhwq8} that
\bea\lbl{fviri27}
\frac{2}{p}\cdot \Big[\tr(\hat{\R}^2)-\frac{p^2}{m}\Big]=2+\frac{4}{m}V_p-\frac{2}{m} \to 2
\eea
in probability.  It follows from this and \eqref{fewuoefw90} that
\beaa
T_{SD}=\frac{\sqrt{2p}}{\sqrt{p(2+4m^{-1}V_p-2m^{-1})}}\cdot \Big[U_p+\frac{1}{m}U_p+\frac{\sqrt{p}}{\sqrt{2}(m-2)}\Big].
\eeaa
By \eqref{dskhwq8} and the Slutsky lemma, $\frac{\sqrt{2p}}{\sqrt{p(2+4m^{-1}V_p-2m^{-1})}}\to 1$ in probability. Therefore, by \eqref{GMBueffet},
\beaa
U_p+\frac{1}{m}U_p+\frac{\sqrt{p}}{\sqrt{2}(m-2)}\to
\begin{cases}
\xi_0, & \text{if $p/n^2\to 0$};\\
\xi_0+\sqrt{h/2}, & \text{if $p/n^2\to h$};\\
\infty, & \text{if $p/n^2\to \infty$}
\end{cases}
\eeaa
in distribution, where $\xi_0\sim N(0, 1)$. So we get part (ii) of Proposition~\ref{Lemma_Remark_1}.

Finally, as for $T_{p,1},$ in lieu of~\eqref{Statistics1},
\beaa
T_{p,1}=\frac{n\bar{\X}^{T} \hat{\D}^{-1}\bar{\X}-pn(n-3)^{-1}}{\sqrt{ 2 \big|\tr(\hat{\R}^2)-p(p-1)(n-1)^{-1}\big|}}
=\frac{\frac{n}{m}\sqrt{2p}U_p}{\sqrt{2p+\frac{4p}{m}V_p}}
=\frac{\frac{n}{m}U_p}{\sqrt{1+\frac{2}{m}V_p}}
\eeaa
by \eqref{fuhw8} and \eqref{dvhw80}. This implies $T_{p,1} \to N(0, 1)$
as $p\to\infty$ by \eqref{GMBueffet} and \eqref{dskhwq8}, regardless of the speeds of $n$ and $p$ going to infinity.
\hfill$\square$

\medskip

\noindent\textbf{Proof of Theorem~\ref{Theorem2}}.
Recall
\beaa
&& \bar{\bd{X}}_i=\frac{1}{n_i}\sum_{j=1}^{n_i}\bd{X}_{ij}\ \mbox{for}\ i=1, 2;\\
&& \hat{\S}=\frac{1}{n_1+n_2}\Big[\sum_{j=1}^{n_1}
(\bd{X}_{1j}-\bar{\bd{X}}_1)(\bd{X}_{1j}-\bar{\bd{X}}_1)^T+\sum_{j=1}^{n_2}
(\bd{X}_{2j}-\bar{\bd{X}}_2)(\bd{X}_{2j}-\bar{\bd{X}}_2)^T\Big].
\eeaa
The $p\times p$ matrix  $\hat{\D}$ is the diagonal matrix of $\hat{\S}$ and $\hat{\R}=\hat{\D}^{-1/2}\hat{\S}\hat{\D}^{-1/2}$ is the $p\times p$  pooled sample correlation matrix. In particular, all of the diagonal entries of $\hat{\R}$ are $1$. Thus,  $\|\bd{R}\|_F\geq \sqrt{p}$. Now,  $\bar{\bd{X}}_i$ is independent of $\sum_{j=1}^{n_i}(\bd{X}_{ij}-\bar{\bd{X}}_i)(\bd{X}_{ij}-\bar{\bd{X}}_i)^T$ for each $i=1, 2$; see, for example, Theorem 3.1.2 from \cite{muirhead1982aspects}. Second, by assumption, the two samples are independent. Therefore, the three random vectors $\bar{\bd{X}}_1$, $\bar{\bd{X}}_2$ and $\hat{\S}$ are independent. In particular, $\bar{\bd{X}}_1$ and $\bar{\bd{X}}_2$ are independent of $\hat{\D}$. By assumption $\bmu_1=\bmu_2$, we obtain
\beaa
\bar{\X}_1-\bar{\X}_2 \sim N_p\Big(\bd{0}, \Big(\frac{1}{n_1}+\frac{1}{n_2}\Big)\bms\Big).
\eeaa
Hence,
\bea\lbl{wefiufwp9}
\sqrt{\frac{n_1n_2}{n_1+n_2}}\big(\bar{\X}_1-\bar{\X}_2\big) \sim N_p(\bd{0}, \bms).
\eea
Furthermore, for each $i$, $\sum_{j=1}^{n_i}
(\bd{X}_{ij}-\bar{\bd{X}}_i)(\bd{X}_{ij}-\bar{\bd{X}}_i)^T \sim W_p(n_i-1, \bms)$; see, for example, Theorem 3.1.2 from \cite{muirhead1982aspects}. Therefore, $(n_1+n_2)\hat{\S}\sim W_p(n_1+n_2-2, \bms)$  by independence.  Consequently,
\bea\lbl{wefr2098}
\hat{\D}_1:=\frac{n_1+n_2}{n_1+n_2-1}\hat{\D}\ \overset{d}{=}\,\mbox{the diagonal matrix of}\ \frac{1}{m}W_p\big(m-1, \bms\big),
\eea
where $m:=n_1+n_2-1.$ As explained earlier, the left hand side of \eqref{wefiufwp9} is independent of $\hat{\D}_1$. In particular, by replacing ``$\bm{\eta}$" and ``$n$" from \eqref{xiwanga} with ``$\sqrt{\frac{n_1n_2}{n_1+n_2}}\big(\bar{\X}_1-\bar{\X}_2\big)$" and ``$m$", respectively, we obtain
\bea\lbl{tongzijieji}
&&E\Big[\frac{n_1n_2}{n_1+n_2-1}\big(\bar{\X}_1-\bar{\X}_2\big)^T
\hat{\D}^{-1}(\bar{\X}_1-\bar{\X}_2)\Big]\nonumber\\
&=&E\Big[\frac{n_1n_2}{n_1+n_2}\big(\bar{\X}_1-\bar{\X}_2\big)^T
\hat{\D}_1^{-1}(\bar{\X}_1-\bar{\X}_2)\Big]\nonumber\\
&=&\frac{(n_1+n_2-1)p}{n_1+n_2-4}.
\eea
This is the reason how come up with the numerator of $T_{p,2}$ defined in \eqref{Statistics2}.
The assumption $\lim_{p\to\infty}\frac{p}{(n_1+n_2)^{a}}=0$ implies $n_1+n_2\to\infty$ as $p\to\infty.$ Thus we know from assumption (b) that  $\lim_{p\to\infty}\frac{p}{m\|\bd{R}\|_F}= 0$.
Replacing ``$n$" with ``$m$" in Lemma~\ref{Jiaodong}, we have from~\eqref{wefiufwp9} and \eqref{wefr2098} that
\bea\lbl{sdhsiu}
&&\frac{\frac{n_1n_2}{n_1+n_2}\big(\bar{\X}_1-\bar{\X}_2\big)^T
\hat{\D}_1^{-1}(\bar{\X}_1-\bar{\X}_2)-p}
{\sqrt{2\tr(\R^2)}} \nonumber\\
&=& \frac{\frac{n_1n_2}{n_1+n_2}\big(\bar{\X}_1-\bar{\X}_2\big)^T
\D^{-1}(\bar{\X}_1-\bar{\X}_2)-p}{\sqrt{2\,\mbox{tr}(\bd{R}^2)}} + o_p(1),
\eea
where  $\D$ is the diagonal matrix of $\bms$, $\R=\D^{-1/2}\bms\D^{-1/2}$ is the population correlation matrix.  By \eqref{h28p0} and \eqref{wefiufwp9},
\bea\lbl{dhew93}
\frac{n_1n_2}{n_1+n_2}\big(\bar{\X}_1-\bar{\X}_2\big)^T
\D^{-1}(\bar{\X}_1-\bar{\X}_2) \overset{d}{=}\sum_{i=1}^p\lambda_i\xi_i^2,
\eea
where $\lambda_1\geq  \cdots \geq \lambda_p \geq 0$ are the eigenvalues of $\bd{R}$ and $\xi_1, \cdots, \xi_p$ are i.i.d. $N(0, 1).$ Similar to the argument as in {\it Step 2} in the proof of Theorem~\ref{Theorem1}, we see
\beaa
\frac{1}{\|\bd{R}\|_F}\Big(-p+\sum_{i=1}^p\lambda_i\xi_i^2\Big) \to \sqrt{2}\Big(1-\sum_{i=1}^{\infty}\rho_i^2\Big)^{1/2}\xi_0 + \sum_{i=1}^{\infty}\rho_i(\xi_i^2-1)
\eeaa
in distribution, where $\xi_0 \sim N(0, 1)$ and $\xi_0$ is independent of $\xi_1, \xi_2, \cdots$. Combining this with \eqref{sdhsiu} and \eqref{dhew93}, we have
\bea\lbl{do230r8}
&&\frac{\frac{n_1n_2}{n_1+n_2-1}\big(\bar{\X}_1-\bar{\X}_2\big)^T
\hat{\D}^{-1}(\bar{\X}_1-\bar{\X}_2)-p}
{\sqrt{2\tr(\R^2)}} \nonumber\\
&=&\frac{\frac{n_1n_2}{n_1+n_2}\big(\bar{\X}_1-\bar{\X}_2\big)^T
\hat{\D}_1^{-1}(\bar{\X}_1-\bar{\X}_2)-p}
{\sqrt{2\tr(\R^2)}} \nonumber\\
&\to &\Big(1-\sum_{i=1}^{\infty}\rho_i^2\Big)^{1/2}\xi_0 + \frac{1}{\sqrt{2}}\sum_{i=1}^{\infty}\rho_i(\xi_i^2-1)
\eea
in distribution.
Evidently,
\beaa
\frac{1}{\sqrt{2\tr(\R^2)}}\Big[p-\frac{(n_1+n_2-1)p}{n_1+n_2-4}\Big]=\frac{1}{\sqrt{2}}\frac{-3p}{(n_1+n_2-4)\|\bd{R}\|_F} \to 0
\eeaa
by the assumption $\lim_{p\to\infty}\frac{p}{(n_1+n_2)\|\bd{R}\|_F}= 0$. Add the left hand sides of the above two assertions to obtain
\bea\lbl{dh6iei}
&&\frac{\frac{n_1n_2}{n_1+n_2-1}\big(\bar{\X}_1-\bar{\X}_2\big)^T
\hat{\D}^{-1}(\bar{\X}_1-\bar{\X}_2)-\frac{(n_1+n_2-1)p}{n_1+n_2-4}}
{\sqrt{2\tr(\R^2)}}\nonumber\\
&\to & \Big(1-\sum_{i=1}^{\infty}\rho_i^2\Big)^{1/2}\xi_0 + \frac{1}{\sqrt{2}}\sum_{i=1}^{\infty}\rho_i(\xi_i^2-1)
\eea
in distribution. Next we will replace $\tr(\R^2)$ by its ratio-unbiased-estimator.

Recall $(n_1+n_2)\hat{\S}\sim W_p(n_1+n_2-2, \bms)$, and $\hat{\D}$ is the diagonal matrix of $\hat{\S}$ and $\hat{\R}=\hat{\D}^{-1/2}\hat{\S}\hat{\D}^{-1/2}$. Set $\kappa=\frac{n_1+n_2}{n_1+n_2-1}$. Then
\beaa
&&\kappa\hat{\S}\sim \frac{1}{n_1+n_2-1}W_p(n_1+n_2-2, \bms);\\
&&\hat{\R}=(\kappa\hat{\D})^{-1/2}(\kappa\hat{\S})(\kappa\hat{\D})^{-1/2}\ \ \mbox{and}\ \ \kappa\hat{\D}\ \mbox{is the diagonal matrix of}\ \kappa\hat{\S}.
\eeaa
 The essential assumption from Theorem~\ref{Theorem3} is that, in its own notation, $\hat{\bd{R}}$ is the sample correlation matrix obtained from $\hat{\S} \sim \frac{1}{n} W_p(n-1, \bms)$. Replace ``$\hat{\S}$" and ``$n$" from Theorem~\ref{Theorem3} with ``$\kappa\hat{\S}$" and ``$n_1+n_2-1$", respectively. Assumption (b) indicates $\lim_{p\to\infty}\frac{p}{(n_1+n_2-1)\|\bd{R}\|_F}= 0$ and $\lim_{p\to\infty}\frac{p}{(n_1+n_2-1)^{a}}=0$. Then by Theorem~\ref{Theorem3},
\bea\lbl{h2390}
\frac{1}{\mbox{tr}(\bd{R}^2)}\Big[\mbox{tr}(\hat{\bd{R}}^2)-\frac{p(p-1)}{n_1+n_2-2}\Big] \to 1
\eea
in probability.
This together with \eqref{dh6iei} and the Slutsky lemma yields
\beaa
T_{p, 2}&=&\frac{\frac{n_1n_2}{n_1+n_2-1}(\bar{\X}_1-\bar{\X}_2)^{T} \hat{\D}^{-1}(\bar{\X}_1-\bar{\X}_2)-\frac{(n_1+n_2-1)p}{n_1+n_2-4}}
{\sqrt{2\big|\tr(\hat{\R}^2)-\frac{p(p-1)}{n_1+n_2-2}\big|}}\\
&\to &  \Big(1-\sum_{i=1}^{\infty}\rho_i^2\Big)^{1/2}\xi_0 + \frac{1}{\sqrt{2}}\sum_{i=1}^{\infty}\rho_i(\xi_i^2-1)
\eeaa
 in distribution as $p\to\infty$.

Finally, by the same argument as deriving \eqref{fehi32r80}, we see $T_{SD}'=T_{p,2}\cdot[1+ o_p(1)] + o_p(1).$ So the conclusion for $T_{SD}'$ follows from the above and the Slutsky lemma. The proof is completed. \hfill$\square$

\subsection{A Lemma and Verification of
\eqref{aiy327},~\eqref{sfyi8} and \eqref{sdispi0sa}}\lbl{duhwoi9}

The following lemma is used in the discussion after Theorem~\ref{Theorem2}.
\begin{lemma}\lbl{Tuan_student}  Let $\{\bd{R}_p;\, p\geq 1\}$ be non-negative definite matrices whose diagonal entries are all equal to $1$. Let $\lambda_{p,1}\geq  \cdots \geq \lambda_{p,p} \geq 0$ be the eigenvalues of $\bd{R}_p$. Assume condition ``$C4$" stated in Theorem 2 from \cite{zhang2020simple} hold, that is, $\lim_{p\to\infty}\frac{\lambda_{p,i}}{\|\bd{R}_p\|_F}=\rho_i$ for all $i\geq 1$ with $\rho_1>0$ and $\lim_{p\to\infty}\sum_{i=1}^p\frac{\lambda_{p,i}}{\|\bd{R}_p\|_F}=\sum_{i=1}^{\infty}\rho_i<\infty$. Then $\sum_{i=1}^{\infty}\rho_i^2=1$.
\end{lemma}
\noindent\textbf{Proof of Lemma~\ref{Tuan_student}}. Notice $\sum_{i=1}^p\lambda_{p,i}=\mbox{tr}(\bd{R}_p)=p$ because all of the diagonal entries of $\bd{R}_p$ are identical to $1$. The assumption $\lim_{p\to\infty}\sum_{i=1}^p\frac{\lambda_{p,i}}{\|\bd{R}_p\|_F}=\sum_{i=1}^{\infty}\rho_i<\infty$ implies that
\bea\lbl{djewi9}
\frac{p}{\|\bd{R}_p\|_F}\to \sum_{i=1}^{\infty}\rho_i
\eea
as $p\to\infty$.  Set $a_{p,i}=\frac{\lambda_{p,i}}{\|\bd{R}_p\|_F}.$ Then $\lim_{p\to\infty}a_{p,i}=\rho_i$ and $a_{p,1}^2+\cdots + a_{p,p}^2=1$ by the definition of $\|\bd{R}_p\|_F$. For any $K\geq 1$, write
\bea\lbl{238cb}
a_{p,1}^2+\cdots +a_{p,K}^2=1-\sum_{i=K+1}^pa_{p,i}^2.
\eea
We claim that $\lim_{K\to\infty}\limsup_{p\to \infty}\sum_{i=K+1}^pa_{p,i}^2=0$. If this is true, by  letting $p\to \infty$ first and then sending $K\to \infty$, then  $\sum_{i=1}^{\infty}\rho_i^2=1$. We now prove the claim. In fact, write
\beaa
\sum_{i=K+1}^pa_{p,i}=\frac{p}{\|\bd{R}_p\|_F}-\sum_{i=1}^Ka_{p,i}.
\eeaa
For fixed $K\geq 1$, let $p\to\infty$ and use \eqref{djewi9} to have
\beaa
\lim_{p\to\infty}\sum_{i=K+1}^pa_{p,i}=\Big(\sum_{i=1}^{\infty}\rho_i\Big)-\sum_{i=1}^{K}\rho_i
=\sum_{i=K+1}^{\infty}\rho_i.
\eeaa
Use the assumption $\sum_{i=1}^{\infty}\rho_i<\infty$ and let $K\to\infty$ to see $\lim_{K\to\infty}\limsup_{p\to \infty}\sum_{i=K+1}^pa_{p,i}=0$. The claim is then verified since $a_{p,i}\leq 1$ for all $1\leq i\leq p$. \hfill$\square$.

\medskip

In the following the notation ``$A_p\sim B_p$" means that $A_p/B_p\to 1$ as $p\to\infty$.

\noindent\textbf{The Verification of \eqref{aiy327}}. Review Example~\ref{Ex1}. We have
\beaa
\bd{R}=
\begin{pmatrix}
1 & r & r &\cdots & r\\
r & 1 & r &\cdots & r\\
\vdots&&&&\\
r & r & r &\cdots & 1
\end{pmatrix}
_{p\times p}.
\eeaa
Easily, $\lambda_1=1+(p-1)r$ and  $\lambda_2=\cdots = \lambda_p=1-r$. To maintain $\bd{R}$ to be non-negative definite, all eigenvalues have to be non-negative, that is, $-\frac{1}{p-1}\leq r \leq 1.$ Recall $r=r_p$ satisfies that $\lim_{p\to\infty}\sqrt{p}\cdot r=c$. Now we consider three cases: $c=0$, $c\in (0, \infty)$ and $c=\infty$.

{\it Case 1: $c=0$}. In this case $\|\bd{R}\|_F^2=[1+(p-1)r]^2+ (p-1)(1-r)^2\sim p$ as $p\to\infty$, thus $\|\bd{R}\|_F \sim  \sqrt{p}$.  Easily, $\frac{p}{n\|\bd{R}\|_F}= O(\frac{\sqrt{p}}{n})\to 0$ provided $p=o(n^2)$. Also,  $\rho_1=\lim_{p\to\infty}\frac{\lambda_1}{\|\bd{R}\|_F}=0$, hence $\lambda_i=0$ for every $i\geq 1$.
By Theorem~\ref{Theorem1}, both $T_{SD}$ and $T_{p,1}$ converge to $N(0, 1)$ in distribution.

{\it Case 2: $c\in (0, \infty)$}. In this case $\|\bd{R}\|_F^2=[1+(p-1)r]^2+ (p-1)(1-r)^2\sim (c^2+1)p$ as $p\to\infty$, hence $\|\bd{R}\|_F \sim  \sqrt{c^2+1}\cdot\sqrt{p}$. Readily,
\beaa
&&\frac{p}{n\|\bd{R}\|_F}= O\Big(\frac{\sqrt{p}}{n}\Big)\to 0~~~~\mbox{and}~~~~\rho_1=\lim_{p\to\infty}\frac{\lambda_1}{\|\bd{R}\|_{F}}\to \frac{c}{\sqrt{c^2+1}};\\
&&
\rho_i=\lim_{p\to\infty}\frac{\lambda_i}{\|\bd{R}\|_F}=0~~ \mbox{for every}~~ i\geq 2
\eeaa
provided $p=o(n^2)$. By Theorem~\ref{Theorem1}, under condition $p=o(n^2)$, we know both $T_{SD}$ and $T_{p,1}$ converge to
\beaa
 \frac{1}{\sqrt{c^2+1}}\xi_0+\frac{c}{\sqrt{2(c^2+1)}}(\xi_1^2-1)
\eeaa
 in distribution, where $\xi_0$ and $\xi_1$ are i.i.d. $N(0, 1).$

{\it Case 3: $c=\infty$}. In this case $\|\bd{R}\|_F^2=[1+(p-1)r]^2+ (p-1)(1-r)^2\sim (pr)^2$ as $p\to\infty$. Therefore, $\|\bd{R}\|_F \sim  pr$.  Notice $\lim_{p\to\infty}p/(n\|\bd{R}\|_F)=\lim_{p\to\infty}1/(nr_p)= 0$, $\rho_1=\lim_{p\to\infty}\lambda_1/\|\bd{R}\|_F=1$ and $\rho_i=0$ for $i\geq 2$ provided $nr_p\to \infty$.
So, under conditions $nr_p\to \infty$ and $p=o(n^a)$ for some constant $a>0$, we have $T_{p,1}\to  \frac{1}{\sqrt{2}}[\chi^2(1)-1]$. \hfill$\square$

\medskip

\noindent\textbf{The Verification of \eqref{sfyi8}}. Recall $r\in (0, 1)$,  $m=[p^{r}]$ and  $\bd{A}_m$ from  Example~\ref{Ex1}. By definition,
\beaa
\bd{R}=
\begin{pmatrix}
\bd{A}_m & \bd{0}\\
\bd{0} & \bd{I}_{p-m}
\end{pmatrix}
.
\eeaa
Then the largest eigenvalues of $\bd{R}$ is $\lambda_1=1+(m-1)r$ and the rest of them are  $\lambda_2=\cdots =\lambda_m=1-r$ and $\lambda_{m+1}=\cdots=\lambda_p=1$. Then
\beaa
\|\bd{R}\|_{F}^2=[1+(m-1)r]^2 + (m-1)(1-r)^2+(p-m)=m^2r^2+p+O(m).
\eeaa
Thus,
\beaa
\|\bd{R}\|_{F}=
\begin{cases}
\sqrt{p}\cdot[1+o(1)], & \text{if $0< r < \frac{1}{2}$};\\
\sqrt{5p/4}\cdot[1+  o(1)],& \text{if $r = \frac{1}{2}$};\\
rp^{r}\cdot[1+ o(1)], & \text{if $\frac{1}{2}< r < 1$.}
\end{cases}
\eeaa
Also,
\beaa
\rho_1=\lim_{p\to\infty}\frac{\lambda_1}{\|\bd{R}\|_{F}} =
\begin{cases}
0, & \text{if $0< r < \frac{1}{2}$};\\
1/\sqrt{5},& \text{if $r = \frac{1}{2}$};\\
1, & \text{if $\frac{1}{2}< r < 1$.}
\end{cases}
\eeaa
Obviously, $\rho_i=0$ for $i\geq 2.$
At last,
\beaa
\frac{p}{n\|\bd{R}\|_F}=
\begin{cases}
O(\sqrt{p}/n), & \text{if $0< r \leq  \frac{1}{2}$};\\
O(p^{1-r}/n), & \text{if $\frac{1}{2}< r < 1$.}
\end{cases}
\eeaa
Assume $\bmu=\bm 0$. Then, we have from Theorem~\ref{Theorem1} that both $T_{SD}$ and $T_{p,1}$ go to
\beaa
\begin{cases}
N(0, 1), & \text{if $0< r < \frac{1}{2}$};\\
\frac{2}{\sqrt{5}}N(0, 1)+\frac{1}{\sqrt{10}}\cdot [\chi^2(1)-1],& \text{if $r = \frac{1}{2}$};\\
\frac{1}{\sqrt{2}}\cdot [\chi^2(1)-1], & \text{if $\frac{1}{2}< r <1$}
\end{cases}
\eeaa
in distribution under condition $p=o(n^2)$ for $0< r \leq  \frac{1}{2}$ and $p=o(n^{1/(1-r)})$ for $\frac{1}{2}< r <1$, where the random variables $N(0, 1)$ and $\chi^2(1)$ are independent as $r = \frac{1}{2}$. \hfill$\square$

\medskip

\noindent\textbf{The Verification of \eqref{sdispi0sa}}. Recall $m=[\log (p+2)]. $ Given $\tau\geq 0$, we have $p'=p-[\tau\sqrt{p}]$, $\lambda_i=1+\tau 2^{-i}\big(1-2^{-m}\big)^{-1}\sqrt{p}$ for $1\leq i \leq m$ and $\lambda_i=1$ for $i=m+1, \cdots, p'-1$,  $\lambda_{p'}=1+[\tau\sqrt{p}]-\tau\sqrt{p} \in [0, 1]$ and the rest of $\lambda_i$ are zero. Obviously, $\lambda_1+\cdots +\lambda_k\geq k$ for each $i=1, \cdots, p'-1$ and $\lambda_1+\cdots + \lambda_{p'}$ is identical to
\beaa
&& \sum_{i=1}^m\big[1+\tau 2^{-i}\big(1-2^{-m}\big)^{-1}\sqrt{p}\big]
+(p'-m-1)\cdot 1 + \big(1+[\tau\sqrt{p}]-\tau\sqrt{p}\big)\\
&=& (m+\tau\sqrt{p}) + \big(p-[\tau\sqrt{p}]-m-1\big) + 1+[\tau\sqrt{p}]-\tau\sqrt{p}\\
&=& p,
\eeaa
where we use the fact $\sum_{i=1}^m2^{-i}\big(1-2^{-m}\big)^{-1}=1$.
Lemma~\ref{dsih328} tells us that there exists a correlation matrix $\bd{R}$ such that $\bd{R}$ has eigenvalues $\lambda_i,\, 1\leq i \leq p.$ Now,
\beaa
\|\bd{R}\|_F^2
&=&\sum_{i=1}^m\big[1+\tau 2^{-i}\big(1-2^{-m}\big)^{-1}\sqrt{p}\big]^2
+(p'-m-1)\cdot 1^2 + \big(1+[\tau\sqrt{p}]-\tau\sqrt{p}\big)^2\\
&=& \tau^2p\big(1-2^{-m}\big)^{-2}\sum_{i=1}^m\frac{1}{4^i}+p +O\big(\sqrt{p}\big)\\
&=& \Big(\frac{\tau^2}{3}+1\Big)p +o(p)
\eeaa
since $\sum_{i=1}^m\frac{1}{4^i}\to \frac{1}{3}$ as $p\to \infty$. In particular, $\lim_{p\to\infty}\frac{p}{n\|\bd{R}\|_F}= 0$ provided $p=o(n^2)$. Easily, for each $i\geq 1$,
\beaa
\rho_i:=\lim_{p\to\infty}\frac{\lambda_i}{\|\bd{R}\|_F}=\lim_{p\to\infty}\frac{1+\tau 2^{-i}\big(1-2^{-m}\big)^{-1}\sqrt{p}}{\sqrt{(\frac{\tau^2}{3}+1)p+o(p)}}\to \sqrt{\frac{3\tau^2}{\tau^2+3}}\cdot \frac{1}{2^i}
\eeaa
as $p\to\infty$. Now,
\beaa
\sum_{i=1}^{\infty}\rho_i^2=\frac{3\tau^2}{\tau^2+3}\sum_{i=1}^{\infty}\frac{1}{4^i}
=\frac{\tau^2}{\tau^2+3}.
\eeaa
Assuming $\bmu=\bm 0$, by Theorem~\ref{Theorem1}, we see that
$T_{p,1}\to  \sqrt{\frac{3}{\tau^2+3}}\xi_0+\sqrt{\frac{3\tau^2}{2(\tau^2+3)}}\sum_{i=1}^{\infty}\frac{1}{2^i}(\xi_i^2-1)$ in distribution, where $\xi_0, \xi_1, \cdots, \xi_d$ are i.i.d. $N(0, 1).$

\bibliography{mybib}
\bibliographystyle{plainnat}

\end{document}